\documentclass[12pt]{article} 
      
\usepackage{amssymb}    
\usepackage{bbm}   
\usepackage{amsmath}       
\usepackage{latexsym}       
\usepackage{amsthm}       
\usepackage{epic,eepic,latexsym,amsbsy,amsmath,amscd,amssymb,mathtools,enumerate}

\numberwithin{equation}{section}

\input xy
\xyoption{all}

 \usepackage{amsmath,amssymb}
\usepackage{booktabs}
\usepackage{array}
\usepackage{tabularx}
\usepackage{rotating} 

\DeclareMathOperator{\Stab}{Stab}
\DeclareMathOperator{\Orb}{Orb}
\DeclareMathOperator{\Per}{Per}
\DeclareMathOperator{\Sym}{Sym}
\DeclareMathOperator{\Surj}{Surj}
\DeclareMathOperator{\diag}{diag}
\newcommand{\id}{\mathrm{id}}

\theoremstyle{definition} 
       
 \newtheorem{definition}{Definition} [section]       
 \newtheorem{remark}[definition]{Remark}       
 \newtheorem{example}[definition]{Example}

\theoremstyle{plain}       

 \newtheorem{proposition}[definition]{Proposition}       
 \newtheorem{theorem}[definition]{Theorem}       
 \newtheorem{corollary}[definition]{Corollary}       
  \newtheorem{lemma}[definition]{Lemma}

\newcommand{\Aut}{\mbox{\rm Aut}}

\renewcommand{\dim}{\mbox{dim}}

\newcommand{\Ind}{\mbox{\rm Ind}}

\newcommand{\Res}{\mbox{\rm Res}}

\newcommand{\Hom}{\mbox{\rm Hom}}

\begin{document} 

\title{Permutation representations of wreath products of finite groups} 
 \author{Fabio Scarabotti}  
 
\maketitle 

\begin{abstract}
We introduce and study a broad class of permutation actions of wreath products $F\wr G$ of finite groups, constructed from a finite set $V$ equipped with a faithful transitive $G$-action, a surjective map $\alpha\colon X\to A$ between finite sets, and a family of pairwise non-isomorphic transitive $F$-spaces $\{Y_a\}_{a\in A}$. Adopting a coordinate-free approach based on set compositions, we classify the $F\wr G$-orbits, determine their stabilizers and isomorphism classes, and describe the double-coset decompositions via iterated orbit partitions.

Applying Clifford theory and the representation theory of wreath products, we decompose the associated permutation representations into irreducible constituents and construct explicit bases for their isotypic components. We establish necessary and sufficient conditions for multiplicity freeness, obtaining explicit formulas for the corresponding spherical functions and their orthogonality relations, as well as criteria for the resulting Gelfand pairs to be symmetric. Finally, we specialize our general framework to three major settings: Young permutation representations for wreath products $F\wr \mathsf{Sym}(V)$ via semistandard tableaux and Kostka numbers; complete centipede subtrees and $k$-insects in finite rooted $q$-ary trees under the action of $\Aut(\mathcal{T}_{q,h})$; and a new family of finite Gelfand pairs for the generalized symmetric group $\mathsf{C}_n\wr \mathsf{Sym}(V)$ whose spherical functions are expressed in terms of Hahn polynomials.

\footnote{{\it AMS 2020 Math. Subj. Class.} Primary: 20C30, 20E22; Secondary: 05E10, 43A90, 20B25, 20E08.

\it Keywords: Wreath product, permutation representation, homogeneous space, Gelfand pair, spherical function, Young subgroup, generalized symmetric group, rooted tree.}
\end{abstract}
\tableofcontents

\section{Introduction}\label{Secintrod}

Wreath products are a classical construction in the theory of permutation groups and play a fundamental role in the study of groups acting on rooted trees. The representation theory of wreath products is also a classical subject, offering a rich and structured application of Clifford theory and the little group method of Mackey and Wigner; see, for instance, \cite{AdinRoich,book3,book6,JK,KT,MishSri}. The aim of the present paper is to construct a broad class of group actions for wreath products of finite groups and to investigate the corresponding permutation representations and harmonic analysis.

Our construction acts as a bridge between our recent work \cite{Trees2} on the representation theory of automorphism groups of finite rooted trees and earlier papers on wreath products and finite Gelfand pairs, namely \cite{CST2} (written jointly with Tullio Ceccherini-Silberstein and Filippo Tolli) and \cite{ScaTol} (written jointly with Filippo Tolli). Indeed, the framework developed here is sufficiently general to encompass, as special cases, both the exponentiation and composition actions for arbitrary finite groups, the Young permutation representations for groups of the form $F\wr \mathsf{Sym}(V)$, and several natural actions on subtrees and configurations of vertices in rooted trees.

Throughout the paper, we adopt a coordinate-free approach in order to avoid, as far as possible, the notational burden of managing multiple indices. A notable exception occurs in the parametrization of the orbits of stabilizer subgroups (such as in \eqref{decYa} and \eqref{decPsixij}), where numerical indices are deliberately introduced so that the corresponding spherical functions can be expressed as functions of discrete numerical parameters—a feature particularly useful in applications to Markov chains and special functions.

The paper is organized as follows.

In Section \ref{Secprelimin}, we review the necessary preliminaries on group actions, set compositions, wreath products, and their representation theory, fixing our coordinate-free notation for extensions, inflations, and inertia groups, and deriving explicit matrix-coefficient formulas for induced representations (Proposition \ref{PropactionH}).

In Section \ref{Sechomspace}, we introduce our general construction. Starting from a faithful transitive $G$-space $V$, a surjective map $\alpha\colon X\to A$, and a family of mutually non-isomorphic transitive $F$-spaces $\{Y_a\}_{a\in A}$, we define a natural action $\gamma$ of $F\wr G$ on a space $\Gamma$ of compatible pairs $(\xi,\varphi)$. We determine the $F\wr G$-orbits, their stabilizers, and necessary and sufficient conditions for two orbits to be isomorphic as $F\wr G$-spaces (Lemma \ref{Lemmaorb}). For a fixed orbit, we establish a general double-coset decomposition into orbits of the stabilizer (Theorem \ref{TheodecOrbxiphi}) and describe the action inductively via Clifford theory (Proposition \ref{propindact}). We then illustrate how this unified framework recovers and extends a wide range of classical actions, including exponentiation, composition, actions on $k$-permutations and subsets, and periodic actions of cyclic groups (Examples \ref{Example1}--\ref{Exampleequiv}).

In Section \ref{Secpermrep}, we analyze the permutation representation of $F\wr G$ on the space of complex-valued functions $L\bigl(\mathsf{Orb}_{F\wr G}(\xi_0,\varphi_{\xi_0})\bigr)$. Using intertwining operators and Clifford theory, we decompose this permutation representation into irreducible constituents and construct explicit bases for its isotypic components (Theorem \ref{bigTheorem} and Corollary \ref{Cormult}), showing that the multiplicity of each induced representation $\Ind_{F\wr\mathsf{I}(\sigma)}^{F\wr G}(\widetilde{\sigma}\otimes\overline{\theta})$ equals the multiplicity of $\theta$ in an auxiliary permutation representation $L(\mathsf{J}_\sigma)$ of the inertia subgroup $\mathsf{I}(\sigma)$.

In Section \ref{Secmultfreepermrep}, we characterize when the permutation representation of $F\wr G$ is multiplicity free (Proposition \ref{Propmultfree} and Corollary \ref{Cormultfree}). In the multiplicity-free case, we derive explicit formulas for the corresponding spherical functions (Theorem \ref{TheoFsigthet}), express them as constant functions on the double cosets in terms of numerical parameters (Theorem \ref{TheodefFsigthetorb}), establish their orthogonality relations \eqref{orthrel}, and provide a criterion for the associated Gelfand pair to be symmetric \eqref{symGel}.

In Section \ref{SecYoungpermrep}, we specialize our results to the case where $G=\mathsf{Sym}(V)$ is the full symmetric group. Here, orbits and isomorphism classes are governed by set compositions and partition-valued functions (Lemma \ref{Lemmaorbsym}), while the inertia subgroups $\mathsf{I}(\sigma)$ are Young subgroups. By applying the Young rule, we obtain an explicit combinatorial formula for the multiplicities in terms of semistandard Young tableaux and Kostka numbers (Theorem \ref{Bigformult}), and we completely classify all multiplicity-free permutation representations of $F\wr\mathsf{Sym}(V)$ (Theorem \ref{TheoSymmultfree}).

In Section \ref{Sectree}, we apply the results of Section \ref{SecYoungpermrep} to the automorphism group $\Aut(\mathcal{T}_{q,h})\cong \Aut(\mathcal{T}_{q,h-1})\wr\mathsf{Sym}(V)$ of a finite rooted $q$-ary tree of height $h$. Using the language of tree compositions and trees of partitions developed in \cite{Trees2}, we obtain the explicit irreducible decompositions for the permutation representations acting on complete centipede subtrees, and on $k$-insects of arbitrary type in $\mathcal{T}_{q,h}$. 

Finally, in Section \ref{Secgensymm}, inspired by Theorem \ref{TheoSymmultfree}, we construct a new family of finite Gelfand pairs for the generalized symmetric group $\mathsf{C}_n\wr\mathsf{Sym}(V)$, where $n=ab$ with $\gcd(a,b)=1$, arising from the action on two complementary cyclic spaces $\mathbb{Z}_b$ and $\mathbb{Z}_a$ indexed by a subset $Z_0\subseteq V$. We describe the double-coset decomposition (Theorem \ref{decOrbxiphib}), prove multiplicity freeness (Proposition \ref{Propmultfreebis}), and compute the spherical functions explicitly in terms of characters of $\mathsf{C}_n^V$ and Hahn polynomials associated with the Johnson scheme.

\section{Preliminaries}\label{Secprelimin}

In this section, we introduce some basic notions and notation that will be widely used throughout the paper. We refer to \cite{Cameron,Machi} for general introductions to group actions and to \cite{book4,Sternberg} for the representation theory of finite groups.

If $X$ is a finite set, we denote by $\mathsf{Sym}(X)$ the symmetric group over $X$, that is, the group of all permutations of $X$; $\mathsf{Sym}_n$ denotes the symmetric group on $n$ letters. If $V$ is another finite set, we denote by $X^V$ the family of all maps from $V$ to $X$; $\mathsf{Surj}(V,X)$ denotes the set of all surjective maps $\xi\in X^V$. Disjoint unions of sets will be denoted by the symbol $\bigsqcup$, direct products of sets or groups by $\bigtimes$, direct products of an indexed family of groups by $\prod$; $V^n$ will denote the direct product of $n$ copies of the finite set $V$, while $G^n$ denotes that of $n$ copies of the finite group $G$. 
An element $\xi\in\mathsf{Surj}(V,X)$ is equivalent to a set composition of $V$
\begin{equation}\label{Vxisetcomp}
V=\bigsqcup_{x\in X}\xi^{-1}(x),
\end{equation}
that is, $V$ is written as a disjoint union of nonempty subsets; cf. \cite[Example 2.9]{Trees2}. Note that the elements of $X$ are distinguishable (cf. \cite[Section 1.4]{Aigner}), which is usually indicated by indexing them, say $X=\left\{x_1,x_2,\dotsc,x_m\right\}$, so that \eqref{Vxisetcomp} may be considered as an ordered partition; cf. \cite[Section 1.8]{Sagan} and \cite[Section 1.9]{Stanley}.
We prefer to use the coordinate-free form in \eqref{Vxisetcomp}, which we call an $X$-\emph{composition} of $V$; alternatively, for typographical reasons, we will denote it by
\begin{equation}\label{Vxisetcomp2}
\left(\xi^{-1}(x)\mathrel{:}x\in X\right).  
\end{equation} 
We usually assume that $0\in\mathbb{N}$.

If $F$ is a finite group, then the set $F^V$ may be seen as the group product $F^{\lvert V\rvert}$ when endowed with the pointwise product: if $f,f'\in F^V$, then their product is denoted by $f\cdot f'$ and it is defined by setting $\left[f\cdot f'\right](v)\coloneqq f(v)f'(v)$ for every $v\in V$. If $A$ is a finite index set and $Y_a$, $a\in A$, is a family of finite sets, then the direct product $\bigtimes_{a\in A}Y_a$ may be seen as the set of all functions $\alpha\colon A\rightarrow\bigsqcup_{a\in A}Y_a$ such that $\alpha(a)\in Y_a$ for every $a\in A$. We denote by $\widehat{G}$ a fixed set of representatives for the equivalence classes of irreducible $G$-representations. We will use the symbol $\sim$ to denote the equivalence of group representations, the symbol $\cong$ to denote isomorphisms of groups, and $\simeq$ to denote isomorphisms of group actions. 

If a group $G$ acts on a set $V$ and $v\in V$, we usually denote by $gv$ the $g$-image of $v$, by $\mathsf{Orb}_G(v)$ the $G$-orbit containing $v$, and by $\mathsf{Stab}_G(v)$ the stabilizer of $v$ in $G$; if $K\leq G$ and $Z\subset V$, then $KZ\coloneqq\{kz\mathrel{:} k\in K, z\in Z\}$. Sometimes, an action is described as a homomorphism $\eta\colon G\rightarrow\mathsf{Sym}(V)$, so that $gv$ is replaced by $\eta(g)v$. With this notation, if $K\leq G$ and $X\subseteq V$, we say that $X$ is $\eta(K)$-\emph{invariant} if $\eta(k)x\in X$ for every $k\in K$ and $x\in X$; then $\left[\Res^G_K\eta\right]_X$ will denote the action on $X$ given by the restriction of $\eta$ to $K$.
If $\beta\colon G\rightarrow\mathsf{Sym}(W)$ is another action of $G$, then the \emph{diagonal action} of $G$ on $V\times W$ is defined by setting 
\begin{equation}\label{defdiagact}
\left[\mathsf{diag}(\eta,\beta)(g)\right](v,w)\coloneqq (\eta(g)v,\beta(g)w)
\end{equation} 
for all $g\in G$, $v\in V$, and $w\in W$.
If $A$ is a finite index set and we have an action $\eta_a$ of the group $F$ on the set $Y_a$ for every $a\in A$, then $\bigtimes_{a\in A}\eta_a$ denotes the product action of $F^A$ on $\bigtimes_{a\in A}Y_a$, that is if $f\in F^A$ and $\zeta\in\bigtimes_{a\in A}Y_a$, then
\begin{equation}\label{defdirprodact}
\left\{\left[\left(\bigtimes_{a\in A}\eta_a\right)f\right]\zeta\right\}(b)\coloneqq \eta_b(f(b))\zeta(b)
\end{equation}
for every $b\in A$. For this action, we will usually use the less cumbersome pointwise notation: $[f\cdot \zeta](b)\coloneqq f(b)\zeta(b)$.

If $G$ acts on $V$ and $X$ is a set, the \emph{natural action} of $G$ on the function space $X^V$ is given by \begin{equation}\label{defnatact}
g\xi\coloneqq \xi\circ g^{-1}\text{ for }\xi\in X^V, g\in G; \;\text{ that is, } [g\xi](v)=\xi\left(g^{-1}v\right)\text{ for every } v\in V. 
\end{equation}
An orbit $\mathsf{Orb}_G(\xi)$ is also called a \emph{pattern} (see, e.g., \cite{Aigner}). The natural action of $G$ on $\mathsf{Surj}(V,X)$ is isomorphic to the action on the $X$-compositions of $V$: if $\xi$ corresponds to \eqref{Vxisetcomp}, then $g\xi$ corresponds to $\bigsqcup_{x\in X}g\left(\xi^{-1}(x)\right)$. With our coordinate-free notation, the orbit of the $X$-composition \eqref{Vxisetcomp2} 
is 
\begin{equation}\label{OrbGpart}
\mathsf{Orb}_G\left(\xi^{-1}(x) \mathrel{:} x\in X\right)\coloneqq\left\{g\left(\xi^{-1}(x) \mathrel{:} x\in X\right)\mathrel{:}g\in G\right\},
\end{equation}
where $g\left(\xi^{-1}(x) \mathrel{:} x\in X\right)\coloneqq\bigl(g\left(\xi^{-1}(x)\right) \mathrel{:} x\in X\bigr)$ and $g\left(\xi^{-1}(x)\right)$ denotes the $g$-image of the subset $\xi^{-1}(x)$.

We denote by $L(V)$ the vector space of all complex-valued functions $h\colon V\rightarrow\mathbb{C}$, and we use $\lambda$ to denote the left action $\lambda(g)h=h\circ g^{-1}$, usually written in the form $[\lambda(g)h](v)=h\left(g^{-1}v\right)$ for $v\in V$. This is the permutation representation of $G$ on $L(V)$; we use the symbol $\lambda_G$ when we need to specify which group is acting. 

Let $V$ be a finite set, let $G$ be a finite group acting faithfully and transitively on $V$ (so that it may be identified with a subgroup of $\mathsf{Sym}(V)$), and let $F$ be another finite group.
By means of the natural action \eqref{defnatact}, $G$ becomes a group of automorphisms of $F^V$: 
\[
g(f\cdot f')=(f\cdot f')\circ g^{-1}=\left(f\circ g^{-1}\right)\cdot\left(f'\circ g^{-1}\right)=(g f)\cdot (g f')
\]
for all $f,f'\in F^V$, $g\in G$. We have $(g f)^{-1}=g\left(f^{-1}\right)$, which will be denoted simply by $g f^{-1}$.
The \emph{wreath product} of $F$ by $G$ is the set $F\wr G\coloneqq F^V\times G$ with the multiplication defined by 
\begin{equation}\label{multwreath}
(f,g)(f',g')\coloneqq (f\cdot gf',g g');
\end{equation}
see, e.g., \cite{book3}.
The identity of $F\wr G$ is $\left(\mathbbm{1}_F,1_G\right)$, where $\mathbbm{1}_F(v)=1_F$ for all $v\in V$ and $1_G$ is the identity of $G$; the inverse of an element is given by the formula $(f,g)^{-1}=\left(g^{-1}f^{-1},g^{-1}\right)$. It is worth noting that this definition of the wreath product does not require the action to be transitive on $V$; indeed, when dealing with stabilizer subgroups, this assumption is generally dropped.

Suppose that $K\leq G$ and that $\eta\colon K\rightarrow \mathsf{Sym}(X)$ is an action of $K$ on a set $X$. An \emph{extension} of $\eta$ to $G$ is an action $\widetilde{\eta}\colon G\rightarrow\mathsf{Sym}(X)$ whose restriction to $K$ coincides with $\eta$: in symbols, $\widetilde{\eta}(k)=\eta(k)$ for all $k\in K$. This gives rise to an extension of the corresponding permutation representation, that is $\Res^G_K\lambda_G\sim \lambda_K$. On the other hand, suppose that $N\triangleleft G$, denote by $\pi\colon G\rightarrow G/N$ the quotient map and suppose that $\beta\colon G/N\rightarrow \mathsf{Sym}(X)$ is an action of $G/N$ on a set $X$. The \emph{inflation} of $\beta$ to $G$ is simply the composition 
\begin{equation}\label{definflperrep}
\overline{\beta}\coloneqq\beta\circ\pi\colon G\xrightarrow{\pi}G/N\xrightarrow{\beta}\mathsf{Sym}(X).
\end{equation}
This gives rise to the inflation of the corresponding permutation representations: $\overline{\lambda}_{G/N}=\lambda_G$. 

Next, suppose that $K\leq G$ and that $K$ acts transitively on a set $X$. If $H$ is the stabilizer in $K$ of a point of $X$, then we may identify $X$ with $K/H$, endowed with the left action. Under these assumptions, we say that the left action of $G$ on $V\coloneqq G/H$ is \emph{induced} by the action of $K$ on $X$. If $\mathsf{R}$ is a set of representatives for the left $K$-cosets in $G$ (with $1_G\in\mathsf{R}$), then 
\begin{equation}\label{defindact}
V=\bigsqcup_{r\in\mathsf{R}}rX. 
\end{equation}
Moreover, for every choice of $g\in G$ and $r\in\mathsf{R}$ there exists a unique $r'\in\mathsf{R}$ and a unique $k\in K$ such that $gr=r'k$, which implies that $grx=r'kx$, for every $x\in X$; this fact is used in \cite[Definition 3.11]{Machi} to define the induced action, without the assumption of transitivity, by identifying $V$ with $\mathsf{R}\times X$. Clearly, if we denote by $L(X)$ the subspace of $L(V)$ of all complex valued functions supported in $X$, then from \eqref{defindact} it follows that 
\begin{equation}\label{indpermrep}
L(V)=\bigoplus_{r\in\mathsf{R}}\lambda(r)L(X),
\end{equation}
where the subspace $\lambda\left(1_G\right)L(X)$ is $K$-invariant and equivalent to $L(X)$ as a $K$-representation;
this is a concrete realization of the fact that $L(V)\sim\Ind_K^GL(X)$. Under the notation above, we also have a \emph{natural projection} $\pi\colon G/H\rightarrow G/K$, that is, a surjective equivariant map defined by setting $\pi(gH)=gK$, for every $g\in G$.

\begin{example}
A simple yet remarkable example of a natural projection is the following. Let $V$ be a finite set, fix $1\leq k\leq\lvert V\rvert$, and denote by $\mathsf{Per}_k(V)$ the family of all $k$-permutations of distinct elements of $V$ (see, e.g., \cite[Section 1.2]{Sagan}), that is all ordered sequences of $k$ distinct elements $\left(v_1,v_2,\dotsc,v_k\right)$ belonging to $V$. Let $\binom{V}{k}$ denote the set of all $k$-element subsets of $V$. Then the following map is a natural projection with respect to the action of $\mathsf{Sym}(V)$:
\begin{equation}\label{natproj}
\begin{array}{rlcc}
\pi\colon&\mathsf{Per}_k(V)&\longrightarrow&\binom{V}{k}\\
&\left(v_1,v_2,\dotsc,v_k\right)&\longmapsto&\left\{v_1,v_2,\dotsc,v_k\right\}.
\end{array}
\end{equation}
\end{example} 

Now suppose that $N\triangleleft G$ and that $\eta\colon N\rightarrow\mathsf{Sym}(Y)$ is a transitive action of $N$ on a set $Y$. Then, for every $g\in G$, the $g$-\emph{conjugate} of $\eta$ is the action $\!^g\eta\colon N\rightarrow\mathsf{Sym}(Y)$ defined by setting $\!^g\eta(n)\coloneqq\eta\left(g^{-1}ng\right)$ for every $n\in N$. The \emph{inertia group} $\mathsf{I}(\eta)$ of $\eta$ in $G$ is the group of all $g\in G$ such that $\!^g\eta\simeq\eta$. Clearly, $N\leq\mathsf{I}(\eta)$; if $\mathsf{Stab}_N\left(y_0\right)$ is the stabilizer of $y_0\in Y$ with respect to $\eta$, then the stabilizer of $y_0$ with respect to $\!^g\eta$ is $g\mathsf{Stab}_N\left(y_0\right)g^{-1}$, so that $g\in\mathsf{I}(\eta)$ if and only if the two stabilizers are \emph{conjugate in} $N$. Moreover, if $K$ is an intermediate subgroup, that is, $N\leq K\leq G$, and there exists an extension $\widetilde{\eta}$ of $\eta$ to $K$, then $K\leq\mathsf{I}(\eta)$. Indeed, if $k\in K$ and $n\in N$, then
\[
\!^k\eta(n)=\eta\left(k^{-1}nk\right)=\widetilde{\eta}\left(k^{-1}\right)\eta(n)\widetilde{\eta}(k).
\]
\begin{remark}
\begin{enumerate}[(1)]
\item\label{Remexta}
While induction and inflation always exist and are uniquely determined, the extension of a group action does not necessarily exist and, if it does, it is not necessarily unique. For instance, consider a simple non abelian group $G$ with an element $g$ of order two, such as $A_5$ with $g$ being the product of two disjoint transpositions. Then we may consider the Cayley action of the subgroup $C_2\coloneqq\langle g\rangle$; this action cannot be extended to $G$ because it would lead to a stabilizer which is a subgroup of index two of $G$. On the other hand, suppose that $G=K\times K$ with $\lvert K\rvert > 1$, and consider the left Cayley action of the subgroup $K\times\left\{1_K\right\}$. We may consider the action of $K\times K$ on $K\times\left\{1_K\right\}$ given by $\left(k_1,k_2\right)\left(k,1_K\right)\coloneqq\left(k_1kk_2^{-1},1_K\right)$ for all $k,k_1,k_2\in K$, but also another action $\eta$ given by: $\left(k_1,k_2\right)\left(k,1_K\right)\coloneqq\left(k_1k,1_K\right)$. Both of them are extensions of the Cayley action of $K\times\left\{1_K\right\}$ but they are not isomorphic because the stabilizers of $\left(1_K,1_K\right)$, namely $\{(k,k)\mathrel{:}k\in K\}$ and $\left\{1_K\right\}\times K$, are not conjugate in $K\times K$.
\item
Consider the semidirect product $(K\times K)\rtimes C_2$, where $C_2$ acts on  $K\times K$ by switching the coordinates (that  is, the wreath product $K\wr C_2$). Then $K\times\left\{1_K\right\}$ and $\left\{1_K\right\}\times K$ are conjugate in $(K\times K)\rtimes C_2$ but not in $K\times K$. Therefore, if $\eta$ is as in \eqref{Remexta} then its inertia group in $(K\times K)\rtimes C_2$ coincides with $K\times K$.
\end{enumerate}
\end{remark}

Following \cite{book3}, we now recall some basic facts on the representation theory of wreath products; other standard references are \cite{Hu,JK}.
Every irreducible representation of $F^V$ is of the form $\sigma\coloneqq\bigotimes_{v\in V}\sigma_v$, where  $ v\mapsto\sigma_v$ varies over all possible maps $V\rightarrow\widehat{F}$; clearly, $\sigma(f)=\bigotimes_{v\in V}\sigma_v(f(v))$ and distinct maps give rise to inequivalent representations. In other words, $\widehat{F^V}$ can be naturally identified with $\widehat{F}^V$. If $(f,g)\in F\wr G$, then the $(f,g)$-\emph{conjugate} of $\sigma$ is given by
\begin{equation}\label{fgsigmaconj}
\!^{(f,g)}\sigma\coloneqq\bigotimes_{v\in V} \;^{f(v)}\sigma_{g^{-1}v}\sim\bigotimes_{v\in V}\sigma_{g^{-1}v},
\end{equation}
where $\left[\,^{f(v)}\sigma_{g^{-1}v}\right]\left(f_0(v)\right)=\sigma_{g^{-1}v}\bigl(f(v)^{-1}\cdot f_0(v)\cdot f(v)\bigr)$ for all $f_0\in F^V$. This yields a natural action of $F\wr G$ on $\widehat{F^V}$ if we identify equivalent representations; in particular, we have $\!^{(f,g)}\sigma\sim\bigotimes_{v\in V} \sigma_{g^{-1}v}$ because $\!^{f(v)}\sigma_{g^{-1}v}\sim \sigma_{g^{-1}v}$ for each $v\in V$.
Therefore, $F^V$ is the kernel of the conjugation action, while $G$ permutes the factors of the tensor products. The stabilizer of $\sigma$ is called \emph{the inertia group} and is given by $F\wr \mathsf{I}(\sigma)$, where $\mathsf{I}(\sigma)\coloneqq\left\{k\in G\colon\sigma_{k^{-1}v}=\sigma_v\text{ for all }v\in V\right\}$; we also set 
\begin{equation}\label{defgsigma}
\!^g\sigma\coloneqq\bigotimes_{v\in V}\sigma_{g^{-1}v}.
\end{equation}
Clearly, the action of $G$ given in \eqref{defgsigma} is isomorphic to the natural action on $\widehat{F}^V$. 

A generic element of $\widehat{F}$ is denoted by $\sigma_*$, whereas $\sigma_0$ is reserved for the trivial representation of $F$. If $\sigma\in\widehat{F^V}$, then $\mathsf{image}(\sigma)$ is the set of all $\sigma_*\in\widehat{F}$ such that $\sigma_v=\sigma_*$ for at least one $v\in V$, while $\mathsf{supp}(\sigma)\coloneqq\{v\in V\colon \sigma_v\neq\sigma_0\}$.
The trivial representation of $F^V$ is denoted by $\sigma^0$; in our notation, $\mathsf{supp}\left(\sigma^0\right)=\emptyset$ and $\mathsf{image}\left(\sigma^0\right)=\left\{\sigma_0\right\}$. Note also that every $\sigma\in\widehat{F^V}$ determines a set composition 
\begin{equation}\label{VsIsigma}
V=\bigsqcup_{\sigma_*\in\mathsf{image}(\sigma)}V_{\sigma_*},\qquad\text{where }\quad V_{\sigma_*}\coloneqq\{v\in V\mathrel{:}\sigma_v=\sigma_*\},
\end{equation}
and that $\mathsf{I}(\sigma)=\left\{k\in G\mathrel{:}kV_{\sigma_*}=V_{\sigma_*}\text{ for every } \sigma_*\in\mathsf{image}(\sigma)\right\}$. Moreover, if $V_{\sigma_0}\neq\emptyset$ and $\mathsf{supp}(\sigma)\neq V$, then $V=V_{\sigma_0}\sqcup\mathsf{supp}(\sigma)$.

We denote by $U_{\sigma_*}$ the representation space of $\sigma_*\in\widehat{F}$, so that
\begin{equation}\label{defUsigma}
U_\sigma=\bigotimes_{v\in V}U_{\sigma_v}
\end{equation}
is the representation space of $\bigotimes_{v\in V}\sigma_v$. Then, for all $\bigotimes_{v\in V}u_v\in U_\sigma$, we set 
\[
\lambda\left(\mathbbm{1}_F,g\right)\left(\bigotimes_{v\in V}u_v\right)\coloneqq\bigotimes_{v\in V}u_{g^{-1}v}\in U_{\,^g\sigma},
\]
so that 
\begin{equation}\label{actgsig}
\!^g\sigma(f)\lambda\left(\mathbbm{1}_F,g\right)\left(\bigotimes_{v\in V}u_v\right)=\bigotimes_{v\in V}\sigma_{g^{-1}v}(f(v))u_{g^{-1}v}.
\end{equation}
Moreover, every $\sigma\in\widehat{F^V}$ has an \emph{extension} $\widetilde{\sigma}$ to its inertia group $F\wr \mathsf{I}(\sigma)$, that is, $\widetilde{\sigma}\in\widehat{F\wr \mathsf{I}(\sigma)}$ and $\Res^{F\wr \mathsf{I}\left(\sigma\right)}_{F^V}\widetilde{\sigma}=\sigma$. It is given by setting
\begin{equation}\label{defextsig}
\widetilde{\sigma}(f,k)\coloneqq\;^k\sigma(f)\lambda\left(\mathbbm{1}_F,k\right)\equiv\sigma(f)\lambda\left(\mathbbm{1}_F,k\right)
\end{equation}
for all $(f,k)\in F\wr \mathsf{I}(\sigma)$, which means that
\begin{equation}\label{defextsig2}
\widetilde{\sigma}(f,k)\left(\bigotimes_{v\in V}u_v\right)=\bigotimes_{v\in V}\sigma_v(f(v))u_{k^{-1}v}
\end{equation}
for all $\bigotimes_{v\in V}u_v\in U_\sigma$. On the other hand, if $\theta\in\widehat{\mathsf{I}(\sigma)}$, then its \emph{inflation} $\overline{\theta}$ to $F\wr \mathsf{I}(\sigma)$ is given by: $\overline{\theta}(f,k)\coloneqq\theta(k)$ for all $(f,k)\in F\wr \mathsf{I}(\sigma)$. Note also that the extension is not unique, because it is defined up to the tensor product with the inflation of a one-dimensional representation of $\mathsf{I}(\sigma)$; however, we will always use \eqref{defextsig}. Finally, let $\Sigma$ be a set of representatives for the $F\wr G$-conjugacy classes of $\widehat{F^V}$, that is for the orbits of $G$ on $\widehat{F^V}$. Then
\begin{equation}\label{irrepwreath}
\widehat{F\wr G}=\left\{\Ind^{F\wr G}_{F\wr \mathsf{I}(\sigma)}\left(\widetilde{\sigma}\otimes\overline{\theta}\right)\mathrel{:}\sigma\in\Sigma,\theta\in\widehat{\mathsf{I}(\sigma)}\right\}.
\end{equation}
The group $F\wr G$ acts on these induced representations by means of left translation, denoted by $\lambda$.

We need to establish more explicit formulas for the irreducible representations in \eqref{irrepwreath}.
Suppose that $\mathsf{Q}_\sigma$ is a set of representatives for the left cosets of $\mathsf{I}(\sigma)$ in $G$, that is, $G=\bigsqcup_{q\in\mathsf{Q}_\sigma}q\mathsf{I}(\sigma)$, and assume that $1_G\in\mathsf{Q}_\sigma$. Then we also have 
\begin{equation}\label{defQsigma}
F\wr G=\bigsqcup_{q\in\mathsf{Q}_\sigma}\left(\mathbbm{1}_F,q\right)\left(F\wr \mathsf{I}(\sigma)\right).
\end{equation}
Indeed, if $(f,g)\in F\wr G$ and $g\in q\mathsf{I}(\sigma)$, then $g=qk$ with $k\in\mathsf{I}(\sigma)$ 
and $(f,g)=(f,qk)=\left(\mathbbm{1}_F,q\right)\left(q^{-1}f,k\right)$. Next, we denote by $W_\theta$ the space of $\theta$. Henceforth, for all $q\in\mathsf{Q}_\sigma$ and $u\otimes w\in U_\sigma\otimes W_\theta$, we define $H_{q,u\otimes w}\colon F\wr G\rightarrow U_\sigma\otimes W_\theta$ by setting, for all $(f,g)\in  F\wr G$,
\begin{equation}\label{defHquw}
H_{q,u\otimes w}(f,g)\coloneqq
\begin{cases}
0&\text{ if }g\notin q\mathsf{I}(\sigma),\\
\widetilde{\sigma}\left(g^{-1}f^{-1},k^{-1}\right)u\otimes\theta\left(k^{-1}\right)w&\text{ if }g=qk, k\in \mathsf{I}(\sigma).
\end{cases}
\end{equation}
\begin{proposition}\label{PropactionH}
\begin{enumerate}[(1)]
\item\label{PropactionH1}
The set of all $H_{q,u\otimes w}$ spans the space of $\Ind^{F\wr G}_{F\wr \mathsf{I}(\sigma)}\left(\widetilde{\sigma}\otimes\overline{\theta}\right)$. 
\item
The action of $F\wr G$ on $H_{q,u\otimes w}$ is the following: if $(f,g)\in F\wr G$ and $gq=q'k$, with $q'\in\mathsf{Q}_\sigma$ and $k\in \mathsf{I}(\sigma)$, then
\begin{equation}\label{actionH}
\lambda(f,g)H_{q,u\otimes w}=H_{q',\widetilde{\sigma}\left(q'^{-1}f,k\right)u\otimes\theta(k)w}.
\end{equation}
\end{enumerate}
\end{proposition}
\begin{proof}
(1) If $(f',k')\in F\wr\mathsf{I}(\sigma)$ and $g=qk$, with $k\in \mathsf{I}(\sigma)$, then
\[
\begin{split}
H_{q,u\otimes w}\left((f,g)(f',k')\right)=&H_{q,u\otimes w}\left(f\cdot gf',gk'\right)\\
(\text{by }\eqref{defHquw})\quad=&\widetilde{\sigma}\Bigl(k'^{-1}k^{-1}q^{-1}\left(gf'^{-1}\cdot f^{-1}\right),k'^{-1}k^{-1}\Bigr)u\otimes\theta\left(k'^{-1}k^{-1}\right)w\\
=&\widetilde{\sigma}\Bigl(k'^{-1}f'^{-1}\cdot k'^{-1}g^{-1}f^{-1},k'^{-1}k^{-1}\Bigr)u\otimes\theta\left(k'^{-1}k^{-1}\right)w\\
=&\widetilde{\sigma}\Bigl((f',k')^{-1}\left(g^{-1}f^{-1},k^{-1}\right)\Bigr)u\otimes\theta\left(k'^{-1}k^{-1}\right)w\\
=&\left[\widetilde{\sigma}\left((f',k')^{-1}\right)\otimes\theta\left(k'^{-1}\right)\right]\left(\widetilde{\sigma}\left(g^{-1}f^{-1},k^{-1}\right)u\otimes\theta\left(k^{-1}\right)w\right)\\
(\text{by }\eqref{defHquw})\quad=&\left[\left(\widetilde{\sigma}\otimes\theta\right)\left((f',k')^{-1}\right)\right]H_{q,u\otimes w}(f,g).
\end{split}
\]
This shows that $H_{q,u\otimes w}$ belongs to the space of $\Ind^{F\wr G}_{F\wr \mathsf{I}(\sigma)}\left(\widetilde{\sigma}\otimes\overline{\theta}\right)$.

On the other hand, if $\psi$ belongs to the space of $\Ind^{F\wr G}_{F\wr \mathsf{I}(\sigma)}\left(\widetilde{\sigma}\otimes\overline{\theta}\right)$, it is zero outside $\left(\mathbbm{1}_F,q\right)\left(F\wr \mathsf{I}(\sigma)\right)$ and satisfies $\psi\left(\mathbbm{1}_F,q\right)=u\otimes w$, then, by the definition of induced representations, we have
\[
\begin{split}
\psi(f,qk)=&\psi\left(\left(\mathbbm{1}_F,q\right)\left(q^{-1}f,k\right)\right)=\left[\left(\widetilde{\sigma}\otimes\overline{\theta}\right)\left(k^{-1}q^{-1}f^{-1},k^{-1}\right)\right]\psi\left(\mathbbm{1}_F,q\right)\\
(g=qk)\quad=&\widetilde{\sigma}\left(g^{-1}f^{-1},k^{-1}\right)u\otimes\theta\left(k^{-1}\right)w,
\end{split}
\]
so that $\psi=H_{q,u\otimes w}$. The assertion then follows from standard properties of induced representations; e.g. see \cite[Proposition 11.1.2]{book4}.

\medskip
(2) Indeed, 
\[
\begin{split}
\left[\lambda(f,g)H_{q,u\otimes w}\right]\left(\mathbbm{1}_F,q'\right)=&H_{q,u\otimes w}\left(g^{-1}f^{-1},g^{-1}q'\right)\\
(g^{-1}q'=qk^{-1})\qquad=&H_{q,u\otimes w}\left(g^{-1}f^{-1},qk^{-1}\right)\\
(\text{by }\eqref{defHquw})\qquad=&\widetilde{\sigma}\left(kq^{-1}g^{-1}f,k\right)u\otimes\theta(k)w\\
(kq^{-1}g^{-1}=q'^{-1})\qquad=&\widetilde{\sigma}\left(q'^{-1}f,k\right)u\otimes\theta(k)w\\
(\text{by }\eqref{defHquw}\text{ with } f \text{ replaced by }\mathbbm{1}_F\text{ and }g \text{ by }q' )\qquad=&H_{q',\widetilde{\sigma}\left(q'^{-1}f,k\right)u\otimes\theta(k)w}\left(\mathbbm{1}_F,q'\right).
\end{split}
\]
\end{proof}

Finally, note that, in the notation of \eqref{OrbGpart}, \eqref{VsIsigma} and \eqref{defQsigma}, the set 
$\mathsf{Orb}_G\left(V_{\sigma_*}\mathrel{:}\sigma_*\in\mathsf{image}(\sigma)\right)$
may be parametrized by $\mathsf{Q}_\sigma$: 
\begin{equation}\label{OrbV0Vt}
\mathsf{Orb}_G\left(V_{\sigma_*}\mathrel{:}\sigma_*\in\mathsf{image}(\sigma)\right)=\bigsqcup_{q\in\mathsf{Q}_\sigma}\left\{q\left(V_{\sigma_*}\mathrel{:}\sigma_*\in\mathsf{image}(\sigma)\right)\right\}.
\end{equation}

\section{Group actions and homogeneous spaces}\label{Sechomspace}

In this section, we introduce the basic construction, namely a family of $F\wr G$-spaces.

\subsubsection*{Basic Assumptions} 
We suppose that 
\begin{enumerate}[(a)]
\item $G$ is a finite group acting faithfully and transitively on a finite set $V$;
\item $X$ and $A$ are finite sets and $1\leq\lvert A\rvert\leq\lvert X\rvert\leq\lvert V\rvert$;
\item  with each $a\in A$, we associate a finite non empty set $Y_a$ equipped with a transitive $F$-action $\eta_a$;
\item the sets $Y_a$, $a\in A$, are mutually disjoint;
\item if $a,a'\in A$, $a\neq a'$, then the $F$-actions $\eta_a$ and $\eta_{a'}$ \emph{are not isomorphic}; 
\item  $\alpha\colon X\rightarrow A$ is a \emph{fixed} surjective map.
\end{enumerate} 

We define $\pi\colon\bigsqcup_{a\in A}Y_a\rightarrow A$ by setting $\pi(y)\coloneqq a$ if $y\in Y_a$.
If $\xi\in\mathsf{Surj}(V,X)$, we denote by $\Phi_\xi$ the set of all maps $\varphi\colon V\rightarrow \bigsqcup_{a\in A}Y_a$ such that $\varphi(v)\in Y_{\alpha(\xi(v))}$ for all $v\in V$; this is equivalent to saying that $\pi\circ\varphi=\alpha\circ\xi$, meaning that the following diagram commutes:

\begin{equation}\label{Diagram}
\xymatrix{
V\ar[rr]^{\xi}\ar[rrrd]_{\varphi}&&X\ar[r]^{\alpha}&A\\
&&&\bigsqcup_{a\in A}Y_a\ar[u]_{\pi}\\
}
\end{equation}

\noindent
In other words, $V=\bigsqcup_{a\in A}(\alpha\circ\xi)^{-1}(a)$ is an $A$-composition of $V$, and $\varphi\in\Phi_\xi$ if and only if $\varphi\bigl((\alpha\circ\xi)^{-1}(a)\bigr)\subseteq Y_a$ for every $a\in A$.
Then we set 
\[
\Gamma\coloneqq\left\{(\xi,\varphi)\mathrel{:}\xi\in\mathsf{Surj}(V,X)\text{ and }\varphi\in\Phi_\xi\right\}.
\]

In addition to the natural actions of $G$ on the function spaces $\mathsf{Surj}(V,X)$ and $\Phi_\xi$ (cf. \eqref{defnatact}), the group $F^V$ also acts on each $\Phi_\xi$ in a pointwise way: 
if $\varphi\in\Phi_\xi$ and $f\in F^V$, then we define $f\cdot\varphi\in\Phi_\xi$ by setting $[f\cdot\varphi](v)\coloneqq f(v)\varphi(v)$ for all $v\in V$, where $f(v)\varphi(v)$ indicates the $\eta_{\alpha(\xi(v))}$-action of $f(v)\in F$ on $\varphi(v)\in Y_{\alpha(\xi(v))}$. For every $a\in A$ we fix a point $y_a\in Y_a$ and for all $\xi\in\mathsf{Surj}(V,X)$ we define $\varphi_\xi\in\Phi_\xi$ by setting 
\[
\varphi_\xi(v)\coloneqq y_{\alpha(\xi(v))}
\]
for each $v\in V$. Then $\Phi_\xi=\mathsf{Orb}_{F^V}\left(\varphi_\xi\right)$ and we have the following isomorphism of $F^V$-spaces:
\begin{equation}\label{FVprod}
\begin{array}{cccl}
\mathsf{Orb}_{F^V}\left(\varphi_\xi\right)&\longrightarrow&\bigtimes\limits_{v\in V}Y_{\alpha(\xi(v))}&\simeq\bigtimes\limits_{a\in A}Y_a^{\left\lvert(\alpha\circ\xi)^{-1}(a)\right\rvert}\\
\varphi&\longmapsto&\left(\varphi(v)\right)_{v\in V}.&
\end{array}
\end{equation}
That is, the action of $F^V$ is equivalent to the direct product action $\eta\coloneqq\bigtimes_{v\in V}\eta_{\alpha(\xi(v))}$; cf. \eqref{defdirprodact}.

Then we define an action $\gamma$ of $F\wr G$ on $\Gamma$ by setting, for $(\xi,\varphi)\in\Gamma$ and $(f,g)\in F\wr G$, 
\begin{equation}\label{defactionOm} 
\gamma(f,g)(\xi,\varphi)\coloneqq(g \xi,f\cdot g\varphi).
\end{equation}
We also set $\Phi\coloneqq\bigcup_{\xi\in\mathsf{Surj}(V,X)}\Phi_\xi$ and denote by $\gamma_0$ the action of $F\wr G$ on $\Phi$ given by the second component of $\gamma$: $\gamma_0(f,g)\varphi\coloneqq f\cdot g\varphi$.
The symbol $\cdot$ is used both in \eqref{multwreath} and in \eqref{defactionOm}, but it is easy to distinguish its specific meaning in each of its applications. Note also that we the projection $(\xi, \varphi) \mapsto \varphi$ from $\mathrm{Orb}_{F\wr G}(\xi_0, \varphi_{\xi_0})$ onto $\mathrm{Orb}_{F\wr G}(\varphi_{\xi_0})$ is equivariant and that $\gamma_0$ recovers the same action $\gamma$ precisely when $\alpha \circ \xi$ determines $\xi$ on the orbit (in particular, when $\alpha$ is bijective, see Example \ref{Example5} below), whereas it collapses fibers when $\alpha$ is not injective.

When no confusion can arise, we will omit the symbols $\eta_a$, $\eta$, $\gamma$, and
 $\gamma_0$. Unless otherwise specified, the notations $\mathsf{Orb}$ and $\mathsf{Stab}$ will refer to these actions.

\begin{lemma}
Both $\gamma$ and $\gamma_0$ are actions of $F\wr G$. 
\end{lemma}
\begin{proof}
First, observe that $(f\cdot g\varphi)(v)=f(v)\varphi\left(g^{-1}v\right)\in Y_{\alpha(g\xi(v))}$, so that 
$(g \xi,f\cdot g\varphi)\in\Gamma$. Furthermore, we have $\left(\mathbbm{1}_F,1_G\right)(\xi,\varphi)=(\xi,\varphi)$. Finally, if $(f',g')\in F\wr G$, then 
\[
\begin{split}
[(f,g)(f',g')](\xi,\varphi)=&(f\cdot gf',gg')(\xi,\varphi)=\bigl((gg')\xi,(f\cdot gf')\cdot(gg')\varphi\bigr)\\
=&\bigl(g(g'\xi),f\cdot \bigl(gf' \cdot g(g'\varphi)\bigr)\bigr)=\bigl(g(g'\xi),f\cdot g(f'\cdot g'\varphi)\bigr)\\
=&(f,g)\bigl(g'\xi, f'\cdot g'\varphi\bigr)=(f,g)\bigl((f',g')(\xi,\varphi)\bigr).
\end{split}
\]
\end{proof}

The decomposition of $\Gamma$ into $F^V$-orbits is
\begin{equation}\label{decGamma}
\Gamma=\bigsqcup_{\xi\in\mathsf{Surj}(V,X)}\bigl(\{\xi\}\times\mathsf{Orb}_{F^V}\left(\varphi_\xi\right)\bigr),
\end{equation}
and 
\[
\mathsf{Stab}_{F^V}\left(\varphi_\xi\right)=\prod_{v\in V}\mathsf{Stab}_F\left(y_{\alpha(\xi(v))}\right)\cong \prod_{a\in A}
\left[\mathsf{Stab}_F\left(y_a\right)\right]^{\left\lvert(\alpha\circ\xi)^{-1}(a)\right\rvert}.
\]
Clearly, $\mathsf{Stab}_G(\xi)=\bigcap_{x\in X}\mathsf{Stab}_G\left(\xi^{-1}(x)\right)$.
\begin{lemma}\label{Lemmastab}
The stabilizer of $\left(\xi,\varphi_{\xi}\right)$ in $F\wr G$ is given by
\begin{equation}\label{Lemmastab1}
\mathsf{Stab}_{F\wr G}\left(\xi,\varphi_{\xi}\right)=\left\{(f,g)\in F\wr G\mathrel{:} g\in\mathsf{Stab}_G(\xi)\; \text{\rm and } f\in \mathsf{Stab}_{F^V}\left(\varphi_\xi\right)\right\}.
\end{equation}
\end{lemma}
\begin{proof}
If $(f,g)\in F\wr G$, then 
\[
\begin{split}
(f,g)\left(\xi,\varphi_{\xi}\right)=\left(\xi,\varphi_{\xi}\right)&\Longleftrightarrow g\xi=\xi\text{ and }
f(v)\varphi_{\xi}\left(g^{-1}v\right)=\varphi_{\xi}(v)\\
&\Longleftrightarrow  g\in \mathsf{Stab}_G(\xi)\text{ and } f(v)y_{\alpha(\xi(v))}=y_{\alpha(\xi(v))}\text{ for every } v\in V.
\end{split}
\]
\end{proof}

The restriction of \eqref{decGamma} to an $F\wr G$-orbit yields the following decomposition: 
\begin{equation}\label{FVorbdec}
\mathsf{Orb}_{F\wr G}\left(\xi,\varphi_{\xi}\right)=\bigsqcup_{\xi'\in\mathsf{Orb}_G(\xi)}\bigl(\{\xi'\}\times\mathsf{Orb}_{F^V}\left(\varphi_{\xi'}\right)\bigr).
\end{equation}

We set $\Omega\coloneqq\prod_{a\in A}\mathsf{Sym}\left(\alpha^{-1}(a)\right)$; that is, $\Omega$ consists of all $\omega\in\mathsf{Sym}(X)$ such that $\alpha\circ\omega=\alpha$. The group $\Omega$ acts on $\Gamma$: if $(\xi,\varphi)\in\Gamma$ and $\omega\in\Omega$, then we set $\omega(\xi,\varphi)\coloneqq(\omega\circ\xi,\varphi)$.  Since $\alpha\circ\omega=\alpha$, it follows that $(\omega\circ\xi,\varphi)\in\Gamma$. We say that $\xi,\xi'\in\mathsf{Surj}(V,X)$ are \emph{equivalent} if there exist $g\in G$ and $\omega\in\Omega$ such that $\omega\circ\xi\circ g^{-1}=\xi'$, meaning that $\xi$ and $\xi'$ belong to the same $\Omega\times G$-orbit.

\begin{lemma}\label{Lemmaorb}
Suppose that $(\xi,\varphi),(\xi',\varphi')\in\Gamma$.
\begin{enumerate}[(1)]
\item\label{Lemmaorb1} We have $\mathsf{Orb}_{F\wr G}(\xi,\varphi)=\mathsf{Orb}_{F\wr G}(\xi',\varphi')$ if and only if $\mathsf{Orb}_G(\xi)=\mathsf{Orb}_G(\xi')$. Moreover, as a homogeneous space, $\mathsf{Orb}_{F\wr G}(\xi,\varphi)\simeq (F\wr G)/\mathsf{Stab}_{F\wr G}\left(\xi,\varphi_{\xi}\right)$.
\item\label{Lemmaorb2} If there exists $a_1\in A$ such that the action of $F$ on $Y_{a_1}$ is faithful and $\left\lvert Y_{a_1}\right\rvert\geq2$, then the action of $F\wr G$ on each orbit on $\Gamma$ is faithful. 
\item\label{Lemmaorb3}  The orbits $\mathsf{Orb}_{F\wr G}(\xi,\varphi)$ and $\mathsf{Orb}_{F\wr G}(\xi',\varphi')$ are isomorphic as $F\wr G$-spaces if and only if 
there exists $g\in G$ such that $g\mathsf{Stab}_G(\xi)g^{-1}=\mathsf{Stab}_G(\xi')$ and $\alpha\circ\xi\circ g^{-1}=\alpha\circ\xi'$. In particular, these conditions are verified if $\xi$ and $\xi'$ are equivalent.
\end{enumerate}
\end{lemma}

\begin{proof}
(1) The only if part is obvious from \eqref{defactionOm}. Now suppose that $\mathsf{Orb}_G(\xi)=\mathsf{Orb}_G(\xi')$. Then there exists $g\in G$ such that $g\xi=\xi'$. Next, from the transitivity of $F$ on $ Y_{\alpha(g\xi(v))}\equiv  Y_{\alpha(\xi'(v))}$, it follows that for every $v\in V$ there exists $f(v)\in F$ such that $f(v)\varphi\left(g^{-1}v\right)=\varphi'(v)$. Then \eqref{defactionOm} yields $(f,g)(\xi,\varphi)=(\xi',\varphi')$ and the if part is proved. In particular, $\left(\xi,\varphi_{\xi}\right)\in\mathsf{Orb}_{F\wr G}(\xi,\varphi)$, and the second statement follows.

\medskip
(2) Suppose that $Y_{a_1}$ is a non trivial faithful $F$-space. Fix $\xi'\in\mathsf{Surj}(V,X)$ and for every $\xi\in\mathsf{Orb}_G(\xi')$ denote by $\Phi^0_\xi$ the set of all $\varphi\in\Phi_\xi$ which are constant on each set $\xi^{-1}\left(\alpha^{-1}(a)\right)$ for $a\in A$. That is, a function $\varphi\in\Phi^0_\xi$ is defined similarly to $\varphi_{\xi}$, except that each $y_a$ is replaced with an arbitrary $y\in Y_a$. Then the kernel of the action of $F\wr G$ on $\mathsf{Orb}_{F\wr G}\left(\xi',\varphi_{\xi'}\right)$ is contained in
\[
\bigcap_{\xi\in\mathsf{Orb}_G(\xi'),\varphi\in\Phi^0_\xi}\mathsf{Stab}_{F\wr G}(\xi,\varphi).
\]
Arguing similarly to the proof of Lemma \ref{Lemmastab}, one can prove that if $(f,g)$ is in this intersection and $v\in V$, then $f(v)y=y$ for all $\xi\in\mathsf{Orb}_G(\xi')$ and all $y\in Y_{\alpha(\xi(v))}$. By the transitivity of $G$ on $V$, for every $v\in V$ there exists $\xi\in\mathsf{Orb}_G(\xi')$ such that $\alpha(\xi(v))=a_1$. Then $f(v)=1_F$; since $v$ is arbitrary, we have $f=\mathbbm{1}_F$.

Finally, if $g\in G$, $g\neq 1_G$, then the faithfulness of the action of $G$ on $V$ yields the existence of $v\in V$ such that $gv=v'\neq v$. Certainly, there exist $\xi\in\mathsf{Orb}_G(\xi')$ and $\varphi\in\Phi_\xi$ such that $\varphi(v)\neq \varphi(v')$: it suffices to choose $\xi$ in such a way that $\xi(v)\in\alpha^{-1}\left(a_1\right)$ and observe that $\left\lvert Y_{a_1}\right\rvert\geq2$. Then $\varphi\left(g^{-1}v'\right)\neq\varphi(v')$, so that $\left(\mathbbm{1}_F,g\right)(\xi,\varphi)\neq(\xi,\varphi)$ and the kernel of the action of $F\wr G$ on $\mathsf{Orb}_{F\wr G}\left(\xi',\varphi_{\xi'}\right)$ is trivial.

\medskip
(3) Suppose that there exists $g\in G$ as in the statement. From the definition of $\varphi_{\xi}$ it follows that $\varphi_{\xi'}(v)=y_{\alpha\circ\xi'(v)}=y_{\alpha\circ\xi\left(g^{-1}v\right)}=\varphi_{\xi}\left(g^{-1}v\right)=\left[g\varphi_{\xi}\right](v)$ for every $v\in V$, that is
\begin{equation}\label{gphixi}
\varphi_{\xi'}=g\varphi_{\xi}.
\end{equation}
If $(f',g')\in\mathsf{Stab}_{F\wr G}\left(\xi,\varphi_{\xi}\right)$, then $g'\in\mathsf{Stab}_G(\xi)$ and therefore
\begin{align*}
\left(\mathbbm{1}_F,g\right)(f',g')\left(\mathbbm{1}_F,g\right)^{-1}\left(\xi',\varphi_{\xi'}\right)
&=\left(gf',gg'g^{-1}\right)\left(\xi',g\varphi_{\xi}\right) \qquad &&(\text{by } \eqref{multwreath} \text{ and } \eqref{gphixi})\\
&=\left(gg'g^{-1}\xi',gf'\cdot gg'\varphi_{\xi}\right) \qquad &&(\text{by } \eqref{defactionOm})\\
&=\left(\xi',g\left(f'\cdot g'\varphi_{\xi}\right)\right) \qquad &&(\text{since } gg'g^{-1}\in \mathsf{Stab}_G(\xi'))\\
&=\left(\xi',g\varphi_{\xi}\right) \qquad &&(\text{since } f'\cdot g'\varphi_{\xi}=\varphi_{\xi})\\
&=\left(\xi',\varphi_{\xi'}\right), \qquad &&(\text{by } \eqref{gphixi})
\end{align*}
We deduce that $\left(\mathbbm{1}_F,g\right)\mathsf{Stab}_{F\wr G}\left(\xi,\varphi_{\xi}\right)\left(\mathbbm{1}_F,g\right)^{-1}\subseteq \mathsf{Stab}_{F\wr G}\left(\xi',\varphi_{\xi'}\right)$ and the proof of the converse inclusion is similar. We conclude that $\mathsf{Stab}_{F\wr G}\left(\xi,\varphi_{\xi}\right)$ and $\mathsf{Stab}_{F\wr G}\left(\xi',\varphi_{\xi'}\right)$ are conjugate in $F\wr G$ and the isomorphism between $\mathsf{Orb}_{F\wr G}(\xi,\varphi)$ and $\mathsf{Orb}_{F\wr G}(\xi',\varphi')$ follows from \eqref{Lemmaorb1} and \cite[Proposition 3.1.11]{book}. 

Conversely, suppose that $\mathsf{Orb}_{F\wr G}(\xi,\varphi)$ and $\mathsf{Orb}_{F\wr G}(\xi',\varphi')$ are isomorphic as $F\wr G$-spaces. Then there exists $(f,g)\in F\wr G$ such that $(f,g)\mathsf{Stab}_{F\wr G}\left(\xi,\varphi_{\xi}\right)(f,g)^{-1}=\mathsf{Stab}_{F\wr G}\left(\xi',\varphi_{\xi'}\right)$. If $(f',g')\in \mathsf{Stab}_{F\wr G}\left(\xi,\varphi_{\xi}\right)$, then
\begin{equation}\label{fgfgfg1}
(f,g)(f',g')(f,g)^{-1}=(f\cdot gf'\cdot gg'g^{-1}f^{-1},gg'g^{-1})\in \mathsf{Stab}_{F\wr G}\left(\xi',\varphi_{\xi'}\right).
\end{equation}
In particular, $g'^{-1}\in\mathsf{Stab}_G(\xi)$ if and only if $gg'^{-1}g^{-1}\in \mathsf{Stab}_G(\xi')$, which implies that $g\mathsf{Stab}_G(\xi) g^{-1}=\mathsf{Stab}_G(\xi')$. Moreover, $gg'g^{-1}\varphi_{\xi'}=\varphi_{\xi'}$ and from \eqref{Lemmastab1} and \eqref{fgfgfg1} it follows that for every $v\in V$ we have 
\begin{equation}\label{fgfgfg2}
f(v)f'\left(g^{-1}v\right)\left[gg'g^{-1}f^{-1}\right](v)\in\mathsf{Stab}_F\left(y_{\alpha(\xi'(v))}\right).
\end{equation}
Setting $f'=\mathbbm{1}_F$ in \eqref{fgfgfg2} we obtain
\begin{equation}\label{fgfgfg3}
f(v)\left[gg'g^{-1}f^{-1}\right](v)\in\mathsf{Stab}_F\left(y_{\alpha(\xi'(v))}\right). 
\end{equation}
Hence,
\begin{align*}
f'\left(g^{-1}v\right)
&\in f^{-1}(v) \mathsf{Stab}_F\left(y_{\alpha(\xi'(v))}\right)\left[gg'g^{-1}f\right](v) \qquad &&(\text{by } \eqref{fgfgfg2})\\
&= f^{-1}(v)\Bigl\{\mathsf{Stab}_F\left(y_{\alpha(\xi'(v))}\right)\left[gg'g^{-1}f\right](v)f^{-1}(v)\Bigr\}f(v) &&\\
&=f^{-1}(v)\mathsf{Stab}_F\left(y_{\alpha(\xi'(v))}\right)f(v). \qquad &&(\text{by } \eqref{fgfgfg3})
\end{align*}
Since $f'\left(g^{-1}v\right)$ is an arbitrary element of $\mathsf{Stab}_F\left(y_{\alpha\left(\xi\left(g^{-1}v\right)\right)}\right)$, this implies that 
$\mathsf{Stab}_F\left(y_{\alpha(\xi'(v))}\right)$ and $\mathsf{Stab}_F\left(y_{\alpha\left(\xi\left(g^{-1}v\right)\right)}\right)$ are conjugate in $F$, so that actually $\alpha(\xi'(v))=\alpha\left(\xi\left(g^{-1}v\right)\right)$, because the $Y_a$'s are mutually non-isomorphic. In conclusion, $\alpha\circ\xi'=\alpha\circ\xi\circ g^{-1}$.

Finally, suppose that $\xi$ and $\xi'$ are equivalent, so that there exist $\omega\in \Omega$ and $g\in G$ such that $\omega\circ\xi\circ g^{-1}=\xi'$. Since $\omega$ is a bijection, we have 
\[
\begin{split}
g'\in\mathsf{Stab}_G(\xi)\quad&\Longleftrightarrow \quad\omega\circ\xi\circ g'^{-1}=\omega\circ\xi\quad\Longleftrightarrow \quad\omega\circ\xi\circ g^{-1}gg'^{-1}g^{-1}=\omega\circ\xi\circ g^{-1}\\
&\Longleftrightarrow \quad\xi'\circ gg'^{-1}g^{-1}=\xi'\quad\Longleftrightarrow \quad gg'^{-1}g^{-1}\in\mathsf{Stab}_G(\xi')
\end{split}
\]
so that $g\mathsf{Stab}_G(\xi)g^{-1}=\mathsf{Stab}_G(\xi')$. Moreover, from $\alpha\circ\omega=\alpha$ we deduce that
\[
\alpha\circ\xi\circ g^{-1}=\alpha\circ\omega\circ\xi\circ g^{-1}=\alpha\circ\xi'.
\]
\end{proof}

When $G$ is the symmetric group, the equivalence is also necessary for the isomorphism in Lemma \ref{Lemmaorb}\eqref{Lemmaorb3}: see Lemma \ref{Lemmaorbsym}\eqref{Lemmaorbsym3}. However, in general, this is not the case: see Example \ref{Exampleequiv}.

In the remainder of this section, we fix $\xi_0\in\mathsf{Surj}(V,X)$ and consider the action of $F\wr G$ on the orbit $\mathsf{Orb}_{F\wr G}\left(\xi_0,\varphi_{\xi_0}\right)$.
For all $a,b\in A$, we denote by
\begin{equation}\label{decYa}
Y_a=\bigsqcup_{i=0}^{\ell_{a,b}}\Upsilon_{a,b,i}
\end{equation}
the decomposition of $Y_a$ into the orbits of $\mathsf{Stab}_F\left(y_b\right)$, with $\Upsilon_{a,a,0}=\left\{y_a\right\}$. Then we denote by 
\begin{equation}\label{decOrbxi}
\mathsf{Orb}_G\left(\xi_0\right)=\bigsqcup_{j=0}^n\Delta_j
\end{equation}
the decomposition of $\mathsf{Orb}_G\left(\xi_0\right)$ into the orbits of $\mathsf{Stab}_G\left(\xi_0\right)$, with $\Delta_0=\left\{\xi_0\right\}$. 
For every $0\leq j\leq n$, we set
\begin{equation}\label{defPsixij}
\Psi_j\coloneqq\left\{(\xi,\mathsf{i})\mathrel{:}\xi\in\Delta_j \text{ and }\mathsf{i}\in\mathbb{N}^V,\text{ with }0\leq \mathsf{i}(v)\leq\ell_{\alpha\circ\xi(v),\alpha\circ\xi_0(v)}\text{ for all } v \in V\right\}. 
\end{equation} 
The group $\mathsf{Stab}_G\left(\xi_0\right)$ acts on $\Psi_j$ in the natural way: if $(\xi,\mathsf{i})\in\Psi_j$ and $k\in\mathsf{Stab}_G\left(\xi_0\right)$, then we set $k(\xi,\mathsf{i})\coloneqq(k\xi,k\mathsf{i})$ and clearly $0\leq k\mathsf{i}(v)=\mathsf{i}\left(k^{-1}v\right)\leq\ell_{\alpha\circ\xi\left(k^{-1}v\right),\alpha\circ\xi_0\left(k^{-1}v\right)}=\ell_{\alpha\circ\left(k\xi\right)(v),\alpha\circ\xi_0(v)}$ for all $v\in V$, so that $k(\xi,\mathsf{i})\in\Psi_j$. Finally, we denote by 
\begin{equation}\label{decPsixij}
\Psi_j=\bigsqcup_{r=0}^{s_j}\Psi_{j,r}
\end{equation} 
the decomposition of $\Psi_j$ into the orbits of $\mathsf{Stab}_G\left(\xi_0\right)$; we assume that $\Psi_{0,0} = \{(\xi_0, 0)\}$, where $0$ denotes the function identically $0$.

\begin{theorem}\label{TheodecOrbxiphi}
The decomposition of $\mathsf{Orb}_{F\wr G}\left(\xi_0,\varphi_{\xi_0}\right)$ into orbits of $\mathsf{Stab}_{F\wr G}\left(\xi_0,\varphi_{\xi_0}\right)$ is given by
\begin{equation}\label{decOrbxiphi}
\mathsf{Orb}_{F\wr G}\left(\xi_0,\varphi_{\xi_0}\right)=\bigsqcup_{j=0}^n\;\bigsqcup_{r=0}^{s_j}\digamma_{j,r},
\end{equation} 
where 
\begin{equation}\label{defdigjr}
\digamma_{j,r}\coloneqq\bigsqcup_{(\xi,\mathsf{i})\in\Psi_{j,r}}\left[\{\xi\}\times\left(\bigtimes\limits_{v\in V}\Upsilon_{\alpha\circ\xi(v),\alpha\circ\xi_0(v),\mathsf{i}(v)}\right)\right].
\end{equation} 
 In particular, $\digamma_{0,0} = \{(\xi_0, \varphi_{\xi_0})\}$
\end{theorem}
\begin{proof}
From \eqref{FVprod}, \eqref{FVorbdec} and \eqref{decOrbxi}, we deduce that
\[
\begin{split}
\mathsf{Orb}_{F\wr G}\left(\xi_0,\varphi_{\xi_0}\right)=&\bigsqcup_{j=0}^n\;\bigsqcup_{\xi\in\Delta_j}\left[\{\xi\}\times\left(\bigtimes\limits_{v\in V}Y_{\alpha(\xi(v))}\right)\right]\\
(\text{by }\eqref{decYa})\quad=&\bigsqcup_{j=0}^n\;\bigsqcup_{\xi\in\Delta_j}\left\{\{\xi\}\times\left[\bigtimes\limits_{v\in V}\left(\bigsqcup_{\mathsf{i}(v)=0}^{\ell_{\alpha\circ\xi(v),\alpha\circ\xi_0(v)}}\Upsilon_{\alpha\circ\xi(v),\alpha\circ\xi_0(v),\mathsf{i}(v)}\right)\right]\right\}\\
(\text{by }\eqref{defPsixij})\quad=&\bigsqcup_{j=0}^n\;\bigsqcup_{(\xi,\mathsf{i})\in\Psi_j}\left\{\{\xi\}\times\left[\bigtimes\limits_{v\in V}\left(\Upsilon_{\alpha\circ\xi(v),\alpha\circ\xi_0(v),\mathsf{i}(v)}\right)\right]\right\}\\
(\text{by }\eqref{decPsixij})\quad=&\bigsqcup_{j=0}^n\;\bigsqcup_{r=0}^{s_j}\left\{\bigsqcup_{(\xi,\mathsf{i})\in\Psi_{j,r}}\left[\{\xi\}\times\left(\bigtimes\limits_{v\in V}\Upsilon_{\alpha\circ\xi(v),\alpha\circ\xi_0(v),\mathsf{i}(v)}\right)\right]\right\},\\
\end{split}
\]
and this proves the decomposition \eqref{decOrbxiphi}. 

Now we show that each $\digamma_{j,r}$ is indeed a $\mathsf{Stab}_{F\wr G}\left(\xi_0,\varphi_{\xi_0}\right)$-homogeneous space. Suppose that $\left(\xi',\varphi'\right)\in\digamma_{j,r}$. Then there exists $\mathsf{i}'\in \mathbb{N}^V$ such that $(\xi',\mathsf{i}')\in\Psi_{j,r}$ and $\varphi'(v)\in\Upsilon_{\alpha\circ\xi'(v),\alpha\circ\xi_0(v),\mathsf{i}'(v)}$ for every $v\in V$; cf. \eqref{FVprod}. If $(f,k)\in\mathsf{Stab}_{F\wr G}\left(\xi_0,\varphi_{\xi_0}\right)$, then $k\in\mathsf{Stab}_G\left(\xi_0\right)$, $k\left(\xi',\mathsf{i}'\right)\in\Psi_{j,r}$, and $f(v)\in \mathsf{Stab}_F\left(y_{\alpha\left(\xi_0(v)\right)}\right)$  for every $v\in V$; cf. \eqref{Lemmastab1}. From \eqref{decYa}, it follows that
\[
\begin{split}
\left(f\cdot k\varphi'\right)(v)=&f(v)\varphi'\left(k^{-1}v\right)\in f(v)\Upsilon_{\alpha\circ\xi'\left(k^{-1}v\right),\alpha\circ\xi_0\left(k^{-1}v\right),\mathsf{i}'\left(k^{-1}v\right)}=f(v)\Upsilon_{\alpha\circ\left(k\xi'\right)(v),\alpha\circ\xi_0(v),\left(k\mathsf{i}'\right)(v)}\\
=&\Upsilon_{\alpha\circ\left(k\xi'\right)(v),\alpha\circ\xi_0(v),\left(k\mathsf{i}'\right)(v)},
\end{split}
\]
and therefore $(f,k)\left(\xi',\varphi'\right)\in\digamma_{j,r}$. That is, $\digamma_{j,r}$ is $\mathsf{Stab}_{F\wr G}\left(\xi_0,\varphi_{\xi_0}\right)$-invariant.

Finally, suppose that we also have $\left(\xi'',\varphi''\right)\in\digamma_{j,r}$, so that there exists $\mathsf{i}''\in \mathbb{N}^V$ such that $(\xi'',\mathsf{i}'')\in\Psi_{j,r}$ and $\varphi''(v)\in\Upsilon_{\alpha\circ\xi''(v),\alpha\circ\xi_0(v),\mathsf{i}''(v)}$ for every $v\in V$. Then there exists $k\in\mathsf{Stab}_G\left(\xi_0\right)$ such that $k\left(\xi',\mathsf{i}'\right)=\left(\xi'',\mathsf{i}''\right)$. Therefore, for every $v\in V$ we have
\[
\left(k\varphi'\right)(v)=\varphi'\left(k^{-1}v\right)\in\Upsilon_{\alpha\circ\left(k\xi'\right)(v),\alpha\circ\xi_0(v),\left(k\mathsf{i}'\right)(v)}=\Upsilon_{\alpha\circ\xi''(v),\alpha\circ\xi_0(v),\mathsf{i}''(v)},
\]
so that there exists $f(v)\in\mathsf{Stab}_F\left(y_{\alpha\left(\xi_0(v)\right)}\right)$ such that $f(v)\left(k\varphi'\right)(v)=\varphi''(v)$. In conclusion, $(f,k)\left(\xi',\varphi'\right)=\left(\xi'',\varphi''\right)$ and the action of $\mathsf{Stab}_{F\wr G}\left(\xi_0,\varphi_{\xi_0}\right)$ on $\digamma_{j,r}$ is transitive.
\end{proof}

Now we apply the notions of extension, inflation, induction, and inertia group described in Section \ref{Secprelimin} to the action of $F\wr G$ on $\mathsf{Orb}_{F\wr G}\left(\xi_0,\varphi_{\xi_0}\right)$.

\begin{proposition}\label{propindact}
Denote by $\eta$ the action of $F^V$ on $\mathsf{Orb}_{F^V}\left(\varphi_{\xi_0}\right)$ (cf. \eqref{defdirprodact} and \eqref{FVprod}) and by $\beta$ the action of $\mathsf{Stab}_G\left(\alpha\circ\xi_0\right)$ on $\mathsf{Orb}_{\mathsf{Stab}_G\left(\alpha\circ\xi_0\right)}\left(\xi_0\right)$.
\begin{enumerate}[(1)]
\item\label{propindact1}
The set $\mathsf{Orb}_{F^V}\left(F\wr\varphi_{\xi_0}\right)$ is $\gamma_0\left(\mathsf{Stab}_G\left(\alpha\circ\xi_0\right)\right)$-invariant, and 
\[
\widetilde{\eta}\coloneqq\left[\Res^{F\wr G}_{F\wr\mathsf{Stab}_G\left(\alpha\circ\xi_0\right)}\gamma_0\right]_{\mathsf{Orb}_{F^V}\left(\varphi_{\xi_0}\right)}
\]
is an extension of $\eta$ to $F\wr\mathsf{Stab}_G\left(\alpha\circ\xi_0\right)$.
\item\label{propindact2}
The inertia group of $\eta$ in $F\wr G$ is $F\wr \mathsf{Stab}_G\left(\alpha\circ\xi_0\right)$.
\item\label{propindact3}
We have the direct product decomposition
\begin{equation}\label{OrStabOr}
\mathsf{Orb}_{F\wr\mathsf{Stab}_G\left(\alpha\circ\xi_0\right)}\left(\xi_0,\varphi_{\xi_0}\right)=\mathsf{Orb}_{\mathsf{Stab}_G\left(\alpha\circ\xi_0\right)}\left(\xi_0\right)\bigtimes\mathsf{Orb}_{F^V}\left(\varphi_{\xi_0}\right),
\end{equation}
and the action of $F\wr\mathsf{Stab}_G\left(\alpha\circ\xi_0\right)$ on this space is isomorphic to the diagonal action $\mathsf{diag}\left(\overline{\beta},\widetilde{\eta}\right)$; cf. \eqref{defdiagact} and \eqref{definflperrep}. 
\item\label{propindact4}
The action of $F\wr G$ on $\mathsf{Orb}_{F\wr G}\left(\xi_0,\varphi_{\xi_0}\right)$ is isomorphic to the action induced by $\mathsf{diag}\left(\overline{\beta},\widetilde{\eta}\right)$.
\end{enumerate}
\end{proposition}
\begin{proof}
(1)
If $k\in\mathsf{Stab}_G\left(\alpha\circ\xi_0\right)$ and $\varphi\in\mathsf{Orb}_{F^V}\left(\varphi_{\xi_0}\right)$, then 
\[
[k\varphi](v)=\varphi\left(k^{-1}v\right)\in Y_{\alpha\left(\xi_0\left(k^{-1}v\right)\right)}=Y_{\left[k\left(\alpha\circ\xi_0\right)\right](v)}= Y_{\alpha\left(\xi_0(v)\right)}\quad\text{ for all }v\in V,
\]
which implies that $k\varphi\in\mathsf{Orb}_{F^V}\left(\varphi_{\xi_0}\right)$. Moreover, if $f\in F^V$, then $\gamma_0(f,k)\varphi=f\cdot k\varphi\in\mathsf{Orb}_{F^V}\left(\varphi_{\xi_0}\right)$, yielding an extension of the action of $F^V$.

\medskip
(2)
From \eqref{propindact}, it follows that $F\wr\mathrm{Stab}_G(\alpha\circ\xi_0)$ is contained in the inertia group of $\eta$. On the other hand, by an argument similar to the one used for \eqref{fgsigmaconj} (cf. \cite[Lemma 2.4.1]{book3}), it is straightforward to check that $\!^{(f,g)}\eta\simeq\bigtimes_{v\in V}\eta_{\alpha\circ\xi_0\left(g^{-1}v\right)}$. 
If $g\in G\setminus\mathsf{Stab}_G\left(\alpha\circ\xi_0\right)$, then there exists $v\in V$ such that $\alpha\circ\xi_0\left(g^{-1}v\right)\neq \alpha\circ\xi_0(v)$, implying that $Y_{\alpha\circ\xi_0\left(g^{-1}v\right)}$ and $Y_{\alpha\left(\xi_0(v)\right)}$ are not isomorphic as $F$-spaces; therefore $\!^{(f,g)}\eta$ and $\eta$ are not isomorphic, because two direct product subgroups of $F^V$ are conjugate in $F^V$ if and only if their factors are conjugate coordinatewise in $F$.

\medskip
(3)
Clearly, $\Phi_{k\xi_0}=\Phi_{\xi_0}=\mathsf{Orb}_{F^V}\left(\varphi_{\xi_0}\right)$ for every $k\in\mathsf{Stab}_G\left(\alpha\circ\xi_0\right)$, so that, by replacing $\mathsf{Orb}_G\left(\xi_0\right)$ with $\mathsf{Orb}_{\mathsf{Stab}_G\left(\alpha\circ\xi_0\right)}\left(\xi_0\right)$ in \eqref{FVorbdec}, we obtain \eqref{OrStabOr}. Moreover, if $\xi\in\mathsf{Orb}_{\mathsf{Stab}_G\left(\alpha\circ\xi_0\right)}\left(\xi_0\right)$, $\varphi\in\mathsf{Orb}_{F^V}\left(\varphi_{\xi_0}\right)$ and $(f,k)\in F\wr\mathsf{Stab}_G\left(\alpha\circ\xi_0\right)$, then 
\[
\gamma(f,k)(\xi,\varphi)=(k\xi,f\cdot k\varphi)=\left(\overline{\beta}(f,k)\xi,\widetilde{\eta}(f,k)\varphi\right).
\]

\medskip
(4)
We give a first proof, which is instructive because it relies on the explicit decomposition of the homogeneous space.
Denote by $\mathsf{R}$ a set of representatives for the left $\mathsf{Stab}_G\left(\alpha\circ\xi_0\right)$-cosets in $G$ (with $1_G\in\mathsf{R}$). Then 
\begin{equation}\label{OrStabOrStab}
\mathsf{Orb}_{F\wr G}\left(\xi_0,\varphi_{\xi_0}\right)=\bigsqcup_{r\in\mathsf{R}}\gamma\left(\mathbbm{1}_F,r\right)\mathsf{Orb}_{F\wr \mathsf{Stab}_G\left(\alpha\circ\xi_0\right)}\left(\xi_0,\varphi_{\xi_0}\right),
\end{equation}
and the subset $\mathsf{Orb}_{F\wr \mathsf{Stab}_G\left(\alpha\circ\xi_0\right)}\left(\xi_0,\varphi_{\xi_0}\right)$ is clearly $\gamma\left(F\wr\mathsf{Stab}_G\left(\alpha\circ\xi_0\right)\right)$-invariant. Hence, the statement follows from \eqref{defindact}. 

Alternatively, from $\mathsf{Stab}_{\mathsf{Stab}_G\left(\alpha\circ\xi_0\right)}\left(\xi_0\right)=\mathsf{Stab}_G\left(\xi_0\right)$ and \eqref{Lemmastab1} it follows that
\[
\mathsf{Stab}_{F\wr\mathsf{Stab}_G\left(\alpha\circ\xi_0\right)}\left(\xi_0,\varphi_{\xi_0}\right)=\mathsf{Stab}_{F\wr G}\left(\xi_0,\varphi_{\xi_0}\right),
\]
and the result follows immediately.
\end{proof}

\begin{example}\label{Example1}
Suppose that $\lvert X\rvert=\lvert A\rvert=1$ and denote by $Y$ the unique $F$-space involved. Then $\Gamma\equiv Y^V$ because in a pair $(\xi,\varphi)\in\Gamma$ the function $\xi$ is trivial (as is $\alpha$) and $\varphi$ is an arbitrary map from $V$ to $Y$. Moreover, the action \eqref{defactionOm} reduces to $(f,g)\varphi=f\cdot g\varphi$, that is, $\gamma\simeq\gamma_0$, which corresponds to the \emph{exponentiation action} and is transitive; cf. \cite[Section 2.1.2]{book3}, \cite[4.1.28]{JK}. In particular, with the notation in Proposition \ref{propindact},
$\beta$ and its inflation are trivial, $\widetilde{\eta}=\gamma_0$, the inertia group is the whole $F\wr G$, and there is no induction. If we fix $y_0\in Y$ and set $\varphi_0(v)=y_0$ for all $v\in V$, then $\mathsf{Stab}_{F\wr G}\left(\varphi_0\right)=\mathsf{Stab}_F\left(y_0\right)\wr G$, so that 
\[
Y^V\simeq (F\wr G)/\left(\mathsf{Stab}_F\left(y_0\right)\wr G\right).
\]
Now \eqref{decYa} reduces to $Y=\bigsqcup_{i=0}^\ell\Upsilon_i$, the decomposition of $Y$ into the orbits of $\mathsf{Stab}_F\left(y_0\right)$. Moreover, \eqref{defPsixij} becomes $\Psi\coloneqq\left\{\mathsf{i}\in\mathbb{N}^V\mathrel{:}0\leq \mathsf{i}(v)\leq \ell\text{ for all }v\in V\right\}$, while \eqref{decPsixij} is $\Psi=\bigsqcup_{r=0}^s\Psi_r$, the decomposition of $\Psi$ into the orbits of $G$. Finally, \eqref{decOrbxiphi} now is the decomposition of $Y^V$ into the orbits of $\mathsf{Stab}_F\left(y_0\right)\wr G$ and it is given by
\begin{equation}\label{YVdigr}
Y^V=\bigsqcup_{r=0}^s\digamma_r,\quad\text{ where }\quad \digamma_r=\bigsqcup_{\mathsf{i}\in\Psi_r}\left\{\varphi\in Y^V\mathrel{:}\varphi(v)\in\Upsilon_{\mathsf{i}(v)}\text{ for all } v \in V\right\}.
\end{equation}
\end{example}

\begin{example}\label{Example2}
Now suppose that $X=\{x_0,x_1\}$, $A=\{a_0,a_1\}$, $Y=Y_{a_0}$ is an arbitrary $F$-space, $\lvert Y_{a_1}\rvert=1$ and $\alpha(x_i)=a_i$, $i=0,1$. For every $v\in V$ define the function $\xi_v\in\mathsf{Surj}(V,X)$ by setting $\xi_v(v)=x_0$ and $\xi_v(v')=x_1$ if $v'\in V$ with $v'\neq v$. Then $g\xi_v=\xi_{gv}$, so that the $G$-action on the space of these functions is isomorphic to the action of $G$ on $V$. Moreover, if $\varphi\in\Phi_{\xi_v}$, then $\varphi(v)\in Y$ whereas if $v'\neq v$, then $\varphi(v')$ is in the trivial space $Y_{a_1}$. Denote by $\Gamma_0$ the set of all pairs $\left(\xi_v,\varphi\right)$ as above and consider the bijection $T\colon\Gamma_0\rightarrow V\times Y$ given by: $T\left(\xi_v,\varphi\right)=(v,\varphi(v))$, for all $v\in V$, $\varphi\in\Phi_{\xi_v}$. The \emph{composition action} $\gamma'$ of $F\wr G$ on $V\times Y$ is defined by setting
\begin{equation}\label{defcompaction}
\gamma'(f,g)(v,y)\coloneqq(gv,f(gv)y)
\end{equation} 
for all $(v,y)\in V\times Y$; see, e.g., \cite[Section 2.1.2]{book3}, \cite[4.1.18]{JK}. Now we show that $T$ is an isomorphism of group actions
\[
T\left[\gamma(f,g)(\xi_v,\varphi)\right]=T\left(g\xi_v,f\cdot g\varphi\right)=\left(gv,f(gv)\varphi(v)\right)=\gamma'(f,g)T\left(\xi_v,\varphi\right).
\]
Clearly, $T\left(\left\{\xi_v\right\}\times\Phi_{\xi_v}\right)=\{v\}\times Y$.
We fix $v_0\in V$, $y_0\in Y_{a_0}$ and define $\varphi_0$ by setting $\varphi_0\left(v_0\right)\coloneqq y_0$ and $\varphi_0(v)\in Y_{a_1}$ for all $v\in V, v\neq v_0$. Then 
\[
\mathsf{Stab}_{F\wr G}\left(\xi_{v_0},\varphi_0\right)=\left\{(f,k)\in F\wr G\mathrel{:} k\in\mathsf{Stab}_G\left(v_0\right), f\left(v_0\right)\in\mathsf{Stab}_F\left(y_0\right)\right\}\equiv\mathsf{Stab}_{F\wr G}\left(v_0,y_0\right).
\]

From now on, we identify $\Gamma_0$ with $V\times Y$. For $a=b=a_0$ we use $Y=\bigsqcup_{i=0}^\ell\Upsilon_i$, as in Example \ref{Example1}, while for $a=a_0, b=a_1$ the right-hand side of \eqref{decYa} reduces to $Y$ and $\ell_{a,b}=0$.
Now \eqref{decOrbxi} may be written in the form $V\simeq\mathsf{Orb}_G\left(v_0\right)=\bigsqcup_{j=0}^n\Delta_j$, where $\Delta_0=\left\{v_0\right\}$; it is the decomposition of $V$ into the orbits of $\mathsf{Stab}_G\left(v_0\right)$. Moreover, \eqref{defPsixij} becomes
\[
\Psi_0\coloneqq \left\{v_0\right\}\times\{0,1,\dotsc,\ell\},\qquad \Psi_j\coloneqq\Delta_j\times\{0\}, \quad\text{ for }j=1,2,\dotsc,n,
\]
and in \eqref{decPsixij} we have simply $\Psi_{0,r}=\left\{v_0\right\}\times\{r\}$, $r=0,1,\dotsc,\ell$ and $\Psi_{j,0}=\Delta_j\times\{0\}$, $j=1,2,\dotsc,n$. Therefore, in the present case, \eqref{decOrbxiphi} has the form 
\begin{equation}\label{orbdeccomp}
V\times Y=\left[\bigsqcup_{i=0}^\ell\left(\left\{v_0\right\}\times\Upsilon_i\right)\right]\bigsqcup\left[\bigsqcup_{j=1}^n\left(\Delta_j\times Y\right)\right];
\end{equation}
see, e.g., \cite[Proposition 3.1.2(ii)]{book3}. Finally, referring to Proposition \ref{propindact}, now the right-hand side of \eqref{OrStabOr} is just $\left\{v_0\right\}\times Y$; $\eta$ is the action of $F^V$ on $Y$ (given by $\eta(f)y=f\left(v_0\right)y$); the $F\wr\mathsf{Stab}_G\left(v_0\right)$-invariance of $\left\{v_0\right\}\times Y$ yields the extension of $\eta$; the inertia group of $\eta$ is $F\wr\mathsf{Stab}_G\left(v_0\right)$; $\beta$ and its inflation in Proposition \ref{propindact}\eqref{propindact3} are trivial; \eqref{OrStabOrStab} in the present setting coincides with
\[
V\times Y\simeq\bigsqcup_{v\in V}\left(\{v\}\times Y\right).
\]
That is, the action of $F\wr G$ on $V\times Y$ is induced by the action of $F\wr\mathsf{Stab}_G\left(v_0\right)$ on $\left\{v_0\right\}\times Y$; cf. Proposition \ref{propindact}\eqref{propindact4}.
\end{example}

\begin{example}\label{Example21}
Now suppose that $\lvert V\rvert\geq 3$, $X=\left\{x_0,x_1,x_2\right\}$, $A=\left\{a_0,a_1\right\}$, $Y=Y_{a_0}$ is an arbitrary $F$-space, $\left\lvert Y_{a_1}\right\rvert=1$ and $\alpha\left(x_0\right)=a_0$, $\alpha\left(x_1\right)=\alpha\left(x_2\right)=a_1$. Recalling the notation in \eqref{natproj}, for all $(v,w)\in\mathsf{Per}_2(V)$ we define the function $\xi_{(v,w)}\in\mathsf{Surj}(V,X)$ by setting $\xi_{(v,w)}(v)=x_0$, $\xi_{(v,w)}(w)=x_1$ and $\xi_{(v,w)}(v')=x_2$ if $v'\in V$, $v'\neq v,w$. Then $g\xi_{(v,w)}=\xi_{(gv,gw)}$, so that the $G$-action on the space of these functions is isomorphic to the action of $G$ on $\mathsf{Per}_2(V)$. Moreover, if $\varphi\in\Phi_{\xi_{(v,w)}}$, then $\varphi(v)\in Y$ whereas $\varphi(v')$ lies in the trivial space $Y_{a_1}$ whenever $v'\neq v$. We denote by $\Gamma_0$ the set of all pairs $\left(\xi_{(v,w)},\varphi\right)$ as above; arguing similarly to Example \ref{Example2}, it is straightforward to prove that the action of $F\wr G$ on $\Gamma_0$ is isomorphic to the action on $\mathsf{Per}_2(V)\times Y$ given by 
\begin{equation}\label{DefactPerY}
(f,g)(v,w,y)=(gv,gw,f(gv)y).
\end{equation}
On the other hand, consider the diagonal action of $F\wr G$ on $V\times Y\times V$, where $F\wr G$ acts on $V\times Y$ via the composition action in Example \ref{Example2} and on the final $V$ factor via the $G$-action. Then \eqref{DefactPerY} is isomorphic to the restriction of this diagonal action to the subset $\{(v,y,w)\in V\times Y\times V\mathrel{:}v\neq w\}$. 

Now we slightly modify the notation of Example \ref{Example2}. We fix $v_0\in V$ and $y_0\in Y_{a_0}$, and assume that $V=\bigsqcup_{j=0}^n\Delta_{v_0,j}$ is the decomposition of $V$ into the orbits of $\mathsf{Stab}_G\left(v_0\right)$, with $\Delta_{v_0,0}=\{v_0\}$. Next, for every $v\in V$ we choose $g\in G$ such that $gv_0=v$ and set $\Delta_{v,j}\coloneqq g\Delta_{v_0,j}$. Clearly, this definition does not depend on the choice of $g$, and since $\mathsf{Stab}_G(v)=g\mathsf{Stab}_G\left(v_0\right)g^{-1}$, we obtain that $V=\bigsqcup_{j=0}^n\Delta_{v,j}$ is the decomposition of $V$ into the orbits of $\mathsf{Stab}_G(v)$, with $\Delta_{v,0}=\{v\}$.
Hence,
\[
\mathsf{Per}_2(V)=\bigsqcup_{v\in V}\bigsqcup_{j=1}^n\left(\{v\}\times\Delta_{v,j}\right)=\bigsqcup_{j=1}^n\left[\bigsqcup_{v\in V}\left(\{v\}\times\Delta_{v,j}\right)\right].
\]
Moreover, if we fix $v_j\in\Delta_{v_0,j}$ we also have 
\begin{equation}\label{orbv0vj}
\mathsf{Orb}_G\left(v_0,v_j\right)=\bigsqcup_{v\in V}\left(\{v\}\times\Delta_{v,j}\right)
\end{equation}
for $j=1,2,\dotsc,n$, and these are the $G$-orbits on $\mathsf{Per}_2(V)$. It readily follows that
\begin{equation}\label{decorbPer2VY}
\mathsf{Per}_2(V)\times Y=\bigsqcup_{j=1}^n\mathsf{Orb}_{F\wr G}\left(v_0,v_j,y_0\right)
\end{equation}
is the decomposition of $\mathsf{Per}_2(V)\times Y$ into the orbits of $F\wr G$; see, e.g., \cite[Exercise 10.4.15]{book4}. In particular, if the group $G$ is transitive on $\mathsf{Per}_2(V)$, that is the action of $G$ on $V$ is $2$-transitive \cite{Cameron,DM}, then $n=1$, $\Delta_{v,1}=V\setminus\{v\}$, and $F\wr G$ is transitive on the whole space $\mathsf{Per}_2(V)\times Y$.

In what follows, we fix $1\leq j\leq n$ and we examine the action of $F\wr G$ on a single orbit in \eqref{decorbPer2VY}. That is, from now on, the function $\xi_0$ in the general theory (cf. \eqref{decOrbxi}) is $\xi_{v_0,v_j}$. First, note that $\mathsf{Stab}_G\left(v_0,v_j\right)=\mathsf{Stab}_G\left(v_0\right)\cap\mathsf{Stab}_G\left(v_j\right)$. Next, we set 
\[
\mathcal{I}_j\coloneqq\left\{(r,s)\mathrel{:}0\leq r,s\leq n\text{ and }\Delta_{v_0,r}\cap\Delta_{v_j,s}\neq\emptyset\right\}. 
\]
For instance, $(j,0)\in\mathcal{I}_j$, because $v_j\in\Delta_{v_0,j}$, but $(r,0)\notin\mathcal{I}_j$ for $r\neq j$, because $v_j\notin\Delta_{v_0,r}$ in this case. Then we can refine \eqref{orbv0vj} by writing it in the form
\[
\mathsf{Orb}_G\left(v_0,v_j\right)=\bigsqcup_{(r,s)\in\mathcal{I}_j}\;\bigsqcup_{v\in \Delta_{v_0,r}\cap\Delta_{v_j,s}}\left(\{v\}\times\Delta_{v,j}\right),
\]
where, for a fixed pair $(r,s)\in\mathcal{I}_j$, the disjoint union over $v$ is $\mathsf{Stab}_G\left(v_0,v_j\right)$-invariant. Therefore, we can suppose that
\[
\bigsqcup_{v\in\Delta_{v_0,r}\cap\Delta_{v_j,s}}\left(\{v\}\times\Delta_{v,j}\right)=\bigsqcup_{t=0}^{m_{r,s}}\Delta_{r,s,t}
\]
is the decomposition of the left-hand side into the $\mathsf{Stab}_G\left(v_0,v_j\right)$-orbits; in particular, if $v_0\in\Delta_{v_j,s_0}$, then $\left(0,s_0\right)\in\mathcal{I}_j$, 
\begin{equation}\label{v0Deltav0}
\left\{v_0\right\}\times\Delta_{v_0,j}=\bigsqcup_{t=0}^{m_{0,s_0}}\Delta_{0,s_0,t}
\end{equation} 
and we can suppose that $\Delta_{0,s_0,0}=\left\{\left(v_0,v_j\right)\right\}$. In conclusion, 
\[
\mathsf{Orb}_G\left(v_0,v_j\right)=\bigsqcup_{(r,s)\in\mathcal{I}_j}\;\bigsqcup_{t=0}^{m_{r,s}}\Delta_{r,s,t}
\]
is the decomposition of $\mathsf{Orb}_G\left(v_0,v_j\right)$ into the $\mathsf{Stab}_G\left(v_0,v_j\right)$-orbits. 

As in Example \ref{Example2}, for $a=b=a_0$ we use $Y=\bigsqcup_{i=0}^\ell\Upsilon_i$, while for $a=a_0, b=a_1$ the right-hand side of \eqref{decYa} reduces to $Y$ and $\ell_{a,b}=0$.
Note that
\[
\mathsf{Stab}_{F\wr G}\left(v_0,v_j,y_0\right)=\left\{(f,k)\in F\wr G\mathrel{:}k\in\mathsf{Stab}_G\left(v_0,v_j\right), f\left(v_0\right)\in\mathsf{Stab}_F\left(y_0\right)\right\}.
\]
In the present setting, \eqref{defPsixij} becomes
\begin{align*}
&\Psi_{0,s_0,t}\coloneqq \Delta_{0,s_0,t}\times\{0,1,\dotsc,\ell\}&&\text{for }t=0,1,\dotsc,m_{0,s_0},\\
&\Psi_{r,s,t}\coloneqq\Delta_{r,s,t}\times\{0\}&&\text{for }(r,s)\neq\left(0,s_0\right), 0\leq t\leq m_{r,s}.
\end{align*}
Indeed, now $\ell_{\alpha\circ\xi(v),\alpha\circ\xi_0(v)}$ is equal to $\ell_{\alpha\circ\xi_{(v',w)}(v),\alpha\circ\xi_{(v_0,v_j)}(v)}$ and we have $\alpha\circ\xi_{(v',w)}(v)=a_0=\alpha\circ\xi_{\left(v_0,v_j\right)}(v)$ if and only if $v'=v=v_0$. This leads to the orbits in \eqref{v0Deltav0}, and in a pair $\left(\xi_{\left(v_0,w\right)},\mathsf{i}\right)$, with $w\in\Delta_{v_0,j}$, we have $0\leq\mathsf{i}\left(v_0\right)\leq\ell$ and $\mathsf{i}(v)=0$ for $v\neq v_0$. This also shows that the action of $\mathsf{Stab}_{F\wr G}\left(v_0,v_j\right)$ on $\mathsf{i}$ is trivial, that is in \eqref{decPsixij} we have $\Psi_{0,s_0,t,r'}=\Delta_{0,s_0,t}\times\{r'\}$, $r'=0,1,\dotsc,\ell$. Clearly, $\Psi_{r,s,t,0} = \Delta_{r,s,t} \times \{0\}$, for $(r,s)\neq\left(0,s_0\right),$ $0\leq t\leq m_{r,s}$. We deduce that, in the present case, \eqref{decOrbxiphi} has the form 
\[
\mathsf{Orb}_{F\wr G}\left(v_0,v_j,y_0\right)=\left[\bigsqcup_{t=0}^{m_{0,s_0}}\bigsqcup_{i=0}^\ell\left(\Delta_{0,s_0,t}\times\Upsilon_i\right)\right]\bigsqcup\left[\bigsqcup_{\substack{(r,s)\in\mathcal{I}_j\\ (r,s)\neq\left(0,s_0\right)}}\;\bigsqcup_{t=0}^{m_{r,s}}\left(\Delta_{r,s,t}\times Y\right)\right].
\]

Finally, note that now $\alpha\circ\xi_{\left(v_0,v_j\right)}=\alpha\circ\xi_{v_0}$ (the latter is as in Example \ref{Example2}), so that $\mathsf{Stab}_G\left(\alpha\circ\xi_{\left(v_0,v_j\right)}\right)=\mathsf{Stab}_G\left(v_0\right)$ and \eqref{OrStabOr} in the present setting is
\[
\mathsf{Orb}_{F\wr\mathsf{Stab}_G\left(v_0\right)}\left(v_0,v_j,y_0\right)=\left(\left\{v_0\right\}\times\Delta_{v_0,j}\right)\times Y.
\] 
In particular, referring to Proposition \ref{propindact}, $\beta$ is the action of $\mathsf{Stab}_G\left(v_0\right)$ on $\left\{v_0\right\}\times\Delta_{v_0,j}$, $\eta$ is the action of $F^V$ on $Y$ (given by $\eta(f)y=f\left(v_0\right)y$), the $F\wr\mathsf{Stab}_G\left(v_0\right)$-invariance of $\left(\left\{v_0\right\}\times\Delta_{v_0,j}\right)\times Y$ yields the extension of $\eta$, the inertia group of $\eta$ is $F\wr\mathsf{Stab}_G\left(v_0\right)$
and \eqref{OrStabOrStab} coincides with
\begin{equation}\label{orbv0vjy0}
\mathsf{Orb}_{F\wr G}\left(v_0,v_j,y_0\right)\simeq\bigsqcup_{v\in V}\left[\left(\left\{v\right\}\times\Delta_{v,j}\right)\times Y\right].
\end{equation}
That is, the action of $F\wr G$ on $\mathsf{Orb}_{F\wr G}\left(v_0,v_j,y_0\right)$ is induced by the action of $F\wr\mathsf{Stab}_G\left(v_0\right)$ on $\left(\left\{v_0\right\}\times\Delta_{v_0,j}\right)\times Y$; cf. Proposition \ref{propindact}\eqref{propindact4}.
\end{example}

\begin{example}\label{Example3}
As in Example \ref{Example2}, suppose again that $X=\left\{x_0,x_1\right\}$, $A=\left\{a_0,a_1\right\}$, $Y=Y_{a_0}$ is an arbitrary $F$-space, $\left\lvert Y_{a_1}\right\rvert=1$ and $\alpha\left(x_i\right)=a_i$, $i=0,1$. To each $\emptyset\subsetneqq Z\subsetneqq V$ we associate the function $\xi_Z\in\mathsf{Surj}(V,X)$ defined by setting $\xi_Z(v)\coloneqq x_0$ if $v\in Z$, $\xi_Z(v)\coloneqq x_1$ if $v\notin Z$. Clearly, $g\xi_Z=\xi_{gZ}$ and $\mathsf{Stab}_G\left(\xi_Z\right)=\mathsf{Stab}_G(Z)$. Now fix a subset $\emptyset\subsetneqq Z_0\subsetneqq V$ and consider the natural, transitive action of $F\wr G$ on $\bigsqcup\limits_{Z\in\mathsf{Orb}_G\left(Z_0\right)}\left(\{Z\}\times Y^Z\right)$ given by
\[
(f,g)(Z,\phi)=\left(g Z,f\cdot g\phi\right)
\]
for all $(f,g)\in F\wr G$, $Z\in\mathsf{Orb}_G\left(Z_0\right)$ and $\phi\in Y^Z$. Next, we fix an element $y_0\in Y$ so that $\varphi_{\xi_{Z_0}}(v)=y_0$ for all $v\in Z_0$ and $\varphi_{\xi_{Z_0}}(v)\in Y_{a_1}$ for all $v\notin Z_0$. Then \eqref{FVorbdec} becomes
\[
\mathsf{Orb}_{F\wr G}\left(\xi_{Z_0},\varphi_{\xi_{Z_0}}\right)=\bigsqcup_{Z\in\mathsf{Orb}_G\left(Z_0\right)}\bigl(\left\{\xi_Z\right\}\times\mathsf{Orb}_{F^V}\left(\varphi_{\xi_Z}\right)\bigr).
\]
Consequently, the map
\begin{equation}\label{isomJohn}
\begin{array}{ccc}
\bigsqcup\limits_{Z\in\mathsf{Orb}_G\left(Z_0\right)}\bigl(\left\{\xi_Z\right\}\times\mathsf{Orb}_{F^V}\left(\varphi_{\xi_Z}\right)\bigr)&\longrightarrow&\bigsqcup\limits_{Z\in\mathsf{Orb}_G\left(Z_0\right)}\left(\{Z\}\times Y^Z\right)\\
\left(\xi_Z,\varphi\right)&\longmapsto&\left(Z,\varphi\rvert_Z\right)
\end{array}
\end{equation}
is an isomorphism of $F\wr G$-spaces. Moreover, \eqref{Lemmastab1} becomes 
\[
\begin{split}
\mathsf{Stab}_{F\wr G}\left(\xi_{Z_0},\varphi_{\xi_{Z_0}}\right)=&\left\{(f,k)\in F\wr G\mathrel{:} k\in\mathsf{Stab}_G\left(Z_0\right), f\left(v\right)\in\mathsf{Stab}_F\left(y_0\right)\text{ for all }v\in Z_0\right\}\\
\equiv&\mathsf{Stab}_{F\wr G}\left(Z_0,\phi_{Z_0}\right),
\end{split}
\]
where $\phi_{Z_0}\in Y^{Z_0}$ is defined by setting $\phi_{Z_0}(v)\coloneqq y_0$ for all $v\in Z_0$. 

From now on, we use the right-hand side of \eqref{isomJohn} to represent our homogeneous space.  As in Example \ref{Example2}, we use $Y=\bigsqcup_{i=0}^\ell\Upsilon_i$ for $a=b=a_0$, while for $a=a_0$ and $b=a_1$ the right-hand side of \eqref{decYa} reduces to $Y$. Now  \eqref{decOrbxi} may be written in the form $\mathsf{Orb}_G\left(Z_0\right)=\bigsqcup_{j=0}^n\Delta_j$, where $\Delta_0=\left\{Z_0\right\}$. Moreover, \eqref{defPsixij} becomes
\[
\Psi_j\coloneqq\bigl\{(Z,\mathsf{i})\mathrel{:}Z\in\Delta_j \text{ and }\mathsf{i}\in\mathbb{N}^Z, 0\leq \mathsf{i}(v)\leq \ell \text{ if }v\in Z\cap Z_0,\; \mathsf{i}(v)=0  \text{ if }v\in Z\setminus Z_0\bigr\}
\]
for $j=0,1,\dotsc,n$;
note that $\left\lvert Z\cap Z_0\right\rvert$ has the same value for all $Z\in\Delta_j$ and that, if $G=\mathsf{Sym}(V)$, then $Z,Z'$ belong to the same $\mathsf{Stab}_G\left(Z_0\right)$-orbit (where $\mathsf{Stab}_G\left(Z_0\right)\cong\mathsf{Sym}_{\left\lvert Z_0\right\rvert}\times\mathsf{Sym}_{\left\lvert V\setminus Z_0\right\rvert}$)
if and only if $\left\lvert Z\cap Z_0\right\rvert=\left\lvert Z'\cap Z_0\right\rvert$. Then we suppose that $\Psi_j=\bigsqcup_{r=0}^{s_j}\Psi_{j,r}$ is the decomposition of $\Psi_j$ into the orbits of $\mathsf{Stab}_G\left(Z_0\right)$. If we set
\[
\digamma_{j,r}\coloneqq\bigsqcup_{(Z,\mathsf{i})\in\Psi_{j,r}} \left\{(Z,\phi)\mathrel{:}\phi\in Y^Z\text{ and } \phi(v)\in\Upsilon_{\mathsf{i}(v)}\text{ for all }v\in Z\cap Z_0\right\},
\]
then \eqref{decOrbxiphi} takes the form 
\[
\mathsf{Orb}_{F\wr G}\left(Z_0,\phi_{Z_0}\right)=\bigsqcup_{j=0}^n\bigsqcup_{r=0}^{s_j}\digamma_{j,r}.
\]

Finally, note that now \eqref{OrStabOr} is simply $\left\{Z_0\right\}\times Y^{Z_0}$; this set is $F\wr\mathsf{Stab}_G\left(Z_0\right)$-invariant, which yields the extension in Proposition \ref{propindact}\eqref{propindact1} (the inertia group in Proposition \ref{propindact}\eqref{propindact2} is clearly $F\wr\mathsf{Stab}_G\left(Z_0\right)$); $\beta$ and its inflation in Proposition \ref{propindact}\eqref{propindact3} are trivial while \eqref{OrStabOrStab} in the present setting coincides with \eqref{isomJohn}.
That is, the action of $F\wr G$ on $\mathsf{Orb}_{F\wr G}\left(Z_0,\phi_{Z_0}\right)$ is induced by the action of $F\wr\mathsf{Stab}_G\left(Z_0\right)$ on $\left\{Z_0\right\}\times Y^{Z_0}$; see Proposition \ref{propindact}\eqref{propindact4}.
\end{example}

\begin{example}\label{Example4}
Suppose that $X=V$, $\lvert A\rvert=1$ (so that $\alpha$ is trivial), and $Y=F$, with the left Cayley action. It follows that $\mathsf{Surj}(V,X)=\mathsf{Sym}(V)$ and if $\xi\in\mathsf{Sym}(V)$ and $g\in G$, then 
\[
\begin{split}
g\xi=\xi&\;\Longleftrightarrow\; \xi\left(g^{-1}v\right)=\xi(v)\text{ for all }v\in V\;\Longleftrightarrow\;g^{-1}v=\xi^{-1}(\xi(v))=v\text{ for all }v\in V\\
&\;\Longleftrightarrow\; g=1_G\quad (\text{by virtue of the faithfulness hypothesis})
\end{split}
\]
meaning that the action of $G$ on $\mathsf{Sym}(V)$ is \emph{free} (or \emph{semiregular}): for every $\xi\in\mathsf{Sym}(V)$, $\mathsf{Stab}_G(\xi)$ is trivial and
the action of $G$ on $\mathsf{Orb}_G(\xi)$ is isomorphic to its left Cayley action. Moreover, $\Phi_{\xi}=F^V$ for every $\xi$ and we may set $\varphi_{\xi}\coloneqq\mathbbm{1}_F$. For $(f,g)\in F\wr G$ we have $f\cdot g\varphi_{\xi}=f$ and therefore
\[
(f,g)\left(\xi,\varphi_{\xi}\right)=\left(g\xi,f\cdot g\varphi_{\xi}\right)=(g\xi,f).
\]
We can conclude that the action of $F\wr G$ on each orbit $\mathsf{Orb}_{F\wr G}\left(\xi_0,\varphi_{\xi_0}\right)$ is isomorphic to its left Cayley action. In particular, referring to Proposition \ref{propindact}, \eqref{OrStabOr} now becomes $\mathsf{Orb}_{F\wr G}\left(\xi_0,\varphi_{\xi_0}\right)=\mathsf{Orb}_G\left(\xi_0\right)\times\mathsf{Orb}_{F^V}\left(\varphi_{\xi_0}\right)\simeq G\times F^V$, $\beta$ and $\eta$ are the left Cayley actions of $G$ and $F^V$, respectively, the inertia group is the whole $F\wr G$, $\widetilde{\eta}=\gamma_0$ and there is no induction.
\end{example}

\begin{example}\label{Example5}
Now suppose that $X=A$ and $\alpha$ is the identity map. Then in a pair $(\xi,\varphi)\in\Gamma$ the second component reduces to a function $\varphi\colon V\rightarrow \bigsqcup_{a\in A}Y_a$ such that $\pi\circ\varphi$ is surjective (cf. \eqref{Diagram}), and $\xi$ is determined by $\varphi$: the unique $\xi$ such that $(\xi,\varphi)\in \Gamma$ is precisely $\xi\coloneqq\pi\circ\varphi$. Therefore, we may view $\Gamma$ simply as the set of all functions $\varphi$ as above and clearly $\gamma\simeq\gamma_0$. Moreover, we set $\Gamma_0\coloneqq\mathsf{Surj}\left(V,\left\{y_a\mathrel{:}a\in A\right\}\right)\subseteq\Gamma$ and if $\varphi\in\Gamma$,
then we define $\varphi^\sharp\in\Gamma_0$ by setting $\varphi^\sharp(v)\coloneqq y_a$ whenever $\varphi(v)\in Y_a$. Clearly, $\Gamma_0$ is $G$-invariant, $\varphi$ and $\varphi^\sharp$ determine the same $\xi$, we have $g\left(\varphi^\sharp\right)=(g\varphi)^\sharp$ when $g\in G$ and $\mathsf{Orb}_{F\wr G}(\varphi)=\mathsf{Orb}_{F\wr G}\left(\varphi^\sharp\right)$; in other words, in the present setting, $\varphi^\sharp$ takes the place of $\varphi_\xi$.

Now we fix $\varphi_0\in\Gamma_0$; by Lemma \ref{Lemmaorb}\eqref{Lemmaorb3}, $\mathsf{Orb}_{F\wr G}(\varphi)$ and $\mathsf{Orb}_{F\wr G}\left(\varphi_0\right)$ are isomorphic as $F\wr G$-spaces if and only if $\varphi^\sharp$ and $\varphi_0$ belong to the same $G$-orbit. The stabilizer \eqref{Lemmastab1} now is
\[
\mathsf{Stab}_{F\wr G}\left(\varphi_0\right)=\left\{(f,g)\in F\wr G\mathrel{:} g\in \mathsf{Stab}_G\left(\varphi_0\right)\; \text{\rm and } f(v)\in \mathsf{Stab}_F\left(\varphi_0(v)\right)\text{ for all }v\in V\right\}.
\]
Moreover, \eqref{decOrbxi} becomes $\mathsf{Orb}_G\left(\varphi_0\right)=\bigsqcup_{j=0}^n\Delta_j$,
the decomposition of $\mathsf{Orb}_G\left(\varphi_0\right)$ into the orbits of $\mathsf{Stab}_G\left(\varphi_0\right)$, with $\Delta_0=\{\varphi_0\}$. The decompositions \eqref{decYa} do not admit any particular simplification in this case, but  we use the notation $\ell_{y_a,y_b}$ in place of $\ell_{a,b}$. Since $\mathsf{Orb}_G\left(\varphi_0\right)\subseteq\Gamma_0$, \eqref{defPsixij} becomes
\[
\Psi_j\coloneqq\bigl\{(\varphi,\mathsf{i})\mathrel{:}\varphi\in\Delta_j \text{ and }\mathsf{i}\in\mathbb{N}^V,\text{ with }0\leq \mathsf{i}(v)\leq \ell_{\varphi(v),\varphi_0(v)}\text{ for all } v \in V\bigr\}. 
\]
Similarly, \eqref{decPsixij} may be written in the form $\Psi_j=\bigsqcup_{r=0}^{s_j}\Psi_{j,r}$,
and this is the decomposition of $\Psi_j$ into the orbits of $\mathsf{Stab}_G\left(\varphi_0\right)$.
It follows that, in the present setting, \eqref{decOrbxiphi} is
\begin{equation}\label{decOrbphi0}
\mathsf{Orb}_{F\wr G}\left(\varphi_0\right)=\bigsqcup_{j=0}^n\;\bigsqcup_{r=0}^{s_j}\digamma_{j,r},
\end{equation}
where
\[
\digamma_{j,r}\coloneqq\bigsqcup_{(\varphi,\mathsf{i})\in\Psi_{j,r}}\left(\bigtimes_{v\in V}\Upsilon_{\varphi(v),\varphi_0(v),\mathsf{i}(v)}\right),
\]
and it is the decomposition into the orbits of $\mathsf{Stab}_{F\wr G}\left(\varphi_0\right)$. 

Finally, regarding Proposition \ref{propindact}, in the present setting $\beta$ and its inflation are both trivial, the inertia group is $F\wr \mathsf{Stab}_G\left(\varphi_0\right)$, \eqref{OrStabOr} reduces to $\mathsf{Orb}_{F\wr\mathsf{Stab}_G\left(\varphi_0\right)}\left(\varphi_0\right)=\mathsf{Orb}_{F^V}\left(\varphi_0\right)$ (the first factor of its right-hand side is trivial). Moreover, the inductive step coincides with \eqref{FVorbdec}, which in the present setting is 
\begin{equation}\label{OrbFGphi0}
\mathsf{Orb}_{F\wr G}\left(\varphi_0\right)=\bigsqcup_{\varphi\in\mathsf{Orb}_G\left(\varphi_0\right)}\mathsf{Orb}_{F^V}(\varphi).
\end{equation}
Therefore, the action of $F\wr G$ on $\mathsf{Orb}_{F\wr G}\left(\varphi_0\right)$ is induced by the action of $F\wr\mathsf{Stab}_G\left(\varphi_0\right)$ on $\mathsf{Orb}_{F^V}\left(\varphi_0\right)$.
\end{example}

\begin{example}\label{Example6}
Suppose now that $X$ is an arbitrary set and $\lvert A\rvert=1$. Then $\Phi_\xi=Y^V$ for all $\xi\in\mathsf{Surj}(V,X)$ and we have the direct product decomposition
\begin{equation}\label{OrStabOrAt}
\mathsf{Orb}_{F\wr G}(\xi,\varphi)=\mathsf{Orb}_G(\xi)\times Y^V,
\end{equation}
where the second factor is endowed with the exponentiation action from Example \ref{Example1}. If $\varphi_0$ is as in that example and we fix $\xi_0\in\mathsf{Surj}(V,X)$, then $\mathsf{Stab}_{F\wr G}\left(\xi_0,\varphi_0\right)=\mathsf{Stab}_F\left(y_0\right)\wr\mathsf{Stab}_G\left(\xi_0\right)$. 
Using \eqref{decOrbxi}, we note that \eqref{defPsixij} now reduces to $\Delta_j\times\Psi$, where $\Psi$ is as in Example \ref{Example1}. Then \eqref{decPsixij} can be written as: $\Delta_j\times\Psi=\bigsqcup_{r=0}^{s_j}\Psi_{j,r}$, \eqref{decOrbxiphi} reduces to
\[
\mathsf{Orb}_G\left(\xi_0\right)\times Y^V=\bigsqcup_{j=0}^n\;\bigsqcup_{r=0}^{s_j}\digamma_{j,r},
\]
where
\[
\digamma_{j,r}\coloneqq\bigsqcup_{(\xi,\mathsf{i})\in\Psi_{j,r}}\left[\{\xi\}\times\left(\bigtimes\limits_{v\in V}\Upsilon_{\mathsf{i}(v)}\right)\right],
\]
which yields the decomposition into the orbits of $\mathsf{Stab}_F\left(y_0\right)\wr\mathsf{Stab}_G\left(\xi_0\right)$. Finally, regarding Proposition \ref{propindact}, in the present setting $\mathsf{Stab}_G(\alpha\circ\xi_0)=G$, $\beta$ is the action of $G$ on $\mathsf{Orb}_G\left(\xi_0\right)$, the inertia group is the whole $F\wr G$, \eqref{OrStabOr} coincides with \eqref{OrStabOrAt}, $\widetilde{\eta}$ is the exponentiation action, the induction is trivial, and therefore the action of $F\wr G$ on $\mathsf{Orb}_{F\wr G}\left(\xi_0,\varphi_0\right)$ is isomorphic to the diagonal action $\mathsf{diag}\left(\overline{\beta},\widetilde{\eta}\right)$.
\end{example}

\begin{example}\label{Example7}
Suppose now that $X=\left\{x_0\right\}\sqcup X'$, $A=\left\{a_0\right\}\sqcup A'$,  $\left\lvert Y_{a_0}\right\rvert=1$, while the $Y_a, a\in A'$ are arbitrary $F$-spaces, with $\alpha\left(x_0\right)=a_0$ and $\alpha(X')=A'$. Now we define $\mathsf{Per}_{X'}(V)$ as the set of all functions $\xi\in\mathsf{Surj}(V,X)$ such that: $\left\lvert\xi^{-1}(x)\right\rvert=1$ for all $x\in X'$. In other words, $\xi\rvert_{\xi^{-1}\left(X'\right)}$ is injective. If we use a numbering for the elements of $X$, say $X=\left\{x_0,x_1,\dotsc, x_k\right\}$, then $\mathsf{Per}_{X'}(V)\simeq\mathsf{Per}_k(V)$ (cf. \eqref{natproj})
via the map $\mathsf{Per}_{X'}(V)\ni\xi\mapsto\left(v_1,v_2,\dotsc,v_k\right)$, where $v_j$ is the unique element in $\xi^{-1}\left(x_j\right)$.
We prefer to use our coordinate-free notation and study the action of $F\wr G$ on homogeneous spaces of the form $\mathsf{Orb}_{F\wr G}\left(\xi,\varphi_\xi\right)$, where $\xi\in\mathsf{Per}_{X'}(V)$. Note that if $\varphi\in\Phi_\xi$ and $v\in(\alpha\circ\xi)^{-1}\left(a_0\right)$, then $\varphi(v)\in Y_{a_0}$. Moreover, if $a\in A'$, then $\left\lvert\left(\alpha\circ\xi\right)^{-1}(a)\right\rvert=\left\lvert\alpha^{-1}(a)\right\rvert$. Therefore, \eqref{FVprod} becomes $\Phi_\xi\simeq\bigtimes_{a\in A'}Y_a^{\left\lvert\alpha^{-1}(a)\right\rvert}$. 

Explicit formulations of \eqref{defdigjr} depend strongly on the group $G$ and on its action on $V$, so that a general analysis is omitted.
Finally, in Proposition \ref{propindact} we have 
\[
\mathsf{Stab}_G(\alpha\circ\xi)=\bigcap_{a\in A'}\mathsf{Stab}_G\left((\alpha\circ\xi)^{-1}(a)\right),
\]
and, in general, this group is larger than $\mathsf{Stab}_G(\xi)$ (for instance, when $G=\mathsf{Sym}(V)$ and $\lvert A'\rvert<\lvert X'\rvert$) and if this is the case then the action $\beta$ is non-trivial. Similarly, in general, $\mathsf{Stab}_G(\alpha\circ\xi)$ is smaller than $G$ (for instance, when $G=\mathsf{Sym}(V)$) and if this is the case then the induction in Proposition \ref{propindact} is non-trivial.
\end{example}

\begin{remark}\label{Remark14}
Note that Example \ref{Example2} is a special case of Example \ref{Example3} (when $\lvert Z_0\rvert=1$) but it deserves a separate treatment (cf. the applications in Example \ref{Example21} and in \eqref{DecLSl1h}). Similarly, Example \ref{Example1} might be considered a special case of Example \ref{Example3} (when $Z_0=V$) but requiring that $Z_0\subsetneqq V$ we have excluded it, since the independent analysis of the case $Z_0=V$ is essential. Finally, if $\lvert A'\rvert=1$ in Example \ref{Example7}, then \eqref{natproj} yields a natural projection onto the homogeneous space in Example \ref{Example3}, given by $(\xi,\varphi)\mapsto\left(\xi^{-1}(X'),\varphi\rvert_{\xi^{-1}(X')}\right)$.
\end{remark}

\begin{remark}\label{Remark15}
In the general case, the actions $\gamma$ and $\gamma_0$ are both imprimitive. Indeed, 
\[
\mathsf{Orb}_{F\wr G}\left(\xi_0,\varphi_{\xi_0}\right)=\bigsqcup_{\xi\in\mathsf{Orb}_G\left(\xi_0\right)}\left(\{\xi\}\times\mathsf{Orb}_{F^V}\left(\varphi_\xi\right)\right),\quad\mathsf{Orb}_{F\wr G}\left(\varphi_{\xi_0}\right)=\!\!\!\bigsqcup_{\alpha\circ\xi\in\mathsf{Orb}_G\left(\alpha\circ\xi_0\right)}\!\!\!\mathsf{Orb}_{F^V}\left(\varphi_\xi\right), 
\]
are decompositions into $\gamma$-blocks and $\gamma_0$-blocks, respectively. These blocks are trivial in three cases: 

\medskip

(1) when $\left\lvert\mathsf{Orb}_{F^V}\left(\varphi_\xi\right)\right\rvert=1$, which holds in the trivial case $\left\lvert Y_a\right\rvert=1$ for every $a\in A$; 

\medskip

(2) for $\gamma$, when $\left\lvert\mathsf{Orb}_G\left(\xi_0\right)\right\rvert=1$, which holds if $\lvert X\rvert=1$, that is, precisely in the setting of Example \ref{Example1};

\medskip

(3) for $\gamma_0$, when $\left\lvert\mathsf{Orb}_G\left(\alpha\circ\xi_0\right)\right\rvert=1$, which holds if $\lvert A\rvert=1$, that is, precisely in the setting of Example \ref{Example6}.

\medskip

In case (2), there is a characterization of the primitive actions: with the notation in Example \ref{Example1}, assuming that the actions of $G$ on $V$ and of $F$ on $Y$ are both transitive and faithful, then the action of $F\wr G$ on $Y^V$ is primitive if and only if the action of $F$ on $Y$ is primitive but not regular; see \cite[Theorem 2.2.6]{book3} and \cite{DM}.
\end{remark}

\begin{example}\label{Exampleequiv}
Now we suppose that $G=\mathsf{C}_n$, the cyclic group of order $n$. We write it multiplicatively and we denote by $\epsilon$ its generator. Moreover, now $V$ is the ring $\mathbb{Z}_n$ of the integers modulo $n$ and the action of $\mathsf{C}_n$ on $\mathbb{Z}_n$ is given by $\epsilon^k(m)\coloneqq m+k\pmod n$ for all $m,k\in\mathbb{Z}_n$.
We denote by $\mathsf{D}(n)$ the set of all divisors of $n$ (positive integers that divide $n$); for each $k\in\mathsf{D}(n)$ there exists exactly one subgroup of $\mathsf{C}_n$ of order $k$, it is $\mathsf{C}_k\coloneqq\left\langle\epsilon^{n/k}\right\rangle$ and these constitute  all the subgroups of $\mathsf{C}_n$; see, e.g., \cite[Chapter 1]{book4}.

For $k\in\mathsf{D}(n)$, a function $\xi\in\mathsf{Surj}\left(\mathbb{Z}_n,X\right)$ is $k$-\emph{periodic} if it is $\mathsf{C}_{n/k}$-invariant, that is $\xi(m+k)=\xi(m)$ for all $m\in\mathbb{Z}_n$. The \emph{period} of a $\xi\in\mathsf{Surj}\left(\mathbb{Z}_n,X\right)$ is the smallest $k\in\mathsf{D}(n)$ such that $\xi$ is $k$-periodic. If $\xi\in\mathsf{Surj}\left(\mathbb{Z}_n,X\right)$ has no such period, that is $\mathrm{Stab}_{C_n}(\xi) = \{1\}$, it is called \emph{aperiodic}.
Then the orbits $\mathsf{Orb}_{F\wr\mathsf{C}_n}(\xi,\varphi)$ and $\mathsf{Orb}_{F\wr\mathsf{C}_n}(\xi',\varphi')$ are isomorphic as $F\wr\mathsf{C}_n$-spaces if and only if $\xi$ and $\xi'$ have the same period and $\alpha\circ\xi$ and $\alpha\circ\xi'$ determine the same pattern (i.e., they belong to the same $\mathsf{C}_n$-orbit). This follows from Lemma \ref{Lemmaorb}\eqref{Lemmaorb3}, just taking into account that two subgroups of $\mathsf{C}_n$ are conjugate if and only if they coincide.

Now suppose that $n\geq 5$, $X=\left\{x_0,x_1,x_2\right\}$, $A=\left\{a_0\right\}$ and $\alpha\left(x_i\right)=a_0$ for $i=0,1,2$. Define $\xi,\xi'\in\mathsf{Surj}\left(\mathbb{Z}_n,X\right)$ by requiring that
\begin{equation}\label{formExampleequiv}
\begin{split}
&\left\lvert\xi^{-1}\left(x_1\right)\right\rvert=\left\lvert\xi^{-1}\left(x_2\right)\right\rvert=1;\qquad
\left\lvert\xi'^{-1}\left(x_1\right)\right\rvert=2\quad\text{and}\quad\left\lvert\xi'^{-1}\left(x_2\right)\right\rvert=1,\\
&\text{ so that}\quad\left\lvert\xi^{-1}\left(x_0\right)\right\rvert=n-2\quad\text{and}\quad\left\lvert\xi'^{-1}\left(x_0\right)\right\rvert=n-3.
\end{split}
\end{equation}
Then both $\xi$ and $\xi'$ are aperiodic (because $\left\lvert\xi^{-1}\left(x_2\right)\right\rvert=\left\lvert\xi'^{-1}\left(x_2\right)\right\rvert=1$). Moreover $\alpha\circ\xi=\alpha\circ\xi'$, so that the conditions for the isomorphism of orbits described above are verified. On the other hand, by \eqref{formExampleequiv}, $\xi$ and $\xi'$ are \emph{not} equivalent, because the multisets of fiber cardinalities $\{n-2, 1, 1\}$ and $\{n-3, 2, 1\}$ are distinct for $n \ge 5$ (which also explains the exact threshold $n \ge 5$). 
\end{example}

\begin{sidewaystable}[p]
\centering
\caption{Summary of the examples in Section \ref{Sechomspace}, classifying the data $(X,A,\alpha,\{Y_a\}_{a\in A})$, the resulting homogeneous space, the inertia group from Proposition \ref{propindact}, and the nontriviality of $\beta$ and of the induction step.}
\label{tab:sec3-examples}
\scriptsize
\renewcommand{\arraystretch}{1.28}
\setlength{\tabcolsep}{3.5pt}
\begin{tabularx}{\textwidth}{%
  >{\raggedright\arraybackslash}p{2.3cm}%
  >{\raggedright\arraybackslash}p{2.0cm}%
  >{\raggedright\arraybackslash}p{1.9cm}%
  >{\raggedright\arraybackslash}p{2.5cm}%
  >{\raggedright\arraybackslash\hsize=0.80\hsize}X%
  >{\raggedright\arraybackslash\hsize=1.40\hsize}X%
  >{\raggedright\arraybackslash\hsize=0.85\hsize}X%
  >{\raggedright\arraybackslash\hsize=0.95\hsize}X%
}
\toprule
\textbf{Example} &
\textbf{$X$} &
\textbf{$A$} &
\textbf{$\alpha$} &
\textbf{$Y_a$} &
\textbf{Homogeneous space} &
\textbf{Inertia group (Prop. \ref{propindact})} &
\textbf{$\beta$ / induction} \\
\midrule

\textbf{\ref{Example1}} (Exponentiation) &
$|X|=1$ &
$|A|=1$ &
trivial &
single arbitrary $Y$ &
$Y^V \simeq (F \wr G)/(\Stab_F(y_0) \wr G)$ ($\gamma \simeq \gamma_0$, transitive) &
$F \wr G$ &
trivial $\beta$ / trivial induction \\
\addlinespace

\textbf{\ref{Example2}} (Composition) &
$\{x_0, x_1\}$ &
$\{a_0, a_1\}$ &
$\alpha(x_i)=a_i$ (bijective) &
$Y_{a_0}=Y$ arbitrary, $|Y_{a_1}|=1$ &
$\Gamma_0 \simeq V \times Y$ (with $\xi_v(v)=x_0$) &
$F \wr \Stab_G(v_0)$ &
trivial $\beta$ / nontrivial induction \\
\addlinespace

\textbf{\ref{Example21}} ($\Per_2(V)\times Y$) &
$\{x_0, x_1, x_2\}$ ($|V|\ge 3$) &
$\{a_0, a_1\}$ &
$\alpha(x_0)=a_0$, $\alpha(x_1)=\alpha(x_2)=a_1$ &
$Y_{a_0}=Y$ arbitrary, $|Y_{a_1}|=1$ &
$\Orb_{F\wr G}(v_0,v_j,y_0) \subseteq \Per_2(V)\times Y$ &
$F \wr \Stab_G(v_0)$ &
nontrivial $\beta$ (on $\{v_0\}\times\Delta_{v_0,j}$, if $|\Delta_{v_0,j}|>1$) / nontrivial induction \\
\addlinespace

\textbf{\ref{Example3}} (Subsets $Z_0$) &
$\{x_0, x_1\}$ &
$\{a_0, a_1\}$ &
$\alpha(x_i)=a_i$ (bijective) &
$Y_{a_0}=Y$ arbitrary, $|Y_{a_1}|=1$ &
$\bigsqcup_{Z\in\Orb_G(Z_0)}(\{Z\}\times Y^Z)$ ($\emptyset\subsetneq Z_0\subsetneq V$) &
$F \wr \Stab_G(Z_0)$ &
trivial $\beta$ / nontrivial induction \\
\addlinespace

\textbf{\ref{Example4}} (Left Cayley) &
$V$ &
$|A|=1$ &
trivial &
single $Y=F$ (left Cayley) &
$\Orb_{F\wr G}(\xi_0,\varphi_{\xi_0}) \simeq G \times F^V \simeq F \wr G$ ($\xi_0\in\Sym(V)$) &
$F \wr G$ &
nontrivial $\beta$ (Cayley of $G$, if $|G|>1$) / trivial induction \\
\addlinespace

\textbf{\ref{Example5}} (Case $\alpha=\id$) &
$X=A$ &
arbitrary ($1\le |A|\le |V|$) &
$\id_A$ &
$\{Y_a\}_{a\in A}$ arbitrary, pairwise non-isomorphic &
$\Orb_{F\wr G}(\varphi_0)$ with $\varphi_0\in\Surj(V,\{y_a\})$ ($\gamma\simeq\gamma_0$) &
$F \wr \Stab_G(\varphi_0)$ &
trivial $\beta$ / nontrivial induction (for $|A|\ge 2$) \\
\addlinespace

\textbf{\ref{Example6}} (Diagonal / product) &
arbitrary ($1\le |X|\le |V|$) &
$|A|=1$ &
trivial &
single arbitrary $Y$ &
$\Orb_G(\xi_0) \times Y^V$ (action $\diag(\overline{\beta},\tilde{\eta})$) &
$F \wr G$ &
nontrivial $\beta$ (on $\Orb_G(\xi_0)$, if $|X|\ge 2$) / trivial induction \\
\addlinespace

\textbf{\ref{Example7}} ($\Per_{X'}(V)$) &
$\{x_0\}\sqcup X'$ &
$\{a_0\}\sqcup A'$ &
$\alpha(x_0)=a_0$, $\alpha(X')=A'$ &
$|Y_{a_0}|=1$, $\{Y_a\}_{a\in A'}$ arbitrary &
$\Orb_{F\wr G}(\xi,\varphi_\xi)$ with $\xi\in\Per_{X'}(V)$ ($\Phi_\xi\simeq\times_{a\in A'}Y_a^{|\alpha^{-1}(a)|}$) &
$F \wr \bigcap_{a\in A'}\Stab_G((\alpha\circ\xi)^{-1}(a))$ &
both $\beta$ and induction nontrivial in general (e.g., $G=\Sym(V)$, $|A'|<|X'|$) \\
\addlinespace

\textbf{\ref{Exampleequiv}} ($G=C_n$ on $\mathbb{Z}_n$, gen.) &
arbitrary &
arbitrary &
arbitrary surjection &
$\{Y_a\}_{a\in A}$ arbitrary &
$\Orb_{F\wr C_n}(\xi,\varphi)$ (classified by period of $\xi$ and $C_n$-orbit of $\alpha\circ\xi$) &
$F \wr \Stab_{C_n}(\alpha\circ\xi)$ ($C_{n/k'}$ for period $k'$, or $\{1\}$) &
both nontrivial in general ($\beta$ trivial iff $\Stab_{C_n}(\xi)=\Stab_{C_n}(\alpha\circ\xi)$; ind.\ trivial iff $\alpha\circ\xi$ constant) \\
\addlinespace

\textbf{\ref{Exampleequiv}} (counter\-example, $n\ge 5$) &
$\{x_0, x_1, x_2\}$ &
$\{a_0\}$ &
constant &
single arbitrary $Y$ &
$\Orb_{F\wr C_n}(\xi,\varphi)\simeq\Orb_{F\wr C_n}(\xi',\varphi')\simeq C_n\times Y^{\mathbb{Z}_n}$ ($\xi,\xi'$ aperiodic, fibre sizes $\{n-2,1,1\}$ and $\{n-3,2,1\}$, inequivalent) &
$F \wr C_n$ &
nontrivial $\beta$ (Cayley of $C_n$) / trivial induction \\

\bottomrule
\end{tabularx}
\end{sidewaystable}

\begin{remark}[Structural overview of the examples]\label{rem:sec3-summary}
Table~\ref{tab:sec3-examples} summarizes the examples in this section and illustrates all four possible regimes arising from Proposition \ref{propindact} with respect to the nontriviality of $\beta$ (which occurs precisely when $\Stab_G(\xi_0) \subsetneq \Stab_G(\alpha\circ\xi_0)$, and hence requires $\alpha$ to be non-injective) and of the induction step (which occurs when $\Stab_G(\alpha\circ\xi_0) \subsetneq G$, i.e., when $|A|\ge 2$):
\begin{itemize}
  \item \emph{Trivial $\beta$, trivial induction:} Example \ref{Example1} (the prototypical case with $|X|=|A|=1$).
	  \item \emph{Trivial $\beta$, nontrivial induction:} Examples \ref{Example2}, \ref{Example3}, and \ref{Example5} (all having $\alpha$ bijective and $|A|\ge 2$).
  \item \emph{Nontrivial $\beta$, trivial induction:} Examples \ref{Example4} and \ref{Example6}, as well as the counterexample in Example \ref{Exampleequiv} (all having $|A|=1$ and $|X|\ge 2$).
  \item \emph{Nontrivial $\beta$, nontrivial induction:} Examples \ref{Example21} and \ref{Example7}, together with the general case of Example \ref{Exampleequiv} ($|A|\ge 2$ and $\alpha$ non-injective).
\end{itemize}
Finally, Remarks \ref{Remark14} and \ref{Remark15} serve as structural bridges across Table~\ref{tab:sec3-examples}: Remark \ref{Remark14} highlights the hierarchical relationships among the examples (Example \ref{Example2} is the $|Z_0|=1$ specialization of Example \ref{Example3}, Example \ref{Example1} formally corresponds to $Z_0=V$, and for $|A'|=1$ Example \ref{Example7} admits a natural projection onto Example \ref{Example3}), whereas Remark \ref{Remark15} characterizes the imprimitivity of $\gamma$ and $\gamma_0$, isolating the cases with trivial block systems (in particular, Example \ref{Example1} for $\gamma$ and Example \ref{Example6} for $\gamma_0$).
\end{remark}

\section{The permutation representations}\label{Secpermrep}

In the present section, we again fix an element $\xi_0\in\mathsf{Surj}(V,X)$ and study the permutation representation $\lambda$ of $F\wr G$ on $L\bigl(\mathsf{Orb}_{F\wr G}\left(\xi_0,\varphi_{\xi_0}\right)\bigr)$. For future reference, it is convenient to express $\lambda$ explicitly: if $(f,g)\in F\wr G$, $H\in L\bigl(\mathsf{Orb}_{F\wr G}\left(\xi_0,\varphi_{\xi_0}\right)\bigr)$ and $(\xi,\varphi)\in\mathsf{Orb}_{F\wr G}\left(\xi_0,\varphi_{\xi_0}\right)$, then
\begin{equation} \label{explpermrep}
\left[\lambda(f,g)H\right](\xi,\varphi)=H\left(g^{-1}\xi,g^{-1}\left(f^{-1}\cdot\varphi\right)\right).
\end{equation}
First, we provide an intermediate result, which is a straightforward consequence of \eqref{propindact3} and \eqref{propindact4} in Proposition \ref{propindact}.

\begin{proposition}\label{Propindinflext}
We have the following equivalence of $F\wr \mathsf{Stab}_G\left(\alpha\circ\xi_0\right)$-permutation representations
\[
L\bigl(\mathsf{Orb}_{F\wr \mathsf{Stab}_G(\alpha\circ\xi_0)}\left(\xi_0,\varphi_{\xi_0}\right)\bigr)\sim\overline{L}\Bigl(\mathsf{Orb}_{\mathsf{Stab}_G\left(\alpha\circ\xi_0\right)}\left(\xi_0\right)\Bigr)\bigotimes\widetilde{L}\Bigl(\mathsf{Orb}_{F^V}\left(\varphi_{\xi_0}\right)\Bigr),
\]
where $\overline{L}$ denotes the inflation from $\mathsf{Stab}_G(\alpha\circ\xi_0)$ to $F\wr\mathsf{Stab}_G(\alpha\circ\xi_0)$, while $\widetilde{L}$ denotes the extension from $F^V$ to $F\wr \mathsf{Stab}_G(\alpha\circ\xi_0)$ associated with $\tilde{\eta}$ in Proposition \ref{propindact}\eqref{propindact1}. Moreover, the permutation representation of $F\wr G$ on $L\bigl(\mathsf{Orb}_{F\wr G}\left(\xi_0,\varphi_{\xi_0}\right)\bigr)$ is equivalent to
\[
\Ind_{F\wr\mathsf{Stab}_G(\alpha\circ\xi_0)}^{F\wr G}\left[\overline{L}\Bigl(\mathsf{Orb}_{\mathsf{Stab}_G(\alpha\circ\xi_0)}(\xi_0)\Bigr)\bigotimes\widetilde{L}\Bigl(\mathsf{Orb}_{F^V}\left(\varphi_{\xi_0}\right)\Bigr)\right].
\]
\end{proposition}

From now on, we use all the notions and notation of the representation theory of wreath products established in Section \ref{Secprelimin}.
For every $a\in A$, let $\widehat{F}_a$ be the subset of $\widehat{F}$ consisting of all irreducible representations contained in the permutation representation of $F$ on $L\left(Y_a\right)$; for $\sigma_*\in\widehat{F}_a$, we denote by $m_{a,\sigma_*}$ the multiplicity of $U_{\sigma_*}$ in $L\left(Y_a\right)$. Then, for every $a\in A$,
\begin{equation}\label{decLYa}
L\left(Y_a\right)=\bigoplus_{\sigma_*\in\widehat{F}_a}m_{a,\sigma_*}U_{\sigma_*}
\end{equation}
is the decomposition of $L\left(Y_a\right)$ into isotypic components. If $\xi\in\mathsf{Surj}(V,X)$ and $\sigma=\bigotimes_{v\in V}\sigma_v\in\widehat{F^V}$, then we write $\sigma\prec\xi$ to denote that $\sigma_v\in\widehat{F}_{\alpha(\xi(v))}$ for every $v\in V$; from \eqref{defgsigma} it follows that
 if $\sigma\prec\xi$ and $g\in G$, then $\!^g\sigma\prec g\xi$.
From \eqref{defUsigma}, \eqref{FVprod} and \eqref{decLYa} we deduce that the decomposition of $L\bigl(\mathsf{Orb}_{F^V}\left(\varphi_{\xi}\right)\bigr)$ into isotypic $F^V$-components is
\begin{equation}\label{decYV}
L\bigl(\mathsf{Orb}_{F^V}\left(\varphi_{\xi}\right)\bigr)=\bigotimes_{v\in V} L\bigl(Y_{\alpha(\xi(v))}\bigr)=\bigoplus_{\sigma\prec\xi}\left(\prod_{v\in V}m_{\alpha(\xi(v)),\sigma_v}\right)U_\sigma.
\end{equation}

In what follows, $\lambda$ will be used to denote the left action of groups such as $F\wr G$, $F^V$ or $F$ on the various function spaces or tensor products involved; the specific group, space and action will be clear from the context.
Moreover, from \eqref{FVorbdec} we deduce the following decomposition into $F^V$-invariant permutation representations
\begin{equation}\label{decFVinv}
L\bigl(\mathsf{Orb}_{F\wr G}\left(\xi_0,\varphi_{\xi_0}\right)\bigr)=\bigoplus_{\xi\in\mathsf{Orb}_G\left(\xi_0\right)}L\bigl(\mathsf{Orb}_{F^V}\left(\varphi_{\xi}\right)\bigr).
\end{equation}
\begin{proposition}\label{Propacctens2}
Assume that $h_v\in L\left(Y_{\alpha(\xi(v))}\right)$, so that $\bigotimes_{v\in V}h_v\in L\bigl(\mathsf{Orb}_{F^V}\left(\varphi_{\xi}\right)\bigr)$ (cf. \eqref{decYV}), and let $(f,g)\in F\wr G$. Then
\begin{equation}\label{acttens2}
\lambda(f,g)\left(\bigotimes_{v\in V}h_v\right)=\bigotimes_{v\in V}\lambda(f(v))h_{g^{-1}v}\in L\bigl(\mathsf{Orb}_{F^V}\left(\varphi_{g\xi}\right)\bigr).
\end{equation}
\end{proposition}
\begin{proof}
First, note that if $\varphi\in\mathsf{Orb}_{F^V}\left(\varphi_{\xi}\right)$, then
\begin{equation}\label{hvtens}
\left(\bigotimes_{v\in V}h_v\right)(\varphi)=\prod_{v\in V}h_v(\varphi(v)).
\end{equation}
It follows that, for every $\varphi \in \mathsf{Orb}_{F^V}(\varphi_{g\xi})$ we have
\begin{align*}
\left[\lambda(f,g)\bigotimes_{v\in V}h_v\right](\varphi)=&\left[\bigotimes_{v\in V}h_v\right]\left(g^{-1}f^{-1}\cdot g^{-1}\varphi\right)&&(\text{by }\eqref{explpermrep})\\
=&\prod_{v\in V}h_v\left(f(gv)^{-1}\varphi(gv)\right)&&(\text{by }\eqref{hvtens})\\
=&\left[\bigotimes_{v\in V}\lambda(f(v))h_{g^{-1}v}\right](\varphi).&&(\text{by }v\rightarrow g^{-1}v\text{ and }\eqref{hvtens})
\end{align*}
\end{proof}

Now, following the notation in \eqref{decLYa} and \eqref{decYV}, for each $a\in A$ and $\sigma_*\in\widehat{F}_a$ we fix a basis $T_{a,\sigma_*,j}$, $1\leq j\leq m_{a,\sigma_*}$, for the space $\Hom_F\left(U_{\sigma_*},L\left(Y_a\right)\right)$. This means that, for every $f\in F$ and $u\in U_{\sigma_*}$ we have
\begin{equation}\label{decULFVbis}
\lambda(f)T_{a,\sigma_*,j}u=T_{a,\sigma_*,j}\sigma_*(f)u.
\end{equation}
If $\sigma\prec\xi$, then we denote by $\mathsf{J}_{\xi,\sigma}$ the set of all functions $\mathsf{j}\in\mathbb{N}^V$ such that
\begin{equation}\label{condj} 
1\leq\mathsf{j}(v)\leq m_{\alpha(\xi(v)),\sigma_v}
\end{equation}
for every $v\in V$.
In particular, if $k\in\mathsf{I}(\sigma)$, then
\begin{equation}\label{condjbis} 
k\mathsf{J}_{\xi,\sigma}=\mathsf{J}_{k\xi,\sigma},
\end{equation}  
because $1\leq[k\mathsf{j}](v)=\mathsf{j}\left(k^{-1}v\right)\leq m_{\alpha\left(\xi\left(k^{-1}v\right)\right),\sigma_v}= m_{\alpha((k\xi)(v)), \sigma_v}$ (by $\sigma_{k^{-1}v} = \sigma_v$ for $k \in \mathsf{I}(\sigma)$).
We define a linear operator $T_{\xi,\sigma,\mathsf{j}}\colon U_\sigma\rightarrow L\left(\mathsf{Orb}_{F^V}\left(\varphi_{\xi}\right)\right)$ (cf. \eqref{decYV}) by setting
\begin{equation}\label{defTxisigj} 
T_{\xi,\sigma,\mathsf{j}}\coloneqq\bigotimes_{v\in V}T_{\alpha(\xi(v)),\sigma_v,\mathsf{j}(v)}.
\end{equation}
Then $\left\{T_{\xi,\sigma,\mathsf{j}}\mathrel{:}\mathsf{j}\in \mathsf{J}_{\xi,\sigma}\right\}$ is a basis for $\Hom_{F^V}\left(U_\sigma,L\bigl(\mathsf{Orb}_{F^V}\left(\varphi_{\xi}\right)\bigr)\right)$. In other words,
\begin{equation}\label{decULFV}
\left(\prod_{v\in V}m_{\alpha(\xi(v)),\sigma_v}\right)U_\sigma=\bigoplus_{\mathsf{j}\in \mathsf{J}_{\xi,\sigma}}T_{\xi,\sigma,\mathsf{j}}U_\sigma
\end{equation}
is a decomposition of the $\sigma$-isotypic component in $L\bigl(\mathsf{Orb}_{F^V}\left(\varphi_{\xi}\right)\bigr)$; cf. \eqref{decYV}. Henceforth, for any $\sigma\in\widehat{F^V}$ such that $U_\sigma$ is contained in one of the summands on the right-hand side of \eqref{decFVinv} (meaning that $\sigma\prec\xi$ for some $\xi\in\mathsf{Orb}_G\left(\xi_0\right)$) we set 
\begin{equation}\label{defJsigma}
\mathsf{J}_\sigma\coloneqq\left\{(\xi,\mathsf{j})\mathrel{:}\xi\in\mathsf{Orb}_G\left(\xi_0\right),\xi\succ\sigma,\mathsf{j}\in \mathsf{J}_{\xi,\sigma}\right\}.
\end{equation}
It follows that $\left\{T_{\xi,\sigma,\mathsf{j}}\mathrel{:}(\xi,\mathsf{j})\in \mathsf{J}_\sigma\right\}$ is a basis for $\Hom_{F^V}\left(U_\sigma,L\bigl(\mathsf{Orb}_{F\wr G}\left(\xi_0,\varphi_{\xi_0}\right)\bigr)\right)$. 
Moreover, the set $\mathsf{J}_\sigma$ is $\mathsf{I}(\sigma)$-invariant: if  $(\xi,\mathsf{j})\in\mathsf{J}_\sigma$ and $k\in \mathsf{I}(\sigma)$, then $\sigma=\!^k\sigma$, so that $k\xi\in\mathsf{Orb}_G\left(\xi_0\right)$, $\sigma\prec k\xi$ and $k\mathsf{j}\in\mathsf{J}_{k\xi,\sigma}$ (because if $\mathsf{j}$ satisfies \eqref{condj}, then $1\leq k\mathsf{j}(v)\leq m_{\alpha(k\xi(v)),\sigma_v}$). Therefore, $(k\xi,k\mathsf{j})\in\mathsf{J}_\sigma$. In particular, 
\begin{equation}\label{defJsig0old}
\mathsf{J}_\sigma\supseteq\left\{(\xi,\mathsf{j})\mathrel{:}\xi\in\mathsf{Orb}_{\mathsf{I}(\sigma)}(\xi_0),\xi\succ\sigma,\mathsf{j}\in \mathsf{J}_{\xi,\sigma}\right\}.
\end{equation} 
We also need to consider the permutation representation of $\mathsf{I}(\sigma)$ on $L\left(\mathsf{J}_\sigma\right)$, and we denote by $\widehat{\mathsf{I}(\sigma)}_0$ the set of all $\theta\in\widehat{\mathsf{I}(\sigma)}$ contained in $L\left(\mathsf{J}_\sigma\right)$. 

\begin{remark}\label{Remark}
\begin{enumerate}[(1)]
\item\label{Remark1}
In general, the action of $\mathsf{I}(\sigma)$ on $\mathsf{J}_\sigma$ is not transitive. For instance, if $m_{a,\sigma_*}>1$ for some $a\in A$ (cf. \eqref{decLYa}) and $\sigma_*$ is one of the factors of $\sigma$, then the set of all $(\xi,\mathsf{j})\in\mathsf{J}_\sigma$ such that $\mathsf{j}(v)=1$ for every $v\in V$ constitutes a proper orbit.
\item The inclusion in \eqref{defJsig0old} may be proper; see Example \ref{Example4bis}.
\end{enumerate}
\end{remark}

Now we study the action of $F\wr G$ and of $F\wr \mathsf{I}(\sigma)$ on the decompositions \eqref{decULFV}.
\begin{proposition}
For $(f,g)\in F\wr G$ and $(\xi,\mathsf{j})\in \mathsf{J}_\sigma$,
we have
\begin{equation}\label{lambsigm2}
\lambda(f,g)T_{\xi,\sigma,\mathsf{j}}=T_{g\xi,\;^g\sigma,g\mathsf{j}}\;^g\sigma(f)\lambda\left(\mathbbm{1}_F,g\right).
\end{equation}
In particular, if $(f,k)\in F\wr \mathsf{I}(\sigma)$, then 
\begin{equation}\label{lambsigm}
\lambda(f,k)T_{\xi,\sigma,\mathsf{j}}=T_{k\xi,\sigma,k\mathsf{j}}\widetilde{\sigma}(f,k).
\end{equation}
\end{proposition}
\begin{proof}
First, note that if $\bigotimes_{v\in V}u_v\in U_\sigma$ and $\varphi\in\mathsf{Orb}_{F^V}(\varphi_\xi)$, then from \eqref{defTxisigj} it follows that
\begin{equation}\label{Txiuv}
\left[T_{\xi,\sigma,\mathsf{j}}\left(\bigotimes_{v\in V}u_v\right)\right](\varphi)=\prod_{v\in V}\left[T_{\alpha(\xi(v)),\sigma_v,\mathsf{j}(v)}u_v\right](\varphi(v)).
\end{equation}
Therefore, we obtain
\begin{align*}
\lambda(f,g)T_{\xi,\sigma,\mathsf{j}}\left(\bigotimes_{v\in V}u_v\right)=&\lambda(f,g)\bigotimes_{v\in V}T_{\alpha(\xi(v)),\sigma_v,\mathsf{j}(v)}u_v&&(\text{by }\eqref{Txiuv})\\
=&\bigotimes_{v\in V}\lambda(f(v))T_{\alpha(g\xi(v)),\sigma_{g^{-1}v},g\mathsf{j}(v)}u_{g^{-1}v}&&(\text{by }\eqref{acttens2})\\
=&\bigotimes_{v\in V}T_{\alpha(g\xi(v)),\sigma_{g^{-1}v},g\mathsf{j}(v)}\sigma_{g^{-1}v}(f(v))u_{g^{-1}v}&&(\text{by }\eqref{decULFVbis})\\
=&\left(\bigotimes_{v\in V}T_{\alpha(g\xi(v)),\sigma_{g^{-1}v},g\mathsf{j}(v)}\right)\left(\bigotimes_{v\in V}\sigma_{g^{-1}v}(f(v))u_{g^{-1}v}\right)&&\\
=&T_{g\xi,\;^g\sigma,g\mathsf{j}}\;^g\sigma(f)\lambda\left(\mathbbm{1}_F,g\right)\left(\bigotimes_{v\in V}u_v\right).&&(\text{by }\eqref{actgsig})
\end{align*}
Then \eqref{lambsigm} follows from \eqref{defextsig}.
\end{proof}

From \eqref{actgsig} and \eqref{lambsigm2} we deduce that
\begin{equation}\label{lambUsigm}
\lambda(f,g)T_{\xi,\sigma,\mathsf{j}}U_\sigma=T_{g\xi,\;^g\sigma,g\mathsf{j}}\;U_{\;^g\sigma}.
\end{equation}

We define $\Sigma_0$ as the subset of $\Sigma$ in \eqref{irrepwreath} consisting of the representatives of the conjugacy classes of those $\sigma\in\widehat{F^V}$ such that $U_\sigma$ is contained in $\Res_{F^V}^{F\wr G}L\bigl(\mathsf{Orb}_{F\wr G}\left(\xi_0,\varphi_{\xi_0}\right)\bigr)$. By \eqref{lambUsigm}, we may assume that $\Sigma_0\subseteq \bigl\{\sigma\in\widehat{F^V}\mathrel{:}\sigma\prec\xi_0\bigr\}$. With this assumption, if $\sigma\in\Sigma_0$, then the condition $\xi\succ\sigma$ on the right-hand side of \eqref{defJsig0old} may be omitted.

We are now ready to prove the fundamental results of this section. Applying Clifford theory (cf. \eqref{irrepwreath}), we show how to obtain
the decomposition of the permutation representation of $F\wr G$ on $L\bigl(\mathsf{Orb}_{F\wr G}\left(\xi_0,\varphi_{\xi_0}\right)\bigr)$ from the decomposition of $L\left(\mathsf{J}_\sigma\right)$ into irreducible $\mathsf{I}(\sigma)$-representations. We recall that $N\in\Hom_{\mathsf{I}(\sigma)}\bigl(W_\theta,L\left(\mathsf{J}_\sigma\right)\bigr)$ means that 
\begin{equation}\label{Nthetak}
[N\theta(k)w](\xi,\mathsf{j})=[Nw]\left(k^{-1}\xi,k^{-1}\mathsf{j}\right)
\end{equation}
for all $w\in W_\theta$ and $k\in\mathsf{I}(\sigma)$.

\begin{lemma}\label{Theointoper}
Assume the notation in \eqref{irrepwreath} and \eqref{defHquw}, and suppose that $\sigma\in\Sigma_0$, $\theta\in\widehat{\mathsf{I}(\sigma)}_0$ and $N\in\Hom_{\mathsf{I}(\sigma)}\bigl(W_\theta,L\left(\mathsf{J}_\sigma\right)\bigr)$. For every $H_{q,u\otimes w}$ define $ N^\flat\left(H_{q,u\otimes w}\right)\in L\bigl(\mathsf{Orb}_{F\wr G}\left(\xi_0,\varphi_{\xi_0}\right)\bigr)$ by setting
\begin{equation}\label{defNflat}
\left[N^\flat\left(H_{q,u\otimes w}\right)\right](\xi,\varphi)\coloneqq
\begin{cases}
\sum_{\mathsf{j}\in\mathsf{J}_{q^{-1}\xi,\sigma}}[Nw]\left(q^{-1}\xi,\mathsf{j}\right)\left[T_{q^{-1}\xi,\sigma,\mathsf{j}}u\right]\left(q^{-1}\varphi\right)&\text{\rm if } q^{-1}\xi\succ\sigma,\\
0&\text{\rm otherwise.}
\end{cases}
\end{equation}
for all $(\xi,\varphi)\in\mathsf{Orb}_{F\wr G}\left(\xi_0,\varphi_{\xi_0}\right)$. Then $N^\flat\in\Hom_{F\wr G}\left(\Ind^{F\wr G}_{F\wr \mathsf{I}(\sigma)}\left(\widetilde{\sigma}\otimes\overline{\theta}\right),L\bigl(\mathsf{Orb}_{F\wr G}\left(\xi_0,\varphi_{\xi_0}\right)\bigr)\right)$. 
\end{lemma}
\begin{proof}
If $(f,g)\in F\wr G$, $gq=q'k$ with $q'\in\mathsf{Q}_\sigma$, $k\in\mathsf{I}(\sigma)$ and $\xi\in q'\mathsf{Orb}_{\mathsf{I}(\sigma)}\left(\xi_0\right)$, then
\[
q'^{-1}\xi \succ \sigma \iff k^{-1}q'^{-1}\xi \succ {}^{k^{-1}}\sigma = \sigma \iff q^{-1}g^{-1}\xi \succ \sigma.
\]
Therefore,
\[ 
\begin{split}
&\left[N^\flat\lambda(f,g)H_{q,u\otimes w}\right](\xi,\varphi)=\left[N^\flat H_{q',\widetilde{\sigma}\left(q'^{-1}f,k\right)u\otimes\theta(k)w}\right](\xi,\varphi)\quad(\text{by }\eqref{actionH})\\
&\quad=\sum_{\mathsf{j}\in\mathsf{J}_{q'^{-1}\xi,\sigma}}[N\theta(k)w]\left(q'^{-1}\xi,\mathsf{j}\right)\left[T_{q'^{-1}\xi,\sigma,\mathsf{j}}\widetilde{\sigma}\left(q'^{-1}f,k\right)u\right]\left(q'^{-1}\varphi\right)\quad(\text{by }\eqref{defNflat})\\
&\quad=\!\!\!\!\!\!\sum_{\mathsf{j}\in\mathsf{J}_{q'^{-1}\xi,\sigma}}\!\!\!\!\![Nw]\!\!\left(k^{-1}q'^{-1}\xi,k^{-1}\mathsf{j}\right)\left[\lambda\left(q'^{-1}f,k\right)\!\left(T_{k^{-1}q'^{-1}\xi,\sigma,k^{-1}\mathsf{j}}u\right)\right]\!\!\left(q'^{-1}\varphi\right)\;(\text{by } \eqref{lambsigm} \text{ and }\eqref{Nthetak})\\
&\quad=\sum_{\mathsf{j}\in\mathsf{J}_{q'^{-1}\xi,\sigma}}[Nw]\left(k^{-1}q'^{-1}\xi,k^{-1}\mathsf{j}\right)\left[T_{k^{-1}q'^{-1}\xi,\sigma,k^{-1}\mathsf{j}}u\right]\left(k^{-1}q'^{-1}f^{-1}\cdot k^{-1}q'^{-1}\varphi\right)\quad(\text{by }\eqref{explpermrep})\\
&\quad=\sum_{\mathsf{j}\in\mathsf{J}_{q^{-1}g^{-1}\xi,\sigma}}[Nw]\left(q^{-1}g^{-1}\xi,\mathsf{j}\right)\left[T_{q^{-1}g^{-1}\xi,\sigma,\mathsf{j}}u\right]\left(q^{-1}g^{-1}\left(f^{-1}\cdot\varphi\right)\right)\quad(\text{by }q'k=gq\text{ and }\eqref{condjbis})\\
&\quad=\left[N^\flat H_{q,u\otimes w}\right]\left(g^{-1}\xi,g^{-1}\left(f^{-1}\cdot\varphi\right)\right)\quad(\text{by }\eqref{defNflat}\text{ and }q^{-1}g^{-1}\xi\succ\sigma)\\
&\quad=\lambda(f,g)\left[N^\flat H_{q,u\otimes w}\right](\xi,\varphi).\quad(\text{by }\eqref{explpermrep})\\
\end{split}
\]
\end{proof}

If $\sigma\in\Sigma_0$ and $\theta\in\widehat{\mathsf{I}(\sigma)}_0$, then we denote by $m_{\sigma,\theta}$ the multiplicity of $\theta$ in $L\left(\mathsf{J}_\sigma\right)$ and let $\left\{N_{\sigma,\theta,l}\mathrel{:}1\leq l\leq m_{\sigma,\theta}\right\}$ denote a basis for $\Hom_{\mathsf{I}(\sigma)}\left(W_\theta,L\left(\mathsf{J}_\sigma\right)\right)$. Therefore,
\begin{equation}\label{decLJ}
L\left(\mathsf{J}_\sigma\right)=\bigoplus_{\theta\in\widehat{\mathsf{I}(\sigma)}_0}\bigoplus_{l=1}^{m_{\sigma,\theta}}N_{\sigma,\theta,l}W_\theta
\end{equation}
is a decomposition into irreducible representations. We recall, however, that the action of $\mathsf{I}(\sigma)$ on $\mathsf{J}_\sigma$ is not necessarily transitive (cf. Remark \ref{Remark}\eqref{Remark1}), so that a preliminary decomposition of $\mathsf{J}_\sigma$ into its $\mathsf{I}(\sigma)$-orbits should be considered.  We also assume that $\left\{w_{\theta,i}\mathrel{:}1\leq i\leq d_\theta\right\}$ is a basis of $W_\theta$. Furthermore, we fix a basis $\left\{u_{\sigma_*,k}\mathrel{:}1\leq k\leq d_{\sigma_*}\right\}$
for every $\sigma_*\in\widehat{F}$ and for $\sigma\in\Sigma_0$ we define $\mathsf{P}_\sigma$ as the set of all functions $p\colon V\rightarrow\mathbb{N}$ such that $1\leq p(v)\leq d_{\sigma_v}$. Then we set $u_{\sigma,p}\coloneqq\bigotimes_{v\in V}u_{\sigma_v,p(v)}$, so that $\left\{u_{\sigma,p}\mathrel{:}p\in\mathsf{P}_\sigma\right\}$ is a basis for $U_\sigma$.

\begin{theorem}\label{bigTheorem}
We have the following decomposition
\begin{equation}\label{decLOxi}
L\bigl(\mathsf{Orb}_{F\wr G}\left(\xi_0,\varphi_{\xi_0}\right)\bigr)=\bigoplus_{\sigma\in\Sigma_0}\bigoplus_{\theta\in\widehat{\mathsf{I}(\sigma)}_0}\bigoplus_{l=1}^{m_{\sigma,\theta}}\Ind_{F\wr\mathsf{I}(\sigma)}^{F\wr G}\left[\bigoplus_{i=1}^{d_\theta}\left(\sum_{(\xi,\mathsf{j})\in\mathsf{J}_\sigma}\left\{\left[N_{\sigma,\theta,l}w_{\theta,i}\right](\xi,\mathsf{j})\right\}T_{\xi,\sigma,\mathsf{j}}\right)U_\sigma\right]
\end{equation}
where each induced representation is irreducible and equivalent to $\Ind^{F\wr G}_{F\wr \mathsf{I}(\sigma)}\left(\widetilde{\sigma}\otimes\overline{\theta}\right)$. A basis for the corresponding isotypic component is given by the complex-valued functions 
\[
E_{l,i,p,q},\quad\text{ for } 1\leq l\leq m_{\sigma,\theta},\; 1\leq i\leq d_\theta,\; p\in\mathsf{P}_\sigma,\; q\in\mathsf{Q}_\sigma,
\]
defined by setting (cf. \eqref{defNflat})
\begin{equation}\label{defElipq}
E_{l,i,p,q}\coloneqq N_{\sigma,\theta,l}^\flat\left(H_{q,u_{\sigma,p}\otimes w_{\theta,i}}\right).
\end{equation}
\end{theorem}
\begin{proof} 
From \eqref{decYV}, \eqref{decULFV} and \eqref{decFVinv}, we deduce that
\begin{equation}\label{decOplus}
\Res_{F^V}^{F\wr G}L\bigl(\mathsf{Orb}_{F\wr G}\left(\xi_0,\varphi_{\xi_0}\right)\bigr)=\bigoplus_{\xi'\in\mathsf{Orb}_G\left(\xi_0\right)}\;\bigoplus_{\sigma'\prec\xi'}\;\bigoplus_{\mathsf{j}'\in\mathsf{J}_{\xi',\sigma'}}T_{\xi',\sigma',\mathsf{j}'}U_{\sigma'}.
\end{equation}
For every $\sigma'$ in \eqref{decOplus} there exists exactly one choice of 
$\sigma\in\Sigma_0$ and of $q\in\mathsf{Q}_\sigma$ such that $\sigma'=\!^q\sigma$; then from \eqref{lambUsigm} we deduce that $T_{\xi',\sigma',\mathsf{j}'}U_{\sigma'}=\lambda\left(\mathbbm{1}_F,q\right)T_{q^{-1}\xi',\sigma,q^{-1}\mathsf{j}'}U_{\sigma}$.
Consequently, setting $(\xi,\mathsf{j})\coloneqq\left(q^{-1}\xi',q^{-1}\mathsf{j}'\right)\in\mathsf{J}_\sigma$ we obtain
\begin{equation}\label{lambdaq}
L\bigl(\mathsf{Orb}_{F\wr G}\left(\xi_0,\varphi_{\xi_0}\right)\bigr)=\bigoplus_{\sigma\in\Sigma_0}\bigoplus_{q\in\mathsf{Q}_\sigma}\lambda\left(\mathbbm{1}_F,q\right)\bigoplus_{(\xi,\mathsf{j})\in\mathsf{J}_\sigma}T_{\xi,\sigma,\mathsf{j}}U_\sigma.
\end{equation}
From \eqref{lambUsigm} it follows that
\begin{equation}\label{lambfkUsigm}
\lambda(f,k)T_{\xi,\sigma,\mathsf{j}}U_\sigma=T_{k\xi,\sigma,k\mathsf{j}}U_\sigma
\end{equation} 
for all $(f,k)\in F\wr\mathsf{I}(\sigma)$.
Taking also into account that $\mathsf{J}_\sigma$ is $\mathsf{I}(\sigma)$-invariant, we 
deduce that the space $\bigoplus_{(\xi,\mathsf{j})\in\mathsf{J}_\sigma}T_{\xi,\sigma,\mathsf{j}}U_\sigma$ is $F\wr\mathsf{I}(\sigma)$-invariant. Then \eqref{lambdaq} implies that 
\begin{equation}\label{lambdaqind}
L\bigl(\mathsf{Orb}_{F\wr G}\left(\xi_0,\varphi_{\xi_0}\right)\bigr)=\bigoplus_{\sigma\in\Sigma_0}\Ind_{F\wr\mathsf{I}(\sigma)}^{F\wr G}\left[\bigoplus_{(\xi,\mathsf{j})\in\mathsf{J}_\sigma}T_{\xi,\sigma,\mathsf{j}}U_\sigma\right].
\end{equation}
On the other hand, from \eqref{lambfkUsigm} it also follows that 
the group $\mathsf{I}(\sigma)$ acts on the direct sum $\bigoplus_{(\xi,\mathsf{j})\in\mathsf{J}_\sigma}T_{\xi,\sigma,\mathsf{j}}U_\sigma$ by permuting the summands. Moreover, the action on these summands is isomorphic to the action on $\mathsf{J}_\sigma$ so that from \eqref{decLJ} we deduce the following direct sum decomposition
\begin{equation}\label{lambdaIsig}
\bigoplus_{(\xi,\mathsf{j})\in\mathsf{J}_\sigma}T_{\xi,\sigma,\mathsf{j}}U_\sigma
=\bigoplus_{\theta\in\widehat{\mathsf{I}(\sigma)}_0}\bigoplus_{l=1}^{m_{\sigma,\theta}}\left[\bigoplus_{i=1}^{d_\theta}\left(\sum_{(\xi,\mathsf{j})\in\mathsf{J}_\sigma}\left\{\left[N_{\sigma,\theta,l}w_{\theta,i}\right](\xi,\mathsf{j})\right\}T_{\xi,\sigma,\mathsf{j}}\right)U_\sigma\right].
\end{equation}
Indeed, the functions $N_{\sigma,\theta,l}w_{\theta,i}$ form a basis for the space $L\left(\mathsf{J}_\sigma\right)$.
Furthermore, if $k\in\mathsf{I}(\sigma)$, then
\[
\begin{split}
\left[N_{\sigma,\theta,l}w_{\theta,i}\right]\left(k^{-1}\xi,k^{-1}\mathsf{j}\right)=&\left[\lambda(k)N_{\sigma,\theta,l}w_{\theta,i}\right](\xi,\mathsf{j})\\
(N_{\sigma,\theta,l}\in\Hom_{\mathsf{I}(\sigma)}\left(W_\theta,L\left(\mathsf{J}_\sigma\right)\right))\qquad=&\left[N_{\sigma,\theta,l}\theta(k)w_{\theta,i}\right](\xi,\mathsf{j})\\
(\theta(k)w_{\theta,i}=\sum_{h=1}^{d_\theta}\gamma_{hi}(k)w_{\theta,h})\qquad=&\sum_{h=1}^{d_\theta}\gamma_{hi}(k)\left[N_{\sigma,\theta,l}w_{\theta,h}\right](\xi,\mathsf{j})
\end{split}
\]
so that, if also $f\in F^V$, $u\in U_\sigma$ and $1\leq i\leq d_\theta$, then from \eqref{lambsigm} we deduce that
\[
\lambda(f,k)\sum_{(\xi,\mathsf{j})\in\mathsf{J}_\sigma}\left\{\left[N_{\sigma,\theta,l}w_{\theta,i}\right](\xi,\mathsf{j})\right\}T_{\xi,\sigma,\mathsf{j}}u
=\sum_{h=1}^{d_\theta}\gamma_{hi}(k)\sum_{(\xi,\mathsf{j})\in\mathsf{J}_\sigma}\left\{\left[N_{\sigma,\theta,l}w_{\theta,h}\right](\xi,\mathsf{j})\right\}T_{\xi,\sigma,\mathsf{j}}\widetilde{\sigma}(f,k)u.
\]
This means that, for each $\theta\in\widehat{\mathsf{I}(\sigma)}_0$ and $1\leq l\leq m_{\sigma,\theta}$, the subspace in the square brackets in \eqref{lambdaIsig} is $F\wr\mathsf{I}(\sigma)$-invariant and equivalent to $\widetilde{\sigma}\otimes\overline{\theta}$ as an $F\wr\mathsf{I}(\sigma)$-representation; therefore, \eqref{decLOxi} follows by invoking \eqref{lambdaqind}. 
Using \eqref{lambdaq}, we may find bases of the subspaces in \eqref{decLOxi} consisting of functions in $L\left(\mathsf{Orb}_{F\wr G}\left(\xi_0,\varphi_{\xi_0}\right)\right)$: if $p\in\mathsf{P}_\sigma$ and $(\xi',\varphi')\in\mathsf{Orb}_{F\wr G}\left(\xi_0,\varphi_{\xi_0}\right)$, with $\xi'\in q\mathsf{Orb}_{\mathsf{I}(\sigma)}\left(\xi_0\right)$, then
\[
\begin{split}
\Biggl[\lambda\left(\mathbbm{1}_F,q\right)\sum_{(\xi,\mathsf{j})\in\mathsf{J}_\sigma}&\left\{\left[N_{\sigma,\theta,l}w_{\theta,i}\right](\xi,\mathsf{j})\right\}T_{\xi,\sigma,\mathsf{j}}u_{\sigma,p}\Biggr](\xi',\varphi')\\
(\text{by }\eqref{explpermrep})&\quad=\Biggl[\sum_{(\xi,\mathsf{j})\in\mathsf{J}_\sigma}\left\{\left[N_{\sigma,\theta,l}w_{\theta,i}\right](\xi,\mathsf{j})\right\}T_{\xi,\sigma,\mathsf{j}}u_{\sigma,p}\Biggr]\left(q^{-1}\xi',q^{-1}\varphi'\right)\\
(\text{by }\eqref{decULFV}\text{ and }\eqref{defJsigma})&\quad=\sum_{\mathsf{j}\in\mathsf{J}_{q^{-1}\xi',\sigma}}\left[N_{\sigma,\theta,l}w_{\theta,i}\right]\left(q^{-1}\xi',\mathsf{j}\right)\left[T_{q^{-1}\xi',\sigma,\mathsf{j}}u_{\sigma,p}\right]\left(q^{-1}\varphi'\right)\\
(\text{by }\eqref{defNflat})&\quad=N_{\sigma,\theta,l}^\flat\left(H_{q,u_{\sigma,p}\otimes w_{\theta,i}}\right)(\xi',\varphi')\\
(\text{by }\eqref{defElipq})&\quad=E_{l,i,p,q}(\xi',\varphi').
\end{split}
\]
Finally, the fact that the functions $E_{l,i,p,q}$ form a basis for the isotypic component follows also from Lemma \ref{Theointoper}, combined with \eqref{lambdaq}, \eqref{lambdaIsig} and the fact that the functions $H_{q,u_{\sigma,p}\otimes w_{\theta,i}}$ form a basis for the space of $\Ind^{F\wr G}_{F\wr \mathsf{I}(\sigma)}\left(\widetilde{\sigma}\otimes\overline{\theta}\right)$.
\end{proof}

\begin{corollary}\label{Cormult}
The operators $N_{\sigma,\theta,l}^\flat$, for $l=1,2,\dotsc, m_{\sigma,\theta}$, form a basis for the space
\[
\Hom_{F\wr G}\left(\Ind^{F\wr G}_{F\wr \mathsf{I}(\sigma)}\left(\widetilde{\sigma}\otimes\overline{\theta}\right),L\bigl(\mathsf{Orb}_{F\wr G}\left(\xi_0,\varphi_{\xi_0}\right)\bigr)\right).
\]
In particular, the multiplicity of $\Ind^{F\wr G}_{F\wr \mathsf{I}(\sigma)}\left(\widetilde{\sigma}\otimes\overline{\theta}\right)$ in $L\bigl(\mathsf{Orb}_{F\wr G}\left(\xi_0,\varphi_{\xi_0}\right)\bigr)$ is equal to the multiplicity of $\theta$ in $L\left(\mathsf{J}_\sigma\right)$. 
\end{corollary}
\begin{proof}
This is an immediate consequence of Lemma \ref{Theointoper} and Theorem \ref{bigTheorem}.
\end{proof}

\begin{example}\label{Example1bis}
In the setting of Example \ref{Example1}, \eqref{decLYa} becomes $L(Y)=\bigoplus_{\sigma_*\in\widehat{F}_0}m_{\sigma_*} U_{\sigma_*}$, where $\widehat{F}_0$ is the set of irreducible $F$-representations contained in $L(Y)$. Here, $\Sigma_0$ is a set of representatives for the orbits of $G$ on $\widehat{F}_0^V$. If $\sigma\in\Sigma_0$, then the set $\mathsf{J}_\sigma$ consists of all $\mathsf{j}\in\mathbb{N}^V$ such that $1\leq \mathsf{j}(v)\leq m_{\sigma_v}$ for all $v\in V$, and the operators $\left\{T_{\sigma,\mathsf{j}}\mathrel{:}\mathsf{j}\in\mathsf{J}_\sigma\right\}$ form a basis for $\Hom_{F^V}\left(U_\sigma,L\bigl(Y^V\bigr)\right)$. Moreover, $\widehat{\mathsf{I}(\sigma)}_0$ is the set of all $\theta\in\widehat{\mathsf{I}(\sigma)}$ contained in $L\left(\mathsf{J}_\sigma\right)$, and \eqref{decLOxi} becomes
\[
L\bigl(Y^V\bigr)=\bigoplus_{\sigma\in\Sigma_0}\bigoplus_{\theta\in\widehat{\mathsf{I}(\sigma)}_0}\bigoplus_{l=1}^{m_{\sigma,\theta}}\Ind_{F\wr\mathsf{I}(\sigma)}^{F\wr G}\left[\bigoplus_{i=1}^{d_\theta}\left(\sum_{\mathsf{j}\in\mathsf{J}_\sigma}\left\{\left[N_{\sigma,\theta,l}w_{\theta,i}\right](\mathsf{j})\right\}T_{\sigma,\mathsf{j}}\right)U_\sigma\right].
\]
In particular, the multiplicity of $\Ind^{F\wr G}_{F\wr \mathsf{I}(\sigma)}\left(\widetilde{\sigma}\otimes\overline{\theta}\right)$ in $L\bigl(Y^V\bigr)$ is equal to the multiplicity of $\theta$ in $L\left(\mathsf{J}_\sigma\right)$ (cf. \cite[Corollary 3.3.12]{book3} and \cite[Corollary 4.9]{ScaTol}). Finally, by virtue of \eqref{defNflat}, \eqref{defElipq} becomes
\[
E_{l,i,p,q}(\varphi)=\sum_{\mathsf{j}\in\mathsf{J}_\sigma}\left\{\left[N_{\sigma,\theta,l}w_{\theta,i}\right](\mathsf{j})\right\}\left[T_{\sigma,\mathsf{j}}u_{\sigma,p}\right]\left(q^{-1}\varphi\right)
\]
for every $\varphi\in Y^V$.
\end{example}

\begin{example}\label{Example2bis}
Now we examine the decomposition of the permutation representation associated with the homogeneous space in Example \ref{Example2}. We assume that the decomposition of $L(Y)$ is as in Example \ref{Example1bis}. Then $\Sigma_0$ is the set of all $\sigma=\bigotimes_{v\in V}\sigma_v$, where $\sigma_{v_0}\in\widehat{F}_0$ while, for all $v\neq v_0$, $\sigma_v$ is the trivial representation. If also $\sigma_{v_0}$ is trivial, then $\sigma=\sigma^0$, the trivial representation of $F^V$. We have to distinguish two cases, according to whether
$\sigma\neq\sigma^0$ or $\sigma=\sigma^0$. 

For $\sigma\neq\sigma^0$, we have $\mathsf{I}(\sigma)=\mathsf{Stab}_G\left(v_0\right)$. Moreover, $\sigma\prec\xi_v$ if and only if $v=v_0$. Hence, $\mathsf{J}_\sigma$ is simply the set of all functions $\mathsf{j}\colon\left\{v_0\right\}\rightarrow\{1,2,\dotsc,m_{\sigma_{v_0}}\}$, and $L\left(\mathsf{J}_\sigma\right)$ is the sum of $m_{\sigma_{v_0}}$ copies of the trivial representation $\theta^0$ of $\mathsf{Stab}_G\left(v_0\right)$; we will denote by $w^0$ a non-zero vector in the space of $\theta^0$. In \eqref{decLJ} we may take $\left[N_{\sigma,\theta^0,l}w^0\right]\left(\xi_{v_0},\mathsf{j}\right)=\delta_{l,\mathsf{j}\left(v_0\right)}$ and $\mathsf{Q}_\sigma=\left\{q_v\mathrel{:}v\in V\right\}$, where $q_vv_0=v$. Moreover, \eqref{Txiuv} yields
\[
\left[T_{\xi_{v_0},\sigma,\mathsf{j}}u_{\sigma,p}\right](\varphi)=\left[T_{a_0,\sigma_{v_0},\mathsf{j}\left(v_0\right)}u_{\sigma_{v_0},p(v_0)}\right]\left(\varphi\left(v_0\right)\right),
\]
where, for $v\neq v_0$ (in which case $Y_a$ is a singleton and $\sigma_v$ is trivial), we assume that $T_{a_1,\sigma^0,1}u_{\sigma^0,1}$ is identically equal to $1$. Therefore, taking \eqref{defNflat} into account and observing that $q_v^{-1}\xi_{v'} \succ \sigma$  if and only if $v=v'$, \eqref{defElipq} becomes
\[
E_{l,1,p,q_v}\left(\xi_{v'},\varphi\right)=\delta_{v,v'}\sum_{\mathsf{j}\left(v_0\right)=1}^{m_{\sigma_{v_0}}}\delta_{l,\mathsf{j}\left(v_0\right)}\left[T_{q_v^{-1}\xi_v,\sigma,\mathsf{j}}u_{\sigma,p}\right]\left(q_v^{-1}\varphi\right)=\delta_{v,v'}\left[T_{a_0,\sigma_{v_0},l}u_{\sigma_{v_0},p(v_0)}\right](\varphi(v)).
\]
From \eqref{explpermrep} we deduce that $\Ind^{F\wr G}_{F\wr \mathsf{I}(\sigma)}\left(\widetilde{\sigma}\otimes\overline{\theta^0}\right)\sim U_\sigma\otimes L(V)$, the latter with the action $(f,g)(u \otimes \delta_v) = \sigma_{v_0}(f(gv))u \otimes \delta_{gv}.$ Then we obtain $m_{\sigma_{v_0}}$ copies of this representation, one for each possible value of $l$ (cf. Corollary \ref{Cormult}).

When $\sigma=\sigma^0$, we have $\mathsf{I}\left(\sigma^0\right)=G$, $\mathsf{J}_{\sigma^0}$ can be identified with $V$, leading us to consider the permutation representation $L(V)=\bigoplus_{\theta\in\widehat{G}_0}m_\theta W_\theta$, where 
$\widehat{G}_0$ denotes the set of all $\theta\in\widehat{G}$ contained in $L(V)$. In particular, with respect to \eqref{irrepwreath}, there is no induction involved and, by Corollary \ref{Cormult},  the decomposition of the permutation representation features $m_\theta$ copies of $U_{\sigma^0}\otimes W_{\overline{\theta}}$.

Finally, under the isomorphism $T$ in Example \ref{Example2}, \eqref{decLOxi} becomes
\begin{equation}\label{deccompact}
L(V\times Y)=\left[\bigoplus_{\sigma\in\Sigma_0\setminus\left\{\sigma^0\right\}}m_\sigma\left(U_\sigma\otimes L(V)\right)\right] \bigoplus \left[\bigoplus_{\theta\in\widehat{G}_0}m_\theta \left(U_{\sigma^0}\otimes W_{\overline{\theta}}\right)\right],
\end{equation}
where $m_\sigma := m_{\sigma_{v_0}}$.
See also \cite[Section 3.1.1]{book3} and \cite[Example 3.13]{ScaTol}.
\end{example}

\begin{example}\label{Example21bis}
Now we examine the decomposition of the permutation representation associated with the homogeneous space in Example \ref{Example21}, namely $\mathsf{Orb}_{F\wr G}\left(v_0,v_j,y_0\right)$. Assume that the decomposition of $L(Y)$ is given as in Example \ref{Example1bis} and note that $\Sigma_0$ is identical to the one in Example \ref{Example2bis}. 

If $\sigma\neq\sigma^0$, then $\mathsf{I}(\sigma)=\mathsf{Stab}_G\left(v_0\right)$; since $\sigma\prec\xi_{(v,w)}$ if and only if $v=v_0$, we have
\[
\mathsf{J}_\sigma\simeq\left(\left\{v_0\right\}\times\Delta_{v_0,j}\right)\times\left\{1,2,\dotsc,m_{\sigma_{v_0}}\right\}\simeq\left(\left\{v_0\right\}\times\left\{1,2,\dotsc,m_{\sigma_{v_0}}\right\}\right)\times\Delta_{v_0,j},
\]
and the action of $\mathsf{Stab}_G\left(v_0\right)$ on the functions $\mathsf{j}\colon\left\{v_0\right\}\rightarrow\{1,2,\dotsc,m_{\sigma_{v_0}}\}$ is trivial. Consequently, $L\left(\mathsf{J}_\sigma\right)$ is the sum of $m_{\sigma_{v_0}}$ copies of the permutation representation of $\mathsf{Stab}_G\left(v_0\right)$ 
on $L\left(\Delta_{v_0,j}\right)$. If $L\left(\Delta_{v_0,j}\right)=\bigoplus_{\theta\in\widehat{\mathsf{I}(\sigma)}_0}m_\theta W_\theta$ is the decomposition into irreducible $\mathsf{Stab}_G\left(v_0\right)$-representations, then the multiplicity of
$\Ind^{F\wr G}_{F\wr \mathsf{I}(\sigma)}\left(\widetilde{\sigma}\otimes\overline{\theta}\right)$ in $L\bigl(\mathsf{Orb}_{F\wr G}\left(v_0,v_j,y_0\right)\bigr)$ is equal to the product $m_{\sigma_{v_0}}m_\theta$.

If $\sigma=\sigma^0$, then $\mathsf{I}\left(\sigma^0\right)=G$ and $\mathsf{J}_{\sigma^0}\simeq\mathsf{Orb}_G\left(v_0,v_j\right)$. Suppose that $L\bigl(\mathsf{Orb}_G\left(v_0,v_j\right)\bigr)=\bigoplus_{\theta\in\widehat{G}_0}m_\theta W_\theta$ is the decomposition of the corresponding permutation representation into irreducible representations, where
$\widehat{G}_0$ denotes the set of all $\theta\in\widehat{G}$ contained in $L\bigl(\mathsf{Orb}_G\left(v_0,v_j\right)\bigr)$. In particular, with respect to \eqref{irrepwreath}, no induction is required. Thus, the decomposition of the permutation representation contains $m_\theta$ copies of $U_{\sigma^0}\otimes W_{\overline{\theta}}$. In conclusion, \eqref{decLOxi} becomes
\begin{equation}\label{DecgenCentip}
\begin{split}
L\bigl(\mathsf{Orb}_{F\wr G}\left(v_0,v_j,y_0\right)\bigr)=&\left[\bigoplus_{\sigma\in\Sigma_0\setminus\left\{\sigma^0\right\}}\bigoplus_{\theta\in\widehat{\mathsf{I}(\sigma)}_0}m_\sigma m_\theta\Ind^{F\wr G}_{F\wr \mathsf{I}(\sigma)}\left(U_{\widetilde{\sigma}}\otimes W_{\overline{\theta}}\right)\right] \bigoplus\\
&\quad\bigoplus \left[\bigoplus_{\theta\in\widehat{G}_0}m_\theta \left(U_{\sigma^0}\otimes W_{\overline{\theta}}\right)\right],
\end{split}
\end{equation}
where $m_\sigma := m_{\sigma_{v_0}}$.
\end{example}

\begin{example}\label{Example3bis}
In the setting of Example \ref{Example3}, we may suppose again that \eqref{decLYa} is simply $L(Y)=\bigoplus_{\sigma_*\in\widehat{F}_0}m_{\sigma_*} U_{\sigma_*}$. Using the isomorphism \eqref{isomJohn}, we replace $\xi_Z$  with $Z$. We 
denote by $\widehat{F^Z}$ the set of all $\sigma\in\widehat{F}_0^V$ such that $\mathsf{supp}(\sigma)\subseteq Z$; this is equivalent to $\sigma\prec\xi_Z$. Let us fix a subset $Z_0\subsetneqq V$, and let $\Sigma_0$ be a set of representatives for the $G$-conjugacy classes of the representations in $\widehat{F^{Z_0}}$. For all $Z\in\mathsf{Orb}_G\left(Z_0\right)$ and $\sigma\in\widehat{F^Z}$, the set $\mathsf{J}_{Z,\sigma}$ consists of all $\mathsf{j}\in\mathbb{N}^V$ such that $1\leq \mathsf{j}(v)\leq m_{\sigma_v}$ if $v\in Z$, while $\mathsf{j}(v)=1$ if $v\notin Z$. The set $\mathsf{J}_\sigma$ in \eqref{defJsigma} consists of all pairs $(Z,\mathsf{j})$ such that $Z\in\mathsf{Orb}_G\left(Z_0\right)$, $Z\supseteq\mathsf{supp}(\sigma)$, and $\mathsf{j}\in\mathsf{J}_{Z,\sigma}$. Note that, if $\sigma\in\Sigma_0$, then the second condition is equivalent to $\mathsf{supp}(\sigma)\subseteq Z\cap Z_0$.
Moreover, $\mathsf{supp}(\sigma)$ is $\mathsf{I}(\sigma)$-invariant and contained in $Z_0$, so that, with the notation in \eqref{VsIsigma},
\begin{equation}\label{OrbIsigZ0}
\mathsf{Orb}_{\mathsf{I}(\sigma)}\left(Z_0\right)=\left\{k\left(Z_0\cap V_{\sigma_0}\right)\sqcup\mathsf{supp}(\sigma)\mathrel{:}k\in\mathsf{I}(\sigma)\right\}\simeq\mathsf{Orb}_{\mathsf{I}(\sigma)}\left(Z_0\cap V_{\sigma_0}\right),
\end{equation}
as $\mathsf{I}(\sigma)$-homogeneous spaces. If we denote by $T_{Z,\sigma,\mathsf{j}}$
the operators in \eqref{defTxisigj}, and take \eqref{isomJohn} into account, then \eqref{decLOxi} becomes 
\[
\bigoplus_{Z\in\mathsf{Orb}_G\left(Z_0\right)}L\bigl(Y^Z\bigr)=\bigoplus_{\sigma\in\Sigma_0}\bigoplus_{\theta\in\widehat{\mathsf{I}(\sigma)}_0}\bigoplus_{l=1}^{m_{\sigma,\theta}}\Ind_{F\wr\mathsf{I}(\sigma)}^{F\wr G}\left[\bigoplus_{i=1}^{d_\theta}\left(\sum_{(Z,\mathsf{j})\in\mathsf{J}_\sigma}\left\{\left[N_{\sigma,\theta,l}w_{\theta,i}\right](Z,\mathsf{j})\right\}T_{Z,\sigma,\mathsf{j}}\right)U_\sigma\right].
\]
In particular, the multiplicity of $\Ind^{F\wr G}_{F\wr \mathsf{I}(\sigma)}\left(\widetilde{\sigma}\otimes\overline{\theta}\right)$ in this decomposition is equal to the multiplicity of $\theta$ in $L\left(\mathsf{J}_\sigma\right)$. 
\end{example}

\begin{example}\label{Example4bis}
In the setting of Example \ref{Example4}, we omit the unique $a\in A$ from the notation for the operators $T_{\sigma_*,j}$, where $\sigma_*\in\widehat{F}$ and $1\leq j\leq d_{\sigma_*}$. These operators form a basis for the space $\Hom_F\left(U_{\sigma_*},L(F)\right)$, thereby providing a decomposition of the left regular representation of $F$. Similarly, we omit the index $\xi$ from the operators in \eqref{defTxisigj}, observing that they provide a decomposition of the group algebra of $F^V$. Moreover, $\mathsf{J}_\sigma'\coloneqq\mathsf{J}_{\xi,\sigma}$ consists of all $\mathsf{j}\in \mathbb{N}^V$ such that $1\leq \mathsf{j}(v)\leq d_{\sigma_v}$ for all $v\in V$ and is thus independent of $\xi$, while $\mathsf{J}_\sigma=\mathsf{Orb}_G\left(\xi_0\right)\times\mathsf{J}_\sigma'$, because $\sigma\prec\xi$ for all $\sigma\in\widehat{F^V}$ and $\xi\in\mathsf{Sym}(V)$. We have $\Sigma_0=\Sigma$ (cf.\ \eqref{irrepwreath}) and the action of $\mathsf{I}(\sigma)$ on $\mathsf{J}_\sigma$ is semiregular; hence, as an $\mathsf{I}(\sigma)$-set, it decomposes into the union of 
\[
\frac{\left\lvert\mathsf{J}_\sigma\right\rvert}{\left\lvert\mathsf{I}(\sigma)\right\rvert}=\frac{\left\lvert\mathsf{Orb}_G\left(\xi_0\right)\right\rvert\cdot\left\lvert\mathsf{J}_\sigma'\right\rvert}{\left\lvert\mathsf{I}(\sigma)\right\rvert}=\frac{\lvert G \rvert\cdot d_\sigma}{\left\lvert\mathsf{I}(\sigma)\right\rvert}=\left\lvert\mathsf{Q}_\sigma\right\rvert\cdot d_\sigma
\]
disjoint copies of the left regular action of $\mathsf{I}(\sigma)$ on itself. From Corollary \ref{Cormult}, we deduce that
$m_{\sigma,\theta}=\left\lvert\mathsf{Q}_\sigma\right\rvert\cdot d_\sigma\cdot d_\theta$, which in turn is equal to $\dim\Ind^{F\wr G}_{F\wr \mathsf{I}(\sigma)}\left(\widetilde{\sigma}\otimes\overline{\theta}\right)$. This agrees with the fact that, in the present case, \eqref{decLOxi} yields a decomposition of the left regular representation of $F\wr G$. Finally, we observe that, in this case, if $\mathsf{I}(\sigma)$ is a proper subgroup of $G$, then the inclusion \eqref{defJsig0old} is strict, its right-hand side being equal to $\mathsf{Orb}_{\mathsf{I}(\sigma)}\left(\xi_0\right)\times\mathsf{J}_\sigma'$. 
\end{example}

\begin{example}\label{Example5bis}
In the setting of Example \ref{Example5}, the decomposition \eqref{OrbFGphi0} leads to
\[
L\bigl(\mathsf{Orb}_{F\wr G}\left(\varphi_0\right)\bigr)\sim \Ind_{F\wr\mathsf{Stab}_G(\varphi_0)}^{F\wr G}\widetilde{L}\Bigl(\mathsf{Orb}_{F^V}\left(\varphi_0\right)\Bigr).
\]
In particular, only the extension factor $\widetilde{L}$ in Proposition \ref{Propindinflext} appears.
The general formulas \eqref{defNflat}, \eqref{decLOxi}, and \eqref{defElipq} do not simplify significantly; however, we indicate some peculiarities and notation adjustments. We write $\sigma\prec\varphi$ to denote that a representation $\sigma\in\widehat{F^V}$ appears in $L\bigl(\mathsf{Orb}_{F^V}(\varphi)\bigr)$.
Fixing $\varphi_0\in\Gamma_0$, we may assume that for every $\sigma\in\Sigma_0$ we have $\sigma\prec\varphi_0$. If $\varphi\in\mathsf{Orb}_G\left(\varphi_0\right)$ and $\sigma\prec\varphi$, then $\mathsf{J}_{\varphi,\sigma}$ (which replaces $\mathsf{J}_{\xi,\sigma}$) is the set of all $\mathsf{j}\in\mathbb{N}^V$ such that $1\leq\mathsf{j}(v)\leq m_{\varphi(v),\sigma_v}$ for all $v\in V$, where the notation $m_{a,\sigma_*}$ in \eqref{decLYa} is replaced by $m_{y_a,\sigma_*}$. Similarly, $\mathsf{J}_\sigma$ may be identified with the set of all pairs $(\varphi,\mathsf{j})$ such that $\varphi\in\mathsf{Orb}_G\left(\varphi_0\right)$, $\sigma\prec\varphi$ and $\mathsf{j}\in\mathsf{J}_{\varphi,\sigma}$. The operators in \eqref{defTxisigj} may now be denoted by $T_{\varphi,\sigma,\mathsf{j}}$. 
\end{example}

\begin{example}\label{Example6bis}
In the setting of Example \ref{Example6}, Proposition \ref{Propindinflext} becomes 
\[
L\bigl(\mathsf{Orb}_{F\wr G}\left(\xi_0,\varphi_0\right)\bigr)\sim\overline{L}\bigl(\mathsf{Orb}_G\left(\xi_0\right)\bigr)\bigotimes\widetilde{L}\left(Y^V\right);
\]
cf. \eqref{OrStabOrAt}. In particular, there is no induction and we have $\sigma\prec\xi$ for every $\sigma\in\widehat{F}_0^V$ and $\xi\in\mathsf{Orb}_G\left(\xi_0\right)$, that is $\prec$ is trivial. Assuming the notation of Example \ref{Example1bis}, we define $\mathsf{J}_\sigma'$ as the set of all $\mathsf{j}\in\mathbb{N}^V$ such that $1\leq \mathsf{j}(v)\leq m_{\sigma_v}$ for all $v\in V$. Then $\mathsf{J}_\sigma=\mathsf{Orb}_G\left(\xi_0\right)\times\mathsf{J}_\sigma'$. 
\end{example}

\section{Multiplicity-free permutation representations and spherical functions}\label{Secmultfreepermrep}

In the present section, we examine when the permutation representation studied in Section \ref{Secpermrep} is multiplicity free and, in this case, we compute the corresponding spherical functions.

\begin{proposition}\label{Propmultfree}
The permutation representation of $F\wr G$ on $L\bigl(\mathsf{Orb}_{F\wr G}\left(\xi_0,\varphi_{\xi_0}\right)\bigr)$ is multiplicity free if and only if both of the following conditions hold
\begin{enumerate}[(1)]
\item\label{Propmultfreeb} for every $\sigma\in\Sigma_0$ we have $\mathsf{J}_\sigma\simeq\mathsf{Orb}_{\mathsf{I}(\sigma)}\left(\xi_0\right)$ as $\mathsf{I}(\sigma)$-homogeneous spaces;
\item\label{Propmultfreec} for every $\sigma\in\Sigma_0$ the permutation representation of $\mathsf{I}(\sigma)$ on $L\bigl(\mathsf{Orb}_{\mathsf{I}(\sigma)}\left(\xi_0\right)\bigr)$ is multiplicity free.
\end{enumerate}
\end{proposition}
\begin{proof}
Assume that $L\bigl(\mathsf{Orb}_{F\wr G}\left(\xi_0,\varphi_{\xi_0}\right)\bigr)$ is multiplicity free. By Corollary \ref{Cormult}, for every $\sigma\in\Sigma_0$ the permutation representation of $\mathsf{I}(\sigma)$ on $L\left(\mathsf{J}_\sigma\right)$ is multiplicity free. If the permutation representation of $F$ on $L\left(Y_a\right)$ is not multiplicity free for some $a\in A$, then, arguing similarly to Remark \ref{Remark}\eqref{Remark1}, we deduce that the action of $\mathsf{I}(\sigma)$ on $\mathsf{J}_\sigma$, for a suitable $\sigma \in \Sigma_0$, is not transitive, which is impossible since each orbit would yield a copy of the trivial representation of $\mathsf{I}(\sigma)$ in $L\left(\mathsf{J}_\sigma\right)$;  
the suitable $\sigma$ has a component $\sigma_v \in \widehat{F}_a$ with multiplicity $\ge 2$, where $v$ satisfies $\alpha(\xi_0(v))=a$.

Similarly, if the inclusion in \eqref{defJsig0old} is strict, then 
$\mathsf{Orb}_{\mathsf{I}(\sigma)}(\xi_0)$ is a proper $\mathsf{I}(\sigma)$-invariant subset of $\{\xi \in \mathsf{Orb}_G(\xi_0) \mathrel{:} \sigma \prec \xi\}$, so that the $\mathsf{I}(\sigma)$-equivariant surjective projection $(\xi,\mathsf{j}) \mapsto \xi$ shows that $\mathsf{J}_\sigma$ has at least two orbits. This is again impossible. We conclude that $\mathsf{J}_\sigma$ must consist of all pairs $(\xi,\mathsf{j})$ such that $\xi\in\mathsf{Orb}_{\mathsf{I}(\sigma)}\left(\xi_0\right)$ and $\mathsf{j}\equiv 1$, and this proves \eqref{Propmultfreeb}.  Finally, from \eqref{Propmultfreeb} it follows that 
 $L\left(\mathsf{J}_\sigma\right)\sim L\bigl(\mathsf{Orb}_{\mathsf{I}(\sigma)}\left(\xi_0\right)\bigr)$ as $\mathsf{I}(\sigma)$-permutation representations, so that the latter must be multiplicity free, and \eqref{Propmultfreec} is also proved.

Conversely, if \eqref{Propmultfreeb} and \eqref{Propmultfreec} hold, then Corollary \ref{Cormult} ensures that $L\bigl(\mathsf{Orb}_{F\wr G}\left(\xi_0,\varphi_{\xi_0}\right)\bigr)$ is multiplicity free.
\end{proof}

\begin{corollary}\label{Cormultfree}
\begin{enumerate}[(1)]
\item\label{Cormultfree1}
The condition in Proposition \ref{Propmultfree}\eqref{Propmultfreeb} is satisfied if and only if both of the following conditions hold
\begin{itemize}
\item for every $a\in A$, the permutation representation of $F$ on $L\left(Y_a\right)$ is multiplicity free;
\item for every $\sigma\in\Sigma_0$, we have equality in \eqref{defJsig0old}.
\end{itemize}
In particular, if the conditions above are satisfied, then we have 
\begin{equation}\label{JsiginCor}
\left\{\xi\in\mathsf{Orb}_G\left(\xi_0\right)\mathrel{:}\sigma\prec\xi\right\}=\mathsf{Orb}_{\mathsf{I}(\sigma)}\left(\xi_0\right),
\end{equation}
and we identify $\mathsf{J}_\sigma$ with this set.
\item\label{Cormultfree2}
If the permutation representation of $F\wr G$ on $L\bigl(\mathsf{Orb}_{F\wr G}\left(\xi_0,\varphi_{\xi_0}\right)\bigr)$ is multiplicity free, then $L\bigl(\mathsf{Orb}_G\left(\xi_0\right)\bigr)$ is also a multiplicity-free permutation representation of $G$.
\end{enumerate}
\end{corollary}
\begin{proof}
Indeed, the ``only if'' part of \eqref{Cormultfree1} follows directly from the proof of Proposition \ref{Propmultfree}\eqref{Propmultfreeb}, since $\mathsf{J}_\sigma\simeq\mathsf{Orb}_{\mathsf{I}(\sigma)}\left(\xi_0\right)$ implies that $\mathsf{I}(\sigma)$ acts transitively on $\mathsf{J}_\sigma$. Conversely, if the two conditions in \eqref{Cormultfree1} hold, then every $(\xi,\mathsf{j})\in\mathsf{J}_\sigma$ satisfies $\xi\in\mathsf{Orb}_{\mathsf{I}(\sigma)}\left(\xi_0\right)$ and $\mathsf{j}\equiv 1$, so that the projection $(\xi,\mathsf{j})\mapsto \xi$ is an isomorphism of $\mathsf{I}(\sigma)$-homogeneous spaces from $\mathsf{J}_\sigma$ onto $\mathsf{Orb}_{\mathsf{I}(\sigma)}\left(\xi_0\right)$. We can prove \eqref{Cormultfree2} by considering the trivial representation $\sigma^0$ of $F^V$, so that $\mathsf{I}\left(\sigma^0\right)=G$, and then applying Proposition \ref{Propmultfree}\eqref{Propmultfreec}.
\end{proof}

Henceforth, we assume that the conditions in Proposition \ref{Propmultfree} hold. By Corollary \ref{Cormultfree}, \eqref{decLYa} becomes
\begin{equation}\label{decLYamf}
L\left(Y_a\right)=\bigoplus_{\sigma_*\in\widehat{F}_a}T_{a,\sigma_*}U_{\sigma_*},
\end{equation}
where the operators $T_{a,\sigma_*}$ do not depend on $j$ because it can only assume the value 1.
If $u_{a,\sigma_*}$ is a unit $\mathsf{Stab}_F\left(y_a\right)$-invariant vector in $U_{\sigma_*}$, then
we may suppose that  $\left[T_{a,\sigma_*}u\right]\left(fy_a\right)=\left\langle u,\sigma_*(f)u_{a,\sigma_*}\right\rangle_{U_{\sigma_*}}$ for all $f\in F$ and $u\in U_{\sigma_*}$ (see \cite[Corollary 4.6.4]{book},  adapted here by omitting the normalization coefficient).  In particular, for $a,b\in A$ and $\sigma_*\in\widehat{F}_a\cap \widehat{F}_b$ we define $\mathcal{S}_{a,b,\sigma_*}\in L\left(Y_a\right)$ by setting
\begin{equation}\label{defSabsig}
\mathcal{S}_{a,b,\sigma_*}\left(fy_a\right)\coloneqq\left\langle u_{b,\sigma_*},\sigma_*(f)u_{a,\sigma_*}\right\rangle_{U_{\sigma_*}}
\end{equation}
for every $f\in F$.

Next, we define the operators in \eqref{defTxisigj} by omitting the index $\mathsf{j}$ (which is identically 1), which yields $T_{\xi,\sigma}=\bigotimes_{v\in V}T_{\alpha(\xi(v)),\sigma_v}$. Finally, for $\sigma\in\Sigma_0$ and $\xi\in\mathsf{J}_\sigma$ (recall \eqref{JsiginCor}), we set $u_\sigma\coloneqq\bigotimes_{v\in V}u_{\alpha\circ\xi_0(v),\sigma_v}$ and define $\mathcal{S}_{\xi,\sigma}\in L\bigl(\mathsf{Orb}_{F^V}\left(\varphi_{\xi}\right)\bigr)$ by setting $\mathcal{S}_{\xi,\sigma}\coloneqq T_{\xi,\sigma}u_\sigma$; 
more explicitly, 
\begin{equation}\label{defSxisigma}
\mathcal{S}_{\xi,\sigma}\left(f\cdot\varphi_\xi\right)=\prod_{v\in V}\left\langle u_{\alpha\circ\xi_0(v),\sigma_v},\sigma_v(f(v))u_{\alpha\circ\xi(v),\sigma_v}\right\rangle_{U_{\sigma_v}}
\end{equation}
for every $f\in F^V$. Alternatively, we have
\begin{equation}\label{defSxisigmabis}
\mathcal{S}_{\xi,\sigma}(\varphi)=\prod_{v\in V}\mathcal{S}_{\alpha\circ\xi(v),\alpha\circ\xi_0(v),\sigma_v}(\varphi(v))
\end{equation}
for every $\varphi\in\mathsf{Orb}_{F^V}\left(\varphi_{\xi}\right)$. Note also that, for $a \neq b$, the individual factor $\mathcal{S}_{a,b,\sigma_*}$ depends, up to a phase factor, on the choice of the unit vectors $u_{a,\sigma_*}$ and $u_{b,\sigma_*}$. However, for every $\xi = k\xi_0 \in \mathsf{Orb}_{\mathsf{I}(\sigma)}(\xi_0)$, we have $\sigma_{k^{-1}v} = \sigma_v$ and $\xi(v) = \xi_0(k^{-1}v)$, so that the multisets $\{(\alpha\circ\xi(v),\sigma_v)\}_{v\in V}$ and $\{(\alpha\circ\xi_0(v),\sigma_v)\}_{v\in V}$ coincide; hence, all phase factors cancel out exactly in the product \eqref{defSxisigma}.

If $\theta\in\widehat{\mathsf{I}(\sigma)}_0$ and $w_{\sigma,\theta}$ is a unit $\mathsf{Stab}_{\mathsf{I}(\sigma)}\left(\xi_0\right)$-invariant vector in $W_\theta$, we may assume that the operator in \eqref{decLJ} is given by $\left[N_{\sigma,\theta}w\right]\left(k\xi_0\right)=\left\langle w,\theta(k)w_{\sigma,\theta}\right\rangle_{W_\theta}$ for all $k\in\mathsf{I}(\sigma)$, $w\in W_\theta$ (omitting $l=1$ and the normalization coefficient). Moreover, the corresponding spherical function is 
\begin{equation}\label{defRsigmatheta}
\mathcal{R}_{\sigma,\theta}\coloneqq N_{\sigma,\theta}w_{\sigma,\theta},
\end{equation}
that is $\mathcal{R}_{\sigma,\theta}\left(k\xi_0\right)=\left\langle w_{\sigma,\theta},\theta(k)w_{\sigma,\theta}\right\rangle_{W_\theta}$ for every $k\in\mathsf{I}(\sigma)$ (see \cite[Corollary 4.6.4]{book}).

\begin{lemma}\label{Lsphfunct}
Suppose that $\sigma \in \Sigma_0$, $\theta \in \widehat{\mathsf{I}(\sigma)}_0$ and $\xi \in \mathsf{Orb}_{\mathsf{I}(\sigma)}(\xi_0)$. 
\begin{enumerate}[(1)]
\item\label{Lsphfunct1}
If $k'\in\mathsf{Stab}_{\mathsf{I}(\sigma)}\left(\xi_0\right)$, then $\lambda\left(k'\right)\mathcal{R}_{\sigma,\theta}=\mathcal{R}_{\sigma,\theta}$ and $\mathcal{S}_{k'^{-1}\xi,\sigma}=\lambda\left(\mathbf{1}_F,k'^{-1}\right)\mathcal{S}_{\xi,\sigma}$.

\item\label{Lsphfunct2}
The function $\mathcal{S}_{\xi,\sigma}$ is $\mathsf{Stab}_{F^V}\left(\varphi_{\xi_0}\right)$-invariant, that is if $f'\in\mathsf{Stab}_{F^V}\left(\varphi_{\xi_0}\right)$, then $\lambda\left(f',1_G\right)\mathcal{S}_{\xi,\sigma}=\mathcal{S}_{\xi,\sigma}$.
\end{enumerate}
\end{lemma}
\begin{proof}
(1)
The function $\mathcal{R}_{\sigma,\theta}$ is $\mathsf{Stab}_{\mathsf{I}(\sigma)}\left(\xi_0\right)$-invariant because it is a spherical function of the Gelfand pair $\left(\mathsf{I}(\sigma),\mathsf{Stab}_{\mathsf{I}(\sigma)}\left(\xi_0\right)\right)$. Furthermore,
\begin{align*}
\mathcal{S}_{k'^{-1}\xi,\sigma}\left(f\cdot\varphi_{k'^{-1}\xi}\right)=&\prod_{v\in V}\left\langle u_{\alpha\circ\xi_0(v),\sigma_v},\sigma_v(f(v))u_{\alpha\circ\xi\left(k'v\right),\sigma_v}\right\rangle_{U_{\sigma_v}}&&(\text{by }\eqref{defSxisigma})\\
=&\prod_{v\in V}\left\langle u_{\alpha\circ\xi_0(v),\sigma_v},\sigma_v\left(f\left(k'^{-1}v\right)\right)u_{\alpha\circ\xi(v),\sigma_v}\right\rangle_{U_{\sigma_v}}&&(k'\in\mathsf{Stab}_{\mathsf{I}(\sigma)}\left(\xi_0\right))\\
=&\mathcal{S}_{\xi,\sigma}\left(k'\left(f\cdot\varphi_{k'^{-1}\xi}\right)\right)&&(k'\varphi_{k'^{-1}\xi}=\varphi_\xi)\\
=&\left[\lambda\left(\mathbf{1}_F,k'^{-1}\right)\mathcal{S}_{\xi,\sigma}\right]\left(f\cdot\varphi_{k'^{-1}\xi}\right).&&
\end{align*}

\medskip
(2)
\[
\begin{split}
\left[\lambda\left(f',1_G\right)\mathcal{S}_{\xi,\sigma}\right]\left(f\cdot\varphi_\xi\right)=&\prod_{v\in V}\left\langle u_{\alpha\circ\xi_0(v),\sigma_v},\sigma_v\left(f'(v)^{-1}f(v)\right)u_{\alpha\circ\xi(v),\sigma_v}\right\rangle_{U_{\sigma_v}}\\
=&\prod_{v\in V}\left\langle\sigma_v\left(f'(v)\right)u_{\alpha\circ\xi_0(v),\sigma_v},\sigma_v(f(v))u_{\alpha\circ\xi(v),\sigma_v}\right\rangle_{U_{\sigma_v}}\\
(f'(v)\in\mathsf{Stab}_F\left(y_{\alpha\circ\xi_0(v)}\right))\quad=&\prod_{v\in V}\left\langle u_{\alpha\circ\xi_0(v),\sigma_v},\sigma_v(f(v))u_{\alpha\circ\xi(v),\sigma_v}\right\rangle_{U_{\sigma_v}}\\
(\text{by } \eqref{defSxisigma})\quad=&\mathcal{S}_{\xi,\sigma}\left(f\cdot\varphi_\xi\right).
\end{split}
\]
\end{proof}

Now we suppose that each $q\in\mathsf{Q}_\sigma$ representing a class in $\mathsf{Stab}_G\left(\xi_0\right)\cdot\mathsf{I}(\sigma)$ belongs to $\mathsf{Stab}_G\left(\xi_0\right)$ and then we denote by $\mathsf{Q}^0_\sigma$ the set of all these representatives. In such a way,
$\mathsf{Q}^0_\sigma$ is the subset of $\mathsf{Q}_\sigma$ defined by the conditions  
\begin{equation}\label{defQ0sig}
\mathsf{Stab}_G\left(\xi_0\right)\cdot\mathsf{I}(\sigma)=\bigsqcup_{q\in\mathsf{Q}^0_\sigma}q\mathsf{I}(\sigma)\qquad\text{ and }\qquad\mathsf{Q}^0_\sigma\subseteq\mathsf{Stab}_G\left(\xi_0\right).
\end{equation}

\begin{lemma}
With the notation in \eqref{decOrbxi} and \eqref{defQ0sig}, we have
\begin{equation}\label{StabOrbIsigma}
\mathsf{Stab}_G\left(\xi_0\right)\mathsf{Orb}_{\mathsf{I}(\sigma)}\left(\xi_0\right)=\bigsqcup_{\Delta_j\cap\mathsf{Orb}_{\mathsf{I}(\sigma)}\left(\xi_0\right)\neq\emptyset}\Delta_j=\bigcup_{q\in\mathsf{Q}^0_\sigma} q\mathsf{Orb}_{\mathsf{I}(\sigma)}\left(\xi_0\right),
\end{equation}
where the second union is in general not disjoint (whenever $\vert{}\mathsf{Q}^0_\sigma\vert{} > 1$, since $\xi_0$ belongs to each subset). Moreover,
\begin{equation}\label{StabStabIsigma}
\mathsf{Stab}_G\left(\xi_0\right)=\bigsqcup_{q\in\mathsf{Q}^0_\sigma} q\mathsf{Stab}_{\mathsf{I}(\sigma)}\left(\xi_0\right).
\end{equation}
\end{lemma}
\begin{proof}
The first set in \eqref{StabOrbIsigma} is clearly $\mathsf{Stab}_G\left(\xi_0\right)$-invariant, so it is the disjoint union of some of the sets $\Delta_j$ in \eqref{decOrbxi}. Moreover, $\mathsf{Stab}_G\left(\xi_0\right)\mathsf{Orb}_{\mathsf{I}(\sigma)}\left(\xi_0\right)\cap\Delta_j\neq\emptyset$ if and only if there exist $h\in\mathsf{Stab}_G\left(\xi_0\right)$ and $k\in\mathsf{I}(\sigma)$ such that $hk\xi_0\in\Delta_j$, which holds if and only if there exists $k\in\mathsf{I}(\sigma)$ such that $k\xi_0\in\Delta_j$ (since $\Delta_j$ is $\mathsf{Stab}_G(\xi_0)$-invariant), or equivalently $\Delta_j\cap\mathsf{Orb}_{\mathsf{I}(\sigma)}\left(\xi_0\right)\neq\emptyset$. The second identity in \eqref{StabOrbIsigma} follows immediately from the definition of 
$\mathsf{Q}^0_\sigma$ in \eqref{defQ0sig}.

To prove \eqref{StabStabIsigma}, observe that for every $h\in\mathsf{Stab}_G\left(\xi_0\right)$, there exist a unique $q\in\mathsf{Q}^0_\sigma$ and a unique $k\in\mathsf{I}(\sigma)$ such that $h=qk$. From the second condition in \eqref{defQ0sig}, it follows that $k=q^{-1}h\in\mathsf{Stab}_G\left(\xi_0\right)$, which implies that $k\in \mathsf{Stab}_{\mathsf{I}(\sigma)}\left(\xi_0\right)$.
\end{proof}

\begin{remark}\label{Remkh}
In general, the product in \eqref{defQ0sig} is not a subgroup. For instance, there may exist $h\in\mathsf{Stab}_G\left(\xi_0\right)$ and $k\in\mathsf{I}(\sigma)$ such that $kh\notin\mathsf{Stab}_G\left(\xi_0\right)\mathsf{I}(\sigma)$; we will show this in Example \ref{Example3ter}. Note also that in this situation, if we choose $q=kh$ as a representative of the left coset containing it, then $q\xi_0=k\xi_0\in\mathsf{Orb}_{\mathsf{I}(\sigma)}\left(\xi_0\right)$. Consequently,  $q\mathsf{Orb}_{\mathsf{I}(\sigma)}\left(\xi_0\right)$ has a non-trivial intersection with the set in \eqref{StabOrbIsigma}, but $q\in\mathsf{Q}_\sigma\setminus\mathsf{Q}^0_\sigma$. On the other hand, in certain cases, the product in \eqref{defQ0sig} is a subgroup; see Example \ref{Example2ter}.
\end{remark}

\begin{theorem}\label{TheoFsigthet}
For $\sigma\in\Sigma_0$ and $\theta\in\widehat{\mathsf{I}(\sigma)}_0$ define $\mathcal{F}_{\sigma,\theta}\in L\bigl(\mathsf{Orb}_{F\wr G}\left(\xi_0,\varphi_{\xi_0}\right)\bigr)$ by setting
\begin{equation}\label{defFsigthet}
\mathcal{F}_{\sigma,\theta}(\xi,\varphi)\coloneqq
\begin{cases}
\frac{1}{\left\lvert\mathsf{Q}^0_\sigma\right\rvert}\!\!\!\!\!\!\!\!\sum\limits_{\substack{q\in\mathsf{Q}^0_\sigma\mathrel{:}\\ q^{-1}\xi\in\mathsf{Orb}_{\mathsf{I}(\sigma)}\left(\xi_0\right)}}\!\!\!\!\!\!\!\!\!\mathcal{R}_{\sigma,\theta}\left(q^{-1}\xi\right)\mathcal{S}_{q^{-1}\xi,\sigma}\left(q^{-1}\varphi\right)&\text{if } \xi\in \mathsf{Stab}_G\left(\xi_0\right)\mathsf{Orb}_{\mathsf{I}(\sigma)}\left(\xi_0\right),\varphi\in\Phi_\xi,\\
\qquad\qquad\qquad\qquad0&\text{otherwise}
\end{cases}
\end{equation}
for all $(\xi,\varphi)\in\mathsf{Orb}_{F\wr G}\left(\xi_0,\varphi_{\xi_0}\right)$.
Then $\mathcal{F}_{\sigma,\theta}$ is the spherical function corresponding to the irreducible representation 
\[
\Ind^{F\wr G}_{F\wr \mathsf{I}(\sigma)}\left(\widetilde{\sigma}\otimes\overline{\theta}\right).
\]
\end{theorem}
\begin{proof}
Fix $\xi\in \mathsf{Stab}_G\left(\xi_0\right)\mathsf{Orb}_{\mathsf{I}(\sigma)}\left(\xi_0\right)$, $\varphi\in\Phi_\xi$, as in \eqref{defFsigthet}, and also $(f,h)\in\mathsf{Stab}_{F\wr G}\left(\xi_0,\varphi_{\xi_0}\right)$. From \eqref{StabStabIsigma} we deduce that for every $q\in \mathsf{Q}^0_\sigma$ there exist a unique $q'\in\mathsf{Q}^0_\sigma$ and a unique $k'_q\in\mathsf{Stab}_{\mathsf{I}(\sigma)}\left(\xi_0\right)$ such that 
\begin{equation}\label{h1qqpkpq}
h^{-1}q=q'k'_q.
\end{equation} 
Moreover, $q\mapsto q'$ is a well-defined bijection from $\mathsf{Q}^0_\sigma$ to $\mathsf{Q}^0_\sigma$, which depends on $h$. Observe that, for every 
$q\in \mathsf{Q}^0_\sigma$ such that  $\xi\in q\mathsf{Orb}_{\mathsf{I}(\sigma)}\left(\xi_0\right)$, there exists $k_q\in\mathsf{I}(\sigma)$ such that $\xi=qk_q\xi_0$. Applying \eqref{h1qqpkpq}, we deduce that
\begin{equation}\label{h1xiqp}
h^{-1}\xi=q'k_q'k_q\xi_0\in q'\mathsf{Orb}_{\mathsf{I}(\sigma)}\left(\xi_0\right).
\end{equation}
Similarly, if $h^{-1}\xi\in q'\mathsf{Orb}_{\mathsf{I}(\sigma)}\left(\xi_0\right)$, by \eqref{h1qqpkpq} we obtain 
\[
q^{-1}\xi=k'^{-1}_qq'^{-1}h^{-1}\xi\in\mathsf{Orb}_{\mathsf{I}(\sigma)}\left(\xi_0\right).
\]
This proves that, under the bijection $q\mapsto q'$, we have
\begin{equation}\label{qp1hxiq1xi}
q'^{-1}h^{-1}\xi\in\mathsf{Orb}_{\mathsf{I}(\sigma)}\left(\xi_0\right)	\;\Longleftrightarrow\; q^{-1}\xi\in\mathsf{Orb}_{\mathsf{I}(\sigma)}\left(\xi_0\right).
\end{equation} 

It follows that  
\[
\begin{split}
&\left[\lambda(f,h)\mathcal{F}_{\sigma,\theta}\right](\xi,\varphi)=\mathcal{F}_{\sigma,\theta}\left(h^{-1}\xi,h^{-1}\left(f^{-1}\cdot\varphi\right)\right)\quad(\text{by }\eqref{explpermrep})\\
&\quad=\frac{1}{\left\lvert\mathsf{Q}^0_\sigma\right\rvert}\!\!\!\!\sum_{\substack{q\in\mathsf{Q}^0_\sigma\mathrel{:}\\ q^{-1}\xi\in\mathsf{Orb}_{\mathsf{I}(\sigma)}\left(\xi_0\right)}}\!\!\!\!\!\!\!\mathcal{R}_{\sigma,\theta}\left(q'^{-1}h^{-1}\xi\right)\mathcal{S}_{q'^{-1}h^{-1}\xi,\sigma}\left(q'^{-1}h^{-1}\left(f^{-1}\cdot\varphi\right)\right)\quad(\text{by }\eqref{defFsigthet}\text{ and }\eqref{qp1hxiq1xi})\\
&\quad=\frac{1}{\left\lvert\mathsf{Q}^0_\sigma\right\rvert}\sum_{\substack{q\in\mathsf{Q}^0_\sigma\mathrel{:}\\ q^{-1}\xi\in\mathsf{Orb}_{\mathsf{I}(\sigma)}\left(\xi_0\right)}}\mathcal{R}_{\sigma,\theta}\left(k'_qk_q\xi_0\right)\mathcal{S}_{k'_qk_q\xi_0,\sigma}\left(k'_qq^{-1}\left(f^{-1}\cdot\varphi\right)\right)\quad(\text{by }\eqref{h1qqpkpq}\text{ and }\eqref{h1xiqp})\\
&\quad=\frac{1}{\left\lvert\mathsf{Q}^0_\sigma\right\rvert}\sum_{\substack{q\in\mathsf{Q}^0_\sigma\mathrel{:}\\ q^{-1}\xi\in\mathsf{Orb}_{\mathsf{I}(\sigma)}\left(\xi_0\right)}}\mathcal{R}_{\sigma,\theta}\left(k_q\xi_0\right)\mathcal{S}_{k\xi_0,\sigma}\left(\left(q^{-1}f^{-1}\right)\cdot\left(q^{-1}\varphi\right)\right)\quad(\text{by Lemma }\ref{Lsphfunct}\eqref{Lsphfunct1})\\
&\quad=\frac{1}{\left\lvert\mathsf{Q}^0_\sigma\right\rvert}\!\!\!\!\sum_{\substack{q\in\mathsf{Q}^0_\sigma\mathrel{:}\\ q^{-1}\xi\in\mathsf{Orb}_{\mathsf{I}(\sigma)}\left(\xi_0\right)}}\!\!\!\!\!\!\!\mathcal{R}_{\sigma,\theta}\left(k_q\xi_0\right)\mathcal{S}_{k\xi_0,\sigma}\left(q^{-1}\varphi\right)\quad(\text{by Lemma }\ref{Lsphfunct}\eqref{Lsphfunct2},\text{ since } q^{-1}f\in\mathsf{Stab}_{F^V}\left(\varphi_{\xi_0}\right))\\
&\quad=\frac{1}{\left\lvert\mathsf{Q}^0_\sigma\right\rvert}\sum_{\substack{q\in\mathsf{Q}^0_\sigma\mathrel{:}\\ q^{-1}\xi\in\mathsf{Orb}_{\mathsf{I}(\sigma)}\left(\xi_0\right)}}\mathcal{R}_{\sigma,\theta}\left(q^{-1}\xi\right)\mathcal{S}_{q^{-1}\xi,\sigma}\left(q^{-1}\varphi\right)\quad(\xi=qk_q\xi_0)\\
&\quad=\mathcal{F}_{\sigma,\theta}(\xi,\varphi)\quad(\text{by }\eqref{defFsigthet}).
\end{split}
\]
Therefore, $\mathcal{F}_{\sigma,\theta}$ is $\mathsf{Stab}_{F\wr G}\left(\xi_0,\varphi_{\xi_0}\right)$-invariant. 

By the second condition in \eqref{defQ0sig}, we have $q^{-1}\xi_0=\xi_0$ and $q^{-1}\varphi_{\xi_0}=\varphi_{\xi_0}$ for every $q\in\mathsf{Q}^0_\sigma$, so that
\[
\mathcal{F}_{\sigma,\theta}\left(\xi_0,\varphi_{\xi_0}\right)=\mathcal{R}_{\sigma,\theta}\left(\xi_0\right)\mathcal{S}_{\xi_0,\sigma}\left(\varphi_{\xi_0}\right)=\left\lVert w_{\sigma,\theta}\right\rVert^2\prod_{v\in V}\left\lVert u_{\alpha\circ\xi_0(v),\sigma_v}\right\rVert^2=1.
\]
Finally, recalling that $u_\sigma=\bigotimes_{v\in V}u_{\alpha\circ\xi_0(v),\sigma_v}$, from \eqref{defNflat} and Corollary \ref{Cormultfree}\eqref{Cormultfree1}, we deduce that
\[
\begin{split}
\left[N_{\sigma,\theta}^\flat\left(\sum_{q\in\mathsf{Q}^0_\sigma}H_{q,u_\sigma\otimes w_{\sigma,\theta}}\right)\right](\xi,\varphi)=&\sum_{\substack{q\in\mathsf{Q}^0_\sigma\mathrel{:}\\ q^{-1}\xi\in\mathsf{Orb}_{\mathsf{I}(\sigma)}\left(\xi_0\right)}}\left[N_{\sigma,\theta}w_{\sigma,\theta}\right]\left(q^{-1}\xi\right)\left[T_{q^{-1}\xi,\sigma}u_\sigma\right]\left(q^{-1}\varphi\right)\\
(\text{by }\eqref{defSxisigma} \text{ and }\eqref{defRsigmatheta})\quad=&\sum_{\substack{q\in\mathsf{Q}^0_\sigma\mathrel{:}\\ q^{-1}\xi\in\mathsf{Orb}_{\mathsf{I}(\sigma)}\left(\xi_0\right)}}\mathcal{R}_{\sigma,\theta}\left(q^{-1}\xi\right)\mathcal{S}_{q^{-1}\xi,\sigma}\left(q^{-1}\varphi\right)\\
(\text{by }\eqref{defFsigthet})\quad=&\left\lvert\mathsf{Q}^0_\sigma\right\rvert \mathcal{F}_{\sigma,\theta}(\xi,\varphi).
\end{split}
\]
Note also that, for $\xi \notin \mathsf{Stab}_G(\xi_0)\mathsf{Orb}_{\mathsf{I}(\sigma)}(\xi_0)$, the $\mathsf{Stab}_G(\xi_0)$-invariance of $\mathsf{Stab}_G(\xi_0)\mathsf{Orb}_{\mathsf{I}(\sigma)}(\xi_0)$ implies that $h^{-1}\xi \notin \mathsf{Stab}_G(\xi_0)\mathsf{Orb}_{\mathsf{I}(\sigma)}(\xi_0)$ (therefore both members of the invariance equation are $0$); analogously, there is not any $q \in \mathsf{Q}^0_\sigma$ with $q^{-1}\xi \in \mathsf{Orb}_{\mathsf{I}(\sigma)}(\xi_0)$, so that the sum in $N^\flat_{\sigma,\theta}$ is empty.
Hence, $\mathcal{F}_{\sigma,\theta}=\frac{1}{\left\lvert\mathsf{Q}^0_\sigma\right\rvert}N_{\sigma,\theta}^\flat\left(\sum_{q\in\mathsf{Q}^0_\sigma}H_{q,u_\sigma\otimes w_{\sigma,\theta}}\right)$, and Lemma
\ref{Theointoper} implies that $\mathcal{F}_{\sigma,\theta}$ belongs to the subspace of $L\bigl(\mathsf{Orb}_{F\wr G}\left(\xi_0,\varphi_{\xi_0}\right)\bigr)$ equivalent to $\Ind^{F\wr G}_{F\wr \mathsf{I}(\sigma)}\left(\widetilde{\sigma}\otimes\overline{\theta}\right)$. This completes the proof.
\end{proof}

\begin{corollary}
For $\xi\in \mathsf{Stab}_G\left(\xi_0\right)\mathsf{Orb}_{\mathsf{I}(\sigma)}\left(\xi_0\right)$ and $\varphi\in\Phi_\xi$, the value of $\mathcal{F}_{\sigma,\theta}(\xi,\varphi)$ can also be written in the form
\begin{equation}\label{defFsigthetbis}
\mathcal{F}_{\sigma,\theta}(\xi,\varphi)=
\frac{1}{\left\lvert\mathsf{Stab}_G\left(\xi_0\right)\right\rvert}\sum\limits_{\substack{h\in\mathsf{Stab}_G\left(\xi_0\right)\mathrel{:}\\ h^{-1}\xi\in\mathsf{Orb}_{\mathsf{I}(\sigma)}\left(\xi_0\right)}}\mathcal{R}_{\sigma,\theta}\left(h^{-1}\xi\right)\mathcal{S}_{h^{-1}\xi,\sigma}\left(h^{-1}\varphi\right).
\end{equation}
\end{corollary}
\begin{proof}
Indeed, \eqref{StabStabIsigma} allows us to write each $h\in\mathsf{Stab}_G\left(\xi_0\right)$ in the form $h=qk$ with $k\in\mathsf{Stab}_{\mathsf{I}(\sigma)}\left(\xi_0\right)$. This rearrangement uses the fact that $\mathcal{R}_{\sigma,\theta}$ is $\mathsf{Stab}_{\mathsf{I}(\sigma)}\left(\xi_0\right)$-invariant, while from Lemma \ref{Lsphfunct}\eqref{Lsphfunct1} we deduce that $\mathcal{S}_{k^{-1}q^{-1}\xi,\sigma}\left(k^{-1}q^{-1}\varphi\right)=\mathcal{S}_{q^{-1}\xi,\sigma}\left(q^{-1}\varphi\right)$. The summation condition $h^{-1}\xi \in \mathsf{Orb}_{\mathsf{I}(\sigma)}\left(\xi_0\right)$ is equivalent to $q^{-1}\xi \in \mathsf{Orb}_{\mathsf{I}(\sigma)}\left(\xi_0\right)$ (since $k \in \mathsf{I}(\sigma)$), and therefore each summand corresponding to an admissible $q$ is repeated exactly $\left\lvert\mathsf{Stab}_{\mathsf{I}(\sigma)}\left(\xi_0\right)\right\rvert$ times.
\end{proof}

Since the spherical functions $\mathcal{F}_{\sigma,\theta}$ are constant on the orbits of $\mathsf{Stab}_{F\wr G}\left(\xi_0,\varphi_{\xi_0}\right)$, our next goal is to determine their expressions as functions of the parameters $j$ and $r$ which are used in \eqref{decOrbxiphi}. Taking into account that the spherical functions are identically zero outside the set in \eqref{StabOrbIsigma}, we assume that the orbits in \eqref{decOrbxi} are numbered in such a way that the first identity in \eqref{StabOrbIsigma} becomes 
\begin{equation}\label{StabOrbIsig0}
\mathsf{Stab}_G\left(\xi_0\right)\mathsf{Orb}_{\mathsf{I}(\sigma)}\left(\xi_0\right)=\bigsqcup_{j=0}^{n^0}\Delta_j,
\end{equation}
for a suitable $n^0\leq n$.
Then note that the function $\mathcal{S}_{a,b,\sigma_*}$ in \eqref{defSabsig} is $\mathsf{Stab}_F\left(y_b\right)$-invariant and belongs to $L\left(Y_a\right)$. This also means that it corresponds to a $\mathsf{Stab}_F\left(y_b\right)$-$\mathsf{Stab}_F\left(y_a\right)$-invariant function defined on $F$; see, e.g., \cite[Exercise 4.2.5]{book} and \cite{Du}. In any case, it is constant on the orbits in \eqref{decYa} and we may set
\begin{equation}\label{defcheckSsing}
\check{\mathcal{S}}_{a,b,\sigma_*}(i)\coloneqq\mathcal{S}_{a,b,\sigma_*}(y)\quad\text{if }y\in\Upsilon_{a,b,i},
\end{equation}
for $i=0,1,\dotsc,\ell_{a,b}$. With the notation in \eqref{defPsixij}, assuming $(\xi,\mathsf{i})\in\bigsqcup_{j=0}^{n^0}\Psi_j$ with $\xi\in\mathsf{Orb}_{\mathsf{I}(\sigma)}(\xi_0)$, and $\varphi\in\bigtimes\limits_{v\in V}\Upsilon_{\alpha\circ\xi(v),\alpha\circ\xi_0(v),\mathsf{i}(v)}$, from \eqref{defSxisigmabis} and \eqref{defcheckSsing} we deduce that
\begin{equation}\label{defcheckSsingb}
\mathcal{S}_{\xi,\sigma}(\varphi)=\prod_{v\in V}\check{\mathcal{S}}_{\alpha\circ\xi(v),\alpha\circ\xi_0(v),\sigma_v}(\mathsf{i}(v))\eqqcolon\check{\mathcal{S}}_{\xi,\sigma}(\mathsf{i}).
\end{equation}
That is, the function $\mathcal{S}_{\xi,\sigma}(\varphi)$ depends only on the index $\mathsf{i}$ corresponding to $\varphi$, so $\check{\mathcal{S}}_{\xi,\sigma}$ is well-defined. Moreover, if $h\in\mathsf{Stab}_G\left(\xi_0\right)$ and $h^{-1}\xi\in\mathsf{Orb}_{\mathsf{I}(\sigma)}(\xi_0)$, then $h^{-1}\varphi\in\bigtimes\limits_{v\in V}\Upsilon_{\alpha\circ\left(h^{-1}\xi\right)(v),\alpha\circ\xi_0(v),\left(h^{-1}\mathsf{i}\right)(v)}$, which implies that
\begin{equation}\label{h1barS}
\mathcal{S}_{h^{-1}\xi,\sigma}\left(h^{-1}\varphi\right)=\check{\mathcal{S}}_{h^{-1}\xi,\sigma}\left(h^{-1}\mathsf{i}\right).
\end{equation}
Similarly, suppose that 
\begin{equation}\label{decPil}
\mathsf{Orb}_{\mathsf{I}(\sigma)}\left(\xi_0\right)=\bigsqcup_{l=0}^{p_\sigma}\Pi_l
\end{equation}
is the decomposition into the orbits of $\mathsf{Stab}_{\mathsf{I}(\sigma)}\left(\xi_0\right)$, with $\Pi_0=\left\{\xi_0\right\}$. 

\begin{lemma}\label{Lemmaajl}
\begin{enumerate}[(1)]
\item\label{Lemmaajl1}
If $0\leq j\leq n^0$, $\xi\in\Delta_j$, and $\Pi_l\subseteq\Delta_j$, then the number $a_{j,l}$ of elements $q\in\mathsf{Q}^0_\sigma$ such that $q^{-1}\xi\in\Pi_l$ is positive and depends only on $j$ and $l$.
\item\label{Lemmaajl2}
If $q,q''\in\mathsf{Q}^0_\sigma$ and $q^{-1}\xi,q''^{-1}\xi\in\Pi_l$, then there exists $k\in\mathsf{Stab}_{\mathsf{I}(\sigma)}\left(\xi_0\right)$ such that $k^{-1}q''^{-1}\xi=q^{-1}\xi$.
\end{enumerate}
\end{lemma}
\begin{proof}
(1)
First, note that there exists $h\in\mathsf{Stab}_G\left(\xi_0\right)$ such that $h^{-1}\xi\in\Pi_l$. From \eqref{StabStabIsigma} we deduce that there exist $q\in\mathsf{Q}^0_\sigma$ and $k\in\mathsf{Stab}_{\mathsf{I}(\sigma)}\left(\xi_0\right)$ such that $h=qk$. Therefore, 
\[
h^{-1}\xi=k^{-1}q^{-1}\xi\in\Pi_l\quad\text{and}\quad k\Pi_l=\Pi_l\qquad\Longrightarrow\qquad q^{-1}\xi\in\Pi_l,
\]
which implies that for all $\xi$ and $l$, there exists at least one $q$ such that $q^{-1}\xi\in\Pi_l$. For any other $\xi'\in\Delta_j$, there exists $h'\in\mathsf{Stab}_G\left(\xi_0\right)$ such that $h'\xi=\xi'$. With this $h'$ in mind, for every $q\in\mathsf{Q}^0_\sigma$ there exist $q'\in\mathsf{Q}^0_\sigma$ and $k'\in\mathsf{Stab}_{\mathsf{I}(\sigma)}\left(\xi_0\right)$ such that $h'q=q'k'$. Then 
\[
k'^{-1}q'^{-1}\xi'=q^{-1}h'^{-1}\xi'=q^{-1}\xi,
\]
yielding
\[
q^{-1}\xi\in\Pi_l \iff q'^{-1}\xi'\in\Pi_l.
\]
Hence, we can invoke the bijection $q\mapsto q'$ in the proof of Theorem \ref{TheoFsigthet} (with $h$ replaced by $h'^{-1}$)  to complete the proof of \eqref{Lemmaajl1}. 

\medskip
(2)
Statement \eqref{Lemmaajl2} is immediate.
\end{proof}

In the following lemma, similarly to \cite[Section 5]{CSSTRam}, we introduce an intermediate set of functions which, in the present case, are sums of spherical functions of $F^V$, while in \cite{CSSTRam} they were sums of characters.

\begin{lemma}\label{LemmabarbarS}
In the notation of \eqref{decPsixij}, \eqref{StabOrbIsig0}, and \eqref{decPil}, and taking into account Lemma \ref{Lemmaajl}, suppose that $0\leq j\leq n^0$, $0\leq r\leq s_j$, $(\xi,\mathsf{i})\in\Psi_{j,r}$, $\Pi_l\subseteq\Delta_j$, and $\varphi\in\bigtimes\limits_{v\in V}\Upsilon_{\alpha\circ\xi(v),\alpha\circ\xi_0(v),\mathsf{i}(v)}$. Then
\begin{equation}\label{barbarcalS}
\sum_{\substack{q\in\mathsf{Q}^0_\sigma\mathrel{:}\\ q^{-1}\xi\in\Pi_l}}\mathcal{S}_{q^{-1}\xi,\sigma}\left(q^{-1}\varphi\right)=\frac{1}{\left\lvert\mathsf{Stab}_{\mathsf{I}(\sigma)}\left(\xi_0\right)\right\rvert}\sum_{\xi'\in\Pi_l}\;\sum_{\substack{h\in\mathsf{Stab}_G\left(\xi_0\right)\mathrel{:}\\ h^{-1}\xi=\xi'}}\check{\mathcal{S}}_{\xi',\sigma}\left(h^{-1}\mathsf{i}\right)\eqqcolon\check{\check{\mathcal{S}}}_{j,l,\sigma}(r).
\end{equation}
In particular, the values of these sums depend only on $j,l,r$. Consequently, $\check{\check{\mathcal{S}}}_{j,l,\sigma}(r)$ is well-defined.
\end{lemma}
\begin{proof}
First, we observe that, by \eqref{StabStabIsigma}, every $h\in\mathsf{Stab}_G\left(\xi_0\right)$ may be uniquely written in the form $h=qk$, with $q\in\mathsf{Q}^0_\sigma$ and $k\in\mathsf{Stab}_{\mathsf{I}(\sigma)}\left(\xi_0\right)$. Moreover,
\begin{equation}\label{h1xiPil}
h^{-1}\xi\in\Pi_l\quad \Longleftrightarrow\quad q^{-1}\xi\in\Pi_l
\end{equation} 
because $k\in\mathsf{Stab}_{\mathsf{I}(\sigma)}(\xi_0)$ stabilizes $\Pi_l$.
Henceforth, by Lemma \ref{Lsphfunct}\eqref{Lsphfunct1}, we have
\[
\begin{split}
\sum_{\substack{q\in\mathsf{Q}^0_\sigma\mathrel{:}\\ q^{-1}\xi\in\Pi_l}}\mathcal{S}_{q^{-1}\xi,\sigma}\left(q^{-1}\varphi\right)=&
\frac{1}{\left\lvert\mathsf{Stab}_{\mathsf{I}(\sigma)}\left(\xi_0\right)\right\rvert}\sum_{\substack{q\in\mathsf{Q}^0_\sigma\mathrel{:}\\ q^{-1}\xi\in\Pi_l}}\sum_{k\in\mathsf{Stab}_{\mathsf{I}(\sigma)}\left(\xi_0\right)}\mathcal{S}_{k^{-1}q^{-1}\xi,\sigma}\left(k^{-1}q^{-1}\varphi\right)\\
(\text{by }\eqref{StabStabIsigma}, \eqref{h1barS}\text{ and }\eqref{h1xiPil})\quad=&
\frac{1}{\left\lvert\mathsf{Stab}_{\mathsf{I}(\sigma)}\left(\xi_0\right)\right\rvert}\sum_{\substack{h\in\mathsf{Stab}_G\left(\xi_0\right)\mathrel{:}\\ h^{-1}\xi\in\Pi_l}}\check{\mathcal{S}}_{h^{-1}\xi,\sigma}\left(h^{-1}\mathsf{i}\right)\\
(\xi'= h^{-1}\xi)\quad=&
\frac{1}{\left\lvert\mathsf{Stab}_{\mathsf{I}(\sigma)}\left(\xi_0\right)\right\rvert}\sum_{\xi'\in\Pi_l}\;\sum_{\substack{h\in\mathsf{Stab}_G\left(\xi_0\right)\mathrel{:}\\ h^{-1}\xi=\xi'}}\check{\mathcal{S}}_{\xi',\sigma}\left(h^{-1}\mathsf{i}\right).\\
\end{split}
\]
To show that the last double sum depends only on $j,l,r$, suppose also that $\left(\xi_1,\mathsf{i}_1\right)\in\Psi_{j,r}$. Then there exists $h_1\in\mathsf{Stab}_G\left(\xi_0\right)$ such that $h_1^{-1}(\xi,\mathsf{i})=\left(\xi_1,\mathsf{i}_1\right)$. Consequently, for $h'\in\mathsf{Stab}_G\left(\xi_0\right)$ and $\xi'\in\Pi_l$, we have $h'^{-1}\xi_1=\xi'$ if and only if the element $h\coloneqq h_1h'$ satisfies $h^{-1}\xi=\xi'$. Moreover, if $h$ is as above, then $h'^{-1}\mathsf{i}_1=h^{-1}h_1\mathsf{i}_1=h^{-1}\mathsf{i}$. Hence,
\[
\sum_{\xi'\in\Pi_l}\;\sum_{\substack{h'\in\mathsf{Stab}_G\left(\xi_0\right)\mathrel{:}\\ h'^{-1}\xi_1=\xi'}}\check{\mathcal{S}}_{\xi',\sigma}\left(h'^{-1}\mathsf{i}_1\right)=\sum_{\xi'\in\Pi_l}\;\sum_{\substack{h\in\mathsf{Stab}_G\left(\xi_0\right)\mathrel{:}\\ h^{-1}\xi=\xi'}}\check{\mathcal{S}}_{\xi',\sigma}\left(h^{-1}\mathsf{i}\right).
\]
\end{proof}

Since $\mathcal{R}_{\sigma,\theta}$ is $\mathsf{Stab}_{\mathsf{I}(\sigma)}\left(\xi_0\right)$-invariant and taking \eqref{decPil} into account, we may define the function $\check{\mathcal{R}}_{\sigma,\theta}$ by setting $\check{\mathcal{R}}_{\sigma,\theta}(l)\coloneqq\mathcal{R}_{\sigma,\theta}(\xi)$ if $\xi\in\Pi_l$. If $\Delta_j\cap\mathsf{Orb}_{\mathsf{I}(\sigma)}\left(\xi_0\right)\neq\emptyset$, then $\Delta_j\cap\mathsf{Orb}_{\mathsf{I}(\sigma)}\left(\xi_0\right)=\bigsqcup_{\Pi_l\subseteq\Delta_j}\Pi_l$, because $\Pi_l\cap\Delta_j\neq\emptyset$ implies that $\Pi_l\subseteq\Delta_j$. Therefore, 
\begin{equation}\label{defbarR}
\xi\in\Delta_j\cap\mathsf{Orb}_{\mathsf{I}(\sigma)}\left(\xi_0\right)\quad\Longrightarrow\quad\mathcal{R}_{\sigma,\theta}(\xi)=\check{\mathcal{R}}_{\sigma,\theta}(l), \text{ if }\xi\in\Pi_l\subseteq\Delta_j.
\end{equation}

Now we express the spherical functions $\mathcal{F}_{\sigma,\theta}$ in \eqref{defFsigthet} in terms of the parameters $j,r$, that is as constant functions on the orbits in \eqref{defdigjr}.
\begin{theorem}\label{TheodefFsigthetorb}
If  $0\leq j\leq n^0$ and $0\leq r\leq s_j$, then the value of $\mathcal{F}_{\sigma,\theta}$ at $(\xi,\varphi)\in\digamma_{j,r}$ depends only on $j$ and $r$ and is given by
\begin{equation}\label{defFsigthetorb}
\mathcal{F}_{\sigma,\theta}(\xi,\varphi)=\frac{1}{\left\lvert\mathsf{Q}^0_\sigma\right\rvert}\sum_{\substack{0\leq l\leq p_\sigma\mathrel{:}\\\Pi_l\subseteq\Delta_j}}\check{\mathcal{R}}_{\sigma,\theta}(l)\check{\check{\mathcal{S}}}_{j,l,\sigma}(r)\eqqcolon\check{\mathcal{F}}_{\sigma,\theta}(j,r).
\end{equation}
\end{theorem}
\begin{proof}
First note that 
\begin{equation}\label{qQ0xiPi}
\begin{split}
&\xi\in\Delta_j,\;q\in\mathsf{Q}^0_\sigma\subseteq\mathsf{Stab}_G(\xi_0)\quad\text{and}\quad q^{-1}\xi\in\mathsf{Orb}_{\mathsf{I}(\sigma)}\left(\xi_0\right)\\
 &\qquad\qquad\qquad\qquad\Longrightarrow\quad q^{-1}\xi\in\Delta_j\cap\mathsf{Orb}_{\mathsf{I}(\sigma)}(\xi_0)=\bigsqcup_{\Pi_l\subseteq\Delta_j}\Pi_l.
\end{split}
\end{equation}
Then,
\begin{align*}
\mathcal{F}_{\sigma,\theta}(\xi,\varphi)=&
\frac{1}{\left\lvert\mathsf{Q}^0_\sigma\right\rvert}\sum_{\substack{q\in\mathsf{Q}^0_\sigma\mathrel{:}\\ q^{-1}\xi\in\mathsf{Orb}_{\mathsf{I}(\sigma)}\left(\xi_0\right)}}\mathcal{R}_{\sigma,\theta}\left(q^{-1}\xi\right)\mathcal{S}_{q^{-1}\xi,\sigma}\left(q^{-1}\varphi\right)&&(\text{by }\eqref{defFsigthet})\\
=&
\frac{1}{\left\lvert\mathsf{Q}^0_\sigma\right\rvert}\sum_{\substack{0\leq l\leq p_\sigma\mathrel{:}\\ \Pi_l\subseteq\Delta_j}}\check{\mathcal{R}}_{\sigma,\theta}(l)\sum_{\substack{q\in\mathsf{Q}^0_\sigma\mathrel{:}\\ q^{-1}\xi\in\Pi_l}}\mathcal{S}_{q^{-1}\xi,\sigma}\left(q^{-1}\varphi\right)&&\!\!\!\!\!\!\!\!\!\!\!\!(\text{by Lemma }\ref{Lemmaajl}\eqref{Lemmaajl1}, \eqref{defbarR}\text{ and }\eqref{qQ0xiPi})\\
=&
\frac{1}{\left\lvert\mathsf{Q}^0_\sigma\right\rvert}\sum_{\substack{0\leq l\leq p_\sigma\mathrel{:}\\\Pi_l\subseteq\Delta_j}}\check{\mathcal{R}}_{\sigma,\theta}(l)\check{\check{\mathcal{S}}}_{j,l,\sigma}(r),&&(\text{by Lemma }\ref{LemmabarbarS})
\end{align*}
which completes the proof.
\end{proof}

In order to formulate the orthogonality relations, we set $\check{\mathcal{F}}_{\sigma,\theta}(j,r)\coloneqq 0$ when $\Delta_j\not\subseteq\mathsf{Stab}_G(\xi_0)\mathsf{Orb}_{\mathsf{I}(\sigma)}(\xi_0)$. We also assume that the numbering $\Delta_0,\dotsc,\Delta_n$ does not depend on $\sigma$ (we suspend the assumption in \eqref{StabOrbIsig0}).

\begin{proposition}
The spherical functions $\check{\mathcal{F}}_{\sigma,\theta}$ satisfy the following orthogonality relations
\begin{equation}\label{orthrel}
\sum_{j=0}^n\sum_{r=0}^{s_j}\check{\mathcal{F}}_{\sigma,\theta}(j,r)\overline{\check{\mathcal{F}}_{\sigma',\theta'}(j,r)}\left\lvert\digamma_{j,r}\right\rvert=\delta_{\sigma,\sigma'}\delta_{\theta,\theta'}\frac{\left\lvert\mathsf{Orb}_{F\wr G}\left(\xi_0,\varphi_{\xi_0}\right)\right\rvert}{\left\lvert\mathsf{Q}_\sigma\right\rvert\text{\rm dim}W_\theta\,\text{\rm dim}U_\sigma},
\end{equation}
where 
\[
\begin{split}
&\left\lvert\digamma_{j,r}\right\rvert=\!\!\!\!\sum_{(\xi,\mathsf{i})\in\Psi_{j,r}}\prod_{v\in V}\left\lvert\Upsilon_{\alpha\circ\xi(v),\alpha\circ\xi_0(v),\mathsf{i}(v)}\right\rvert,\quad \left\lvert\mathsf{Orb}_{F\wr G}\left(\xi_0,\varphi_{\xi_0}\right)\right\rvert=\left\lvert\mathsf{Orb}_G\left(\xi_0\right)\right\rvert\prod_{a\in A}\left\lvert Y_a\right\rvert^{\left\lvert\left(\alpha\circ\xi_0\right)^{-1}(a)\right\rvert},\\
&\text{\rm dim}U_\sigma=\prod_{v\in V}\text{\rm dim}U_{\sigma_v}.
\end{split}
\]
\end{proposition}
\begin{proof}
This follows from the general orthogonality relations for spherical functions; see, e.g., \cite[Proposition 4.7.1]{book}. The formula for the cardinality of $\digamma_{j,r}$ follows immediately from the explicit formula \eqref{defdigjr}. Similarly, the cardinality of $\mathsf{Orb}_{F\wr G}\left(\xi_0,\varphi_{\xi_0}\right)$ is obtained via \eqref{FVprod} and \eqref{FVorbdec}. Finally, the dimension of the irreducible representation involved can be computed by means of the general formula \cite[(11.10)]{book4}.
\end{proof}

We recall that $(G,K)$ is a symmetric Gelfand pair if $g^{-1}\in KgK$ for every $g\in G$; see, e.g., \cite[Example 4.3.2]{book}. This condition can be formulated in our setting as follows.
\begin{proposition}
The pair $\bigl(F\wr G,\mathsf{Stab}_{F\wr G}\left(\xi_0,\varphi_{\xi_0}\right)\bigr)$ is a symmetric Gelfand pair if and only if for every $(f,g)\in F\wr G$ there exist $\left(f_1,h_1\right),\left(f_2,h_2\right)\in\mathsf{Stab}_{F\wr G}\left(\xi_0,\varphi_{\xi_0}\right)$ such that
\begin{equation}\label{symGel}
\begin{cases}
h_1gh_2=g^{-1},&\\
f_1(v)f\left(h_1^{-1}v\right)f_2\left(g^{-1}h_1^{-1}v\right)=f(gv)^{-1}&\text{for every }v\in V.
\end{cases}
\end{equation}
Moreover, if the conditions in \eqref{symGel} hold, then $\bigl(G,\mathsf{Stab}_G\left(\xi_0\right)\bigr)$ and $\bigl(F,\mathsf{Stab}_F\left(y_a\right)\bigr)$, $a\in A$, are symmetric Gelfand pairs.
\end{proposition}
\begin{proof}
We have $\left(f_1,h_1\right)(f,g)\left(f_2,h_2\right)=(f,g)^{-1}$ if and only if
\[
\left(f_1\cdot h_1f\cdot h_1gf_2,h_1gh_2\right)=\left(g^{-1}f^{-1},g^{-1}\right),
\]
the explicit form of which is exactly \eqref{symGel}. Next, we observe that the first condition in \eqref{symGel} is equivalent to the statement that $\bigl(G,\mathsf{Stab}_G\left(\xi_0\right)\bigr)$ is a symmetric Gelfand pair. If we choose $g=1_G$ and a function $f$ that is constant on each set $\left(\alpha\circ\xi_0\right)^{-1}(a)$, then the second condition in \eqref{symGel} implies that each pair $\bigl(F,\mathsf{Stab}_F\left(y_a\right)\bigr)$ is symmetric. Indeed, since $h_1\in\mathsf{Stab}_G\left(\xi_0\right)$, we have $\alpha\circ\xi_0\left(h_1^{-1}v\right)=\alpha\circ\xi_0(v)=a$ for every $v\in\left(\alpha\circ\xi_0\right)^{-1}(a)$, so that $f\left(h_1^{-1}v\right)=f(v)$ and $f_1(v),f_2\left(h_1^{-1}v\right)\in\mathsf{Stab}_F\left(y_a\right)$; hence, the second condition in \eqref{symGel} reduces to $f_1(v)f(v)f_2\left(h_1^{-1}v\right)=f(v)^{-1}$, where $f(v)\in F$ can be chosen arbitrarily.
\end{proof}

\begin{example}\label{Example1ter}
Now we discuss Example \ref{Example1bis}. First, from the results in Example \ref{Example1}, we deduce that the left-hand side of \eqref{defQ0sig} is the whole group $G$, so $\mathsf{Q}_\sigma^0=\mathsf{Q}_\sigma$, and also that the set in \eqref{StabOrbIsigma} has cardinality 1. Furthermore, we note that the conditions in Proposition \ref{Propmultfree} are verified if and only if $L(Y)$ is multiplicity free as a permutation representation of $F$. Indeed, if this is the case, then both $\mathsf{J}_\sigma$ and $\mathsf{Orb}_{\mathsf{I}(\sigma)}\left(\xi_0\right)$ are singletons for every $\sigma\in\Sigma_0$.  From now on, we assume that this condition is verified. 

Since the set in \eqref{decPil} has cardinality 1, the functions $\mathcal{R}$ in \eqref{defRsigmatheta} are trivial (identically equal to $1$). The functions $\mathcal{S}$ in \eqref{defSxisigma} have the form $\mathcal{S}_\sigma\left(f\cdot\varphi_0\right)=\prod_{v\in V}\left\langle u_{\sigma_v},\sigma_v(f(v))u_{\sigma_v}\right\rangle_{U_{\sigma_v}}$, $f\in F^V$, where $u_{\sigma_v}$ is a unit $\mathsf{Stab}_F\left(y_0\right)$-invariant vector in $U_{\sigma_v}$. Then the spherical functions in \eqref{defFsigthet} and \eqref{defFsigthetbis} are given, respectively, by
\[
\mathcal{F}_\sigma(\varphi)=
\frac{1}{\left\lvert\mathsf{Q}_\sigma\right\rvert}\sum_{q\in\mathsf{Q}_\sigma}\mathcal{S}_\sigma\left(q^{-1}\varphi\right)=\frac{1}{\left\lvert G\right\rvert}\sum_{g\in G}\mathcal{S}_\sigma\left(g^{-1}\varphi\right)
\]
for every $\sigma\in\Sigma_0$ and $\varphi\in Y^V$;
see also \cite[Section 5.4]{CST2}, \cite[Theorem 3.3.18]{book3}. 

With the notation of Example \ref{Example1} and in \eqref{defFsigthetorb}, we can express the spherical functions as constant functions on the orbits, namely
\begin{align*}
&\check{\mathcal{S}}_{\sigma_*}(i)=\mathcal{S}_{\sigma_*}(y),\;\sigma_*\in\widehat{F}_0, y\in\Upsilon_i, 0\leq i\leq\ell; &&\check{\mathcal{S}}_\sigma(\mathsf{i})=\prod_{v\in V}\check{\mathcal{S}}_{\sigma_v}(\mathsf{i}(v)),\;\sigma\in\Sigma_0, \mathsf{i}\in\Psi;\\
&\check{\check{\mathcal{S}}}_\sigma(r)=\frac{1}{\left\lvert\mathsf{I}(\sigma)\right\rvert}\sum_{g\in G}\check{\mathcal{S}}_\sigma\left(g^{-1}\mathsf{i}\right),\; \mathsf{i}\in\Psi_r, 0\leq r\leq s; &&\check{\mathcal{F}}_\sigma(r)=\frac{1}{\left\lvert\mathsf{Q}_\sigma\right\rvert}\check{\check{\mathcal{S}}}_\sigma(r),\; 0\leq r\leq s.
\end{align*}
Next, we compute the coefficients in the orthogonality relations \eqref{orthrel}.
From \eqref{YVdigr} we deduce that $\left\lvert\digamma_r\right\rvert=\left\lvert\Psi_r\right\rvert\prod_{i=0}^\ell\left\lvert\Upsilon_i\right\rvert^{m_i}$ where $m_i\coloneqq\left\lvert\{v\in V\mathrel{:}\mathsf{i}(v)=i\}\right\rvert$ for $\mathsf{i}\in\Psi_r$ (it is independent of the choice of $\mathsf{i}$). Moreover, it is clear that $\left\lvert\mathsf{Orb}_{F\wr G}\left(\xi_0,\varphi_{\xi_0}\right)\right\rvert=\vert Y\rvert^{\lvert V\rvert}$. 

Finally, the conditions in \eqref{symGel} are satisfied if and only if $\bigl(F,\mathsf{Stab}_F\left(y_0\right)\bigr)$ is a symmetric Gelfand pair. Indeed, if this is the case and we take $h_1=h_2=g^{-1}$, then the first condition in \eqref{symGel} holds, while the second one reduces to the existence of elements $f_1(v),f_2(v)\in\mathsf{Stab}_F\left(y_0\right)$ such that $f(gv)^{-1}=f_1(v)f(gv)f_2(v)$.
\end{example}

\begin{example}\label{Example2ter}
In Example \ref{Example2bis}, the conditions in Proposition \ref{Propmultfree} are satisfied if and only if both $L(Y)$ and $L(V)$ are multiplicity free, as permutation representations of $F$ and $G$, respectively. Indeed, if this is the case, then $\mathsf{J}_\sigma$ is a singleton for every $\sigma\in\Sigma_0\setminus\left\{\sigma^0\right\}$, while $\mathsf{I}\left(\sigma^0\right)=G$ and $\mathsf{J}_{\sigma^0}\simeq V$. This fact agrees with the explicit decomposition \eqref{deccompact}. From now on, we assume that this condition is satisfied. 

If $\sigma\neq\sigma^0$, then $\mathsf{I}(\sigma)=\mathsf{Stab}_G\left(v_0\right)$ and $\mathsf{Q}_\sigma^0=\left\{1_G\right\}$. The spherical functions $\mathcal{R}$ are trivial (identically equal to $1$), while the functions $\mathcal{S}$ have the form 
\[
\mathcal{S}_{v_0,\sigma}\left(fy_0\right)=\left\langle u_{\sigma_{v_0}},\sigma_{v_0}(f)u_{\sigma_{v_0}}\right\rangle_{U_{\sigma_{v_0}}}
\]
for every $f\in F$. Here, $u_{\sigma_{v_0}}$ is a unit $\mathsf{Stab}_F\left(y_0\right)$-invariant vector in $U_{\sigma_{v_0}}$ and we note that, in the present setting, the condition $\sigma\prec\xi$ assumed in \eqref{defSxisigma} is replaced by $v=v_0$. The set in \eqref{StabOrbIsigma} reduces to $\left\{v_0\right\}$ and from \eqref{defFsigthet} we finally deduce that 
\[
\mathcal{F}_{\sigma,\theta^0}\left(v,y\right)=\begin{cases}
\mathcal{S}_{v_0,\sigma}(y)&\text{if }v=v_0,y\in Y,\\
0 &\text{if }v\neq v_0,y\in Y.
\end{cases} 
\]

If $\sigma=\sigma^0$, then $\mathsf{I}\left(\sigma^0\right)=G$, $\mathsf{Q}_{\sigma^0}=\left\{1_G\right\}$, the corresponding functions $\mathcal{S}$ are trivial (identically equal to $1$), the functions $\mathcal{R}$ have the form $\mathcal{R}_\theta\left(gv_0\right)=\left\langle w_\theta,\theta(g)w_\theta\right\rangle_{W_\theta}$ for every $\theta\in\widehat{G}_0$ and $g\in G$, where $w_\theta$ is a unit $\mathsf{Stab}_G\left(v_0\right)$-invariant vector in $W_\theta$. The spherical functions are $\mathcal{F}_{\sigma^0,\theta}(v,y)=\mathcal{R}_\theta(v)$, $v\in V, y\in Y$. 

With the notation in \eqref{orbdeccomp} and in \eqref{defFsigthetorb}, we can express the spherical functions as constant functions on the orbits, namely
\begin{align*}
&\check{\mathcal{S}}_{v_0,\sigma}(i)=\mathcal{S}_{v_0,\sigma}(y),\; \sigma\neq\sigma^0, 0\leq i\leq\ell, y\in\Upsilon_i, &&\check{\mathcal{F}}_{\sigma,\theta^0}(i)=\check{\mathcal{S}}_{v_0,\sigma}(i),\\
&\check{\mathcal{R}}_\theta(j)=\mathcal{R}_\theta(v),\;\theta\in\widehat{G}_0, 0\leq j\leq n, v\in\Delta_{v_0,j}, &&\check{\mathcal{F}}_{\sigma^0,\theta}(j)=\check{\mathcal{R}}_\theta(j).
\end{align*}
Moreover, 
$\check{\mathcal{F}}_{\sigma,\theta^0}$ ($\sigma\neq\sigma^0$) is identically $0$ on the orbits $\Delta_{v_0,j}\times Y$ ($1\le j\le n$).
See also \cite[Section 4.1]{CST2}, \cite[Section 3.1.2]{book3}.
In the present case, the orthogonality relations \eqref{orthrel} are a disguised form of those satisfied by the functions $\mathcal{S}$ and $\mathcal{R}$. 

Finally, from \eqref{symGel} we deduce that the pair $\left(F\wr G,\mathsf{Stab}_{F\wr G}\left(v_0,y_0\right)\right)$ is a symmetric Gelfand pair if and only if both pairs $\left(G,\mathsf{Stab}_G\left(v_0\right)\right)$ and $\left(F,\mathsf{Stab}_F\left(y_0\right)\right)$ are symmetric. Indeed, the first condition in \eqref{symGel} is precisely the symmetry of $\left(G,\mathsf{Stab}_G\left(v_0\right)\right)$. Moreover, for any fixed $g\in G$ and $h_1,h_2\in\mathsf{Stab}_G\left(v_0\right)$ satisfying $h_1gh_2=g^{-1}$, the map $v\mapsto g^{-1}h_1^{-1}v$ is a bijection of $V$; hence, as $v$ varies in $V$, the equations in the second condition of \eqref{symGel} involve pairwise independent pairs of unknowns $\left(f_1(v),f_2\left(g^{-1}h_1^{-1}v\right)\right)$ and can be solved independently for each $v\in V$. Recall that $f_1(v)$ is required to lie in $\mathsf{Stab}_F\left(y_0\right)$ only when $v=v_0$, whereas $f_2\left(g^{-1}h_1^{-1}v\right)$ is required to lie in $\mathsf{Stab}_F\left(y_0\right)$ only when $g^{-1}h_1^{-1}v=v_0$, that is, $v=h_1gv_0$. We now distinguish two cases:
\begin{itemize}
\item If $g\in\mathsf{Stab}_G\left(v_0\right)$, then $h_1gv_0=v_0$ and $h_1^{-1}v_0=g^{-1}h_1^{-1}v_0=gv_0=v_0$. For $v=v_0$, the second condition in \eqref{symGel} reduces to $f_1\left(v_0\right)f\left(v_0\right)f_2\left(v_0\right)=f\left(v_0\right)^{-1}$ with $f_1\left(v_0\right),f_2\left(v_0\right)\in\mathsf{Stab}_F\left(y_0\right)$, which holds for every $f\left(v_0\right)\in F$ if and only if $\left(F,\mathsf{Stab}_F\left(y_0\right)\right)$ is symmetric. For every $v\neq v_0$, both $f_1(v)$ and $f_2\left(g^{-1}h_1^{-1}v\right)$ can be chosen arbitrarily in $F$, so the equation is trivially satisfied.
\item If $g\notin\mathsf{Stab}_G\left(v_0\right)$, then $h_1gv_0\neq v_0$ (since $h_1\in\mathsf{Stab}_G\left(v_0\right)$). Thus the two constrained vertices $v=v_0$ and $v=h_1gv_0$ are distinct, so that for every $v\in V$ at least one of the elements $f_1(v)$ and $f_2\left(g^{-1}h_1^{-1}v\right)$ can be chosen arbitrarily in $F$ (namely, $f_2\left(g^{-1}v_0\right)$ for $v=v_0$, $f_1\left(h_1gv_0\right)$ for $v=h_1gv_0$, and both for all other $v$). Therefore, the second condition in \eqref{symGel} is trivially satisfied for every $v\in V$.
\end{itemize}
See \cite[Proposition 4.3]{CST2} for a different proof.
\end{example}

\begin{example}\label{Example3ter}
From the description of $\mathsf{J}_\sigma$ for $\sigma\in\Sigma_0$ in Example \ref{Example3bis}, it follows that we have equality in \eqref{defJsig0old} if and only if the following condition is satisfied: for every $Z\in\mathsf{Orb}_G\left(Z_0\right)$ such that $\mathsf{supp}(\sigma)\subseteq Z\cap Z_0$ there exists $k\in\mathsf{I}(\sigma)$ such that $kZ_0=Z$. 
If we use the notation in \eqref{VsIsigma}, this is equivalent to the statement that there exists $k\in G$ such that
\begin{equation}\label{kV0V0kVs}
\begin{array}{ll}
kV_{\sigma_*}=V_{\sigma_*}\; \text{ for every }\; \sigma_*\in\mathsf{image}(\sigma), &\qquad(k\in\mathsf{I}(\sigma)),\\ 
\quad&\\
k\left(Z_0\setminus\mathsf{supp}(\sigma)\right)=Z\setminus\mathsf{supp}(\sigma),  &\qquad(kZ_0=Z).
\end{array}
\end{equation}
If these conditions hold for every $\sigma\in\Sigma_0$, then from Corollary \ref{Cormultfree}\eqref{Cormultfree1}, we deduce that the condition \eqref{Propmultfreeb} in Proposition \ref{Propmultfree} is satisfied if and only if $L(Y)$ is multiplicity free. Moreover, the condition \eqref{Propmultfreec} in Proposition \ref{Propmultfree} now reduces to the statement that, for every $\sigma\in\Sigma_0$, the permutation representation of $\mathsf{I}(\sigma)$ on $L\bigl(\mathsf{Orb}_{\mathsf{I}(\sigma)}\left(Z_0\right)\bigr)$ is multiplicity free. In conclusion, the permutation representation in the present example is multiplicity free if and only if, for every $\sigma\in\Sigma_0$, \eqref{kV0V0kVs} holds and both $L(Y)$ and each $L\bigl(\mathsf{Orb}_{\mathsf{I}(\sigma)}\left(Z_0\right)\bigr)$ are multiplicity free. Henceforth, we assume that these conditions are satisfied and we note that \eqref{JsiginCor} becomes $\mathsf{J}_\sigma=\mathsf{Orb}_{\mathsf{I}(\sigma)}\left(Z_0\right)=\left\{Z\in\mathsf{Orb}_G\left(Z_0\right)\mathrel{:}\mathsf{supp}(\sigma)\subseteq Z\right\}$.

We define $\mathsf{Orb}_{\mathrm{I}(\sigma)}(Z_0 \cap V_{\sigma_0})^0$ directly as a subset of $\mathsf{Orb}_{\mathrm{I}(\sigma)}(Z_0 \cap V_{\sigma_0})$ obtained by choosing exactly one $Z'$ for each $\mathsf{Stab}_G(Z_0)$-orbit on $\mathsf{Stab}_G(Z_0)\mathsf{Orb}_{\mathrm{I}(\sigma)}(Z_0)$ (i.e., a single $Z'$ for each orbit of the form $\mathsf{Orb}_{\mathsf{Stab}_G(Z_0)}(Z'\sqcup\mathsf{supp}(\sigma))$). 
For instance, for $G=\mathsf{Sym}(V)$ it reduces to choosing one $Z'$ for each value of $\vert{}Z'\cap Z_0\vert{}$.
Then we can write
\[
\mathsf{Stab}_G\left(Z_0\right)\mathsf{Orb}_{\mathsf{I}(\sigma)}\left(Z_0\cap V_{\sigma_0}\right)=\bigsqcup_{Z'\in\mathsf{Orb}_{\mathsf{I}(\sigma)}\left(Z_0\cap V_{\sigma_0}\right)^0}\mathsf{Orb}_{\mathsf{Stab}_G\left(Z_0\right)}\left(Z'\right).
\]
We note that, for $h\in\mathsf{Stab}_G\left(Z_0\right)$ and $k\in\mathsf{I}(\sigma)$, we have 
\[
hkZ_0=hk\left[\left(Z_0\cap V_{\sigma_0}\right)\sqcup\mathsf{supp}(\sigma)\right]=hk\left(Z_0\cap V_{\sigma_0}\right)\sqcup h\mathsf{supp}(\sigma).
\]
We deduce the following coordinate-free form of \eqref{StabOrbIsig0}
\begin{equation}\label{qZ0V0Vtbis1}
\mathsf{Stab}_G\left(Z_0\right)\mathsf{Orb}_{\mathsf{I}(\sigma)}\left(Z_0\right)=\bigsqcup_{Z'\in\mathsf{Orb}_{\mathsf{I}(\sigma)}\left(Z_0\cap V_{\sigma_0}\right)^0}\mathsf{Orb}_{\mathsf{Stab}_G\left(Z_0\right)}\left(Z'\sqcup\mathsf{supp}(\sigma)\right).
\end{equation}

On the other hand, from \eqref{defQ0sig} we deduce that
\[
\mathsf{Stab}_G\left(Z_0\right)\mathsf{Orb}_{\mathsf{I}(\sigma)}\left(Z_0\cap V_{\sigma_0}\right)=\bigcup_{q\in\mathsf{Q}_\sigma^0}q\mathsf{Orb}_{\mathsf{I}(\sigma)}\left(Z_0\cap V_{\sigma_0}\right).
\]
Then, taking \eqref{OrbIsigZ0} into account, the following version of the second equality in \eqref{StabOrbIsigma} can be derived
\begin{equation}\label{qZ0V0Vtbis2}
\mathsf{Stab}_G\left(Z_0\right)\mathsf{Orb}_{\mathsf{I}(\sigma)}\left(Z_0\right)=\bigcup_{q\in\mathsf{Q}_\sigma^0}q\left(\;\bigsqcup_{Z'\in\mathsf{Orb}_{\mathsf{I}(\sigma)}\left(Z_0\cap V_{\sigma_0}\right)}\left(Z'\sqcup\mathsf{supp}(\sigma)\right)\right).
\end{equation} 
We also observe that, assuming the notation in \eqref{OrbV0Vt}, from \eqref{StabStabIsigma} and the first condition in \eqref{kV0V0kVs} it follows that the map 
\begin{equation}\label{qV1Vtbijec}
\begin{array}{ccc}
\mathsf{Q}_\sigma^0&\longrightarrow&\mathsf{Orb}_{\mathsf{Stab}_G\left(Z_0\right)}\bigl(V_{\sigma_*}\mathrel{:}\sigma_*\in\mathsf{image}(\sigma)\setminus\left\{\sigma_0\right\}\bigr)\\
\quad&&\\
q&\longmapsto& \bigl(qV_{\sigma_*}\mathrel{:}\sigma_*\in\mathsf{image}(\sigma)\setminus\left\{\sigma_0\right\}\bigr)
\end{array}
\end{equation}
is a bijection. Recall also that $\mathsf{supp}(\sigma)\subseteq Z_0$.

In the setting of this example, it is straightforward to provide an instance of the phenomenon described in Remark \ref{Remkh}. Suppose that $Z_0\neq V$ and that $\mathsf{supp}(\sigma)$
and $V_{\sigma_0}\cap Z_0$ are both non-trivial, and that there exist $h\in\mathsf{Stab}_G\left(Z_0\right)$ such that $(h\mathsf{supp}(\sigma))\cap\left(V_{\sigma_0}\cap Z_0\right)\neq\emptyset$ and $k\in\mathsf{I}(\sigma)$ such that 
$k\left[(h\mathsf{supp}(\sigma))\cap\left(V_{\sigma_0}\cap Z_0\right)\right]\setminus Z_0\neq\emptyset$. Then $kh\mathsf{supp}(\sigma)\nsubseteq Z_0$, which implies that $kh\notin\mathsf{Stab}_G\left(Z_0\right)\cdot\mathsf{I}(\sigma)$. Elements $h$ and $k$ with these properties certainly exist if $G=\mathsf{Sym}(V)$. 

Next, we examine the condition on $q$ in \eqref{defFsigthet}, which in the present setting has the form $q\in\mathsf{Q}_\sigma^0\mathrel{:}q^{-1}Z\in\mathsf{Orb}_{\mathsf{I}(\sigma)}\left(Z_0\right)$, where $Z\in\mathsf{Stab}_G\left(Z_0\right)\mathsf{Orb}_{\mathsf{I}(\sigma)}\left(Z_0\right)$. From \eqref{qZ0V0Vtbis2} we deduce that
\begin{equation}\label{qZminusV1Vtter}
\begin{split}
q^{-1}Z\in\mathsf{Orb}_{\mathsf{I}(\sigma)}\left(Z_0\right)\quad&\Longleftrightarrow\quad Z=qZ'\sqcup q\mathsf{supp}(\sigma)\;\text{ with }\;Z'\in\mathsf{Orb}_{\mathsf{I}(\sigma)}\left(Z_0\cap V_{\sigma_0}\right)\\
\Longleftrightarrow\quad &Z\supseteq q\mathsf{supp}(\sigma)\;\text{ and }\;Z\setminus q\mathsf{supp}(\sigma) \in q\mathsf{Orb}_{\mathsf{I}(\sigma)}\left(Z_0\cap V_{\sigma_0}\right).
\end{split}
\end{equation}

We turn now to the functions in \eqref{defSxisigma}. For every $\sigma_*\in\mathsf{image}(\sigma)\setminus\left\{\sigma_0\right\}$ we set $\mathcal{S}_{\sigma_*}(fy_0)\coloneqq\left\langle u_{\sigma_*},\sigma_*(f)u_{\sigma_*}\right\rangle_{U_{\sigma_*}}$ for all $f\in F$, where $u_{\sigma_*}$ is a unit $\mathsf{Stab}_F\left(y_0\right)$-invariant vector in $U_{\sigma_*}$, so that $\mathcal{S}_{\sigma_*}$ is also $\mathsf{Stab}_F\left(y_0\right)$ bi-invariant. Then we define the functions  $\mathcal{S}$ in the following way
\begin{equation}\label{Ssigmaphip}
\mathcal{S}_\sigma\left(\phi'\right)\coloneqq\prod_{\sigma_*\in\mathsf{image}(\sigma)\setminus\left\{\sigma_0\right\}}\;\prod_{v\in V_{\sigma_*}}\mathcal{S}_{\sigma_*}\left(\phi'(v)\right)\quad\text{for all }  \phi'\in Y^{\mathsf{supp}(\sigma)}.
\end{equation}
Here we omit the dependence on $Z\in\mathsf{J}_\sigma$ because outside $\mathsf{supp}(\sigma)$ we have the trivial representation which yields only the constant value $1$. Note also that, taking \eqref{StabStabIsigma} into account, for all $\phi\in Y^Z$ with $Z\supseteq q\mathsf{supp}(\sigma)$ and $q\in\mathsf{Q}_\sigma^0$, we have $\left.\left(q^{-1}\phi\right)\right\rvert_{\mathsf{supp}(\sigma)}=q^{-1}\left(\left.\phi\right\rvert_{q\mathsf{supp}(\sigma)}\right)$, and therefore
\begin{equation}\label{Ssigmaphiorb}
\mathcal{S}_\sigma\left(\left.\left(q^{-1}\phi\right)\right\rvert_{\mathsf{supp}(\sigma)}\right)=\mathcal{S}_\sigma\left(q^{-1}\left(\left.\phi\right\rvert_{q\mathsf{supp}(\sigma)}\right)\right).
\end{equation}

Similarly, recalling \eqref{OrbIsigZ0}, we may think of the functions $\mathcal{R}$ in \eqref{defRsigmatheta} as being defined on $\mathsf{Orb}_{\mathsf{I}(\sigma)}\left(Z_0\cap V_{\sigma_0}\right)$, that is, we set $\mathcal{R}_{\sigma,\theta}\left(k\left(Z_0\cap V_{\sigma_0}\right)\right)\coloneqq\left\langle w_\theta,\theta(k)w_\theta\right\rangle_{W_\theta}$ for every $k\in\mathsf{I}(\sigma)$. With the notation in \eqref{qZminusV1Vtter}, we can write $\mathcal{R}_{\sigma,\theta}$ as a function of $q^{-1}Z\setminus\mathsf{supp}(\sigma)=Z'\in\mathsf{Orb}_{\mathsf{I}(\sigma)}\left(Z_0\cap V_{\sigma_0}\right)$, hence in the following form
\begin{equation}\label{qZminusV1Vtbis}
\mathcal{R}_{\sigma,\theta}\left(q^{-1}Z\setminus\mathsf{supp}(\sigma)\right).
\end{equation}

Taking \eqref{qV1Vtbijec}, \eqref{Ssigmaphip}, \eqref{Ssigmaphiorb} and \eqref{qZminusV1Vtbis} into account, we can now write the spherical functions in \eqref{defFsigthet}. For all $Z\in\mathsf{Stab}_G\left(Z_0\right)\mathsf{Orb}_{\mathsf{I}(\sigma)}\left(Z_0\right)$, $\phi\in Y^Z$, we have
\[
\begin{split}
\mathcal{F}_{\sigma,\theta}(Z,\phi)=&\frac{1}{\left\lvert\mathsf{Orb}_{\mathsf{Stab}_G\left(Z_0\right)}\bigl(V_{\sigma_*}\mathrel{:}\sigma_*\in\mathsf{image}(\sigma)\setminus\left\{\sigma_0\right\}\bigr)\right\rvert}\cdot\\
&\qquad\cdot\sum_{\substack{q\in \mathsf{Q}_\sigma^0\mathrel{:}\\ q^{-1}Z\in\mathsf{Orb}_{\mathsf{I}(\sigma)}\left(Z_0\right)}}\!\!\!\!\!\!\!\mathcal{R}_{\sigma,\theta}\left(q^{-1}Z\setminus\mathsf{supp}(\sigma)\right)\mathcal{S}_\sigma\left(\left.\left(q^{-1}\phi\right)\right\rvert_{\mathsf{supp}(\sigma)}\right)\\
=&\frac{1}{\left\lvert\mathsf{Stab}_G\left(Z_0\right)\right\rvert}\sum_{\substack{h\in\mathsf{Stab}_G\left(Z_0\right)\mathrel{:}\\ h^{-1}Z\in\mathsf{Orb}_{\mathsf{I}(\sigma)}\left(Z_0\right)}}\mathcal{R}_{\sigma,\theta}\left(h^{-1}Z\setminus\mathsf{supp}(\sigma)\right)\mathcal{S}_\sigma\left(\left.\left(h^{-1}\phi\right)\right\rvert_{\mathsf{supp}(\sigma)}\right).
\end{split}
\]
The second sum follows immediately from the first since \eqref{StabStabIsigma} allows us to set $h=qk$ with $k\in\mathsf{Stab}_{\mathsf{I}(\sigma)}\left(Z_0\right)$. This rearrangement uses the fact that $\mathcal{R}_{\sigma,\theta}$ is $\mathsf{Stab}_{\mathsf{I}(\sigma)}\left(Z_0\right)$-invariant, $\mathcal{S}_\sigma$ is $\mathsf{I}(\sigma)$-invariant, and $\mathsf{supp}(\sigma)$ is stabilized by $\mathsf{Stab}_{\mathsf{I}(\sigma)}\left(Z_0\right)$. 

Next, we show how to express the spherical functions as constant functions on the orbits. We use the notation in Example \ref{Example3}, taking into account that now \eqref{StabOrbIsig0} is replaced by \eqref{qZ0V0Vtbis1}. In particular, $j$ is replaced by $Z'$ and the condition 
$0\leq j\leq n^0$ by $Z'\in\mathsf{Orb}_{\mathsf{I}(\sigma)}\left(Z_0\cap V_{\sigma_0}\right)^0$.
If $\sigma\in\Sigma_0$, then we set
\[
\begin{split}
&\check{\mathcal{S}}_{\sigma_v}(i)\coloneqq\mathcal{S}_{\sigma_v}(y)\quad\text{ for }\;v\in\mathsf{supp}(\sigma),\; 0\leq i\leq\ell,\; y\in\Upsilon_i,\\
&\check{\mathcal{S}}_\sigma\left(\mathsf{i}'\right)\coloneqq\prod_{v\in \mathsf{supp}(\sigma)}\check{\mathcal{S}}_{\sigma_v}\left(\mathsf{i}'(v)\right)\quad\text{ for }\; \mathsf{i}'\in\{0,\dotsc,\ell\}^{\mathsf{supp}(\sigma)},\\
\end{split}
\]
so that
\[
\check{\check{\mathcal{S}}}_{Z',l,\sigma}(r)=\frac{1}{\left\lvert\mathsf{Stab}_{\mathsf{I}(\sigma)}\left(Z_0\cap V_{\sigma_0}\right)\right\rvert}\sum_{Z''\in\Pi_l}\sum_{\substack{h\in\mathsf{Stab}_G\left(Z_0\right)\mathrel{:}\\ h^{-1}Z=Z''}}\check{\mathcal{S}}_\sigma\left(\left.\left(h^{-1}\mathsf{i}\right)\right\rvert_{\mathsf{supp}(\sigma)}\right)
\]
for $Z'\in\mathsf{Orb}_{\mathsf{I}(\sigma)}\left(Z_0\cap V_{\sigma_0}\right)^0$, $0\leq l\leq p_\sigma$, $\Pi_l\subseteq\mathsf{Orb}_{\mathsf{Stab}_G(Z_0)}(Z'\sqcup\mathsf{supp}(\sigma))$, $0\leq r\leq s_{Z'}$ and any choice of $(Z,\mathsf{i})\in\Psi_{Z',r}$. Then $\check{\mathcal{R}}_{\sigma,\theta}$ and $\check{\mathcal{F}}_{\sigma,\theta}$ are exactly as in the general formulas.
Using the results in Example \ref{Example3}, we easily obtain the coefficients in the orthogonality relations \eqref{orthrel}
\[
\left\lvert\digamma_{Z',r}\right\rvert=\sum_{(Z,\mathsf{i})\in\Psi_{Z',r}}\left[\left(\prod_{v\in Z\cap Z_0}\left\lvert\Upsilon_{\mathsf{i}(v)}\right\rvert\right)\left\lvert Y\right\rvert^{\left\lvert Z\setminus Z_0\right\rvert}\right]
\]
and $\left\lvert\mathsf{Orb}_{F\wr G}\left(Z_0,\phi_{Z_0}\right)\right\rvert=\left\lvert\mathsf{Orb}_G\left(Z_0\right)\right\rvert\lvert Y\rvert^{\left\lvert Z_0\right\rvert}$.

Finally, we determine when $\left(F\wr G,\mathsf{Stab}_{F\wr G}\left(Z_0,\phi_{Z_0}\right)\right)$ is a symmetric Gelfand pair. The first equation in \eqref{symGel} is satisfied if and only if the pairs $\left(G,\mathsf{Stab}_G\left(Z_0\right)\right)$ and $\left(F,\mathsf{Stab}_F\left(y_0\right)\right)$ are both symmetric. Assuming these conditions, $\left(F\wr G,\mathsf{Stab}_{F\wr G}\left(Z_0,\phi_{Z_0}\right)\right)$ is symmetric if and only if the following condition holds true.

\medskip
($\ast$)  For every $g\in G\setminus\mathsf{Stab}_G\left(Z_0\right)$ and every $f\in F^V$, there exist $(f_1,h_1),(f_2,h_2)\in\mathsf{Stab}_{F\wr G}\left(Z_0,\phi_{Z_0}\right)$ with $h_1gh_2=g^{-1}$ such that the second equation of \eqref{symGel} holds simultaneously for every $v\in Z_0\cap g^{-1}Z_0$.

\medskip

Note that the symmetry of $\left(F,\mathsf{Stab}_F\left(y_0\right)\right)$ is automatically guaranteed by $(\ast)$ whenever there exists $g \in G \setminus \mathsf{Stab}_G(Z_0)$ such that $Z_0 \cap g^{-1}Z_0 \neq \emptyset$. In extreme cases where $Z_0 \cap g^{-1}Z_0 = \emptyset$ for every $g \in G \setminus \mathsf{Stab}_G(Z_0)$ (such as when $Z_0 = V$, $\vert{}Z_0\vert{} = 1$, or when $\mathsf{Orb}_G(Z_0)$ consists of pairwise disjoint blocks), condition $(\ast)$ is vacuously true, and symmetry of both component pairs must be explicitly required.

\medskip
\noindent
To prove that the second condition in \eqref{symGel} holds in full generality, it suffices to note the following elementary facts.
\begin{itemize}
\item From ($\ast$) it follows that $\left(F,\mathsf{Stab}_F\left(y_0\right)\right)$ is symmetric. It suffices to observe that we can choose $f$ constant on $V$ (or on $Z_0$, given that $h_1^{-1}v, gv\in Z_0$), so that the equality $f(h_1^{-1}v)=f(gv)$ holds automatically for any $h_1\in\mathsf{Stab}_G(Z_0)$, and that $g^{-1}h_1^{-1}v=h_2gv\in Z_0$, because $h_2\in\mathsf{Stab}_G\left(Z_0\right)$.
\item On the other hand, if $g\in\mathsf{Stab}_G\left(Z_0\right)$ and $v\in Z_0$, then also $gv\in Z_0$ (and vice versa). We may take $h_1=h_2=g^{-1}$, so that the second equality in \eqref{symGel} reduces to $f_1(v)f(gv)f_2(v)=f(gv)^{-1}$, and $f_1(v),f_2(v)\in\mathsf{Stab}_F\left(y_0\right)$ exist by the symmetry of $\left(F,\mathsf{Stab}_F\left(y_0\right)\right)$.
\item Similarly, if $g\in\mathsf{Stab}_G\left(Z_0\right)$ and $v\notin Z_0$, then the value of $f_1(v)$ can be chosen arbitrarily in $F$ and the second equality in \eqref{symGel} is trivially satisfied.
\item Finally, if $g\notin\mathsf{Stab}_G\left(Z_0\right)$ and at least one of the elements $v$ or $gv$ is not in $Z_0$, then the corresponding values of $f_1$ or $f_2$ can be chosen arbitrarily in $F$ (noting that the map $v \mapsto g^{-1}h_1^{-1}v$ is bijective, so that the choices for different $v \in V$ do not interfere with each other).
\end{itemize}
\end{example}

\begin{example}\label{Example5ter}
In Example \ref{Example5bis}, the conditions in Proposition \ref{Propmultfree} are satisfied if and only if
the following conditions hold:
\begin{itemize}
\item each permutation representation $L\left(Y_a\right)$ for $a\in A$ is multiplicity free;
\item for every $\sigma\in\Sigma_0$, we have $\left\{\varphi\in\mathsf{Orb}_G\left(\varphi_0\right)\mathrel{:}\sigma\prec\varphi\right\}=\mathsf{Orb}_{\mathsf{I}(\sigma)}\left(\varphi_0\right)$;
\item for every $\sigma\in\Sigma_0$, the permutation representation of $\mathsf{I}(\sigma)$ on $L\bigl(\mathsf{Orb}_{\mathsf{I}(\sigma)}\left(\varphi_0\right)\bigr)$ is multiplicity free.
\end{itemize}

We assume that the conditions above are satisfied and recall that, with the notation of Example \ref{Example5}, $\varphi_\xi$ is replaced by $\varphi^\sharp$. Furthermore, we slightly change the notation in \eqref{defSxisigma} by denoting the unit $\mathsf{Stab}_F\left(y_a\right)$-invariant vectors by $u_{y_a,\sigma_*}$ and setting
\begin{equation}\label{Sxisigmaphi}
\mathcal{S}_{\varphi^\sharp,\sigma}\left(f\cdot\varphi^\sharp\right)=\prod_{v\in V}\left\langle u_{\varphi_0(v),\sigma_v},\sigma_v(f(v))u_{\varphi^\sharp(v),\sigma_v}\right\rangle_{U_{\sigma_v}}
\end{equation}
for $\sigma\in\Sigma_0$, $\varphi^\sharp\in\mathsf{Orb}_{\mathsf{I}(\sigma)}\left(\varphi_0\right)$, and $f\in F^V$.
Clearly, $\mathcal{S}_{\varphi^\sharp,\sigma}$ is a $\mathsf{Stab}_{F^V}\left(\varphi_0\right)$-invariant function in $L\bigl(\mathsf{Orb}_{F^V}\left(\varphi^\sharp\right)\bigr)$.
Then we can write \eqref{defFsigthet} in the form 
\begin{equation}\label{Fsigthetphi}
\mathcal{F}_{\sigma,\theta}(\varphi)=
\frac{1}{\left\lvert\mathsf{Q}^0_\sigma\right\rvert}\sum_{\substack{q\in\mathsf{Q}^0_\sigma\mathrel{:}\\ q^{-1}\varphi^\sharp\in\mathsf{Orb}_{\mathsf{I}(\sigma)}\left(\varphi_0\right)}}\mathcal{R}_{\sigma,\theta}\left(q^{-1}\varphi^\sharp\right)\mathcal{S}_{q^{-1}\varphi^\sharp,\sigma}\left(q^{-1}\varphi\right)
\end{equation}
if $\varphi^\sharp\in \mathsf{Stab}_G\left(\varphi_0\right)\mathsf{Orb}_{\mathsf{I}(\sigma)}\left(\varphi_0\right)$, and $\mathcal{F}_{\sigma,\theta}(\varphi)=0$ otherwise. 
\end{example}

\begin{example}\label{Example6ter}
In the setting of Example \ref{Example6bis}, from the characterization of $\mathsf{J}_\sigma$ and Corollary \ref{Cormultfree}, we deduce that we have a multiplicity free permutation representation if and only if the following conditions are satisfied:
\begin{itemize} 
\item the permutation representation of $F$ on $L(Y)$ is multiplicity free;
\item for every $\sigma\in\Sigma_0$, the group $\mathsf{I}(\sigma)$ is transitive on $\mathsf{Orb}_G\left(\xi_0\right)$;
\item for every $\sigma\in\Sigma_0$, the permutation representation of $\mathsf{I}(\sigma)$ on $L\bigl(\mathsf{Orb}_G\left(\xi_0\right)\bigr)$ is multiplicity free. 
\end{itemize}
From now on, we assume that these conditions are satisfied.

The functions $\mathcal{S}$ have the form $\mathcal{S}_\sigma\left(f\cdot\varphi_0\right)=\prod_{v\in V}\left\langle u_{\sigma_v},\sigma_v(f(v))u_{\sigma_v}\right\rangle_{U_{\sigma_v}}$ for $f\in F^V$, where $u_{\sigma_v}$ is a unit $\mathsf{Stab}_F\left(y_0\right)$-invariant vector in $U_{\sigma_v}$; the set in \eqref{StabOrbIsigma} is the whole $\mathsf{Orb}_G\left(\xi_0\right)$ and $\mathcal{R}_{\sigma,\theta}$ is a $\mathsf{Stab}_{\mathsf{I}(\sigma)}\left(\xi_0\right)$-invariant function defined on $\mathsf{Orb}_G\left(\xi_0\right)\equiv\mathsf{Orb}_{\mathsf{I}(\sigma)}\left(\xi_0\right)$. Because $\mathsf{I}(\sigma)$ is transitive on $\mathsf{Orb}_G\left(\xi_0\right)$, one actually has $\mathsf{Stab}_G\left(\xi_0\right)\mathsf{I}(\sigma)=G$, and therefore $\mathsf{Q}^0_\sigma=\mathsf{Q}_\sigma$. The condition $q^{-1}\xi\in\mathsf{Orb}_{\mathsf{I}(\sigma)}\left(\xi_0\right)$ in \eqref{defFsigthet} is satisfied by all $q\in\mathsf{Q}_\sigma$ and the spherical functions have the form
\[
\mathcal{F}_{\sigma,\theta}(\xi,\varphi)=
\frac{1}{\left\lvert\mathsf{Q}_\sigma\right\rvert}\sum_{q\in\mathsf{Q}_\sigma}\mathcal{R}_{\sigma,\theta}\left(q^{-1}\xi\right)\mathcal{S}_\sigma\left(q^{-1}\varphi\right)
\]
for all $(\xi,\varphi)\in\mathsf{Orb}_G\left(\xi_0\right)\times Y^V$. 

Now $n^0=n$ and $\mathsf{Orb}_G\left(\xi_0\right)=\bigsqcup_{l=0}^{p_\sigma}\Pi_l$ is the decomposition into the orbits of $\mathsf{Stab}_{\mathsf{I}(\sigma)}\left(\xi_0\right)$. Furthermore, with the notation of Example \ref{Example6}, for $\sigma_*\in\widehat{F}_0$ and $\sigma\in\Sigma_0$ we have
\[
\check{\mathcal{S}}_{\sigma_*}(i)=\mathcal{S}_{\sigma_*}(y),\; y\in\Upsilon_i, \;0\leq i\leq\ell; \qquad\qquad\check{\mathcal{S}}_\sigma(\mathsf{i})=\prod_{v\in V}\check{\mathcal{S}}_{\sigma_v}(\mathsf{i}(v)), \;\mathsf{i}\in\Psi.
\]
Therefore, by \eqref{barbarcalS} and observing that both $\mathcal{S}$ and $\check{\mathcal{S}}$ do not depend on $\xi$ or $\xi'$, for $0\leq j\leq n$, $0\leq r\leq s_j$, $(\xi,\mathsf{i})\in\Psi_{j,r}$, and $\Pi_l\subseteq\Delta_j$ we obtain
\[
\check{\check{\mathcal{S}}}_{j,l,\sigma}(r)
=\frac{1}{\left\lvert\mathsf{Stab}_{\mathsf{I}(\sigma)}\left(\xi_0\right)\right\rvert}\sum_{\xi'\in\Pi_l}\;\sum_{\substack{h\in\mathsf{Stab}_G\left(\xi_0\right)\mathrel{:}\\ h^{-1}\xi=\xi'}}\check{\mathcal{S}}_\sigma\left(h^{-1}\mathsf{i}\right)
=\frac{1}{\left\lvert\mathsf{Stab}_{\mathsf{I}(\sigma)}\left(\xi_0\right)\right\rvert}\sum_{\substack{h\in\mathsf{Stab}_G\left(\xi_0\right)\mathrel{:}\\ h^{-1}\xi\in\Pi_l}}\check{\mathcal{S}}_\sigma\left(h^{-1}\mathsf{i}\right).
\]
Finally, $\check{\mathcal{F}}_{\sigma,\theta}$ is as in \eqref{defFsigthetorb}.
\end{example}

\section{Young permutation representations for wreath products $F\wr \mathsf{Sym}(V)$}\label{SecYoungpermrep}

In this section, assuming the notions and notation in Section \ref{Sechomspace}, we examine the case $G=\mathsf{Sym}(V)$. 
A map $\xi\in\mathsf{Surj}(V,X)$ determines the set composition $V=\bigsqcup_{x\in X}\xi^{-1}(x)$; the \emph{type} of $\xi$ is the function $\left\lvert\xi^{-1}\right\rvert\in\mathbb{N}^X$ defined by setting 
\begin{equation}\label{deftype}
\left\lvert\xi^{-1}\right\rvert(x)\coloneqq\left\lvert\xi^{-1}(x)\right\rvert
\end{equation}
for all $x\in X$, meaning that its values are the cardinalities of the blocks of the set composition. It also determines a \emph{partition-valued function} $\Lambda_\xi$ defined on $A$ by letting $\Lambda_\xi(a)$ be the integer partition whose parts are $\left\lvert\xi^{-1}(x)\right\rvert$ for $x\in \alpha^{-1}(a)$. That is, $\Lambda_\xi(a)=\left(\lambda_1,\lambda_2,\dotsc,\lambda_k\right)$ means that $\alpha^{-1}(a)=\{x_1,x_2,\dotsc,x_k\}$, $\lambda_j=\left\lvert\xi^{-1}\right\rvert\left(x_j\right)$, and the elements in $\alpha^{-1}(a)$ are numbered in such a way that $\left\lvert\xi^{-1}\right\rvert\left(x_1\right)\geq\left\lvert\xi^{-1}\right\rvert\left(x_2\right)\geq\dotsb\geq\left\lvert\xi^{-1}\right\rvert\left(x_k\right)$. 

\begin{lemma}\label{Lemmaorbsym}
Suppose that $(\xi,\varphi),(\xi',\varphi')\in\Gamma$.
\begin{enumerate}[(1)]
\item\label{Lemmaorbsym1} We have $\mathsf{Orb}_{F\wr\mathsf{Sym}(V)}(\xi,\varphi)=\mathsf{Orb}_{F\wr\mathsf{Sym}(V)}(\xi',\varphi')$ if and only if $\left\lvert\xi^{-1}\right\rvert=\left\lvert\xi'^{-1}\right\rvert$.
\item\label{Lemmaorbsym2} We have the isomorphisms
\begin{equation}\label{orbsym2}
\begin{split}
\mathsf{Stab}_{F\wr\mathsf{Sym}(V)}\left(\xi,\varphi_{\xi}\right)&\;\cong\;\prod_{x\in X}\left[\mathsf{Stab}_F\left(y_{\alpha(x)}\right)\wr \mathsf{Sym}\left(\xi^{-1}(x)\right)\right]\\
&\;\cong\;\prod_{a\in A}\left\{\mathsf{Stab}_F\left(y_a\right)\wr\left[\prod_{x\in\alpha^{-1}(a)} \mathsf{Sym}\left(\xi^{-1}(x)\right)\right]\right\}.
\end{split}
\end{equation}
\item\label{Lemmaorbsym3} The following facts are equivalent:
\begin{itemize} 
\item $\mathsf{Orb}_{F\wr\mathsf{Sym}(V)}(\xi,\varphi)$ and $\mathsf{Orb}_{F\wr\mathsf{Sym}(V)}(\xi',\varphi')$ are isomorphic as $F\wr\mathsf{Sym}(V)$-spaces;
\item $\Lambda_\xi=\Lambda_{\xi'}$;
\item $\xi$ and $\xi'$ are equivalent.
\end{itemize}
\end{enumerate}
\end{lemma}
\begin{proof}
(1)
It follows from Lemma \ref{Lemmaorb}\eqref{Lemmaorb1}. Just note that $\xi'\in\mathsf{Orb}_{\mathsf{Sym}(V)}(\xi)$ if and only if
there exists $\pi\in\mathsf{Sym}(V)$ such that $\pi\left(\xi^{-1}(x)\right)=\xi'^{-1}(x)$ for all $x\in X$. Such a permutation $\pi$ exists if and only if $\left\lvert\xi^{-1}\right\rvert=\left\lvert\xi'^{-1}\right\rvert$. 

\medskip
(2)
It follows from Lemma \ref{Lemmastab}: we have $\pi\xi=\xi$, $\pi\in\mathsf{Sym}(V)$, if and only if $\pi\left(\xi^{-1}(x)\right)=\xi^{-1}(x)$ for every $x\in X$, so that
$\mathsf{Stab}_{\mathsf{Sym}(V)}(\xi)=\prod_{x\in X}\mathsf{Sym}\left(\xi^{-1}(x)\right)$. Note that the second isomorphism in \eqref{orbsym2} follows immediately by grouping $x\in X$ according to $a=\alpha(x)$ and letting $\prod_{x\in\alpha^{-1}(a)}\mathsf{Sym}(\xi^{-1}(x))$ act naturally on $\xi^{-1}(\alpha^{-1}(a))$.

\medskip
(3)
First, note that $\Lambda_\xi=\Lambda_{\xi'}$ if and only if there exists $\omega\in \Omega$ such that $\left\lvert\left(\omega\circ\xi\right)^{-1}\right\rvert=\left\lvert\xi'^{-1}\right\rvert$, because $\omega$ is a suitable reordering of the elements within each set $\alpha^{-1}(a)$. Arguing similarly to the proof of \eqref{Lemmaorbsym1}, we deduce that this happens if and only if there exists $\pi\in\mathsf{Sym}(V)$ such that $\omega\circ\xi\circ\pi^{-1}=\xi'$. That is, $\Lambda_\xi=\Lambda_{\xi'}$ if and only if $\xi$ and $\xi'$ are equivalent and, if this is the case, then we may invoke Lemma \ref{Lemmaorb}\eqref{Lemmaorb3} to conclude that $\mathsf{Orb}_{F\wr\mathsf{Sym}(V)}(\xi,\varphi)$ and $\mathsf{Orb}_{F\wr\mathsf{Sym}(V)}(\xi',\varphi')$ are isomorphic.

Conversely, suppose that $\mathsf{Orb}_{F\wr\mathsf{Sym}(V)}(\xi,\varphi)$ and $\mathsf{Orb}_{F\wr\mathsf{Sym}(V)}(\xi',\varphi')$ are isomorphic. By Lemma \ref{Lemmaorb}\eqref{Lemmaorb3} there exists $\pi\in\mathsf{Sym}(V)$ such that $\pi\mathsf{Stab}_{\mathsf{Sym}(V)}(\xi)\pi^{-1}=\mathsf{Stab}_{\mathsf{Sym}(V)}(\xi')$ and $\alpha\circ\xi\circ\pi^{-1}=\alpha\circ\xi'$. Since 
\[
\pi\mathsf{Stab}_{\mathsf{Sym}(V)}(\xi)\pi^{-1}=\pi\left(\prod_{x\in X}\mathsf{Sym}\left(\xi^{-1}(x)\right)\right)\pi^{-1}=\prod_{x\in X}\mathsf{Sym}\bigl(\pi\left(\xi^{-1}(x)\right)\bigr),
\] 
we deduce that there exists $\omega\in\mathsf{Sym}(X)$ such that $\pi\bigl(\xi^{-1}\left(\omega^{-1}x\right)\bigr)=\xi'^{-1}(x)$ for all $x\in X$. This is valid even when blocks of size $1$ occur, because every block $\xi^{-1}(x)$ is nonempty ($\xi,\xi'\in\mathsf{Surj}(V,X)$), so the blocks $\{\pi(\xi^{-1}(x))\}_{x\in X}$ and $\{\xi'^{-1}(x)\}_{x\in X}$ are precisely the orbits of these two equal subgroups on $V$ (with size-$1$ blocks being fixed points). That is, $\omega\circ\xi\circ\pi^{-1}=\xi'$. Then 
\[
\alpha\circ\omega\circ\xi\circ\pi^{-1}=\alpha\circ\xi'=\alpha\circ\xi\circ\pi^{-1}\quad\Longrightarrow\quad\alpha\circ\omega\circ\xi=\alpha\circ\xi
\]
which in turn forces $\alpha\circ\omega=\alpha$ (because $\xi$ is surjective) so that $\omega\in\Omega$. Therefore, $\xi$ and $\xi'$ are equivalent. 
\end{proof}

For every $\sigma\in\Sigma_0$, $\mathsf{I}(\sigma)$ is a Young subgroup: with the notation in \eqref{VsIsigma}, we have 
\begin{equation}\label{IsigmaSym}
\mathsf{I}(\sigma)=\prod_{\sigma_*\in\mathsf{image}(\sigma)}\mathsf{Sym}\left(V_{\sigma_*}\right). 
\end{equation}
Then we introduce some specific notation.
\begin{itemize}
\item
For every choice of $\xi\in\mathsf{Surj}(V,X)$, $\mathsf{j}\in\mathbb{N}^V$ and $\sigma\in\Sigma_0$, 
\begin{itemize}
\item
we denote by $\mathsf{supp}(\xi\wedge\mathsf{j}\wedge\sigma)$ the set of all triples $\left(x,r,\sigma_*\right)\in X\times\mathbb{N}\times\mathsf{image}(\sigma)$ such that $\xi^{-1}(x)\cap\mathsf{j}^{-1}(r)\cap V_{\sigma_*}\neq\emptyset$;
\item 
we denote by $\xi\wedge\mathsf{j}\wedge\sigma$ the set composition 
\[
V=\bigsqcup\limits_{\left(x,r,\sigma_*\right)\in\mathsf{supp}(\xi\wedge\mathsf{j}\wedge\sigma)}\left[\xi^{-1}(x)\cap\mathsf{j}^{-1}(r)\cap V_{\sigma_*}\right];
\]
\item
we denote by $\lvert\xi\wedge\mathsf{j}\wedge\sigma\rvert$ the positive integer-valued function 
\[
\mathsf{supp}(\xi\wedge\mathsf{j}\wedge\sigma)\ni
\left(x,r,\sigma_*\right)\longmapsto\lvert\xi\wedge\mathsf{j}\wedge\sigma\rvert\left(x,r,\sigma_*\right)\coloneqq\left\lvert\xi^{-1}(x)\cap\mathsf{j}^{-1}(r)\cap V_{\sigma_*}\right\rvert.
\]
\end{itemize}
\item
For every choice of $\sigma\in\Sigma_0$, we denote by $\mathsf{cosupp}(\sigma)$ the set of all triples $\left(x,r,\sigma_*\right)\in X\times\mathbb{N}\times\mathsf{image}(\sigma)$ such that $\sigma_*\in\widehat{F}_{\alpha(x)}$ and $1\leq r\leq m_{\alpha(x),\sigma_*}$.
\item
For every choice of $\xi\in\mathsf{Surj}(V,X)$ and $\sigma\in\Sigma_0$, 
\begin{itemize}
\item
we denote by $\xi^{-1}\wedge\sigma$ the family of all set compositions $V=\bigsqcup_{\left(x,r,\sigma_*\right)\in\mathsf{supp}(\mathsf{c})}\mathsf{c}\left(x,r,\sigma_*\right)$, where $\mathsf{supp}(\mathsf{c})\subseteq X\times\mathbb{N}\times\mathsf{image}(\sigma)$ is the index set of $\mathsf{c}$, such that 
\begin{itemize}
\item
$\mathsf{supp}(\mathsf{c})\subseteq\mathsf{cosupp}(\sigma)$,
\item
 $\sum_{r,\sigma_*\mathrel{:} \left(x,r,\sigma_*\right)\in\mathsf{supp}(\mathsf{c})}\left\lvert \mathsf{c}\left(x,r,\sigma_*\right)\right\rvert=\left\lvert\xi^{-1}(x)\right\rvert$ for all $x\in X$,
\item
$\bigsqcup_{x,r\mathrel{:} \left(x,r,\sigma_*\right)\in\mathsf{supp}(\mathsf{c})} \mathsf{c}\left(x,r,\sigma_*\right)=V_{\sigma_*}$ for all $\sigma_*\in\mathsf{image}(\sigma)$; 
\end{itemize}
\item 
we denote by $\left\lvert\xi^{-1}\wedge\sigma\right\rvert$ the set of all positive integer-valued functions $\tau\colon\mathsf{supp}(\tau)\ni\left(x,r,\sigma_*\right)\mapsto\tau\left(x,r,\sigma_*\right)\in\mathbb{N}$, where $\mathsf{supp}(\tau)\subseteq X\times\mathbb{N}\times\mathsf{image}(\sigma)$ is the domain of $\tau$, such that
\begin{itemize}
\item 
$\mathsf{supp}(\tau)\subseteq\mathsf{cosupp}(\sigma)$,
\item 
$\sum_{r,\sigma_*\mathrel{:} \left(x,r,\sigma_*\right)\in\mathsf{supp}(\tau)}\tau\left(x,r,\sigma_*\right)=\left\lvert\xi^{-1}(x)\right\rvert$ for all $x\in X$,
\item
\begin{equation}\label{recalleq}
\sum_{x,r\mathrel{:} \left(x,r,\sigma_*\right)\in\mathsf{supp}(\tau)} \tau\left(x,r,\sigma_*\right)=\left\lvert V_{\sigma_*}\right\rvert\text{ for all }\sigma_*\in\mathsf{image}(\sigma); 
\end{equation}
\end{itemize}
\item
for $\tau\in\left\lvert\xi^{-1}\wedge\sigma\right\rvert$ and $\sigma_*\in\mathsf{image}(\sigma)$, we denote by $\tau^{\left(\sigma_*\right)}$ the numerical partition of $\left\lvert V_{\sigma_*}\right\rvert$ obtained by arranging the positive integers $\tau\left(x,r,\sigma_*\right)$ ($\left(x,r,\sigma_*\right)\in\mathsf{supp}(\tau)$) in weakly decreasing order (cf.\ \eqref{recalleq});
\item
if $\tau,\tau'\in\left\lvert\xi^{-1}\wedge\sigma\right\rvert$, then we say that they are equivalent, and we write $\tau\approx\tau'$, if $\tau^{\left(\sigma_*\right)}=\tau'^{\left(\sigma_*\right)}$ for all $\sigma_*\in\mathsf{image}(\sigma)$ (equivalently, the corresponding Young subgroups of $\mathsf{Sym}\left(V_{\sigma_*}\right)$ are conjugate for all $\sigma_*\in\mathsf{image}(\sigma)$). In particular, if $\tau\approx\tau'$, then 
\begin{equation}\label{deftaueqtaup}
\prod_{x,r\mathrel{:}\left(x,r,\sigma_*\right)\in\mathsf{supp}(\tau)}\mathsf{Sym}_{\tau\left(x,r,\sigma_*\right)}\cong \prod_{x,r\mathrel{:}\left(x,r,\sigma_*\right)\in\mathsf{supp}(\tau')}\mathsf{Sym}_{\tau'\left(x,r,\sigma_*\right)}
\end{equation}
for all $\sigma_*\in\mathsf{image}(\sigma)$. We denote by 
$\left\lvert\xi^{-1}\wedge\sigma\right\rvert/\approx$ the corresponding quotient set and by $[\tau]$ the equivalence class containing $\tau$.
\end{itemize}
\end{itemize}

Regarding the definition of $\approx$, if $\tau\in\left\lvert\xi^{-1}\wedge\sigma\right\rvert$, $x\in X$, and $\sigma_*\in\mathsf{image}(\sigma)$ satisfy $m_{\alpha(x),\sigma_*}\geq 2$, with $\left(x,1,\sigma_*\right)\in\mathsf{supp}(\tau)$ and $\left(x,2,\sigma_*\right)\notin\mathsf{supp}(\tau)$, then we may define $\tau'\in\left\lvert\xi^{-1}\wedge\sigma\right\rvert$ with $\mathsf{supp}(\tau')=\bigl(\mathsf{supp}(\tau)\setminus\left\{\left(x,1,\sigma_*\right)\right\}\bigr)\cup\left\{\left(x,2,\sigma_*\right)\right\}$ by setting $\tau'\left(x,2,\sigma_*\right)=\tau\left(x,1,\sigma_*\right)$ and $\tau'\left(x',r,\sigma_\dag\right)=\tau\left(x',r,\sigma_\dag\right)$ for all $\left(x',r,\sigma_\dag\right)\in\mathsf{supp}(\tau)\setminus\left\{\left(x,1,\sigma_*\right)\right\}$. Then, clearly $\mathsf{supp}(\tau)\neq\mathsf{supp}(\tau')$ but $\tau\approx\tau'$.

From now on, we fix $\xi_0\in\mathsf{Surj}(V,X)$.
\begin{proposition}\label{ProporbstabSym}
\begin{enumerate}[(1)]
\item
For every $\sigma\in\Sigma_0$ we have
\begin{equation}\label{DecJsigSym}
\begin{split}
\mathsf{J}_\sigma=&
\left\{(\xi,\mathsf{j})\in\mathsf{Surj}(V,X)\times\mathbb{N}^V\mathrel{:}\xi\wedge\mathsf{j}\wedge\sigma\in\xi_0^{-1}\wedge\sigma\right\}\\
=&\bigsqcup_{\tau\in\left\lvert\xi_0^{-1}\wedge\sigma\right\rvert}\left\{(\xi,\mathsf{j})\in\mathsf{Surj}(V,X)\times\mathbb{N}^V\mathrel{:}\lvert\xi\wedge\mathsf{j}\wedge\sigma\rvert=\tau\right\}.
\end{split}
\end{equation}
The last expression is also the decomposition into the orbits of $\mathsf{I}(\sigma)$ on $\mathsf{J}_\sigma$. 
\item
If $(\xi,\mathsf{j})\in\mathsf{J}_\sigma$, then 
\[
\begin{split}
\mathsf{Stab}_{\mathsf{I}(\sigma)}(\xi,\mathsf{j})=&\prod_{\sigma_*\in\mathsf{image}(\sigma)}\left(\prod_{x,r\mathrel{:}\left(x,r,\sigma_*\right)\in\mathsf{supp}(\xi\wedge\mathsf{j}\wedge\sigma)}\mathsf{Sym}\left(\xi^{-1}(x)\cap\mathsf{j}^{-1}(r)\cap V_{\sigma_*}\right)\right)\\
\cong&\prod_{\sigma_*\in\mathsf{image}(\sigma)}\left(\prod_{x,r\mathrel{:}\left(x,r,\sigma_*\right)\in\mathsf{supp}(\xi\wedge\mathsf{j}\wedge\sigma)}\mathsf{Sym}_{\lvert\xi\wedge\mathsf{j}\wedge\sigma\rvert\left(x,r,\sigma_*\right)}\right).
\end{split}
\]
\end{enumerate}
\end{proposition}
\begin{proof}
(1) 
First, note that $\mathsf{supp}(\xi\wedge\mathsf{j}\wedge\sigma)\subseteq\mathsf{cosupp}(\sigma)$ is equivalent to the occurrence of both the following conditions: $\sigma\prec\xi$ and $1\leq \mathsf{j}(v)\leq m_{\alpha(\xi(v)),\sigma_v}$ for all $v\in V$.
Since $\sigma\in\Sigma_0$ we have also $\sigma\prec\xi_0$. Moreover, the condition
\[
\sum_{r,\sigma_*\mathrel{:} \left(x,r,\sigma_*\right)\in\mathsf{supp}(\xi\wedge\mathsf{j}\wedge\sigma)}\left\lvert (\xi\wedge\mathsf{j}\wedge\sigma)\left(x,r,\sigma_*\right)\right\rvert=\left\lvert\xi^{-1}_0(x)\right\rvert\quad\text{for all}\quad x\in X
\]
is equivalent to the existence of an element $\pi\in\mathsf{Sym}(V)$ such that $\pi\xi_0=\xi$. Therefore, from \eqref{condj} and \eqref{defJsigma} we deduce the first equality in the statement. In other words, the map
\[
\mathsf{J}_\sigma\ni(\xi,\mathsf{j})\longmapsto\xi\wedge\mathsf{j}\wedge\sigma\in\xi_0^{-1}\wedge\sigma
\]
is a bijection. Note also that $\bigsqcup_{x,r}(\xi\wedge\mathsf{j}\wedge\sigma)(x,r,\sigma_*)=V_{\sigma_*}$ holds automatically, that $(\xi,\mathsf{j})$ is uniquely recovered from the index $(x,r,\sigma_*)$ of the unique block containing $v$ (giving injectivity), and that any $\mathsf{c}\in\xi_0^{-1}\wedge\sigma$ defines such a pair $(\xi,\mathsf{j})$ with $\left\lvert\xi^{-1}(x)\right\rvert=\left\lvert\xi_0^{-1}(x)\right\rvert\geq 1$ (giving surjectivity).
The second equality in the statement follows by observing that the map
\[
\xi_0^{-1}\wedge\sigma\ni\xi\wedge\mathsf{j}\wedge\sigma\longmapsto\lvert\xi\wedge\mathsf{j}\wedge\sigma\rvert\in\left\lvert\xi_0^{-1}\wedge\sigma\right\rvert
\]
is surjective. Indeed, given $\tau\in\left\lvert\xi_0^{-1}\wedge\sigma\right\rvert$, by \eqref{recalleq} we can choose a function $\xi\colon V\to X$ such that 
\[
\left\lvert\xi^{-1}(x)\cap V_{\sigma_*}\right\rvert=\sum_{r\mathrel{:} \left(x,r,\sigma_*\right)\in\mathsf{supp}(\tau)}\tau\left(x,r,\sigma_*\right)
\]
for all $x\in X$ and $\sigma_*\in\mathsf{image}(\sigma)$. Summing over $\sigma_*\in\mathsf{image}(\sigma)$, the second condition in the definition of $\left\lvert\xi_0^{-1}\wedge\sigma\right\rvert$ guarantees $\left\lvert\xi^{-1}(x)\right\rvert=\left\lvert\xi_0^{-1}(x)\right\rvert\geq 1$ for all $x\in X$, so $\xi\in\mathsf{Surj}(V,X)$ and $\xi\in\mathsf{Orb}_{\mathsf{Sym}(V)}(\xi_0)$. Next, we can define $\mathsf{j}\in\mathbb{N}^V$ in such a way that $\lvert\xi\wedge\mathsf{j}\wedge\sigma\rvert=\tau$. In particular, $\mathsf{supp}(\xi\wedge\mathsf{j}\wedge\sigma)=\mathsf{supp}(\tau)\subseteq\mathsf{cosupp}(\sigma)$ ensures that \eqref{condj} holds. Finally, every $\tau\in\left\lvert\xi_0^{-1}\wedge\sigma\right\rvert$ determines an orbit of $\mathsf{I}(\sigma)$ because this group is precisely the stabilizer of each block of the set composition $V=\bigsqcup_{\sigma_*\in\mathsf{image}(\sigma)}V_{\sigma_*}$.

\medskip
(2)
An element $\pi\in\mathsf{I}(\sigma)$ stabilizes $(\xi,\mathsf{j})$ if and only if it stabilizes each block of the set composition 
$V=\bigsqcup_{\left(x,r,\sigma_*\right)\in\mathsf{supp}(\xi\wedge\mathsf{j}\wedge\sigma)}\xi^{-1}(x)\cap\mathsf{j}^{-1}(r)\cap V_{\sigma_*}$.
\end{proof}

Now we recall some basic facts on the representation theory of the symmetric group; see, e.g., \cite{book2,JK,Mayer,SagSym}.
If $\mu$ and $\lambda$ are numerical partitions of the positive integer $n$ with $\lambda\trianglelefteq\mu$ (where $\trianglelefteq$ is the majorization order), we denote by $\mathsf{STab}(\mu,\lambda)$ the set of semistandard tableaux of shape $\mu$ and weight (or content) $\lambda$; its cardinality is equal to the multiplicity of $S^\mu$, the irreducible representation canonically associated with $\mu$, in the permutation module $M^\lambda$. Similarly, an irreducible representation of $\mathsf{I}(\sigma)$ has the form $\bigotimes_{\sigma_*\in\mathsf{image}(\sigma)}S^{\mu^{\left(\sigma_*\right)}}$, where $\mu^{\left(\sigma_*\right)}$ is a partition of $\left\lvert V_{\sigma_*}\right\rvert$. To unify and simplify the notation, the space of $\sigma$ will be denoted by $\bigotimes_{\sigma_*\in\mathsf{image}(\sigma)}S^{\sigma_*}$.

\begin{theorem}\label{ThmmultSym}
For every $\sigma\in\Sigma_0$ and every family of partitions $\mu^{\left(\sigma_*\right)}$ of $\left\lvert V_{\sigma_*}\right\rvert$ ($\sigma_*\in\mathsf{image}(\sigma)$), the multiplicity of the irreducible $\mathsf{I}(\sigma)$-representation $\bigotimes_{\sigma_*\in\mathsf{image}(\sigma)}S^{\mu^{\left(\sigma_*\right)}}$ in $L\left(\mathsf{J}_\sigma\right)$ is equal to
\begin{equation}\label{Bigformult}
\sum_{\substack{[\tau]\in\left\lvert\xi_0^{-1}\wedge\sigma\right\rvert/\approx\mathrel{:}\\ \tau^{\left(\sigma_*\right)}\trianglelefteq\mu^{\left(\sigma_*\right)}\; \text{\rm for }\sigma_*\in\mathsf{image}(\sigma)}}\lvert[\tau]\rvert\prod_{\sigma_*\in\mathsf{image}(\sigma)}\left\lvert\mathsf{STab}\left(\mu^{\left(\sigma_*\right)},\tau^{\left(\sigma_*\right)}\right)\right\rvert.
\end{equation}
\end{theorem}
\begin{proof}
From Proposition \ref{ProporbstabSym} it follows that the permutation representation of the group \eqref{IsigmaSym} on the orbit of $\mathsf{J}_\sigma$ corresponding to $\tau\in\left\lvert\xi_0^{-1}\wedge\sigma\right\rvert$ is equivalent to the tensor product of permutation representations $\bigotimes_{\sigma_*\in\mathsf{image}(\sigma)}M^{\tau^{\left(\sigma_*\right)}}$. Then an application of the Young Rule to each factor of the tensor product in \eqref{IsigmaSym} leads to the following intermediate result for the multiplicity in the statement:
\[
\sum_{\substack{\tau\in\left\lvert\xi_0^{-1}\wedge\sigma\right\rvert\mathrel{:}\\ \tau^{\left(\sigma_*\right)}\trianglelefteq\mu^{\left(\sigma_*\right)}\; \text{\rm for }\sigma_*\in\mathsf{image}(\sigma)}}\prod_{\sigma_*\in\mathsf{image}(\sigma)}\left\lvert\mathsf{STab}\left(\mu^{\left(\sigma_*\right)},\tau^{\left(\sigma_*\right)}\right)\right\rvert.
\]
Finally, by definition of $\approx$, equivalent elements $\tau\in\left\lvert\xi_0^{-1}\wedge\sigma\right\rvert$ determine the same partitions $\tau^{\left(\sigma_*\right)}$ (for all $\sigma_*\in\mathsf{image}(\sigma)$) and therefore yield the same summand above; grouping the sum over $\left\lvert\xi_0^{-1}\wedge\sigma\right\rvert/\approx$ proves \eqref{Bigformult}.
\end{proof}

\begin{corollary}
For every $\sigma\in\Sigma_0$ and every family of partitions $\mu^{\left(\sigma_*\right)}$ of $\left\lvert V_{\sigma_*}\right\rvert$ ($\sigma_*\in\mathsf{image}(\sigma)$), the multiplicity of the irreducible $F\wr\mathsf{Sym}(V)$-representation 
\[
\Ind^{F\wr\mathsf{Sym}(V)}_{F\wr\mathsf{I}(\sigma)}\left[\widetilde{\left(\bigotimes_{\sigma_*\in\mathsf{image}(\sigma)}S^{\sigma_*}\right)} \otimes\overline{\left(\bigotimes_{\sigma_*\in\mathsf{image}(\sigma)}S^{\mu^{\left(\sigma_*\right)}}\right)}\right]
\]
in the permutation representation $L\left(\mathsf{Orb}_{F\wr\mathsf{Sym}(V)}\left(\xi_0,\varphi_{\xi_0}\right)\right)$ is equal to \eqref{Bigformult}.
\end{corollary}
\begin{proof}
Immediate consequence of Theorem \ref{ThmmultSym} and Corollary \ref{Cormult}.
\end{proof}

\begin{example}\label{Exsymmultfr}
Now we examine the following case: for every $a\in A$ the $F$-permutation representation $L\left(Y_a\right)$ decomposes without multiplicity. As a consequence, the indices $\mathsf{j}$ and $r$ disappear. We briefly write the corresponding changes in the definitions. If $\xi\in\mathsf{Surj}(V,X)$ and $\sigma\in\Sigma_0$, then:
\begin{itemize}
\item
$\mathsf{supp}(\xi\wedge\sigma)$ is the set of all pairs $\left(x,\sigma_*\right)\in X\times\mathsf{image}(\sigma)$ such that $\xi^{-1}(x)\cap V_{\sigma_*}\neq\emptyset$, $\xi\wedge\sigma$ is the set composition 
\[
V=\bigsqcup_{\left(x,\sigma_*\right)\in\mathsf{supp}(\xi\wedge\sigma)}\left[\xi^{-1}(x)\cap V_{\sigma_*}\right], 
\]
and $\lvert\xi\wedge\sigma\rvert$ is the positive integer-valued function $\mathsf{supp}(\xi\wedge\sigma)\ni\left(x,\sigma_*\right)\mapsto\left\lvert\xi^{-1}(x)\cap V_{\sigma_*}\right\rvert$;
\item
$\mathsf{cosupp}(\sigma)$ is the set of all pairs $\left(x,\sigma_*\right)\in X\times\mathsf{image}(\sigma)$ such that $\sigma_*\in\widehat{F}_{\alpha(x)}$; we have $\sigma\prec\xi$ if and only if $\mathsf{supp}(\xi\wedge\sigma)\subseteq\mathsf{cosupp}(\sigma)$; 
\item
$\xi^{-1}\wedge\sigma$ is the family of all set compositions $V=\bigsqcup_{\left(x,\sigma_*\right)\in\mathsf{supp}(\mathsf{c})}\mathsf{c}\left(x,\sigma_*\right)$ such that $\mathsf{supp}(\mathsf{c})\subseteq\mathsf{cosupp}(\sigma)$, $\sum_{\sigma_*\mathrel{:}\left(x,\sigma_*\right)\in\mathsf{supp}(\mathsf{c})}\left\lvert\mathsf{c}\left(x,\sigma_*\right)\right\rvert=\left\lvert\xi^{-1}(x)\right\rvert$ for all $x\in X$, and $\bigsqcup_{x\mathrel{:} \left(x,\sigma_*\right)\in\mathsf{supp}(\mathsf{c})} \mathsf{c}\left(x,\sigma_*\right)=V_{\sigma_*}$ for all $\sigma_*\in\mathsf{image}(\sigma)$; $\left\lvert\xi^{-1}\wedge\sigma\right\rvert$ is the set of all positive integer-valued functions $\tau\colon\mathsf{supp}(\tau)\ni\left(x,\sigma_*\right)\mapsto\tau\left(x,\sigma_*\right)\in\mathbb{N}$ such that $\mathsf{supp}(\tau)\subseteq\mathsf{cosupp}(\sigma)$, $\sum_{\sigma_*\mathrel{:}\left(x,\sigma_*\right)\in\mathsf{supp}(\tau)}\tau\left(x,\sigma_*\right)=\left\lvert\xi^{-1}(x)\right\rvert$ for all $x\in X$, and $\sum_{x\mathrel{:} \left(x,\sigma_*\right)\in\mathsf{supp}(\tau)} \tau\left(x,\sigma_*\right)=\left\lvert V_{\sigma_*}\right\rvert$ for all $\sigma_*\in\mathsf{image}(\sigma)$. 
\end{itemize}

Then we fix $\xi_0\in\mathsf{Surj}(V,X)$ and \eqref{DecJsigSym} becomes
\[
\begin{split}
\mathsf{J}_\sigma=\left\{\xi\in\mathsf{Orb}_{\mathsf{Sym}(V)}\left(\xi_0\right)\mathrel{:}\sigma\prec\xi\right\}=&
\left\{\xi\in\mathsf{Surj}(V,X)\mathrel{:}\xi\wedge\sigma\in\xi_0^{-1}\wedge\sigma\right\}\\
=&\bigsqcup_{\tau\in\left\lvert\xi_0^{-1}\wedge\sigma\right\rvert}\left\{\xi\in\mathsf{Surj}(V,X)\mathrel{:}\lvert\xi\wedge\sigma\rvert=\tau\right\}.
\end{split}
\]
The last expression is also the decomposition into the orbits of $\mathsf{I}(\sigma)$; that is, every $\tau\in\left\lvert\xi_0^{-1}\wedge\sigma\right\rvert$ determines an orbit of $\mathsf{I}(\sigma)$ on $\mathsf{J}_\sigma$. Similarly, if $\xi\in\mathsf{J}_\sigma$, then 
\[
\begin{split}
\mathsf{Stab}_{\mathsf{I}(\sigma)}(\xi)=&\prod_{\sigma_*\in\mathsf{image}(\sigma)}\left(\prod_{x\mathrel{:}\left(x,\sigma_*\right)\in\mathsf{supp}(\xi\wedge\sigma)}\mathsf{Sym}\left(\xi^{-1}(x)\cap V_{\sigma_*}\right)\right)\\
\cong&\prod_{\sigma_*\in\mathsf{image}(\sigma)}\left(\prod_{x\mathrel{:}
\left(x,\sigma_*\right)\in\mathsf{supp}(\xi\wedge\sigma)}\mathsf{Sym}_{\lvert\xi\wedge\sigma\rvert\left(x,\sigma_*\right)}\right).
\end{split}
\]

For $\tau\in\left\lvert\xi_0^{-1}\wedge\sigma\right\rvert$ and $\sigma_*\in\mathsf{image}(\sigma)$, we denote by $\tau^{\left(\sigma_*\right)}$ the numerical partition of $\left\lvert V_{\sigma_*}\right\rvert$ obtained by arranging the positive integers $\tau\left(x,\sigma_*\right)$ ($\left(x,\sigma_*\right)\in\mathsf{supp}(\tau)$) in weakly decreasing order. We then define the equivalence relation $\approx$ on $\left\lvert\xi_0^{-1}\wedge\sigma\right\rvert$ by setting $\tau\approx\tau'$ if $\tau^{\left(\sigma_*\right)}=\tau'^{\left(\sigma_*\right)}$ for all $\sigma_*\in\mathsf{image}(\sigma)$ (equivalently, \eqref{deftaueqtaup} holds with the index $r$ deleted). Also in this case distinct $\tau$ may be equivalent: for instance, if $x,x'\in X$ are distinct, $\sigma_*,\sigma_\dagger\in\mathsf{image}(\sigma)$ are distinct, and $\tau\left(x,\sigma_*\right)=1=\tau\left(x',\sigma_\dagger\right)$, $\tau\left(x,\sigma_\dagger\right)=2=\tau\left(x',\sigma_*\right)$, then we may define $\tau'\in\left\lvert\xi_0^{-1}\wedge\sigma\right\rvert$ with $\mathsf{supp}(\tau')=\mathsf{supp}(\tau)$ by swapping these values, namely $\tau'\left(x,\sigma_*\right)=2=\tau'\left(x',\sigma_\dagger\right)$ and $\tau'\left(x,\sigma_\dagger\right)=1=\tau'\left(x',\sigma_*\right)$, while leaving all other values unchanged. Moreover, it will be convenient to extend each $\tau\in\left\lvert\xi_0^{-1}\wedge\sigma\right\rvert$ to $X\times\mathsf{image}(\sigma)$ by setting $\tau\left(x,\sigma_*\right)\coloneqq 0$ whenever $\left(x,\sigma_*\right)\notin\mathsf{supp}(\tau)$. Finally, for every $\sigma\in\Sigma_0$ and every family of partitions $\mu^{\left(\sigma_*\right)}$ of $\left\lvert V_{\sigma_*}\right\rvert$ ($\sigma_*\in\mathsf{image}(\sigma)$), the multiplicity of $\bigotimes_{\sigma_*\in\mathsf{image}(\sigma)}S^{\mu^{\left(\sigma_*\right)}}$ in $L\left(\mathsf{J}_\sigma\right)$ is equal to
\begin{equation}\label{Bigformult2}
\sum_{\substack{[\tau]\in\left\lvert\xi_0^{-1}\wedge\sigma\right\rvert/\approx\mathrel{:}\\ \tau^{\left(\sigma_*\right)}\trianglelefteq\mu^{\left(\sigma_*\right)}\; \text{\rm for }\sigma_*\in\mathsf{image}(\sigma)}}\lvert[\tau]\rvert\prod_{\sigma_*\in\mathsf{image}(\sigma)}\left\lvert\mathsf{STab}\left(\mu^{\left(\sigma_*\right)},\tau^{\left(\sigma_*\right)}\right)\right\rvert.
\end{equation}
\end{example}

\begin{theorem}\label{TheoSymmultfree}
The permutation representation of $F\wr\mathsf{Sym}(V)$ on $\mathsf{Orb}_{F\wr\mathsf{Sym}(V)}(\xi_0,\varphi_{\xi_0})$ is multiplicity free if and only if the following conditions are satisfied:
\begin{enumerate}[(1)]
\item\label{TheoSymmultfree1}
for every $a\in A$ the $F$-permutation representation $L\left(Y_a\right)$ decomposes without multiplicity, that is, we are in the framework of Example \ref{Exsymmultfr};
\item\label{TheoSymmultfree2}
at least one of the following conditions holds:
\begin{itemize}
\item
$\lvert X\rvert=1$;
\item
$\lvert X\rvert=2$, say $X=\{x,y\}$, and $\left\lvert\widehat{F}_{\alpha(x)}\cap\widehat{F}_{\alpha(y)}\right\rvert=1$ (in particular, if $\alpha(x)=\alpha(y)=a$, then $\left\lvert\widehat{F}_a\right\rvert=1$).
\end{itemize}
\end{enumerate}
\end{theorem}

\begin{proof}
By Corollary \ref{Cormultfree}\eqref{Cormultfree1}, condition \eqref{TheoSymmultfree1} is necessary for the permutation representation to be multiplicity free. Assuming \eqref{TheoSymmultfree1}, the multiplicity of $\bigotimes_{\sigma_*\in\mathsf{image}(\sigma)} S^{\mu^{\left(\sigma_*\right)}}$ in $L\left(\mathsf{J}_\sigma\right)$ is given by \eqref{Bigformult2}.

We first show that $\lvert X\rvert\leq 2$ is necessary. Since $L\left(Y_a\right)$ is a permutation representation, the trivial representation $\sigma^0$ belongs to $\widehat{F}_a$ for every $a\in A$. Consider $\sigma\in\Sigma_0$ with $\mathsf{image}(\sigma)=\{\sigma^0\}$ and $V_{\sigma^0}=V$. Because $\xi_0\in\mathsf{Surj}(V,X)$, the unique element $\tau\in\left\lvert\xi_0^{-1}\wedge\sigma\right\rvert$ is given by $\tau(x,\sigma^0)=\left\lvert\xi_0^{-1}(x)\right\rvert\geq 1$ for all $x\in X$, so $\tau^{(\sigma^0)}$ is a numerical composition of $\lvert V\rvert$ with $k=\lvert X\rvert$ positive parts. Let $\lambda=(\lambda_1,\dotsc,\lambda_k)$ be the partition obtained by arranging the parts of $\tau^{(\sigma^0)}$ in weakly decreasing order. If $k\geq 3$, set $\mu=(\lambda_1+1,\lambda_2,\dotsc,\lambda_{k-1},\lambda_k-1)$ (omitting the last entry if $\lambda_k=1$). Then $\mu$ is a partition of $\lvert V\rvert$ satisfying $\lambda\trianglelefteq\mu$, and there exist at least two distinct semistandard Young tableaux of shape $\mu$ and content $\lambda$:
\begin{itemize}
\item
$T_1$, whose first row consists of $\lambda_1$ entries equal to $1$ followed by one $k$, whose $i$-th row ($2\leq i\leq k-1$) consists of $\lambda_i$ entries equal to $i$, and (if $\lambda_k> 1$) whose $k$-th row consists of $\lambda_k-1$ entries equal to $k$;
\item
$T_2$, whose first row consists of $\lambda_1$ entries equal to $1$ followed by one $2$, whose $i$-th row ($2\leq i\leq k-1$) consists of $\lambda_i-1$ entries equal to $i$ followed by one $i+1$, and (if $\lambda_k> 1$) whose $k$-th row consists of $\lambda_k-1$ entries equal to $k$.
\end{itemize}
Since $k\geq 3$, the rightmost entry of the first row distinguishes $T_1$ from $T_2$, so $\left\lvert\mathsf{STab}\left(\mu,\tau^{(\sigma^0)}\right)\right\rvert\geq 2$. By \eqref{Bigformult2}, the representation is not multiplicity free. Hence $\lvert X\rvert\leq 2$.

If $\lvert X\rvert=1$, then we are in the case of Example \ref{Example1ter}, which yields a multiplicity free permutation representation.

Assume henceforth that $\lvert X\rvert=2$, $X=\{x,y\}$, and set $a=\alpha(x)$ and $b=\alpha(y)$ (where $a$ and $b$ need not be distinct). Notice that $\sigma^0\in\widehat{F}_a\cap\widehat{F}_b$, so $\left\lvert\widehat{F}_a\cap\widehat{F}_b\right\rvert\geq 1$.

To prove necessity, suppose that $\left\lvert\widehat{F}_a\cap\widehat{F}_b\right\rvert\geq 2$, and choose $\sigma_1\in\left(\widehat{F}_a\cap\widehat{F}_b\right)\setminus\{\sigma^0\}$. Let $m=\left\lvert\xi_0^{-1}(x)\right\rvert\geq 1$ and $n=\left\lvert\xi_0^{-1}(y)\right\rvert\geq 1$, so that $\lvert V\rvert=m+n\geq 2$. Pick $\sigma\in\Sigma_0$ with $\mathsf{image}(\sigma)=\{\sigma^0,\sigma_1\}$, $\left\lvert V_{\sigma^0}\right\rvert=\lvert V\rvert-1$, and $\left\lvert V_{\sigma_1}\right\rvert=1$. Since $\sigma^0,\sigma_1\in\widehat{F}_a\cap\widehat{F}_b$, we have $\mathsf{cosupp}(\sigma)=X\times\{\sigma^0,\sigma_1\}$. Define $\tau_1,\tau_2\in\left\lvert\xi_0^{-1}\wedge\sigma\right\rvert$ on their respective supports by
\begin{align*}
\tau_1\left(x,\sigma^0\right)&=m, & \tau_1\left(y,\sigma_1\right)&=1, & \tau_1\left(y,\sigma^0\right)&=n-1 \quad\text{if } n> 1,\\
\tau_2\left(y,\sigma^0\right)&=n, & \tau_2\left(x,\sigma_1\right)&=1, & \tau_2\left(x,\sigma^0\right)&=m-1 \quad\text{if } m> 1.
\end{align*}
Then $\tau_1$ and $\tau_2$ are distinct elements of $\left\lvert\xi_0^{-1}\wedge\sigma\right\rvert$, since $\left(y,\sigma_1\right)\in\mathsf{supp}(\tau_1)\setminus\mathsf{supp}(\tau_2)$. Choosing the one-row partitions $\mu^{(\sigma^0)}=(\lvert V\rvert-1)$ and $\mu^{(\sigma_1)}=(1)$, every $\tau\in\left\lvert\xi_0^{-1}\wedge\sigma\right\rvert$ satisfies $\tau^{(\sigma_*)}\trianglelefteq\mu^{(\sigma_*)}$ and $\left\lvert\mathsf{STab}\left(\mu^{(\sigma_*)},\tau^{(\sigma_*)}\right)\right\rvert=1$ for both $\sigma_*\in\{\sigma^0,\sigma_1\}$. Therefore \eqref{Bigformult2} yields
\[
\sum_{[\tau]\in\left\lvert\xi_0^{-1}\wedge\sigma\right\rvert/\approx}\lvert[\tau]\rvert=\left\lvert\left\lvert\xi_0^{-1}\wedge\sigma\right\rvert\right\rvert\geq 2,
\]
contradicting multiplicity freeness. Thus $\left\lvert\widehat{F}_a\cap\widehat{F}_b\right\rvert=1$.

Conversely, suppose that $\widehat{F}_a\cap\widehat{F}_b=\{\sigma^0\}$, and fix $\sigma\in\Sigma_0$ such that $\left\lvert\xi_0^{-1}\wedge\sigma\right\rvert\neq\emptyset$. Let $\tau\in\left\lvert\xi_0^{-1}\wedge\sigma\right\rvert$. Because $\mathsf{supp}(\tau)\subseteq\mathsf{cosupp}(\sigma)$ and $\left\lvert V_{\sigma_*}\right\rvert\geq 1$ for all $\sigma_*\in\mathsf{image}(\sigma)$, we have $\mathsf{image}(\sigma)\subseteq\widehat{F}_a\cup\widehat{F}_b$, and:
\begin{itemize}
\item
if $\sigma_*\in\mathsf{image}(\sigma)\cap\left(\widehat{F}_a\setminus\{\sigma^0\}\right)$, then $\left(y,\sigma_*\right)\notin\mathsf{cosupp}(\sigma)$, so $\left(x,\sigma_*\right)\in\mathsf{supp}(\tau)$ and $\tau\left(x,\sigma_*\right)=\left\lvert V_{\sigma_*}\right\rvert$;
\item
if $\sigma_*\in\mathsf{image}(\sigma)\cap\left(\widehat{F}_b\setminus\{\sigma^0\}\right)$, then $\left(x,\sigma_*\right)\notin\mathsf{cosupp}(\sigma)$, so $\left(y,\sigma_*\right)\in\mathsf{supp}(\tau)$ and $\tau\left(y,\sigma_*\right)=\left\lvert V_{\sigma_*}\right\rvert$;
\item
if $\sigma^0\in\mathsf{image}(\sigma)$, then the row-sum constraints uniquely determine the values of $\tau$ at $\sigma^0$, namely
\begin{align*}
\tau\left(x,\sigma^0\right)&=\left\lvert\xi_0^{-1}(x)\right\rvert-\sum_{\sigma_*\in\mathsf{image}(\sigma)\cap\left(\widehat{F}_a\setminus\{\sigma^0\}\right)}\left\lvert V_{\sigma_*}\right\rvert,\\
\tau\left(y,\sigma^0\right)&=\left\lvert\xi_0^{-1}(y)\right\rvert-\sum_{\sigma_*\in\mathsf{image}(\sigma)\cap\left(\widehat{F}_b\setminus\{\sigma^0\}\right)}\left\lvert V_{\sigma_*}\right\rvert,
\end{align*}
with $\left(x,\sigma^0\right)$ (respectively, $\left(y,\sigma^0\right)$) belonging to $\mathsf{supp}(\tau)$ precisely when the corresponding right-hand side is strictly positive.
\end{itemize}
Consequently, $\left\lvert\left\lvert\xi_0^{-1}\wedge\sigma\right\rvert\right\rvert=1$, so the sum in \eqref{Bigformult2} contains at most one term, with $\lvert[\tau]\rvert=1$. Furthermore, for this unique $\tau$, every $\tau^{(\sigma_*)}$ has at most two parts (in fact, one part for $\sigma_*\neq\sigma^0$ and at most two parts for $\sigma_*=\sigma^0$). Whenever $\tau^{(\sigma_*)}\trianglelefteq\mu^{(\sigma_*)}$, the partition $\mu^{(\sigma_*)}$ also has at most two parts, and by Young's rule $\left\lvert\mathsf{STab}\left(\mu^{(\sigma_*)},\tau^{(\sigma_*)}\right)\right\rvert=1$. Hence the multiplicity \eqref{Bigformult2} is at most $1$ for every $\sigma\in\Sigma_0$ and every choice of $\left(\mu^{(\sigma_*)}\right)_{\sigma_*\in\mathsf{image}(\sigma)}$, so the permutation representation is multiplicity free.
\end{proof}

\begin{example}\label{Example7quater}
Now we assume all the hypotheses and notation in Example \ref{Example7} and in Example \ref{Exsymmultfr}. Recall that $X=\{x_0\}\sqcup X'$ with $\lvert V\rvert-\lvert X'\rvert\geq 1$, and that $\widehat{F}_{\alpha(x_0)}=\{\sigma_0\}$. If $\xi\in\mathsf{Per}_{X'}(V)$ and $\sigma\in\Sigma_0$, then:
\begin{itemize} 
\item
whenever $\sigma\prec\xi$ (that is, $\xi\in\mathsf{J}_\sigma$), we have
\[
\mathsf{supp}(\xi\wedge\sigma)=\left\{\left(x_0,\sigma_0\right)\right\}\sqcup\left\{\left(x,\sigma_{\xi^{-1}(x)}\right)\mathrel{:}x\in X'\right\},
\]
$\lvert\xi\wedge\sigma\rvert\left(x_0,\sigma_0\right)=\lvert V\rvert-\lvert X'\rvert$, and $\lvert\xi\wedge\sigma\rvert\left(x,\sigma_*\right)=1$ for every $\left(x,\sigma_*\right)\in\mathsf{supp}(\xi\wedge\sigma)$ with $x\in X'$ (where we slightly abuse notation by writing $\sigma_{\xi^{-1}(x)}$ for $\sigma_v$ when $\xi^{-1}(x)=\{v\}$);
\item
$\mathsf{cosupp}(\sigma)=\left\{\left(x_0,\sigma_0\right)\right\}\sqcup\bigsqcup_{x\in X'}\Bigl(\left\{x\right\}\times\left(\widehat{F}_{\alpha(x)}\cap\mathsf{image}(\sigma)\right)\Bigr)$;
\item
an element $\tau\in\left\lvert\xi^{-1}\wedge\sigma\right\rvert$ is uniquely determined by a function $\mathsf{s}\mathrel{:} X'\rightarrow\mathsf{image}(\sigma)$ satisfying $\mathsf{s}(x)\in\widehat{F}_{\alpha(x)}$ for all $x\in X'$ and
\begin{itemize} 
\item
$\mathsf{supp}(\tau)=\left\{\left(x_0,\sigma_0\right)\right\}\sqcup\{(x,\mathsf{s}(x))\mathrel{:}x\in X'\}$, 
\item
$\lvert\{x\in X'\mathrel{:}\mathsf{s}(x)=\sigma_0\}\rvert=\left\lvert V_{\sigma_0}\cap \xi_0^{-1}(X')\right\rvert$, while
$\lvert\{x\in X'\mathrel{:}\mathsf{s}(x)=\sigma_*\}\rvert=\left\lvert V_{\sigma_*}\right\rvert$ for all $\sigma_*\in\mathsf{image}(\sigma)\setminus\left\{\sigma_0\right\}$, 
\item
$\tau\left(x_0,\sigma_0\right)=\lvert V\rvert-\lvert X'\rvert$ and $\tau(x,\mathsf{s}(x))=1$ for every $x\in X'$.
\end{itemize}
\end{itemize}

Next, we fix $\xi_0\in\mathsf{Per}_{X'}(V)$. Since $\sigma\in\Sigma_0$ implies $\sigma\prec\xi_0$, we have $\xi_0^{-1}(x_0)\subseteq V_{\sigma_0}$ and $\left\lvert\xi_0^{-1}(x_0)\right\rvert=\lvert V\rvert-\lvert X'\rvert\geq 1$, so that $\left\lvert\xi_0^{-1}\wedge\sigma\right\rvert\neq\emptyset$ (as $\lvert\xi_0\wedge\sigma\rvert\in\left\lvert\xi_0^{-1}\wedge\sigma\right\rvert$). Moreover,
\[
\mathsf{J}_\sigma
=\bigsqcup_{\tau\in\left\lvert\xi_0^{-1}\wedge\sigma\right\rvert}\left\{\xi\in\mathsf{Per}_{X'}(V)\mathrel{:}\lvert\xi\wedge\sigma\rvert=\tau\right\}
\]
is the decomposition into the orbits of $\mathsf{I}(\sigma)$. Similarly, if $\xi\in\mathsf{J}_\sigma$, then 
\[
\begin{split}
\mathsf{Stab}_{\mathsf{I}(\sigma)}(\xi)
\;&\cong\;\prod_{\sigma_*\in\mathsf{image}(\sigma)}\left(\prod_{x\mathrel{:}
\left(x,\sigma_*\right)\in\mathsf{supp}(\xi\wedge\sigma)}\mathsf{Sym}_{\lvert\xi\wedge\sigma\rvert\left(x,\sigma_*\right)}\right)\\
&\cong\;\mathsf{Sym}_{\lvert V\rvert-\lvert X'\rvert}\times\left(\mathsf{Sym}_1\right)^{\lvert X'\rvert}.
\end{split}
\]
Furthermore, for every $\sigma_*\in\mathsf{image}(\sigma)\setminus\left\{\sigma_0\right\}$ we have $\tau^{\left(\sigma_*\right)}=\left(1^{\left\lvert V_{\sigma_*}\right\rvert}\right)=(1,1,\dotsc,1)$ ($\left\lvert V_{\sigma_*}\right\rvert$ times), while 
\begin{equation}\label{tauex7quater}
 \tau^{\left(\sigma_0\right)}=\left(\lvert V\rvert-\lvert X'\rvert,1^{\left\lvert\xi_0^{-1}( X')\cap V_{\sigma_0}\right\rvert}\right)
\end{equation}
(note that $\lvert V\rvert-\lvert X'\rvert\geq 1$, so the parts in \eqref{tauex7quater} are indeed in weakly decreasing order). Since these partitions do not depend on the choice of $\tau\in\left\lvert\xi_0^{-1}\wedge\sigma\right\rvert$, the equivalence relation $\approx$ is trivial: all elements of $\left\lvert\xi_0^{-1}\wedge\sigma\right\rvert$ belong to a single equivalence class.

If $\mu\vdash n$, then we denote by $\mathsf{Tab}(\mu)$ the set of all standard tableaux of shape $\mu$; its cardinality equals the dimension of the irreducible representation $S^\mu$ and therefore also its multiplicity in the regular representation of $\mathsf{Sym}_n$. For $\sigma_*\in\mathsf{image}(\sigma)\setminus\left\{\sigma_0\right\}$, the permutation representation associated with $\tau^{\left(\sigma_*\right)}=\left(1^{\left\lvert V_{\sigma_*}\right\rvert}\right)$ coincides with the regular representation of $\mathsf{Sym}\left(V_{\sigma_*}\right)$, and $\tau^{\left(\sigma_*\right)}\trianglelefteq\mu^{\left(\sigma_*\right)}$ holds for every partition $\mu^{\left(\sigma_*\right)}$ of $\left\lvert V_{\sigma_*}\right\rvert$. On the other hand, $\left\lvert\mathsf{STab}\left(\mu^{\left(\sigma_0\right)},\tau^{\left(\sigma_0\right)}\right)\right\rvert\geq 1$ if and only if $\tau^{\left(\sigma_0\right)}\trianglelefteq\mu^{\left(\sigma_0\right)}$. Finally, let $\mathcal{S}(\sigma)$ denote the set of all functions $\mathsf{s}\mathrel{:} X'\rightarrow\mathsf{image}(\sigma)$ satisfying $\mathsf{s}(x)\in\widehat{F}_{\alpha(x)}$ for all $x\in X'$ and the fiber-size conditions above; by the third bullet, $\left\lvert\left\lvert\xi_0^{-1}\wedge\sigma\right\rvert\right\rvert=\lvert\mathcal{S}(\sigma)\rvert\geq 1$. From \eqref{Bigformult2}, we deduce that $\bigotimes_{\sigma_*\in\mathsf{image}(\sigma)}S^{\mu^{\left(\sigma_*\right)}}$ is contained in $L\left(\mathsf{J}_\sigma\right)$ if and only if $\tau^{\left(\sigma_0\right)}\trianglelefteq\mu^{\left(\sigma_0\right)}$, with multiplicity
\[
\lvert\mathcal{S}(\sigma)\rvert\cdot\left\lvert\mathsf{STab}\left(\mu^{\left(\sigma_0\right)},\tau^{\left(\sigma_0\right)}\right)\right\rvert\cdot\prod_{\sigma_*\in\mathsf{image}(\sigma)\setminus\left\{\sigma_0\right\}}\left\lvert\mathsf{Tab}\left(\mu^{\left(\sigma_*\right)}\right)\right\rvert.
\]
In particular, if $\widehat{F}_{\alpha(x)}$ is independent of $x\in X'$ (for instance, when $\alpha$ is constant on $X'$, so that $Y_{\alpha(x)}=Y_a$ for all $x\in X'$), then $\mathsf{image}(\sigma)\subseteq\widehat{F}_a$ and the condition $\mathsf{s}(x)\in\widehat{F}_{\alpha(x)}$ is automatically satisfied; writing $\mathsf{image}(\sigma)=\left\{\sigma_0,\sigma_1,\dotsc,\sigma_d\right\}$, $\lvert\mathcal{S}(\sigma)\rvert$ is then the number of ordered partitions of $X'$ into $d+1$ (possibly empty) subsets of cardinalities $\left\lvert\xi_0^{-1}( X')\cap V_{\sigma_0}\right\rvert,\left\lvert V_{\sigma_1}\right\rvert,\dotsc,\left\lvert V_{\sigma_d}\right\rvert$, and the multiplicity becomes
\begin{equation}\label{abstmultkins}
\binom{\left\lvert X'\right\rvert}{\left\lvert\xi_0^{-1}( X')\cap V_{\sigma_0}\right\rvert,\left\lvert V_{\sigma_1}\right\rvert,\dotsc,\left\lvert V_{\sigma_d}\right\rvert}\cdot\left\lvert\mathsf{STab}\left(\mu^{\left(\sigma_0\right)},\tau^{\left(\sigma_0\right)}\right)\right\rvert\cdot\prod_{\sigma_*\in\mathsf{image}(\sigma)\setminus\left\{\sigma_0\right\}}\left\lvert\mathsf{Tab}\left(\mu^{\left(\sigma_*\right)}\right)\right\rvert.
\end{equation}
\end{example}

\begin{example}\label{Example21quater}
Now we examine the decomposition of the permutation representation in Example \ref{Example21bis} when $G=\mathsf{Sym}(V)$. Assume that $m=\lvert V\rvert\geq 4$ (with obvious modifications for $m=2,3$). We have $\mathsf{Stab}_G\left(v_0\right)=\mathsf{Sym}\left(\left\{v_0\right\}\right)\times\mathsf{Sym}\left(V\setminus\left\{v_0\right\}\right)$, $n=1$, and $\Delta_{v_0,1}=V\setminus\left\{v_0\right\}$; moreover, $G$ is transitive on $\mathsf{Per}_2(V)$, so the homogeneous space is $\mathsf{Per}_2(V)\times Y$. Under the action of $G\cong\mathsf{Sym}_m$, the space $L\left(\mathsf{Per}_2(V)\right)$ is isomorphic to the permutation module $M^{(m-2,1,1)}$, so by Young's rule we have
\[
L\left(\mathsf{J}_{\sigma^0}\right)=L\left(\mathsf{Per}_2(V)\right)=S^{m}\oplus 2S^{m-1,1}\oplus S^{m-2,2}\oplus S^{m-2,1,1}.
\]
Similarly, for $\sigma\neq\sigma^0$, under the action of $\mathsf{Stab}_G\left(v_0\right)\cong\mathsf{Sym}_{m-1}$ we have $L\left(\Delta_{v_0,1}\right)\cong M^{(m-2,1)}$, whence
\[
L\left(\mathsf{J}_\sigma\right)=L\left(\Delta_{v_0,1}\right)=S^{m-1}\oplus S^{m-2,1}. 
\] 
Since $m_{\sigma^0}=1$, in \eqref{DecgenCentip} we therefore obtain $m_\sigma m_\theta=m_\sigma$ if $\sigma\neq\sigma^0$ and $W_\theta\in\left\{S^{m-1},S^{m-2,1}\right\}$, whereas $m_{\sigma^0} m_\theta=1$ if $W_\theta\in\left\{S^m,S^{m-2,2},S^{m-2,1,1}\right\}$, $m_{\sigma^0} m_\theta=2$ if $W_\theta=S^{m-1,1}$, and $m_\sigma m_\theta=0$ in all remaining cases.
\end{example}

\section{Centipedes and $k$-insects in $q$-ary trees}\label{Sectree}

Set $q\coloneqq\lvert V\rvert$;
the finite, rooted $q$-ary tree of height $h$ is the tree with vertex set
\[
\mathcal{T}_{q,h}\coloneqq \emptyset\sqcup\bigsqcup_{j=1}^h V^j
\]
where two vertices $(v_1,\dotsc,v_r)$ and $(w_1,\dotsc,w_s)$ are adjacent if $|r-s| = 1$
and $v_j= w_j$ for $j=1,2,\ldots, \min\{r,s\}$; the root is $\emptyset$.
The first level decomposition of the tree is
\begin{equation}\label{firstlevdec}
\mathcal{T}_{q,h}=\emptyset\sqcup V\sqcup \left[\bigsqcup_{v\in V}\left(\{v\}\times\left(\mathcal{T}_{q,h-1}\setminus\{\emptyset\}\right) \right)\right];
\end{equation}
the corresponding first level decomposition of the automorphism group is
\begin{equation}\label{firstlevautdec}
\Aut\left(\mathcal{T}_{q,h}\right)\cong\Aut\left(\mathcal{T}_{q,h-1}\right)\wr \mathsf{Sym}(V)
\end{equation}
(see, e.g., \cite{book3,OOR,Trees2}).
More explicitly, the action of $\Aut\left(\mathcal{T}_{q,h-1}\right)\wr \mathsf{Sym}(V)$ on $\mathcal{T}_{q,h}$ is given by \eqref{defcompaction}: for $(f,\pi)\in\Aut\left(\mathcal{T}_{q,h-1}\right)\wr \mathsf{Sym}(V)$ we have
\begin{equation}\label{firstlevaction}
(f,\pi)(v,w)= (\pi v, f(\pi v)w),
\end{equation}
for all $v\in V$ and $w\in\mathcal{T}_{q,h-1}$. Since $f(v)\in\Aut\left(\mathcal{T}_{q,h-1}\right)$, \eqref{firstlevaction} defines an automorphism of the tree.

Assume that $h\geq 1$ and $q\geq 2$. A \emph{complete centipede subtree} of $\mathcal{T}_{q,h}$ is a subtree with $2h+1$ vertices $v_0,v_1,\dotsc, v_h$, $w_1,\dotsc,w_h$ such that $v_0,v_1,\dotsc, v_h$ is a \emph{path} from the root $v_0=\emptyset$ to a leaf $v_h$ and $w_j$ is a child of $v_{j-1}$ distinct from $v_j$ (a \emph{pendant}).

\begin{picture}(400,85)
\thicklines
\put(190,75){${\scriptstyle \emptyset=v_0}$}
\put(200,70){\circle*{2}}
\drawline(200,70)(180,50)
\put(180,50){\circle*{2}}
\put(170,53){${\scriptstyle v_1}$}

\drawline(200,70)(220,50)
\put(220,50){\circle*{2}}
\put(222,50){${\scriptstyle w_1}$}

\drawline(180,50)(200,30)
\put(200,30){\circle*{2}}
\put(202,30){${\scriptstyle w_2}$}

\dottedline{3}(180,50)(160,30)
\put(160,30){\circle*{2}}
\put(140,33){${\scriptstyle v_{h-1}}$}

\drawline(160,30)(180,10)
\put(180,10){\circle*{2}}
\put(182,10){${\scriptstyle {w_h}}$}

\drawline(160,30)(140,10)
\put(140,10){\circle*{2}}
\put(130,13){${\scriptstyle v_h}$}

\end{picture}

In the following picture, a complete centipede subtree of $\mathcal{T}_{3,3}$ is shown.

\begin{picture}(400,140)
\put(200,110){\circle*{4}}
\put(188,115){${\scriptstyle v_0=\emptyset}$}
\drawline(200,110)(80,80)
\drawline(200,110)(200,80)
\drawline(200,110)(320,80)
\put(80,80){\circle*{4}}
\put(67,80){${\scriptstyle v_1}$}
\put(200,80){\circle*{2}}
\put(320,80){\circle*{4}}
\put(323,80){${\scriptstyle w_1}$}

\drawline(80,80)(45,50)
\drawline(80,80)(80,50)
\drawline(80,80)(115,50)
\put(45,50){\circle*{2}}
\put(80,50){\circle*{4}}
\put(67,50){${\scriptstyle v_2}$}
\put(115,50){\circle*{4}}
\put(119,50){${\scriptstyle w_2}$}

\drawline(200,80)(165,50)
\drawline(200,80)(200,50)
\drawline(200,80)(235,50)
\put(165,50){\circle*{2}}
\put(200,50){\circle*{2}}
\put(235,50){\circle*{2}}

\drawline(320,80)(285,50)
\drawline(320,80)(320,50)
\drawline(320,80)(355,50)
\put(285,50){\circle*{2}}
\put(320,50){\circle*{2}}
\put(355,50){\circle*{2}}

\drawline(45,50)(35,20)
\drawline(45,50)(45,20)
\drawline(45,50)(55,20)
\put(35,20){\circle*{2}}
\put(45,20){\circle*{2}}
\put(55,20){\circle*{2}}

\drawline(80,50)(70,20)
\drawline(80,50)(80,20)
\drawline(80,50)(90,20)
\put(70,20){\circle*{4}}
\put(65,10){${\scriptstyle v_3}$}
\put(80,20){\circle*{2}}
\put(90,20){\circle*{4}}
\put(90,10){${\scriptstyle w_3}$}

\drawline(115,50)(105,20)
\drawline(115,50)(115,20)
\drawline(115,50)(125,20)
\put(105,20){\circle*{2}}
\put(115,20){\circle*{2}}
\put(125,20){\circle*{2}}

\drawline(165,50)(155,20)
\drawline(165,50)(165,20)
\drawline(165,50)(175,20)
\put(155,20){\circle*{2}}
\put(165,20){\circle*{2}}
\put(175,20){\circle*{2}}

\drawline(200,50)(190,20)
\drawline(200,50)(200,20)
\drawline(200,50)(210,20)
\put(190,20){\circle*{2}}
\put(200,20){\circle*{2}}
\put(210,20){\circle*{2}}

\drawline(235,50)(225,20)
\drawline(235,50)(235,20)
\drawline(235,50)(245,20)
\put(225,20){\circle*{2}}
\put(235,20){\circle*{2}}
\put(245,20){\circle*{2}}

\drawline(285,50)(275,20)
\drawline(285,50)(285,20)
\drawline(285,50)(295,20)
\put(275,20){\circle*{2}}
\put(285,20){\circle*{2}}
\put(295,20){\circle*{2}}

\drawline(320,50)(310,20)
\drawline(320,50)(320,20)
\drawline(320,50)(330,20)
\put(310,20){\circle*{2}}
\put(320,20){\circle*{2}}
\put(330,20){\circle*{2}}

\drawline(355,50)(345,20)
\drawline(355,50)(355,20)
\drawline(355,50)(365,20)
\put(345,20){\circle*{2}}
\put(355,20){\circle*{2}}
\put(365,20){\circle*{2}}

\Thicklines
\drawline(200,110)(80,80)
\drawline(200,110)(320,80)

\drawline(80,80)(80,50)
\drawline(80,80)(115,50)

\drawline(80,50)(70,20)
\drawline(80,50)(90,20)
\end{picture}

Denote by $\mathsf{Cen}_{q,h}$ the set of all complete centipede subtrees of $\mathcal{T}_{q,h}$. Then $\mathsf{Cen}_{q,h}$ is a homogeneous $\Aut\left(\mathcal{T}_{q,h}\right)$-space. Now we describe this space as an application of the construction in Example \ref{Example21} with $G=\mathsf{Sym}(V)$; see also Example \ref{Example21quater}. We take $Y=\mathsf{Cen}_{q,h-1}$ (where $\mathsf{Cen}_{q,0}$ is a singleton), and if $(v_1,w_1,y)\in\mathsf{Per}_2(V)\times\mathsf{Cen}_{q,h-1}$, then we think of $y$ as a centipede in the subtree rooted at $v_1$; namely, the vertices of $y$ are $v_1,\dotsc,v_h$ (the path) and $w_2,\dotsc,w_h$ (the pendants). Then $(v_1,w_1,y)$ may be identified with the centipede in $\mathsf{Cen}_{q,h}$ with vertex set $v_0\coloneqq\emptyset,v_1,\dotsc, v_h$, $w_1,\dotsc,w_h$.

Now we want to describe the homogeneous spaces in this section, in particular $\mathsf{Cen}_{q,h}$, in the framework of \emph{tree compositions}, which we developed in \cite{Trees2}. We recall that a tree composition of $\mathcal{T}_{q,h}$ is a surjective homomorphism $\phi\colon\mathcal{T}_{q,h}\rightarrow\mathcal{P}$, where $\mathcal{P}$ is a finite tree of height $h$, that satisfies the following condition: if $u,v\in\mathcal{T}_{q,h}$ and $\phi(u)=\phi(v)$, then there exists $\tau\in\Aut\left(\mathcal{T}_{q,h}\right)$ such that $\tau(u)=v$ and $\phi(w)=(\phi\circ\tau)(w)$ for every descendant $w$ of $u$. In the setting of the present section, tree compositions that describe the various homogeneous spaces may be obtained by iterating the basic construction at the beginning of Section \ref{Sechomspace}. In particular, for $\mathsf{Cen}_{q,h}$ we construct a tree $\mathcal{P}_{q,h}$ in the following way. The tree has three vertices at the first level, denoted by $x_1,y_1,z_1$ (for $q=2$ only $x_1,y_1$). They correspond to $x_0,x_1,x_2$ in Example \ref{Example21}, respectively. The subtrees rooted at $y_1$ and $z_1$ are paths of length $h-1$, whereas the subtree rooted at $x_1$ is constructed by iteration: if we have constructed the tree below $x_1$ up to level $l-1$ ($2\leq l\leq h$), then at level $l$ we attach three children to $x_{l-1}$, namely $x_l,y_l,z_l$ (only $x_l,y_l$ if $q=2$). Below $y_l$ and $z_l$ there are paths of length $h-l$, whereas below $x_l$ we iterate the construction. Finally, the tree composition $\phi\colon\mathcal{T}_{q,h}\rightarrow\mathcal{P}_{q,h}$ corresponding to a complete centipede in $\mathsf{Cen}_{q,h}$ maps $\emptyset$ to the root of $\mathcal{P}_{q,h}$ and is defined on the remaining vertices by setting
\begin{equation}\label{defphicaterpillar}
\begin{split}
\phi\left(v_l\right)=x_l\quad&\text{for }l=1,2,\dotsc,h;\\
\phi\left(w_l\right)=y_l\quad&\text{for }l=1,2,\dotsc,h;\\
\phi\left(u\right)=z_l\;\quad&\text{for all children }u\text{ of }v_{l-1}\text{ different from }v_l\text{ and }w_l,\\
&\quad\text{where } l=1,2,\dotsc,h,
\end{split}
\end{equation}
and mapping all descendants of $w_l$ (respectively, of each such child $u$) level-by-level onto the unique path of length $h-l$ below $y_l$ (respectively, below $z_l$) when $l<h$.
This is a picture of $\mathcal{P}_{3,3}$.

\begin{picture}(400,140)
\put(200,110){\circle*{2}}
\put(195,115){${\emptyset}$}
\drawline(200,110)(80,80)
\drawline(200,110)(200,80)
\drawline(200,110)(320,80)
\put(80,80){\circle*{2}}
\put(67,80){${\scriptstyle x_1}$}
\put(200,80){\circle*{2}}
\put(203,80){${\scriptstyle y_1}$}
\put(320,80){\circle*{2}}
\put(323,80){${\scriptstyle z_1}$}

\drawline(80,80)(45,50)
\drawline(80,80)(80,50)
\drawline(80,80)(115,50)
\put(45,50){\circle*{2}}
\put(30,50){${\scriptstyle x_2}$}
\put(80,50){\circle*{2}}
\put(67,50){${\scriptstyle y_2}$}
\put(115,50){\circle*{2}}
\put(119,50){${\scriptstyle z_2}$}

\drawline(200,80)(200,50)
\put(200,50){\circle*{2}}

\drawline(320,80)(320,50)
\put(320,50){\circle*{2}}

\drawline(45,50)(35,20)
\drawline(45,50)(45,20)
\drawline(45,50)(55,20)
\put(35,20){\circle*{2}}
\put(28,12){${\scriptstyle x_3}$}
\put(45,20){\circle*{2}}
\put(40,12){${\scriptstyle y_3}$}
\put(55,20){\circle*{2}}
\put(53,12){${\scriptstyle z_3}$}

\drawline(80,50)(80,20)
\put(80,20){\circle*{2}}

\drawline(115,50)(115,20)
\put(115,20){\circle*{2}}
\drawline(200,50)(200,20)
\put(200,20){\circle*{2}}
\drawline(320,50)(320,20)
\put(320,20){\circle*{2}}

\end{picture}

To describe the decomposition of $L\left(\mathsf{Cen}_{q,h}\right)$ into irreducible representations, we use the parametrization of $\widehat{\Aut\left(\mathcal{T}_{q,h}\right)}$ by means of \emph{trees of partitions}, which we introduced in \cite{Trees2}. We just need some very particular cases, which we describe directly. 
For $h=1$ and $q\geq 4$ we have $\mathsf{Cen}_{q,1}=\mathsf{Per}_2(V)$, so that $L\left(\mathsf{Cen}_{q,1}\right)=S^q\oplus 2S^{q-1,1}\oplus S^{q-2,2}\oplus S^{q-2,1,1}$, and we parametrize these four representations by means of the tree $\mathcal{T}_{1,1}$, whose unique vertex at level $1$ is labeled by the corresponding partition.

\begin{picture}(400,55)
\put(55,25){${\scriptstyle S^q}$}
\put(80,40){\circle*{2}}
\put(78,45){${\scriptstyle\emptyset}$}
\thicklines
\drawline(80,40)(80,20)
\put(80,20){\circle*{2}}
\put(75,10){${\scriptstyle (q)}$}

\put(125,25){${\scriptstyle S^{q-1,1}}$}
\put(180,40){\circle*{2}}
\put(178,45){${\scriptstyle\emptyset}$}
\thicklines
\drawline(180,40)(180,20)
\put(180,20){\circle*{2}}
\put(165,10){${\scriptstyle (q-1,1)}$}

\put(230,25){${\scriptstyle S^{q-2,2}}$}
\put(280,40){\circle*{2}}
\put(278,45){${\scriptstyle\emptyset}$}
\thicklines
\drawline(280,40)(280,20)
\put(280,20){\circle*{2}}
\put(265,10){${\scriptstyle (q-2,2)}$}

\put(335,25){${\scriptstyle S^{q-2,1,1}}$}
\put(390,40){\circle*{2}}
\put(388,45){${\scriptstyle\emptyset}$}
\thicklines
\drawline(390,40)(390,20)
\put(390,20){\circle*{2}}
\put(375,10){${\scriptstyle (q-2,1,1)}$}

\end{picture}

\noindent
In order to avoid too many particular cases, from now on we assume that $q\geq 4$.  

Next, we define inductively the required trees of partitions in the following way. There is a non-negative integer $l\leq h$, which we call the \emph{depth} of the tree of partitions. For $l\geq 1$, there are two families of special vertices, denoted by $a_0\coloneqq\emptyset,a_1,a_2,\dotsc, a_{l-1}$ and $b_1,b_2,\dotsc,b_l$, with the following properties.
\begin{itemize}
\item
The trivial representation $S((q))$ is represented separately by the tree of depth $l=0$, which reduces to a single path with all vertices at levels $1,2,\dotsc,h$ labeled by $(q)$.
\item
For $l=1$, the tree reduces to a single path whose first-level vertex $b_1$ is labeled by one of the three non-trivial partitions $\lambda^{(1)}\in\{(q-1,1),(q-2,2),(q-2,1,1)\}$, while all vertices below $b_1$ (when $h\geq 2$) are labeled by $(q)$.
\item 
For $j=1,2,\dotsc,l-1$, the vertex $b_j$ is labeled by a partition $\lambda^{(j)}\in\{(q-1),(q-2,1)\}$, and the vertices below it form a path labeled by the trivial partition $(q)$.
\item 
For $j=0,1,\dotsc,l-2$, the vertex $a_j$ has two children, which are precisely $a_{j+1}$ and $b_{j+1}$, and it is labeled by the trivial partition $(1)$.
\item
The vertices below $a_{l-1}$ form a path: its unique child is $b_l$, the vertex $b_l$ is labeled by a partition $\lambda^{(l)}\in\{(q-1,1),(q-2,2),(q-2,1,1)\}$, and all vertices below $b_l$ (when $l<h$) are labeled by $(q)$.
\end{itemize}
The following picture illustrates the procedure for $l\geq 2$. 

\begin{equation}\label{Caterpillar}
\begin{picture}(400,150)
\thicklines
\put(200,145){${\scriptstyle a_0=\emptyset}$}
\put(200,140){\circle*{2}}
\drawline(200,140)(180,120)
\put(180,120){\circle*{2}}
\put(170,123){${\scriptstyle (1)}$}
\put(184,118){${\scriptstyle a_1}$}

\drawline(200,140)(220,120)
\put(220,120){\circle*{2}}
\put(222,120){${\scriptstyle \lambda^{(1)}}$}
\put(212,116){${\scriptstyle b_1}$}
\drawline(220,120)(240,100)
\put(240,100){\circle*{2}}
\put(242,100){${\scriptstyle (q)}$}
\dottedline{3}(240,100)(320,20)
\put(320,20){\circle*{2}}
\put(320,23){${\scriptstyle (q)}$}

\drawline(180,120)(200,100)
\put(200,100){\circle*{2}}
\put(202,100){${\scriptstyle \lambda^{(2)}}$}
\put(192,96){${\scriptstyle b_2}$}
\drawline(200,100)(220,80)
\put(220,80){\circle*{2}}
\put(222,80){${\scriptstyle (q)}$}
\dottedline{3}(220,80)(280,20)
\put(280,20){\circle*{2}}
\put(280,23){${\scriptstyle (q)}$}

\dottedline{3}(180,120)(160,100)
\put(160,100){\circle*{2}}
\put(150,103){${\scriptstyle (1)}$}
\put(164,98){${\scriptstyle a_{l-2}}$}

\drawline(160,100)(180,80)
\put(180,80){\circle*{2}}
\put(182,80){${\scriptstyle \lambda^{(l-1)}}$}
\put(164,76){${\scriptstyle b_{l-1}}$}
\drawline(180,80)(200,60)
\put(200,60){\circle*{2}}
\put(202,60){${\scriptstyle (q)}$}
\dottedline{3}(200,60)(240,20)
\put(240,20){\circle*{2}}
\put(240,23){${\scriptstyle (q)}$}

\drawline(160,100)(140,80)
\put(140,80){\circle*{2}}
\put(144,78){${\scriptstyle a_{l-1}}$}
\put(130,83){${\scriptstyle (1)}$}

\drawline(140,80)(120,60)
\put(120,60){\circle*{2}}
\put(110,62){${\scriptstyle \lambda^{(l)}}$}
\put(124,56){${\scriptstyle b_l}$}
\put(60,60){${\scriptstyle  \text{\rm level }l\;\rightarrow}$}

\drawline(120,60)(100,40)
\put(100,40){\circle*{2}}
\put(90,43){${\scriptstyle (q)}$}

\dottedline{3}(100,40)(80,20)
\put(80,20){\circle*{2}}
\put(70,23){${\scriptstyle (q)}$}
\put(20,20){${\scriptstyle  \text{\rm level }h\;\rightarrow}$}
\end{picture}
\end{equation}

For $l\geq 2$, the tree in \eqref{Caterpillar} corresponds to a representation in the first summand of \eqref{DecgenCentip}. Namely, the tree below $a_1$ represents the \emph{non-trivial} irreducible representation $\sigma_{v_0}$ of $\Aut\left(\mathcal{T}_{q,h-1}\right)$, whereas $b_1$ is labeled by the partition representing $\theta$, that is, $(q-1)$ or $(q-2,1)$; cf.\ Example \ref{Example21quater}. For $l=1$, it reduces to a single path whose first-level vertex is labeled by one of the three non-trivial partitions $(q-1,1)$, $(q-2,2)$, $(q-2,1,1)$, whereas its descendants (when $h\geq 2$) are labeled by $(q)$; it corresponds to a representation in the second summand of \eqref{DecgenCentip}. We also have the case $l=0$, in which the corresponding tree reduces to a single path with all vertices below the root labeled by $(q)$. For $l\geq 1$, we denote by $S\left(\lambda^{(1)},\dotsc,\lambda^{(l)}\right)$ the corresponding irreducible representation, and by $S((q))$ the trivial representation ($l=0$). 
Clearly, $\dim S((q))=1$, whereas for $l\geq 1$
\[
\dim S\left(\lambda^{(1)},\dotsc,\lambda^{(l)}\right)=q^{l-1}\prod_{j=1}^l\dim S^{\lambda^{(j)}}.
\]
By induction on $h$, from Example \ref{Example21quater} we deduce that $S((q))$ has multiplicity $1$, while for $l\geq 1$ the multiplicity of $S\left(\lambda^{(1)},\dotsc,\lambda^{(l)}\right)$ in $L\left(\mathsf{Cen}_{q,h}\right)$ is $2$ if $\lambda^{(l)}=(q-1,1)$, and $1$ otherwise. In conclusion, the decomposition is
\[
\begin{split}
L\left(\mathsf{Cen}_{q,h}\right)=&S((q))\bigoplus\left[\bigoplus_{l=1}^h\bigoplus_{\substack{\lambda^{(1)},\dotsc,\lambda^{(l-1)}\in\\
\{(q-1),(q-2,1)\}}}\; \bigoplus_{\substack{\lambda^{(l)}\in\\
\{(q-2,2),(q-2,1,1)\}}}S\left(\lambda^{(1)},\dotsc,\lambda^{(l)}\right)\right]\\
&\bigoplus\left[\bigoplus_{l=1}^h\bigoplus_{\substack{\lambda^{(1)},\dotsc,\lambda^{(l-1)}\in\\
\{(q-1),(q-2,1)\}}}2S\left(\lambda^{(1)},\dotsc,\lambda^{(l-1)},(q-1,1)\right)\right].
\end{split}
\]

\medskip

For $1\leq k\leq q$, a $k$-\emph{insect} of type $1\leq l_1\leq\dotsb\leq l_k\leq h$ in $\mathcal{T}_{q,h}$ is an ordered sequence $\left(w_1,\dotsc, w_k\right)$ of $k$ distinct vertices of $\mathcal{T}_{q,h}$ such that
\begin{itemize}
\item the ancestors $v_1,\dotsc, v_k$ at the first level of the vertices $w_1,\dotsc, w_k$ are $k$ distinct vertices;
\item $w_j$ is a vertex at level $l_j$, $j=1,2,\dotsc,k$.
\end{itemize}
In the following picture, a $3$-insect of type $(1,2,3)$ in $\mathcal{T}_{3,3}$ is shown.\\

\begin{picture}(400,130)
\put(200,120){\circle*{2}}
\put(198,125){${\scriptstyle\emptyset}$}
\drawline(200,120)(80,90)
\drawline(200,120)(200,90)
\drawline(200,120)(320,90)
\put(80,90){\circle*{4}}
\put(50,90){${\scriptstyle w_1=v_1}$}
\put(200,90){\circle*{2}}
\put(185,90){${\scriptstyle v_2}$}
\put(320,90){\circle*{2}}
\put(325,90){${\scriptstyle v_3}$}

\drawline(80,90)(45,60)
\drawline(80,90)(80,60)
\drawline(80,90)(115,60)
\put(45,60){\circle*{2}}
\put(80,60){\circle*{2}}
\put(115,60){\circle*{2}}

\drawline(200,90)(165,60)
\drawline(200,90)(200,60)
\drawline(200,90)(235,60)
\put(165,60){\circle*{4}}
\put(150,60){${\scriptstyle w_2}$}
\put(200,60){\circle*{2}}
\put(235,60){\circle*{2}}

\drawline(320,90)(285,60)
\drawline(320,90)(320,60)
\drawline(320,90)(355,60)
\put(285,60){\circle*{2}}
\put(320,60){\circle*{2}}
\put(355,60){\circle*{2}}

\drawline(45,60)(35,30)
\drawline(45,60)(45,30)
\drawline(45,60)(55,30)
\put(35,30){\circle*{2}}
\put(45,30){\circle*{2}}
\put(55,30){\circle*{2}}

\drawline(80,60)(70,30)
\drawline(80,60)(80,30)
\drawline(80,60)(90,30)
\put(70,30){\circle*{2}}
\put(80,30){\circle*{2}}
\put(90,30){\circle*{2}}

\drawline(115,60)(105,30)
\drawline(115,60)(115,30)
\drawline(115,60)(125,30)
\put(105,30){\circle*{2}}
\put(115,30){\circle*{2}}
\put(125,30){\circle*{2}}

\drawline(165,60)(155,30)
\drawline(165,60)(165,30)
\drawline(165,60)(175,30)
\put(155,30){\circle*{2}}
\put(165,30){\circle*{2}}
\put(175,30){\circle*{2}}

\drawline(200,60)(190,30)
\drawline(200,60)(200,30)
\drawline(200,60)(210,30)
\put(190,30){\circle*{2}}
\put(200,30){\circle*{2}}
\put(210,30){\circle*{2}}

\drawline(235,60)(225,30)
\drawline(235,60)(235,30)
\drawline(235,60)(245,30)
\put(225,30){\circle*{2}}
\put(235,30){\circle*{2}}
\put(245,30){\circle*{2}}

\drawline(285,60)(275,30)
\drawline(285,60)(285,30)
\drawline(285,60)(295,30)
\put(275,30){\circle*{2}}
\put(285,30){\circle*{4}}
\put(282,20){${\scriptstyle w_3}$}
\put(295,30){\circle*{2}}

\drawline(320,60)(310,30)
\drawline(320,60)(320,30)
\drawline(320,60)(330,30)
\put(310,30){\circle*{2}}
\put(320,30){\circle*{2}}
\put(330,30){\circle*{2}}

\drawline(355,60)(345,30)
\drawline(355,60)(355,30)
\drawline(355,60)(365,30)
\put(345,30){\circle*{2}}
\put(355,30){\circle*{2}}
\put(365,30){\circle*{2}}

\Thicklines
\drawline(200,120)(80,90)
\drawline(200,120)(200,90)
\drawline(200,120)(320,90)

\drawline(200,90)(165,60)

\drawline(320,90)(285,60)
\drawline(285,60)(285,30)

\end{picture}

We denote by $\mathsf{Ins}\left(l_1,\dotsc,l_k;h\right)$ the set of all $k$-insects of type $l_1,\dotsc,l_k$ in $\mathcal{T}_{q,h}$. 

Now we describe the $k$-insects by means of tree compositions, as in \eqref{defphicaterpillar}; in order to avoid special cases, in this construction we assume that $k<q$. The first level of the quotient tree $\mathcal{P}_{q,h;l_1,\dotsc,l_k}$ is precisely the set $X$ in Example \ref{Example7}, and we denote the first-level vertices by $\{0,(1,1),(1,2),\dotsc,(1,k)\}$, so that $X'=\{(1,1),(1,2),\dotsc,(1,k)\}$. The subtree rooted at $0$ is a path of length $h-1$. Then, for $j=1,2,\dotsc,k$, in the subtree rooted at $(1,j)$:
\begin{itemize}
\item
all non-leaf vertices have one child, except for those denoted by $(i,j)$ for $i=1,2,\dotsc,l_j-1$, which have two children;
\item 
for $i=1,2,\dotsc,l_j-1$, $(i,j)$ is at level $i$ and $(i+1,j)$ is one of its children.
\end{itemize}
In other words, the vertices $(1,j),(2,j),\dotsc,(l_j,j)$ form a \emph{principal path} (extended below $(l_j,j)$ to level $h$), and, down to $(l_j-1,j)$, a \emph{secondary path} is attached to each $(i,j)$. Then the tree composition $\phi\colon\mathcal{T}_{q,h}\rightarrow\mathcal{P}_{q,h;l_1,\dotsc,l_k}$ corresponding to a $k$-insect $\left(w_1,\dotsc,w_k\right)\in\mathsf{Ins}\left(l_1,\dotsc,l_k;h\right)$ is defined as follows: $\phi(\emptyset)$ is the root of $\mathcal{P}_{q,h;l_1,\dotsc,l_k}$; all non-root vertices not in the subtrees rooted at $v_1,v_2,\dotsc,v_k$ are mapped level-by-level onto the path rooted at $0$; for $j=1,2,\dotsc,k$, the path from $v_j$ to $w_j$ is mapped onto $(1,j),(2,j),\dotsc,(l_j,j)$ (with $\phi\left(v_j\right)=(1,j)$ and $\phi\left(w_j\right)=\left(l_j,j\right)$), while all remaining vertices below $v_j$ are mapped level-by-level onto the corresponding secondary paths. This is a picture of $\mathcal{P}_{3,3;1,2,3}$. 

\begin{picture}(400,130)
\put(200,120){\circle*{2}}
\put(198,125){${\scriptstyle\emptyset}$}
\drawline(200,120)(80,90)
\drawline(200,120)(200,90)
\drawline(200,120)(320,90)
\put(80,90){\circle*{4}}
\put(55,90){${\scriptstyle (1,1)}$}
\put(200,90){\circle*{2}}
\put(180,90){${\scriptstyle (1,2)}$}
\put(320,90){\circle*{2}}
\put(325,90){${\scriptstyle (1,3)}$}

\drawline(80,90)(80,60)
\put(80,60){\circle*{2}}

\drawline(200,90)(165,60)
\drawline(200,90)(235,60)
\put(165,60){\circle*{4}}
\put(145,60){${\scriptstyle (2,2)}$}
\put(235,60){\circle*{2}}

\drawline(320,90)(285,60)
\drawline(320,90)(355,60)
\put(285,60){\circle*{2}}
\put(355,60){\circle*{2}}


\drawline(80,60)(80,30)
\put(80,30){\circle*{2}}

\drawline(165,60)(165,30)
\put(165,30){\circle*{2}}

\drawline(235,60)(235,30)
\put(235,30){\circle*{2}}

\drawline(285,60)(285,30)
\put(285,30){\circle*{4}}
\put(282,20){${\scriptstyle (3,3)}$}

\drawline(355,60)(355,30)
\put(355,30){\circle*{2}}

\Thicklines
\drawline(200,120)(80,90)
\drawline(200,120)(200,90)
\drawline(200,120)(320,90)

\drawline(200,90)(165,60)

\drawline(320,90)(285,60)
\drawline(285,60)(285,30)

\end{picture}

As above, we assume that $q\geq 4$ and decompose $L\bigl(\mathsf{Ins}\left(l_1,\dotsc,l_k;h\right)\bigr)$ into irreducible $\Aut\left(\mathcal{T}_{q,h}\right)$-representations.
First of all, we recall the decomposition in the case $k=1$ (see, e.g., \cite[Theorem 4.2]{CST2} and \cite[Theorem 3.3.1]{book3}; the original source is \cite{Figa}). We show how to obtain it by means of a repeated application of \eqref{deccompact}, taking also \eqref{firstlevautdec} into account. We denote by $W_{0,h}$ the trivial representation of $\Aut\left(\mathcal{T}_{q,h}\right)$, we set $W_{1,h}\coloneqq W_{0,h-1}\otimes S^{q-1,1}$ (cf.\ the second summand of \eqref{deccompact}) and furthermore $W_{l,h}\coloneqq W_{l-1,h-1}\otimes L(V)$ (cf.\ the first summand of \eqref{deccompact}), for $l=2,3,\dotsc,h$. For $l_1=1$ and every $h\geq 1$, we have $L\bigl(\mathsf{Ins}\left(1;h\right)\bigr)\equiv L(V)=W_{0,h}\oplus W_{1,h}$. For $2\leq l_1\leq h$, we use the fact that $\mathsf{Ins}\left(l_1;h\right)\simeq V\times\mathsf{Ins}\left(l_1-1;h-1\right)$, with the composition action in Example \ref{Example2}; if $L\bigl(\mathsf{Ins}\left(l_1-1;h-1\right)\bigr)=\bigoplus_{l=0}^{l_1-1}W_{l,h-1}$ with respect to the action of $\Aut\left(\mathcal{T}_{q,h-1}\right)$, then by induction
\begin{equation}\label{DecLSl1h}
\begin{split}
L\bigl(\mathsf{Ins}\left(l_1;h\right)\bigr)=&\left[\bigoplus_{l'=1}^{l_1-1}\left(W_{l',h-1}\otimes L(V)\right)\right]\bigoplus\left[\left(W_{0,h-1}\otimes S^{q}\right)\oplus\left(W_{0,h-1}\otimes S^{q-1,1}\right)\right]\\
=&\bigoplus_{l=0}^{l_1}W_{l,h}.
\end{split}
\end{equation}
Note that, for $l=1,2,\dotsc,h$, the tree representing $W_{l,h}$ is \eqref{Caterpillar} with $\lambda^{(j)}=(q-1)$ for $j=1,\dotsc,l-1$ and $\lambda^{(l)}=(q-1,1)$; we denote it by $\mathcal{Q}(l,h)$. In other words, with the notation of the centipede above, 
\begin{equation}\label{Qlhcase}
W_{l,h}=S\left((q-1),\dotsc,(q-1),(q-1,1)\right),
\end{equation}
where $(q-1)$ is repeated $l-1$ times.

In the general case ($1\leq k<q$), we are in the setting of Example \ref{Example7quater}. First of all, note that there exist $1\leq m\leq k$ and indices $1\leq i_1<i_2<\dotsb<i_m= k$ such that 
\[
l_1=l_2=\dotsb= l_{i_1}<l_{i_1+1}=l_{i_1+2}=\dotsb=l_{i_2}<\dotsb<l_{i_{m-1}+1}=l_{i_{m-1}+2}=\dotsb=l_{i_m}=l_k.
\]
We also set $i_0=0$. Then we can suppose that $X=\{0,1,\dotsc,k\}$, $x_0=0$, $A=\left\{0,1,\dotsc, l_k-1\right\}$, we can set $\alpha(0)=0$ and $\alpha(i)=l_{i_j}-1$ if $i_{j-1}<i\leq i_j$ for $i\in X'$, so that $Y_{l_{i_j}-1}=\mathsf{Ins}\left(l_{i_j}-1;h-1\right)$ for $j=1,2,\dotsc,m$ (where $\mathsf{Ins}(0;h-1)$ is a singleton).
We fix a $k$-tuple $\left(v_1,\dotsc,v_k\right)$ of distinct elements of $V$ and define $\xi_0$ by setting $\xi_0\left(v_i\right)\coloneqq i$ for $i=1,2,\dotsc,k$, and $\xi_0(v)\coloneqq 0$ for $v\neq v_1,\dotsc, v_k$. By \eqref{DecLSl1h}, $\Sigma_0$ is parametrized by the set 
\[
\widetilde{\Sigma}_0\coloneqq\left\{\left(l'_1,\dotsc,l_k'\right)\mathrel{:} l'_1\leq\dotsb\leq l_k'\text{ and }0\leq l'_i\leq l_i-1,\; i=1,2,\dotsc,k\right\}.
\]
Furthermore, for every $\left(l'_1,\dotsc,l_k'\right)\in\widetilde{\Sigma}_0$ we define $\sigma_{l'_1,\dotsc,l'_k}$ by setting
$U_{\sigma_{l'_1,\dotsc,l'_k}}\coloneqq\bigotimes_{v\in V}U_{\sigma_v}$, where $U_{\sigma_{v_i}}=W_{l'_i,h-1}$ for $i=1,2,\dotsc,k$ and $U_{\sigma_v}=W_{0,h-1}$ for $v\neq v_1,\dotsc, v_k$. Hence, we may take $\Sigma_0\coloneqq\bigl\{\sigma_{l'_1,\dotsc,l'_k}\mathrel{:}\left(l'_1,\dotsc,l_k'\right)\in\widetilde{\Sigma}_0\bigr\}$. Moreover, there exist $0\leq m'\leq k$ and indices $0\leq i'_0< i'_1<i'_2<\dotsb<i'_{m'}= k$ such that 
\[
0\eqcolon l'_0=\dotsb= l'_{i'_0}<l'_{i'_0+1}=\dotsb= l'_{i'_1}<l'_{i'_1+1}=\dotsb=l'_{i'_2}<\dotsb<l'_{i'_{m'-1}+1}=\dotsb=l'_{i'_{m'}}=l'_k 
\]
(where $m'=0$ and $i'_0=k$ if $l'_k=0$, while $i'_0=0$ if $l'_1>0$), so that, with a slight modification of the notation in \eqref{VsIsigma}, we set $V_s\coloneqq\left\{v_{i'_{s-1}+1},\dotsc,v_{i'_s}\right\}$ for $s=1,2,\dotsc, m'$, whereas $V_0\coloneqq V\setminus\left\{v_{i'_0+1},\dotsc,v_k\right\}$ and $\tau^{(\sigma_0)}=\left(q-k, 1^{i'_0}\right)$; cf.\ \eqref{tauex7quater}. Therefore,
\[
\mathsf{I}(\sigma)=\prod_{s=0}^{m'}\mathsf{Sym}\left(V_s\right).
\]
Now let $\sigma_0\coloneqq W_{0,h-1}$ and $\sigma_s\coloneqq W_{l'_{i'_s},h-1}$ for $s=1,2,\dotsc,m'$. For each $s\in\{1,2,\dotsc,m'\}$, we have $\sigma_s\in\widehat{F}_{\alpha(i)}$ ($i\in X'$) if and only if $l'_{i'_s}\leq l_i-1$; setting $r_s\coloneqq\left\lvert\left\{i\in\{1,2,\dotsc,k\}\mathrel{:} l_i-1<l'_{i'_s}\right\}\right\rvert$, this condition is equivalent to $i\in\{r_s+1,\dotsc,k\}$ (note that $0\leq r_1\leq\dotsb\leq r_{m'}\leq k-1$ and $r_s\leq i'_{s-1}$). Counting the valid maps $\mathsf{s}\mathrel{:} X'\rightarrow\mathsf{image}(\sigma)$ in Example \ref{Example7quater} sequentially from $s=m'$ down to $s=1$ gives $\lvert\mathcal{S}(\sigma)\rvert=\prod_{s=1}^{m'}\binom{i'_s-r_s}{i'_s-i'_{s-1}}$ (which reduces to the multinomial coefficient in \eqref{abstmultkins} when $l_1=\dotsb=l_k$, since $r_1=\dotsb=r_{m'}=0$). Consequently, the irreducible $\Aut\left(\mathcal{T}_{q,h}\right)$-representation
\begin{equation}\label{Wl1mu0}
W_{l'_1,\dotsc,l'_k;\mu^{(0)},\dotsc,\mu^{(m')}}\coloneqq\Ind^{\text{\rm Aut}(\mathcal{T}_{q,h})}_{\text{\rm Aut}(\mathcal{T}_{q,h-1})\wr\mathsf{I}(\sigma)}\left[\widetilde{U}_{\sigma_{l'_1,\dotsc,l'_k}}\otimes\left(\overline{\bigotimes}_{s=0}^{\;m'}S^{\mu^{(s)}}\right)\right]
\end{equation}
is contained in $L\bigl(\mathsf{Ins}\left(l_1,\dotsc,l_k;h\right)\bigr)$ if and only if $\mu^{(s)}\vdash i'_s-i'_{s-1}$ for $s=1,2,\dotsc,m'$, and $\mu^{(0)}\trianglerighteq\tau^{(\sigma_0)}$, with multiplicity equal to
\begin{equation}\label{multinsgen}
\left[\prod_{s=1}^{m'}\binom{i'_s-r_s}{i'_s-i'_{s-1}}\right]\cdot\left\lvert\mathsf{STab}\left(\mu^{(0)},\tau^{(\sigma_0)}\right)\right\rvert\cdot\prod_{s=1}^{m'}\left\lvert\mathsf{Tab}\left(\mu^{(s)}\right)\right\rvert.
\end{equation}
With the notation in \cite{Trees2}, the irreducible representation in \eqref{Wl1mu0} is represented by the following tree of partitions:

\begin{equation}\label{TreeWW}
\begin{picture}(450,85)
\put(200,70){\circle*{5}}
\put(199,75){${\scriptstyle \emptyset}$}

\thicklines
\drawline(200,70)(100,30)
\put(100,30){\circle*{3}}
\put(84,30){${\scriptstyle \mu^{(0)}}$}

\drawline(200,70)(160,30)
\put(160,30){\circle*{3}}
\put(144,30){${\scriptstyle \mu^{(1)}}$}

\put(215,30){\circle*{1}}
\put(225,30){\circle*{1}}
\put(235,30){\circle*{1}}

\drawline(200,70)(300,30)
\put(300,30){\circle*{3}}
\put(305,30){${\scriptstyle \mu^{(m')}}$}

\thinlines
\put(100,30){\line(-1,-1){15}}
\put(100,30){\line(1,-1){15}}
\put(82,5){${\scriptstyle \mathcal{Q}(0,h-1)}$}

\put(160,30){\line(-1,-1){15}}
\put(160,30){\line(1,-1){15}}
\put(140,5){${\scriptstyle \mathcal{Q}\bigl(l'_1,h-1\bigr)}$}

\put(300,30){\line(-1,-1){15}}
\put(300,30){\line(1,-1){15}}
\put(280,5){${\scriptstyle \mathcal{Q}\bigl(l'_k,h-1\bigr)}$}
\end{picture}
\end{equation}

\begin{remark}\label{Remlprime01}
We want to express the representations in \eqref{DecLSl1h} by means of the notation in \eqref{Wl1mu0}; that is, we examine the representations in the case $k=1$, previously considered, from the point of view of the general construction.
If $l'_1=\dotsb=l'_k=0$, then $m'=0$, $i'_0=k$, and $\mathsf{I}(\sigma)=\mathsf{Sym}(V)$, so that if we choose $\mu^{(0)}=(q)$, then \eqref{Wl1mu0} coincides with $W_{0,h}$, the trivial representation of $\Aut\left(\mathcal{T}_{q,h}\right)$ (with multiplicity $1$). If $l'_1=\dotsb=l'_{k-1}=0 < l'_k$ (for $k\geq 2$; for $k=1$ we simply have $0<l'_1$), then $m'=1$, $i'_0=k-1$, $i'_1=k$, $\mu^{(1)}=(1)$, and if we choose $\mu^{(0)}=(q-1)$, then \eqref{Wl1mu0} coincides with $W_{l'_k + 1, h}$; compare
\eqref{Caterpillar} in the $\mathcal{Q}(l,h)$-case (cf.\ \eqref{Qlhcase}) with \eqref{TreeWW} in the present case. Moreover, setting $r_1\coloneqq\left\lvert\left\{i\in\{1,\dotsc,k\}\mathrel{:} l_i-1<l'_k\right\}\right\rvert$, the multiplicity of $W_{l'_k+1,h}$ in $L\bigl(\mathsf{Ins}\left(l_1,\dotsc,l_k;h\right)\bigr)$, computed by means of \eqref{multinsgen}, is equal to
\[
\binom{k-r_1}{1}\left\lvert\mathsf{STab}\bigl((q-1),\left(q-k,1^{k-1}\right)\bigr)\right\rvert=k-r_1
\]
(which equals $k$ whenever $l'_k\leq l_1-1$, and in particular when $l_1=\dotsb=l_k$).
\end{remark}

\begin{example}
Now we examine the case $k=2$ (with $1\leq l_1\leq l_2\leq h$ and $q\geq 4$).
\begin{itemize}
\item If $l_1'=l_2'=0$, then $m'=0$, $i'_0=2$, and $\tau^{(\sigma_0)}=(q-2,1,1)$. Therefore for $\mu^{(0)}$ we have four possibilities, namely $\mu^{(0)}=(q)$, $\mu^{(0)}=(q-1,1)$, $\mu^{(0)}=(q-2,2)$, and $\mu^{(0)}=(q-2,1,1)$. By \eqref{multinsgen}, the multiplicity is $2$ in the second case ($\mu^{(0)}=(q-1,1)$) and $1$ in the remaining three cases (cf.\ Example \ref{Example21quater}). The corresponding irreducible representations are of the form $W_{0,0;\mu^{(0)}}$.

\item If $l_1'=0<l_2'$ (so $1\leq l'_2\leq l_2-1$), then $m'=1$, $i'_0=1$, $i'_1=2$, and $\tau^{(\sigma_0)}=(q-2,1)$. Therefore for $\mu^{(0)}$ we have two possibilities, namely $\mu^{(0)}=(q-1)$ and $\mu^{(0)}=(q-2,1)$, whereas $\mu^{(1)}=(1)$. For each $l'_2$, we obtain two irreducible representations, namely $W_{0,l'_2;(q-1),(1)}$ and $W_{0,l'_2;(q-2,1),(1)}$; by \eqref{multinsgen}, both have multiplicity $2$ if $l'_2\leq l_1-1$ (where $r_1=0$), and multiplicity $1$ if $l_1\leq l'_2\leq l_2-1$ (where $r_1=1$).

\item If $0<l_1'=l_2'$ (so $1\leq l'_1=l'_2\leq l_1-1$), then $m'=1$, $i'_0=0$, $i'_1=2$, and $r_1=0$. Moreover, since $\tau^{(\sigma_0)}=(q-2)$, we have $\mu^{(0)}=(q-2)$, whereas for $\mu^{(1)}$ we have two possibilities, namely $\mu^{(1)}=(2)$ or $\mu^{(1)}=(1,1)$. In both cases, the multiplicity, computed by means of \eqref{multinsgen}, is equal to $1$.

\item If $0<l_1'<l_2'$ (with $l'_1\leq l_1-1$ and $l'_2\leq l_2-1$), then $m'=2$, $i'_0=0$, $i'_1=1$, $i'_2=2$, $\mu^{(0)}=(q-2)$, and $\mu^{(1)}=\mu^{(2)}=(1)$. By \eqref{multinsgen}, the multiplicity is equal to $2$ if $l'_2\leq l_1-1$ (where $r_1=r_2=0$), and equal to $1$ if $l_1\leq l'_2\leq l_2-1$ (where $r_1=0$ and $r_2=1$).
\end{itemize}
In conclusion, we obtain
\[
\begin{split}
L\bigl(\mathsf{Ins}\left(l_1,l_2;h\right)\bigr)=&W_{0,0;(q)}\oplus 2W_{0,0;(q-1,1)}\oplus W_{0,0;(q-2,2)}\oplus W_{0,0;(q-2,1,1)}\\
&\bigoplus \left(\bigoplus_{l'_2=1}^{l_1-1}2W_{0,l'_2;(q-1),(1)}\right)\bigoplus\left(\bigoplus_{l'_2=l_1}^{l_2-1}W_{0,l'_2;(q-1),(1)}\right)\\
&\bigoplus \left(\bigoplus_{l'_2=1}^{l_1-1}2W_{0,l'_2;(q-2,1),(1)}\right)\bigoplus\left(\bigoplus_{l'_2=l_1}^{l_2-1}W_{0,l'_2;(q-2,1),(1)}\right)\\
&\bigoplus\left(\bigoplus_{0<l_1'=l_2'\leq l_1-1}W_{l'_1,l'_1;(q-2),(2)}\right)\bigoplus\left(\bigoplus_{0<l_1'=l_2'\leq l_1-1}W_{l'_1,l'_1;(q-2),(1,1)}\right)\\
&\bigoplus\left(\bigoplus_{0<l_1'<l_2'\leq l_1-1}2W_{l'_1,l'_2;(q-2),(1),(1)}\right)\bigoplus\left(\bigoplus_{\substack{0<l_1'\leq l_1-1\\ l_1\leq l'_2\leq l_2-1}}W_{l'_1,l'_2;(q-2),(1),(1)}\right).
\end{split}
\]
\end{example}

\section{A Gelfand pair with the generalized symmetric group $\mathsf{C}_n\wr\mathsf{Sym}(V)$}\label{Secgensymm}

In this section, inspired by Theorem \ref{TheoSymmultfree}, we construct a family of multiplicity-free permutation representations for the generalized symmetric group $\mathsf{C}_n\wr\mathsf{Sym}(V)$; for related results on Gelfand pairs and spherical functions for wreath products with cyclic groups, see also \cite{Miz2}.
Assuming all the notation in Example \ref{Exampleequiv}, we set 
\[
\omega_n\coloneqq\exp\left(\frac{2\pi i}{n}\right)\qquad\text{and}\qquad\chi_{n,t}(m)\coloneqq\omega_n^{tm}\quad\text{for }t,m\in\mathbb{Z},
\]
so that the latter are the characters of $\mathsf{C}_n$ realized as functions on $\mathbb{Z}_n$. We also regard them as $n$-periodic functions defined on $\mathbb{Z}$. 
Then we set $W_{n,t}\coloneqq\left\langle\chi_{n,t}\right\rangle$, so that $L\left(\mathbb{Z}_n\right)=\bigoplus_{t\in\mathbb{Z}_n}W_{n,t}$ is the decomposition of the left regular representation of $\mathsf{C}_n$. 

More generally, if $b\in\mathsf{D}(n)$, then we consider the action of $\mathsf{C}_n$ on $\mathbb{Z}_b$ given by $\epsilon^k(m)\coloneqq m+k\pmod b$ for all $m\in\mathbb{Z}_b$ and $k\in\mathbb{Z}$. With respect to this action, the stabilizer of $0_{\mathbb{Z}_b}$ is $\mathsf{C}_a\coloneqq\langle \epsilon^b\rangle$, where $a\coloneqq \frac{n}{b}$. Furthermore, assuming $n>1$ and that $a$ and $b$ as above are \emph{coprime}, that is, $\gcd(a,b)=1$ (in particular, $a\neq b$), we set $\mathsf{C}_b\coloneqq\langle \epsilon^a\rangle$, which is the stabilizer of $0_{\mathbb{Z}_a}$ with respect to the action of $\mathsf{C}_n$ on $\mathbb{Z}_a$. We have $\omega_b=\omega_{n}^a$ and $\omega_a=\omega_{n}^b$, so that $\chi_{b,t}=\chi_{n,t}^a=\chi_{n,at}$ and $\chi_{a,u}=\chi_{n,u}^b=\chi_{n,bu}$ for all $t,u\in\mathbb{Z}$. In particular, 
\begin{equation}\label{decLZrZs}
L\left(\mathbb{Z}_a\right)=\bigoplus_{u\in\mathbb{Z}_a}W_{a,u}\qquad\text{and}\qquad L\left(\mathbb{Z}_b\right)=\bigoplus_{t\in\mathbb{Z}_b}W_{b,t}
\end{equation} 
are the decompositions into $\mathsf{C}_n$-irreducible representations. By virtue of the assumption that $a$ and $b$ are coprime, these spaces have only the trivial $\mathsf{C}_n$-representation in common.

We use the notation in Theorem \ref{TheoSymmultfree}, with 
\[
F=\mathsf{C}_n,\quad X=A=\{a,b\},\quad Y_a=\mathbb{Z}_b,\quad\text{and}\quad Y_b=\mathbb{Z}_a,
\]
where $\mathbb{Z}_a$ and $\mathbb{Z}_b$ are regarded as disjoint sets (since $a\neq b$). We think of $a$ and $b$ as generators of the groups $\mathbb{Z}_b$ and $\mathbb{Z}_a$, respectively; see, e.g., \cite[Proposition 1.2.5]{book4}. We are also in the framework of Example \ref{Example5ter} and of the related Example \ref{Example5}, and we use the notation therein; we also use some ideas from Example \ref{Example3ter}, without explicit references. To simplify notation, we set 
\[
Y_0\coloneqq\left\{0_{\mathbb{Z}_a},0_{\mathbb{Z}_b}\right\}
\]
and we assume that $0_{\mathbb{Z}_b}$ and $0_{\mathbb{Z}_a}$ are the base points $y_a$ and $y_b$, respectively; cf.\ Example \ref{Example5}. For every $\psi\in\Gamma$ we set 
\[
\psi^\flat(v)\coloneqq\begin{cases}
a\in\mathbb{Z}_b&\text{if }\psi(v)\in\mathbb{Z}_b,\\
b\in\mathbb{Z}_a&\text{if }\psi(v)\in\mathbb{Z}_a.
\end{cases}
\]
Subsequently, for each $Z\subseteq V$, we define $\varphi_Z\in\Gamma_0$ by setting 
\[
\varphi_Z(v)\coloneqq\begin{cases}
0_{\mathbb{Z}_b}&\text{if }v\in Z,\\
0_{\mathbb{Z}_a}&\text{if }v\in V\setminus Z.
\end{cases}
\]
Clearly, $\Gamma_0$ is the set of all functions from $V$ to $Y_0$, the map $Z\mapsto\varphi_Z$ is a bijection between subsets of $V$ and functions in $\Gamma_0$, and, for every $\varphi\in\Gamma$, we have $\varphi^\sharp=\varphi_Z$, where $Z=\left\{v\in V\mathrel{:}\varphi(v)\in\mathbb{Z}_b\right\}$. Moreover, from \eqref{defactionOm} we deduce that $\pi\varphi_Z=\varphi_{\pi Z}$ for all $Z\subseteq V$ and $\pi\in\mathsf{Sym}(V)$, so that in the decomposition
\[
\Gamma=\bigsqcup_{m=0}^{\lvert V\rvert}\left[\bigsqcup_{Z\subseteq V\mathrel{:}\lvert Z\rvert=m}\mathsf{Orb}_{\mathsf{C}_n^V}\left(\varphi_Z\right)\right]
\]
the sets in square brackets are the orbits of $\mathsf{C}_n\wr\mathsf{Sym}(V)$ on $\Gamma$. Therefore, from now on, we fix a proper nonempty subset $Z_0$ of $V$ and we study the permutation representation of $\mathsf{C}_n\wr\mathsf{Sym}(V)$ on $\mathsf{Orb}_{\mathsf{C}_n\wr\mathsf{Sym}(V)}\left(\varphi_{Z_0}\right)$. In particular, from \eqref{orbsym2} we deduce that 
\[
\mathsf{Stab}_{\mathsf{C}_n\wr\mathsf{Sym}(V)}\left(\varphi_{Z_0}\right)=\left(\mathsf{C}_a\wr\mathsf{Sym}\left(Z_0\right)\right)\bigtimes\left(\mathsf{C}_b\wr\mathsf{Sym}\left(V\setminus Z_0\right)\right). 
\]


Now we write, in the present setting, the decomposition \eqref{decOrbphi0}. We need some preliminary results.

\begin{itemize}
\item
First, \eqref{decOrbxi} is the decomposition of $\mathsf{Orb}_{\mathsf{Sym}(V)}\left(Z_0\right)$ into the orbits of $\mathsf{Stab}_{\mathsf{Sym}(V)}\left(Z_0\right)$, and we can write it in the form 
\begin{equation}\label{decOrbxiGel}
\mathsf{Orb}_{\mathsf{Sym}(V)}\left(Z_0\right)=\bigsqcup_{j=0}^{M}\Delta_j,
\end{equation}
where $M\coloneqq\min\left\{\left\lvert Z_0\right\rvert,\left\lvert V\setminus Z_0\right\rvert\right\}$
and $\Delta_j\coloneqq\left\{Z\in\mathsf{Orb}_{\mathsf{Sym}(V)}\left(Z_0\right)\mathrel{:}\left\lvert Z\setminus Z_0\right\rvert=j\right\}$; in particular, $\Delta_0=\left\{Z_0\right\}$. 

\item
In the decomposition \eqref{decYa}, we have $\ell_{a,b}=\ell_{b,a}=0$, $\ell_{a,a}=b-1$, and $\ell_{b,b}=a-1$. Moreover, $\Upsilon_{a,b,0}=\mathbb{Z}_b$, $\Upsilon_{b,a,0}=\mathbb{Z}_a$, $\Upsilon_{a,a,i}=\{i\}\subseteq\mathbb{Z}_b$ for $i=0,\dotsc,b-1$, and $\Upsilon_{b,b,i}=\{i\}\subseteq\mathbb{Z}_a$ for $i=0,\dotsc,a-1$.

\item
We write \eqref{defPsixij} in the form
\begin{equation}\label{defPsixijGel}
\Psi_j\coloneqq\bigl\{(Z,\mathsf{i})\mathrel{:}Z\in\Delta_j \text{ and }\mathsf{i}\in\mathbb{N}^V,\text{ with }0\leq \mathsf{i}(v)\leq \ell_{\varphi_Z^\flat(v),\varphi_{Z_0}^\flat(v)}\text{ for all } v \in V\bigr\}. 
\end{equation}

\item
For $0\leq j\leq M$, we define $\mathsf{T}(j)$ as the set of all pairs of nonnegative integer-valued functions $(\tau,\tau')$ such that $\tau$ is defined on $\{0,1,\dotsc,b-1\}$ and satisfies $\sum_{i=0}^{b-1}\tau(i)=\left\lvert Z_0\right\rvert-j$, whereas $\tau'$ is defined on $\{0,1,\dotsc,a-1\}$ and satisfies $\sum_{i=0}^{a-1}\tau'(i)=\left\lvert V\setminus Z_0\right\rvert-j$.

\item
For $(\tau,\tau')\in\mathsf{T}(j)$, we define $\Psi_{j,\tau,\tau'}$ as the set of all $(Z,\mathsf{i})\in\Psi_j$ such that $\left\lvert Z\cap Z_0\cap \mathsf{i}^{-1}(i)\right\rvert=\tau(i)$ for $i=0,1,\dotsc,b-1$,
and $\left\lvert(V\setminus Z)\cap\left(V\setminus Z_0\right)\cap\mathsf{i}^{-1}(i)\right\rvert=\tau'(i)$ for $i=0,1,\dotsc,a-1$.

\item
For $0\leq j\leq M$ and $(\tau,\tau')\in\mathsf{T}(j)$, we set
\begin{equation}\label{Orbxiphib}
\digamma_{j,\tau,\tau'}\coloneqq\bigsqcup_{(Z,\mathsf{i})\in\Psi_{j,\tau,\tau'}}\left(\bigtimes_{v\in V}\Upsilon_{\varphi_Z^\flat(v),\varphi_{Z_0}^\flat(v),\mathsf{i}(v)}\right).
\end{equation}

\end{itemize}

\begin{theorem}\label{decOrbxiphib}
In the present setting, \eqref{decOrbxiphi} is
\[
\mathsf{Orb}_{\mathsf{C}_n\wr\mathsf{Sym}(V)}\left(\varphi_{Z_0}\right)=\bigsqcup_{j=0}^{M}\;\bigsqcup_{(\tau,\tau')\in\mathsf{T}(j)}\digamma_{j,\tau,\tau'},
\]
and it is the decomposition into the orbits of $\mathsf{Stab}_{\mathsf{C}_n\wr\mathsf{Sym}(V)}\left(\varphi_{Z_0}\right)$.
\end{theorem}
\begin{proof}
Suppose that $(Z,\mathsf{i})\in\Psi_j$. If $v\in Z_0\cap Z$, then $\varphi_{Z_0}^\flat(v)=\varphi_Z^\flat(v)=a$, so that $\ell_{\varphi_Z^\flat(v),\varphi_{Z_0}^\flat(v)}=b-1$ and $\mathsf{i}\big|_{Z_0\cap Z}$ takes its values in $\{0,1,\dotsc,b-1\}$. Similarly, if $v\in V\setminus\left(Z_0\cup Z\right)=(V\setminus Z)\cap\left(V\setminus Z_0\right)$, then $\varphi_{Z_0}^\flat(v)=\varphi_Z^\flat(v)=b$, so that $\ell_{\varphi_Z^\flat(v),\varphi_{Z_0}^\flat(v)}=a-1$ and $\mathsf{i}\big|_{V\setminus\left(Z_0\cup Z\right)}$ takes its values in $\{0,1,\dotsc,a-1\}$. Finally, if $v\in\left(Z_0\setminus Z\right)\cup\left(Z\setminus Z_0\right)$, then $\varphi_{Z_0}^\flat(v)\neq\varphi_Z^\flat(v)$, so that $\ell_{\varphi_Z^\flat(v),\varphi_{Z_0}^\flat(v)}=0$ and $\mathsf{i}\big|_{\left(Z_0\setminus Z\right)\cup\left(Z\setminus Z_0\right)}$ is identically zero.

Since every $Z\in\Delta_j$ satisfies $\left\lvert Z\right\rvert=\left\lvert Z_0\right\rvert$ and $\left\lvert Z\setminus Z_0\right\rvert=j$, we also have $\left\lvert Z_0\setminus Z\right\rvert=j$, $\left\lvert Z\cap Z_0\right\rvert=\left\lvert Z_0\right\rvert-j$, and $\left\lvert(V\setminus Z)\cap\left(V\setminus Z_0\right)\right\rvert=\left\lvert V\setminus Z_0\right\rvert-j$. Recall that $\mathsf{Stab}_{\mathsf{Sym}(V)}\left(Z_0\right)=\mathsf{Sym}\left(Z_0\right)\bigtimes\mathsf{Sym}\left(V\setminus Z_0\right)$. From the above considerations, it follows that two pairs $(Z,\mathsf{i}),(Z',\mathsf{i}')\in\Psi_j$ lie in the same $\mathsf{Stab}_{\mathsf{Sym}(V)}\left(Z_0\right)$-orbit if and only if
\[
\begin{split}
&\left\lvert Z\cap Z_0\cap \mathsf{i}^{-1}(i)\right\rvert=\left\lvert Z'\cap Z_0\cap (\mathsf{i}')^{-1}(i)\right\rvert\quad\text{for }i=0,1,\dotsc,b-1,\qquad\text{and}\\
&\left\lvert(V\setminus Z)\cap\left(V\setminus Z_0\right)\cap\mathsf{i}^{-1}(i)\right\rvert=\left\lvert(V\setminus Z')\cap\left(V\setminus Z_0\right)\cap(\mathsf{i}')^{-1}(i)\right\rvert\quad\text{for }i=0,1,\dotsc,a-1;
\end{split}
\]
indeed, necessity is immediate because any $\pi\in\mathsf{Stab}_{\mathsf{Sym}(V)}\left(Z_0\right)$ preserves both $Z_0$ and $V\setminus Z_0$, while sufficiency follows by combining a bijection of $Z_0$ mapping $Z_0\setminus Z$ onto $Z_0\setminus Z'$ and each $Z\cap Z_0\cap\mathsf{i}^{-1}(i)$ onto $Z'\cap Z_0\cap(\mathsf{i}')^{-1}(i)$ with a bijection of $V\setminus Z_0$ mapping $Z\setminus Z_0$ onto $Z'\setminus Z_0$ and each $(V\setminus Z)\cap\left(V\setminus Z_0\right)\cap\mathsf{i}^{-1}(i)$ onto $(V\setminus Z')\cap\left(V\setminus Z_0\right)\cap(\mathsf{i}')^{-1}(i)$. Thus $(Z,\mathsf{i})$ and $(Z',\mathsf{i}')$ are in the same $\mathsf{Stab}_{\mathsf{Sym}(V)}\left(Z_0\right)$-orbit if and only if they belong to the same set $\Psi_{j,\tau,\tau'}$. Moreover, for every $(\tau,\tau')\in\mathsf{T}(j)$, the set $\Psi_{j,\tau,\tau'}$ is nonempty: choosing any $Z\in\Delta_j$, we may partition $Z\cap Z_0$ into subsets of cardinalities $\tau(0),\dotsc,\tau(b-1)$ and $(V\setminus Z)\cap\left(V\setminus Z_0\right)$ into subsets of cardinalities $\tau'(0),\dotsc,\tau'(a-1)$, and define $\mathsf{i}\in\mathbb{N}^V$ accordingly (setting $\mathsf{i}\equiv 0$ on $\left(Z_0\setminus Z\right)\cup\left(Z\setminus Z_0\right)$) to obtain $(Z,\mathsf{i})\in\Psi_{j,\tau,\tau'}$. Therefore, \eqref{decPsixij}, which is now the decomposition of $\Psi_j$ into the orbits of $\mathsf{Stab}_{\mathsf{Sym}(V)}\left(Z_0\right)$, becomes 
\begin{equation}\label{decPsijtau} 
\Psi_j=\bigsqcup_{(\tau,\tau')\in\mathsf{T}(j)}\Psi_{j,\tau,\tau'}.
\end{equation}
The statement then follows directly from Theorem \ref{TheodecOrbxiphi} and \eqref{decOrbphi0}.
\end{proof}


For $Z\subseteq V$ and $\psi\in\mathsf{Orb}_{\mathsf{C}_n^V}\left(\varphi_Z\right)$, we define $\chi_\psi\in L\left(\mathsf{Orb}_{\mathsf{C}_n^V}\left(\varphi_Z\right)\right)$ by setting
\begin{equation}\label{Defchipsi}
\chi_\psi(\varphi)\coloneqq\prod_{v\in V}\chi_{n/\psi^\flat(v),\,\psi(v)}(\varphi(v))
=\prod_{v\in Z}\chi_{b,\psi(v)}(\varphi(v))\prod_{v\in V\setminus Z}\chi_{a,\psi(v)}(\varphi(v))
\end{equation}
for every $\varphi\in\mathsf{Orb}_{\mathsf{C}_n^V}\left(\varphi_Z\right)$. Henceforth, we slightly modify the notation in Section \ref{Secprelimin}. 
We set $\mathsf{supp}_a(\psi)\coloneqq\psi^{-1}\left(\mathbb{Z}_b\setminus\left\{0_{\mathbb{Z}_b}\right\}\right)$ and $\mathsf{supp}_b(\psi)\coloneqq\psi^{-1}\left(\mathbb{Z}_a\setminus\left\{0_{\mathbb{Z}_a}\right\}\right)$. Then 
\[
\mathsf{supp}(\psi)\coloneqq\mathsf{supp}_a(\psi)\sqcup\mathsf{supp}_b(\psi)=\mathsf{supp}\left(\chi_\psi\right). 
\]
Similarly, we set $\mathsf{image}(\psi)\coloneqq\psi(V)$, whereas $\mathsf{image}\left(\chi_\psi\right)$ is as in Section \ref{Secprelimin} but, by virtue of \eqref{decLZrZs}, we think of it as a subset of $\left(\mathbb{Z}_a\setminus\left\{0_{\mathbb{Z}_a}\right\}\right)\sqcup\left(\mathbb{Z}_b\setminus\left\{0_{\mathbb{Z}_b}\right\}\right)\sqcup\{0\}$. Therefore,
\[
\mathsf{image}\left(\chi_\psi\right)\setminus\{0\}=\mathsf{image}(\psi)\setminus\left\{0_{\mathbb{Z}_a},0_{\mathbb{Z}_b}\right\};
\]
that is, in $\mathsf{image}\left(\chi_\psi\right)$, the trivial representation in both decompositions in \eqref{decLZrZs} yields the same element, namely $0$, while in $\mathsf{image}(\psi)$ the zeros $0_{\mathbb{Z}_a}$ and $0_{\mathbb{Z}_b}$ are distinct elements.
The distinction between $\mathsf{image}\left(\chi_\psi\right)$ and $\mathsf{image}(\psi)$ is crucial in the description of the Gelfand pair for the inertia group.
With the convention $\psi^{-1}(0)\coloneqq\psi^{-1}(Y_0)$, if $\varphi\in\Gamma_0$, we denote by $\mathsf{supp}(\varphi\wedge\psi)$ the set of all pairs $(y,z)\in Y_0\times\mathsf{image}\left(\chi_\psi\right)$ such that $\varphi^{-1}(y)\cap \psi^{-1}(z)\neq\emptyset$. Then we define $\varphi\wedge\psi$ as the following set composition of $V$:
\begin{equation}\label{Vassetcomp}
V=\bigsqcup_{(y,z)\in\mathsf{supp}(\varphi\wedge\psi)}\left[\varphi^{-1}(y)\cap \psi^{-1}(z)\right],
\end{equation}
where each intersection is nonempty.

From \eqref{IsigmaSym} we deduce that, for every $\psi\in\Gamma$, the inertia group of $\chi_\psi$, which we denote simply by $\mathsf{I}\left(\psi\right)$, is
\begin{equation}\label{Ichipsi}
\mathsf{I}\left(\psi\right)=\mathsf{Sym}\left(\psi^{-1}\left(Y_0\right)\right)\bigtimes\left[\prod_{z\in\mathsf{image}\left(\chi_\psi\right)\setminus\{0\}}\mathsf{Sym}\left(\psi^{-1}(z)\right)\right].
\end{equation}
This is distinct from the stabilizer of $\psi$ in $\mathsf{Sym}(V)$, which is
\begin{equation}\label{Stabpsi}
\begin{split}
\mathsf{Stab}_{\mathsf{Sym}(V)}(\psi)&=\mathsf{Sym}\left(\psi^{-1}\left(0_{\mathbb{Z}_a}\right)\right)\\
&\quad\bigtimes\mathsf{Sym}\left(\psi^{-1}\left(0_{\mathbb{Z}_b}\right)\right)\bigtimes\left[\prod_{z\in\mathsf{image}\left(\chi_\psi\right)\setminus \{0\}}\mathsf{Sym}\left(\psi^{-1}(z)\right)\right].
\end{split}
\end{equation}
Again, the distinction between these two stabilizers is a crucial ingredient in the description of the Gelfand pair for the inertia group.
We slightly modify the relation $\prec$ by defining it in the following equivalent way: we set $\psi\prec Z$ when $\mathsf{supp}_a(\psi)\subseteq Z$ and $\mathsf{supp}_b(\psi)\subseteq V\setminus Z$. With the usual notation, this is equivalent to $\chi_\psi\prec\varphi_Z$; see also Example \ref{Example5bis}. 
Then \eqref{defJsigma} may be written in the form
\begin{equation}\label{descJpsi}
\mathsf{J}_\psi=\left\{Z\in\mathsf{Orb}_{\mathsf{Sym}(V)}\left(Z_0\right)\mathrel{:}\psi\prec Z\right\}.
\end{equation}

From now on, we assume that $2\leq\left\lvert Z_0\right\rvert\leq\left\lvert V\right\rvert -2$. In order to simplify some points, we think of $\Sigma_0$ as a set of functions in $\Gamma$ (using $\psi$ in place of $\chi_\psi$) and also assume that for every $\psi\in\Sigma_0$ we have 
\begin{equation}\label{assumpSigma0}
\psi^{-1}\left(\mathbb{Z}_b\right)= Z_0\qquad\text{and}\qquad\psi^{-1}\left(\mathbb{Z}_a\right)= V\setminus Z_0.
\end{equation}
This is slightly stronger than the requirement $\psi\prec Z_0$ on the function $\psi$, but involves no loss of generality for the character $\chi_\psi$, since changing the value of $\psi(v)$ between $0_{\mathbb{Z}_a}$ and $0_{\mathbb{Z}_b}$ does not alter $\chi_\psi$. Moreover, in order to formulate the following proposition, for each $\psi\in\Sigma_0$ we set 
\begin{equation}\label{defhandm}
h\coloneqq\left\lvert \psi^{-1}\left(Y_0\right)\right\rvert\qquad\text{and}\qquad m\coloneqq\min\left\{\left\lvert \psi^{-1}\left(0_{\mathbb{Z}_a}\right)\right\rvert,\left\lvert \psi^{-1}\left(0_{\mathbb{Z}_b}\right)\right\rvert\right\},
\end{equation}
and for $0\leq k\leq m$ we define $\overline{S}^{h-k,k}$ as the inflation from $\mathsf{Sym}\left(\psi^{-1}\left(Y_0\right)\right)$ to $\mathsf{I}\left(\psi\right)$ of the irreducible representation $S^{h-k,k}$. Note that $h$ and $m$ depend on $\psi$, but we omit any subscript. Henceforth, for every subset $Z\subseteq V$ containing $\mathsf{supp}_a(\psi)$, we set 
\[
Z'\coloneqq Z\setminus\mathsf{supp}_a(\psi)
\]
(so that $Z_0'=\psi^{-1}\left(0_{\mathbb{Z}_b}\right)$ by \eqref{assumpSigma0}), and then 
\[
\Pi_j'\coloneqq\left\{Z'\subseteq \psi^{-1}\left(Y_0\right)\mathrel{:}\left\lvert Z'\right\rvert=\left\lvert Z_0'\right\rvert\text{ and }\left\lvert Z'\setminus Z_0'\right\rvert=j\right\}
\]
for $j=0,1,\dotsc,m$.


\begin{proposition}\label{Propmultfreebis}
\begin{enumerate}[(1)]

\item\label{Propmultfreebis1} The permutation representation of $\mathsf{C}_n\wr\mathsf{Sym}(V)$ on $\mathsf{Orb}_{\mathsf{C}_n\wr\mathsf{Sym}(V)}\left(\varphi_{Z_0}\right)$ is multiplicity free; that is, $\left(\mathsf{C}_n\wr\mathsf{Sym}(V),\mathsf{Stab}_{\mathsf{C}_n\wr\mathsf{Sym}(V)}\left(\varphi_{Z_0}\right)\right)$ is a Gelfand pair.

\item\label{Propmultfreebis2} For every $\psi\in\Sigma_0$, we have $\mathsf{J}_\psi=\mathsf{Orb}_{\mathsf{I}\left(\psi\right)}\left(Z_0\right)$.

\item\label{Propmultfreebis3} For every $\psi\in\Sigma_0$, 
\begin{equation}\label{StabpsiIZ0}
\mathsf{Stab}_{\mathsf{I}\left(\psi\right)}\left(Z_0\right)=\mathsf{Stab}_{\mathsf{Sym}(V)}(\psi),
\end{equation} 
and the action of $\mathsf{I}\left(\psi\right)$ on $\mathsf{Orb}_{\mathsf{I}\left(\psi\right)}\left(Z_0\right)$ is the inflation of the action of $\mathsf{Sym}\left(\psi^{-1}\left(Y_0\right)\right)$ on
\[
\mathsf{Sym}\left(\psi^{-1}\left(Y_0\right)\right)\big/\left[\mathsf{Sym}\left(\psi^{-1}\left(0_{\mathbb{Z}_a}\right)\right)\bigtimes\mathsf{Sym}\left(\psi^{-1}\left(0_{\mathbb{Z}_b}\right)\right)\right].
\]

\item\label{Propmultfreebis4} 
For every $\psi\in\Sigma_0$, the permutation representation of $\mathsf{I}\left(\psi\right)$ on $\mathsf{Orb}_{\mathsf{I}\left(\psi\right)}\left(Z_0\right)$ is multiplicity free. Moreover, the decomposition into irreducible representations is
$L\left(\mathsf{Orb}_{\mathsf{I}\left(\psi\right)}\left(Z_0\right)\right)\sim\bigoplus_{k=0}^m\overline{S}^{h-k,k}$.

\item\label{Propmultfreebis5}
For every $\psi\in\Sigma_0$, the decomposition \eqref{decPil} is now
\begin{equation}\label{decPilGel}
\mathsf{Orb}_{\mathsf{I}\left(\psi\right)}\left(Z_0\right)=\bigsqcup_{j=0}^m\Pi_j,
\end{equation}
where $\Pi_j\coloneqq \left\{Z'\sqcup\mathsf{supp}_a(\psi)\mathrel{:}Z'\in\Pi_j'\right\}$.

\item\label{Propmultfreebis6}
For every $\psi\in\Sigma_0$, the decomposition \eqref{StabOrbIsig0} is now
\begin{equation}\label{OrbIpsiZ0bis}
\mathsf{Stab}_{\mathsf{Sym}(V)}\left(Z_0\right)\mathsf{Orb}_{\mathsf{I}\left(\psi\right)}\left(Z_0\right)=\bigsqcup_{j=0}^m\Delta_j,
\end{equation}
where $\Delta_j$ is as in \eqref{decOrbxiGel}. Moreover, we have $\Pi_j\subseteq\Delta_j$ for $j=0,1,\dotsc,m$.

\item\label{Propmultfreebis7} 
For every $\psi\in\Sigma_0$, with the notation in \eqref{OrbV0Vt} and \eqref{defQ0sig}, the map 
\[
\begin{array}{ccc}
\mathsf{Q}_\psi^0&\longrightarrow&\mathsf{Orb}_{\mathsf{Stab}_{\mathsf{Sym}(V)}\left(Z_0\right)}\left(\psi^{-1}(z)\mathrel{:}z\in\mathsf{image}(\psi)\right)\\[4pt]
\kappa&\longmapsto& \left(\kappa\psi^{-1}(z)\mathrel{:}z\in\mathsf{image}(\psi)\right)
\end{array}
\]
is a bijection. 

\item\label{Propmultfreebis8} 
For every $\psi\in\Sigma_0$, in the present setting, the condition on $q$ in \eqref{defFsigthet} takes the following form: for each $Z\in\mathsf{Stab}_{\mathsf{Sym}(V)}\left(Z_0\right)\mathsf{Orb}_{\mathsf{I}\left(\psi\right)}\left(Z_0\right)$, we must consider the set of all $\kappa\in \mathsf{Q}_\psi^0$ such that
\[
\mathsf{supp}_a(\psi)\subseteq\kappa^{-1}Z\qquad\text{and}\qquad\left[\kappa^{-1}Z\right]'\in\mathsf{Orb}_{\mathsf{I}(\psi)}\left(Z_0'\right).
\]

\end{enumerate}
\end{proposition}
\begin{proof}
(1) Indeed, the last conditions in Theorem \ref{TheoSymmultfree}\eqref{TheoSymmultfree2} are satisfied, because $\left\lvert\widehat{F}_a\cap\widehat{F}_b\right\rvert=1$.

\medskip

(2) From the description of $\mathsf{J}_\psi$ in \eqref{descJpsi}, it follows that the equality in the statement holds (that is, we have equality in \eqref{defJsig0old}) because the following condition is satisfied: for every $Z\in\mathsf{Orb}_{\mathsf{Sym}(V)}\left(Z_0\right)$ such that $\psi\prec Z$, there exists $\varepsilon\in\mathsf{I}\left(\psi\right)$ such that $\varepsilon Z_0=Z$. This is in turn equivalent to the simultaneous occurrence of the following conditions:
\[
\begin{array}{lll}
\varepsilon\left(\psi^{-1}(z)\right)=\psi^{-1}(z) &\text{for each } z\in\mathsf{image}\left(\chi_\psi\right) &\qquad\left(\varepsilon\in\mathsf{I}\left(\psi\right)\right),\\[4pt]
\varepsilon\left(\varphi_{Z_0}^{-1}(y)\right)=\varphi^{-1}_Z(y)&\text{for each } y\in Y_0 &\qquad\left(\varepsilon Z_0=Z\right).
\end{array}
\]
To see that they hold for a suitable $\varepsilon\in\mathsf{Sym}(V)$, we show that the set compositions
\[
V=\bigsqcup_{(y,z)\in\mathsf{supp}\left(\varphi_{Z_0}\wedge\psi\right)}\left[\varphi_{Z_0}^{-1}(y)\cap \psi^{-1}(z)\right]=\bigsqcup_{(y,z)\in\mathsf{supp}\left(\varphi_Z\wedge\psi\right)}\left[\varphi_Z^{-1}(y)\cap \psi^{-1}(z)\right]
\] 
have the same support and the same block cardinalities (and hence the same type; cf.\ \eqref{deftype} and \eqref{Vassetcomp}). Indeed, since $\psi\prec Z_0$ and $\psi\prec Z$, for every $z\in\mathbb{Z}_b\setminus\left\{0_{\mathbb{Z}_b}\right\}$ we have $\psi^{-1}(z)\subseteq\mathsf{supp}_a(\psi)\subseteq Z_0\cap Z$, so that
\[
\varphi_{Z_0}^{-1}\left(0_{\mathbb{Z}_b}\right)\cap\psi^{-1}(z)=\varphi_Z^{-1}\left(0_{\mathbb{Z}_b}\right)\cap\psi^{-1}(z)=\psi^{-1}(z)
\]
and $\varphi_{Z_0}^{-1}\left(0_{\mathbb{Z}_a}\right)\cap\psi^{-1}(z)=\varphi_Z^{-1}\left(0_{\mathbb{Z}_a}\right)\cap\psi^{-1}(z)=\emptyset$; symmetrically, for every $z\in\mathbb{Z}_a\setminus\left\{0_{\mathbb{Z}_a}\right\}$ we have $\psi^{-1}(z)\subseteq\mathsf{supp}_b(\psi)\subseteq (V\setminus Z_0)\cap(V\setminus Z)$, so that
\[
\varphi_{Z_0}^{-1}\left(0_{\mathbb{Z}_a}\right)\cap\psi^{-1}(z)=\varphi_Z^{-1}\left(0_{\mathbb{Z}_a}\right)\cap\psi^{-1}(z)=\psi^{-1}(z)
\]
and $\varphi_{Z_0}^{-1}\left(0_{\mathbb{Z}_b}\right)\cap\psi^{-1}(z)=\varphi_Z^{-1}\left(0_{\mathbb{Z}_b}\right)\cap\psi^{-1}(z)=\emptyset$. For $z=0$, recalling that $\psi^{-1}(0)=\psi^{-1}\left(Y_0\right)$ and $\lvert Z\rvert=\lvert Z_0\rvert$, we obtain
\[
\varphi_{Z_0}^{-1}\left(0_{\mathbb{Z}_b}\right)\cap\psi^{-1}(0)=Z_0\setminus\mathsf{supp}_a(\psi)
\quad\text{and}\quad
\varphi_Z^{-1}\left(0_{\mathbb{Z}_b}\right)\cap\psi^{-1}(0)=Z\setminus\mathsf{supp}_a(\psi),
\]
both of cardinality $\left\lvert Z_0\right\rvert-\left\lvert\mathsf{supp}_a(\psi)\right\rvert$, and similarly
\[
\varphi_{Z_0}^{-1}\left(0_{\mathbb{Z}_a}\right)\cap\psi^{-1}(0)=\left(V\setminus Z_0\right)\setminus\mathsf{supp}_b(\psi)
\quad\text{and}\quad
\varphi_Z^{-1}\left(0_{\mathbb{Z}_a}\right)\cap\psi^{-1}(0)=(V\setminus Z)\setminus\mathsf{supp}_b(\psi),
\]
both of cardinality $\left\lvert V\setminus Z_0\right\rvert-\left\lvert\mathsf{supp}_b(\psi)\right\rvert$. Hence $\mathsf{supp}\left(\varphi_{Z_0}\wedge\psi\right)=\mathsf{supp}\left(\varphi_Z\wedge\psi\right)$ and there exists $\varepsilon\in\mathsf{Sym}(V)$ such that
\[
\varepsilon\left[\varphi_{Z_0}^{-1}(y)\cap \psi^{-1}(z)\right]=\varphi_Z^{-1}(y)\cap \psi^{-1}(z)\quad\text{for every }(y,z)\in\mathsf{supp}\left(\varphi_{Z_0}\wedge\psi\right).
\]
Taking unions over $y\in Y_0$ and over $z\in\mathsf{image}\left(\chi_\psi\right)$ yields $\varepsilon\in\mathsf{I}(\psi)$ and $\varepsilon Z_0=Z$.

\medskip

(3)
By \eqref{assumpSigma0}, we have $Z_0=\mathsf{supp}_a(\psi)\sqcup\psi^{-1}\left(0_{\mathbb{Z}_b}\right)$ and $V\setminus Z_0=\mathsf{supp}_b(\psi)\sqcup\psi^{-1}\left(0_{\mathbb{Z}_a}\right)$. Since every element of $\mathsf{I}\left(\psi\right)$ stabilizes $\mathsf{supp}_a(\psi)$, $\mathsf{supp}_b(\psi)$, and $\psi^{-1}\left(Y_0\right)$, it stabilizes $Z_0$ if and only if it stabilizes $\psi^{-1}\left(0_{\mathbb{Z}_b}\right)$ and $\psi^{-1}\left(0_{\mathbb{Z}_a}\right)$. Comparing \eqref{Ichipsi} and \eqref{Stabpsi} immediately yields $\mathsf{Stab}_{\mathsf{I}\left(\psi\right)}\left(Z_0\right)=\mathsf{Stab}_{\mathsf{Sym}(V)}(\psi)$, and the second assertion follows from \eqref{definflperrep}, \eqref{Ichipsi}, and \eqref{Stabpsi}.

\medskip

(4) This follows from \eqref{Propmultfreebis3} and the Young rule, taking also \eqref{defhandm} into account.

\medskip

(5) These are elementary facts in the combinatorics of the Johnson scheme; see, e.g., \cite[Exercise 3.14.2]{book}. Indeed, we consider the action of $\mathsf{I}\left(\psi\right)$ on the subsets $Z'$, which may be identified with the points of the homogeneous space in \eqref{Propmultfreebis3}, taking into account that $\mathsf{supp}_a(\psi)$ is contained in $Z_0$ and is stabilized by $\mathsf{I}\left(\psi\right)$. More explicitly, for all $\varepsilon\in\mathsf{I}\left(\psi\right)$, we have
\[
\varepsilon\left(Z_0\setminus\mathsf{supp}_a(\psi)\right)=\left(\varepsilon Z_0\right)\setminus\mathsf{supp}_a(\psi).
\]
Therefore, the bijection  
\begin{equation}\label{OrbIpsiZ0}
\mathsf{Orb}_{\mathsf{I}\left(\psi\right)}\left(Z_0\right)\ni Z\;\longmapsto\;Z'\in\mathsf{Orb}_{\mathsf{I}\left(\psi\right)}\left(Z_0'\right)
\end{equation}
establishes an isomorphism of $\mathsf{I}\left(\psi\right)$-spaces. Under this isomorphism, the orbits of $\mathsf{Stab}_{\mathsf{I}\left(\psi\right)}\left(Z_0\right)$ on $\mathsf{Orb}_{\mathsf{I}\left(\psi\right)}\left(Z_0\right)$ correspond to the orbits $\Pi_j'$ ($j=0,1,\dotsc,m$) of $\mathsf{Stab}_{\mathsf{I}\left(\psi\right)}\left(Z_0'\right)$ on $\mathsf{Orb}_{\mathsf{I}\left(\psi\right)}\left(Z_0'\right)$, yielding \eqref{decPilGel}.

\medskip

(6) We first verify that $\Pi_j\subseteq\Delta_j$ for each $j=0,1,\dotsc,m$: if $Z\in\Pi_j$, then $\lvert Z\rvert=\lvert Z'\rvert+\left\lvert\mathsf{supp}_a(\psi)\right\rvert=\left\lvert Z_0'\right\rvert+\left\lvert\mathsf{supp}_a(\psi)\right\rvert=\left\lvert Z_0\right\rvert$ and $Z\setminus Z_0=Z'\setminus Z_0'$, so that $\left\lvert Z\setminus Z_0\right\rvert=j$ and hence $Z\in\Delta_j$. Applying $\mathsf{Stab}_{\mathsf{Sym}(V)}\left(Z_0\right)$ to both sides of \eqref{decPilGel} yields
\[
\mathsf{Stab}_{\mathsf{Sym}(V)}\left(Z_0\right)\mathsf{Orb}_{\mathsf{I}\left(\psi\right)}\left(Z_0\right)=\bigcup_{j=0}^m\mathsf{Stab}_{\mathsf{Sym}(V)}\left(Z_0\right)\Pi_j.
\]
Since each $\Pi_j$ is nonempty and contained in $\Delta_j$, and each $\Delta_j$ is a single orbit of $\mathsf{Stab}_{\mathsf{Sym}(V)}\left(Z_0\right)$, we have $\mathsf{Stab}_{\mathsf{Sym}(V)}\left(Z_0\right)\Pi_j=\Delta_j$ for $j=0,1,\dotsc,m$, and \eqref{OrbIpsiZ0bis} follows because the sets $\Delta_j$ are pairwise disjoint.

\medskip

(7)
This follows from \eqref{StabStabIsigma} and \eqref{StabpsiIZ0}, observing that
\[
\mathsf{Stab}_{\mathsf{I}\left(\psi\right)}\left(Z_0\right)=\mathsf{Stab}_{\mathsf{Sym}(V)}(\psi)=\bigcap_{z\in\mathsf{image}(\psi)}\mathsf{Stab}_{\mathsf{Sym}(V)}\left[\psi^{-1}(z)\right]
\]
is contained in $\mathsf{Stab}_{\mathsf{Sym}(V)}\left(Z_0\right)$ and therefore coincides with the stabilizer in $\mathsf{Stab}_{\mathsf{Sym}(V)}\left(Z_0\right)$ of the tuple $\left(\psi^{-1}(z)\mathrel{:}z\in\mathsf{image}(\psi)\right)$.

\medskip

(8)
The condition on $q$ in \eqref{defFsigthet} is now given by
\[
\kappa\in \mathsf{Q}_\psi^0 \mathrel{:} \kappa^{-1}Z\in\mathsf{Orb}_{\mathsf{I}\left(\psi\right)}\left(Z_0\right) \qquad \text{for } Z\in\mathsf{Stab}_{\mathsf{Sym}(V)}\left(Z_0\right)\mathsf{Orb}_{\mathsf{I}\left(\psi\right)}\left(Z_0\right). 
\]
From \eqref{decPilGel} and \eqref{OrbIpsiZ0}, we deduce that
\[
\kappa^{-1}Z\in\mathsf{Orb}_{\mathsf{I}\left(\psi\right)}\left(Z_0\right) \quad\Longleftrightarrow\quad Z=\kappa\left(Z'\sqcup\mathsf{supp}_a(\psi)\right) \text{ for some } Z'\in\mathsf{Orb}_{\mathsf{I}(\psi)}\left(Z_0'\right),
\]
which means precisely that $\mathsf{supp}_a(\psi)\subseteq\kappa^{-1}Z$ and $\left[\kappa^{-1}Z\right]'\in\mathsf{Orb}_{\mathsf{I}(\psi)}\left(Z_0'\right)$.
\end{proof}

\begin{remark}
On the one hand, taking Proposition \ref{Propmultfree} into account, Proposition \ref{Propmultfreebis}\eqref{Propmultfreebis2} and the multiplicity-freeness in Proposition \ref{Propmultfreebis}\eqref{Propmultfreebis4} follow from Proposition \ref{Propmultfreebis}\eqref{Propmultfreebis1}. This is not surprising, because everything is based on the machinery developed in Section \ref{SecYoungpermrep}. On the other hand, we need Proposition \ref{Propmultfreebis}\eqref{Propmultfreebis3} because it contains the identification of an essential homogeneous space; similarly, we need the explicit decomposition in Proposition \ref{Propmultfreebis}\eqref{Propmultfreebis4}. 
\end{remark}


The invariant vectors in \eqref{defSabsig} are now the characters in $L\left(\mathbb{Z}_a\right)$ and in $L\left(\mathbb{Z}_b\right)$. We use them to give a precise expression for the functions $\mathcal{S}_{\varphi,\sigma}$.

\begin{proposition}\label{Propchifphi}
For every $\psi\in\Sigma_0$, the function $f\mapsto\mathcal{S}_{\varphi_Z,\chi_\psi}\left(f\cdot\varphi_Z\right)$ in \eqref{Sxisigmaphi} does not depend on $Z\in\mathsf{Orb}_{\mathsf{I}\left(\psi\right)}\left(Z_0\right)$ (that is, on $\varphi^\sharp=\varphi_Z$ for $\varphi\in\mathsf{Orb}_{\mathsf{C}_n^V}\left(\varphi_Z\right)$). Moreover, it is given by the character
\[
\mathsf{C}_n^V\ni f\longmapsto\chi_\psi\left(f\cdot\psi^\sharp\right);
\]
equivalently, $\mathcal{S}_{\varphi_Z,\chi_\psi}(\varphi)=\chi_\psi(\varphi)$ for every $\varphi\in\mathsf{Orb}_{\mathsf{C}_n^V}\left(\varphi_Z\right)$.
\end{proposition}
\begin{proof}
First, note that, with the notation in \eqref{Sxisigmaphi}, the function $f\mapsto\mathcal{S}_{\varphi_Z,\chi_\psi}\left(f\cdot\varphi_Z\right)$ is a unit matrix coefficient of the one-dimensional representation $\sigma$ associated with $\chi_\psi$ that is right $\mathsf{Stab}_{\mathsf{C}_n^V}\left(\varphi_Z\right)$-invariant and left $\mathsf{Stab}_{\mathsf{C}_n^V}\left(\varphi_{Z_0}\right)$-invariant; these properties determine it uniquely up to a scalar of modulus $1$, and the scalar is fixed by its value at the identity of $\mathsf{C}_n^V$.

Since $\psi\in\Sigma_0$ and $Z\in\mathsf{Orb}_{\mathsf{I}\left(\psi\right)}\left(Z_0\right)$, by \eqref{assumpSigma0} and Proposition \ref{Propmultfreebis}\eqref{Propmultfreebis2} (see also \eqref{Ichipsi}) we have $\psi^\sharp=\varphi_{Z_0}$, $\mathsf{supp}_a(\psi)\subseteq Z\cap Z_0$, and $\mathsf{supp}_b(\psi)\subseteq V\setminus\left(Z\cup Z_0\right)$, so that $\varphi_Z(v)=\varphi_{Z_0}(v)=\psi^\sharp(v)$ for every $v\in\mathsf{supp}(\psi)$. In particular, in \eqref{defSabsig}, the chosen invariant unit vectors coincide at every $v\in\mathsf{supp}(\psi)$ and correspond to the trivial representation at every $v\notin\mathsf{supp}(\psi)$, which yields $\mathcal{S}_{\varphi_Z,\chi_\psi}\left(\varphi_Z\right)=1$. Moreover, by \eqref{Defchipsi}, the factors corresponding to $v\notin\mathsf{supp}(\psi)$ are identically $1$, so the function $f\mapsto\chi_\psi\left(f\cdot\psi^\sharp\right)$ is a character of $\mathsf{C}_n^V$ (taking the value $1$ at the identity) and satisfies
\[
\chi_\psi\left(f\cdot\varphi_Z\right)=\chi_\psi\left(f\cdot\varphi_{Z_0}\right)=\chi_\psi\left(f\cdot\psi^\sharp\right)
\]
for every $f\in\mathsf{C}_n^V$. Since $\mathsf{C}_n^V$ is abelian, this implies that the function $f\mapsto\chi_\psi\left(f\cdot\psi^\sharp\right)$ is both $\mathsf{Stab}_{\mathsf{C}_n^V}\left(\varphi_Z\right)$-invariant and $\mathsf{Stab}_{\mathsf{C}_n^V}\left(\varphi_{Z_0}\right)$-invariant. Hence $\mathcal{S}_{\varphi_Z,\chi_\psi}\left(f\cdot\varphi_Z\right)=\chi_\psi\left(f\cdot\psi^\sharp\right)=\chi_\psi\left(f\cdot\varphi_Z\right)$ for all $f\in\mathsf{C}_n^V$, which ends the proof.
\end{proof}

Using the notation in Proposition \ref{Propmultfreebis}\eqref{Propmultfreebis4}, for each $\psi\in\Sigma_0$ and $0\leq k\leq m$ we denote by $\theta_k$ the irreducible representation of $\mathsf{I}\left(\psi\right)$ on $\overline{S}^{h-k,k}$. Moreover, we consider the spherical functions in \eqref{defRsigmatheta} as functions defined on $\mathsf{Orb}_{\mathsf{I}\left(\psi\right)}\left(Z_0\right)$; cf.\ Proposition \ref{Propmultfreebis}\eqref{Propmultfreebis3}.
\begin{proposition}
If $\psi\in\Sigma_0$ and $0\leq k\leq m$, then the spherical function associated with the irreducible representation $\Ind^{\mathsf{C}_n\wr\mathsf{Sym}(V)}_{\mathsf{C}_n\wr\mathsf{I}(\psi)}\left(\widetilde{\chi}_\psi\otimes\overline{\theta}_k\right)$
has the following form:
\[
\mathcal{F}_{\psi,\theta_k}(\varphi)=
\frac{1}{\left\lvert\mathsf{Q}^0_\psi\right\rvert}\sum_{\substack{\kappa\in\mathsf{Q}^0_\psi\mathrel{:}\;\mathsf{supp}_a(\psi)\subseteq\kappa^{-1}Z,\\ \left[\kappa^{-1}Z\right]'\in\mathsf{Orb}_{\mathsf{I}(\psi)}\left(Z_0'\right)}}\mathcal{R}_{\psi,\theta_k}\left(\kappa^{-1}Z\right)\chi_\psi\left(\kappa^{-1}\varphi\right)
\]
for every $Z\in\mathsf{Stab}_{\mathsf{Sym}(V)}\left(Z_0\right)\mathsf{Orb}_{\mathsf{I}\left(\psi\right)}\left(Z_0\right)$ and $\varphi\in\mathsf{Orb}_{\mathsf{C}_n^V}\left(\varphi_Z\right)$, whereas $\mathcal{F}_{\psi,\theta_k}(\varphi)=0$ for all $Z\in\mathsf{Orb}_{\mathsf{Sym}(V)}\left(Z_0\right)\setminus\mathsf{Stab}_{\mathsf{Sym}(V)}\left(Z_0\right)\mathsf{Orb}_{\mathsf{I}\left(\psi\right)}\left(Z_0\right)$ and $\varphi\in\mathsf{Orb}_{\mathsf{C}_n^V}\left(\varphi_Z\right)$.
\end{proposition}
\begin{proof}
This follows from \eqref{Fsigthetphi}, Proposition \ref{Propmultfreebis}\eqref{Propmultfreebis8}, and Proposition \ref{Propchifphi}.
\end{proof}


At the end of this section, we represent the spherical functions as functions that are constant on the orbits described in Theorem \ref{decOrbxiphib}. We need some preliminary results.

\begin{itemize}

\item
For every $\psi\in\Sigma_0$, from \eqref{Stabpsi} and Proposition \ref{Propmultfreebis}\eqref{Propmultfreebis3} we deduce that
\[
\left\lvert\mathsf{Stab}_{\mathsf{I}(\psi)}\left(Z_0\right)\right\rvert=\prod_{z\in\mathsf{image}(\psi)}\left\lvert\psi^{-1}(z)\right\rvert!.
\]

\item
For every $\psi\in\Sigma_0$, from \eqref{assumpSigma0} and Proposition \ref{Propmultfreebis}\eqref{Propmultfreebis7} we deduce that
\[
\left\lvert\mathsf{Q}^0_\psi\right\rvert
=\binom{\left\lvert Z_0\right\rvert}{\left(\left\lvert\psi^{-1}(z)\right\rvert\right)_{z\in\mathbb{Z}_b}}
\binom{\left\lvert V\setminus Z_0\right\rvert}{\left(\left\lvert\psi^{-1}(z)\right\rvert\right)_{z\in\mathbb{Z}_a}},
\]
where $\left\lvert\psi^{-1}\left(\mathbb{Z}_b\right)\right\rvert=\left\lvert Z_0\right\rvert$ and $\left\lvert\psi^{-1}\left(\mathbb{Z}_a\right)\right\rvert=\left\lvert V\setminus Z_0\right\rvert$.

\item
From \eqref{Orbxiphib}, for all $0\leq j\leq\min\left\{\left\lvert Z_0\right\rvert,\left\lvert V\setminus Z_0\right\rvert\right\}$ and $(\tau,\tau')\in\mathsf{T}(j)$, we can easily obtain the coefficients in the orthogonality relations \eqref{orthrel}:
\[
\begin{split}
\left\lvert\digamma_{j,\tau,\tau'}\right\rvert
&=\sum_{(Z,\mathsf{i})\in\Psi_{j,\tau,\tau'}}a^jb^j=\left\lvert\Psi_{j,\tau,\tau'}\right\rvert a^jb^j\\
&=\binom{\left\lvert Z_0\right\rvert}{j}\binom{\left\lvert V\setminus Z_0\right\rvert}{j}\binom{\left\lvert Z_0\right\rvert - j}{\tau(0),\dotsc,\tau(b-1)}\binom{\left\lvert V\setminus Z_0\right\rvert -j}{\tau'(0),\dotsc,\tau'(a-1)}(ab)^j
\end{split}
\]
and $\left\lvert\mathsf{Orb}_{\mathsf{C}_n\wr\mathsf{Sym}(V)}\left(\varphi_{Z_0}\right)\right\rvert=\dbinom{\lvert V\rvert}{\left\lvert Z_0\right\rvert}a^{\left\lvert V\setminus Z_0\right\rvert}b^{\left\lvert Z_0\right\rvert}$.

\end{itemize} 

With the notation in \eqref{defPsixijGel} and Proposition \ref{Propmultfreebis}, taking Proposition \ref{Propchifphi} into account, for every $\psi\in\Sigma_0$ and every $(Z,\mathsf{i})\in\bigsqcup_{j=0}^m\Psi_j$ with $Z\in\mathsf{Orb}_{\mathsf{I}(\psi)}\left(Z_0\right)$ (that is, $Z\in\Pi_j$ when $(Z,\mathsf{i})\in\Psi_j$), we define $\check{\chi}_\psi(\mathsf{i})$ by setting
\[
\check{\chi}_\psi(\mathsf{i})\coloneqq\chi_\psi(\varphi),
\]
where $\varphi\in\bigtimes_{v\in V}\Upsilon_{\varphi_Z^\flat(v),\varphi_{Z_0}^\flat(v),\mathsf{i}(v)}$ is arbitrary (since $\psi\prec Z$, the value $\chi_\psi(\varphi)$ does not depend on the choice of $\varphi$). This is \eqref{defcheckSsingb} in the present context. Henceforth, recalling \eqref{barbarcalS} and \eqref{decPsijtau} and taking Proposition \ref{Propmultfreebis} into account, suppose that $0\leq j\leq m$, $(\tau,\tau')\in\mathsf{T}(j)$, and $(Z,\mathsf{i})\in\Psi_{j,\tau,\tau'}$. Then we set
\[
\check{\check{\chi}}_{j,\psi}(\tau,\tau')\coloneqq\frac{1}{\left\lvert\mathsf{Stab}_{\mathsf{I}(\psi)}\left(Z_0\right)\right\rvert}\sum_{Z''\in\Pi_j}\;\sum_{\substack{\zeta\in\mathsf{Stab}_{\mathsf{Sym}(V)}\left(Z_0\right)\mathrel{:}\\ \zeta^{-1}Z=Z''}}\check{\chi}_{\psi}\left(\zeta^{-1}\mathsf{i}\right).
\]
In particular, the value of this sum depends only on $j$, $\tau$, and $\tau'$, and corresponds to \eqref{barbarcalS}; by Proposition \ref{Propmultfreebis}\eqref{Propmultfreebis6}, the index $l$ disappears. Setting $\theta_k\coloneqq\overline{S}^{h-k,k}$ for $0\leq k\leq m$, for all $0\leq j\leq m$ and $(\tau,\tau')\in\mathsf{T}(j)$, the value of $\mathcal{F}_{\psi,\theta_k}$ at $(Z,\varphi)\in\digamma_{j,\tau,\tau'}$ depends only on $j$ and $(\tau,\tau')$ and is given by
\[
\mathcal{F}_{\psi,\theta_k}(Z,\varphi)=\frac{1}{\left\lvert\mathsf{Q}^0_\psi\right\rvert}\check{\mathcal{R}}_{\psi,\theta_k}(j)\check{\check{\chi}}_{j,\psi}(\tau,\tau')\eqqcolon\check{\mathcal{F}}_{\psi,\theta_k}(j,\tau,\tau'),
\]
whereas $\mathcal{F}_{\psi,\theta_k}(Z,\varphi)=0$ for $(Z,\varphi)\in\digamma_{j,\tau,\tau'}$ with $j>m$. This corresponds to \eqref{defFsigthetorb}; by Proposition \ref{Propmultfreebis}\eqref{Propmultfreebis6}, there is no summation. Finally, we point out that the spherical functions $\check{\mathcal{R}}_{\psi,\theta_k}(j)$ can be expressed in terms of Hahn polynomials; see, e.g., \cite[Remark 6.2.4]{book}.

\vspace{1cm}

\noindent
Fabio Scarabotti, Dipartimento SBAI, Sapienza Universit\`a di Roma, via A. Scarpa 8, 00161, Roma, Italy\\
\noindent
Istituto Nazionale di Alta Matematica ``Francesco Severi'', 00185 Roma, Italy\\
\noindent
 {\it e-mail:} {\tt fabio.scarabotti@sbai.uniroma1.it}

\end{document}